\documentclass[12pt]{memo-l}

\usepackage{amsmath,amsthm,amssymb}
\usepackage{lscape}

\newtheorem{thm}{Theorem}[chapter]
\newtheorem{lem}[thm]{Lemma}
\newtheorem{prop}[thm]{Proposition}
\newtheorem{cor}[thm]{Corollary}
\newtheorem{exl}[thm]{Example}
\newtheorem{hypo}[thm]{Hypothesis}

\theoremstyle{dfn}
\newtheorem{dfn}[thm]{Definition}

\theoremstyle{rem}
\newtheorem{rem}[thm]{Remark}

\numberwithin{section}{chapter}
\numberwithin{equation}{chapter}
\newenvironment{prf}{{\sc Proof.}}{\hfill $\diamond$}

\def\char{ \text{\rm char } }

\def\Aut{ \text{\rm Aut} }
\def\Hom{ \text{\rm Hom} }
\def\Irr{ \text{\rm Irr} }
\def\IBr{ \text{\rm IBr} }
\def\Ind{ \text{\rm Ind} }
\def\End{ \text{\rm End} }
\def\Res{ \text{\rm Res} }

\def\Spin{ \text{\rm Spin} }
\def\Stab{ \text{\rm Stab} }

\def\GL{ \text{\rm GL} }
\def\GU{ \text{\rm GU} }
\def\GO{ \text{\rm GO} }
\def\SO{ \text{\rm SO} }
\def\HN{ \text{\rm HN} }
\def\ON{ \text{\rm O'N} }
\def\HS{ \text{\rm HS} }
\def\Suz{ \text{\rm Suz} }
\def\Ru{ \text{\rm Ru} }
\def\He{ \text{\rm He} }
\def\Fi{ \text{\rm Fi} }
\def\Co{ \text{\rm Co} }
\def\Ly{ \text{\rm Ly} }
\def\Th{ \text{\rm Th} }
\def\McL{ \text{\rm McL} }
\def\Sp{ \text{\rm Sp} }
\def\CSp{ \text{\rm CSp} }
\def\PSp{ \text{\rm PSp} }
\def\SL{ \text{\rm SL} }
\def\PSL{ \text{\rm PSL} }
\def\PGL{ \text{\rm PGL} }
\def\SU{ \text{\rm SU} }
\def\PSU{ \text{\rm PSU} }
\def\SO{ \text{\rm SO} }

\def\CSO{ \text{\rm CSO} }

\begin{document}

\frontmatter

\title{Imprimitive irreducible modules for finite quasisimple groups}

%    Remove any unused author tags.

%    author one information
\author{Gerhard Hiss}
\address{Lehrstuhl D f{\"u}r Mathematik, RWTH Aachen University,
52056 Aachen, Germany}
%\curraddr{}
\email{gerhard.hiss@math.rwth-aachen.de}
%\thanks{}

%    author two information
\author{William J.~Husen}
\address{Department of Mathematics, The Ohio State University, 
OH 43210-1174, U.S.A.}
%100 Math Tower, 231 West 18th Avenue, Columbus, OH 43210-1174, U.S.A.}
%\curraddr{}
\email{husen@math.ohio-state.edu}
%\thanks{}

\author{Kay Magaard}
\address{School of Mathematics, University of Birmingham, Edgbaston,
Birmingham B15 2TT, UK}
%\curraddr{}
\email{k.magaard@bham.ac.uk}
%\thanks{}

%    \date is required; it is the date received by the editor.
\date{26.11.2012}

%    The 2010 edition of the Mathematics Subject Classification is
%    now available.  If you are citing a classification from the
%    new scheme, use the following input coding instead.
%\subjclass[2010]{Primary }


%\dedication{Dedication text (use \\[2pt] for line break if necessary)}

\begin{abstract}
Motivated by the maximal subgroup problem of the finite classical groups
we begin the classification of imprimitive irreducible modules of finite 
quasisimple groups over algebraically closed fields~$K$. A module of a 
group~$G$ over~$K$ is {\em imprimitive}, if it is induced from a module 
of a proper subgroup of~$G$.

We obtain our strongest results when $\mbox{\rm char}(K) = 0$, although
much of our analysis carries over into positive characteristic. If~$G$
is a finite quasisimple group of Lie type, we prove that an imprimitive
irreducible $KG$-module is Harish-Chandra induced. This being true for
$\mbox{\rm char}(K)$ different from the defining characteristic of~$G$,
we specialize to the case $\mbox{\rm char}(K) = 0$ and apply Harish-Chandra 
philosophy to classify irreducible Harish-Chandra induced modules in terms 
of Harish-Chandra series, as well as in terms of Lusztig series. 
We determine the asymptotic proportion of the irreducible imprimitive 
$KG$-modules, when~$G$ runs through a series groups of fixed (twisted) Lie 
type. 
One of the surprising outcomes of our investigations is the fact that
these proportions tend to~$1$, if the Lie rank of the groups tends to 
infinity.

For exceptional groups~$G$ of Lie type of small rank, and for sporadic
groups~$G$, we determine all irreducible imprimitive $KG$-modules for 
arbitrary characteristic of~$K$.
\end{abstract}

\maketitle

\tableofcontents

%    Include unnumbered chapters (preface, acknowledgements, etc.) here.
%\include{}

\section*{Acknowledgements}
It is our pleasure to thank the many individuals who have made significant
contributions to this monograph. First of all we thank Thomas Breuer for his
invaluable help with groups and character tables, as well as for computing 
a large number of character tables of maximal subgroups.
Without Thomas our results would not be as complete as they are now.

We also thank Frank L{\"u}beck for providing Lemma~\ref{FranksLemma},
for proving the results of~\cite{LuebeckNegativeCycles} on demand,
for various fruitful discussions on the
subjects of this article, as well as for reading large parts of it.
Thanks are also due to Meinolf Geck for suggesting the argument
in the proof of Lemma~\ref{WeylInd}, as well as to C{\'e}dric
Bonnaf{\'e} for providing Example~\ref{Bonnafe}.

We thank Felix Noeske for computing a basic set modulo~$2$ for the group
$3.\Fi_{24}'$ and we thank Rod Gow for bringing reference \cite{SchurIndicators} 
to our attention.

We are particularly grateful to Gunter Malle for his careful reading of our 
manuscript, for pointing out several errors and for suggesting numerous 
improvements and generalizations, most notably in the formulation of 
Theorem~\ref{FirstLusztigSeriesResultConverse}.

Lastly, the third author would like to thank the NSA for its generous 
support through grant H98230-07-1-0049 and its predecessors.
\mainmatter
%    Include main chapters here.
%\include{}

\chapter{Introduction}

This monograph is a contribution to the classification of the maximal
subgroups of the finite classical groups and the start of the program of 
classifying 
imprimitive modules for finite quasisimple groups. The two programs are 
inextricably linked in a way we now explain. Suppose that we want to study 
the maximal subgroups of a finite classical group~$X$ whose natural module 
is~$V$.

In his 1984 paper on the maximal subgroups of classical groups~\cite{Asch1}
Aschbacher defines eight collections of geometric subgroups $\mathcal{C}_i(X)$,
$1 \leq i \leq 8$. His main theorem asserts that if~$L$ is maximal in
a classical group~$X$ with natural module~$V$, then either~$L$ is an element
of $\mathcal{C}_i(X)$ for some $1 \leq i \leq 8$, or the following are true:
\begin{enumerate}
\item $F^*(L)'$ is quasisimple;
\item $F^*(L)'$ acts absolutely irreducibly on~$V$;
\item the action of~$F^*(L)'$ on~$V$ can not be defined over a smaller field;
\item any bilinear, quadratic or sesquilinear form on~$V$ that is stabilized 
by~$F^*(L)'$ is also stabilized by~$X$.
\end{enumerate}
(Here, $F^*(L)$ denotes the generalized Fitting subgroup of~$L$ and
$F^*(L)'$ its derived subgroup.) 
Asch\-bacher denotes the collection of subgroups of~$X$ which satisfy 
Conditions~$1$ through~$4$ above by $\mathcal{S}$.

In their recent book~\cite{BrayHoltDougal}, Bray, Holt and Roney-Dougal have 
determined the maximal subgroups of the classical groups with $\mbox{\rm dim}(V) 
\leq 12$. This dimension is chosen to complement the work by Kleidman and 
Liebeck in~\cite{kleidman:4}. In this monograph the authors determine the 
maximal members of each collection $\mathcal{C}_i(X)$, and for 
$\mbox{\rm dim}(V) > 12$ they determine those instances where a maximal member 
of $\mathcal{C}_i(X)$ is a maximal subgroup of~$X$. In case a maximal member~$L$ 
of some $\mathcal{C}_i(X)$ fails to be maximal in~$X$, the overgroups of~$L$ are 
also determined. To complete the determination of all maximal subgroups of~$X$ 
one needs to answer the following question. When is a member of~$\mathcal{S}$ 
maximal in~$X$?

We note that if $L \leq X$ is a member of $\mathcal{S}$, then so is $N_X(L)$.
Hence in order to investigate the question of the maximality of~$L$ in~$X$, we
may assume that $L = N_X(L)$. 
Now if $L \leq \tilde{L} \leq X$, then the definition of $\mathcal{S}$ implies 
that $\tilde{L} \not \in \mathcal{C}_1 \cup \mathcal{C}_3 \cup \mathcal{C}_5 
\cup \mathcal{C}_8$. If $\tilde{L} \in \mathcal{C}_4 \cup \mathcal{C}_7$ 
then~$L$ stabilizes a tensor product decomposition of~$V$. Examples of such 
$(L,V)$ can be found in the paper \cite{MagaardTiep} by Magaard and Tiep, and 
in the papers~\cite{BK1,BK2,BK3} by Bessenrodt and Kleshchev.
Examples of $(L,V)$ where $\tilde{L} \in \mathcal{C}_6$ can be found in 
a preprint by Magaard and Tiep \cite{MagaardTiep2}.

The examples $(L,V)$ where $\tilde{L} \in \mathcal{S}$ come in several 
varieties depending on the isomorphism types of~$L$,~$\tilde{L}$ and~$X$. For 
instance by generalizing work of Dynkin, the examples where $L$,~$\tilde{L}$
and~$X$ are of Lie type and of the same characteristic, were classified by 
Seitz and Testerman (see \cite{seitz:0,SeitzTesterman,Testerman}).
Liebeck, Saxl and Seitz classified the examples where~$L$ and~$X$ are of Lie 
type of equal characteristic, but~$\tilde{L}$ is not (see 
\cite{LiebeckSaxlSeitz0}).
A list of possible examples where~$L$ and~$\tilde{L}$ are of Lie type and equal 
characteristic but of characteristic coprime to~$V$ was determined by Seitz 
in~\cite{seitz:1}. In \cite{SchaefferFry} Schaeffer Fry has resolved 
one family of these difficult cases.
In recent work Magaard, R\"ohrle and Testerman produced a possible list of 
examples where~$X$ and~$\tilde{L}$ are of Lie type and of equal characteristic, 
and~$L$ is of Lie type but of characteristic coprime to that of~$X$ (see
\cite{MaRoTe}). Finalizing this list is work in progress by Magaard and 
Testerman.
The case where $F^*(L)'$ is alternating was the subject of Husen's thesis and 
subsequent work. Only one infinite series of examples occurs when~$V$ is the 
$k$th exterior power of the reduced permutation module~$W$ of $F^*(L)'$, where
$2 \leq k \leq \mbox{\rm dim}(W)/2$. In this case $X = \Omega(V)$ (if~$X$ is
quasisimple) and $F^*( \tilde{L} )' = \Omega(W)$. For these results see 
\cite{Hus1,Hus2,Hus3}.

The situation to which the results of this monograph apply is when $\tilde{L} 
\in \mathcal{C}_2$, which implies that~$L$ acts imprimitively on~$V$. 
In particular this implies that~$F^*(L)'$ also acts imprimitively on~$V$. 
As~$F^*(L)'$ is quasisimple, this is precisely the situation we study here.

To introduce our results, let~$K$ be a field and~$G$ a finite group. A 
$KG$-module will always be assumed to be finitely generated with~$G$ acting
from the right.

\begin{dfn}
{\rm
We say that a $KG$-module~$M$ is {\em imprimitive with block stabilizer}~$H$,
if and only if there are subspaces $M_1, M_2, \ldots, M_n$ with $n =
[G\colon\!H] > 1$ such that
\begin{equation}
\label{DecompositionOfM}
M = \bigoplus_{i = 1}^n M_i,
\end{equation}
the~$M_i$ are transitively permuted by the action of~$G$ and~$H$ is the 
stabilizer of~$M_1$, i.e., $H = \{ g \in G \mid M_1g = M_1 \}$.
}
\end{dfn}
In other words a $KG$-module~$M$ is imprimitive with block stabilizer~$H$
if and only if $H$ is a proper subgroup of~$G$ and there is a
$KH$-module~$M_1$ such that $M \cong \Ind_H^G( M_1 ) := 
M_1 \otimes_{KH} KG$.

Let us now describe the principal results of this monograph in more detail.
If a proper subgroup~$H$ of~$G$ is the block stabilizer of some $KG$-module,
then so is any intermediate subgroup $H \leq L \leq G$, since induction is 
transitive. Thus in order to classify the imprimitive irreducible $KG$-modules,
one may restrict attention to the maximal subgroups~$H$ of~$G$ as potential 
block stabilizers. We therefore assume that~$H$ is a maximal subgroup of~$G$ in 
the following. We also assume that~$K$ is algebraically closed. The problem
of classifying the imprimitive irreducible modules of quasisimple groups
over an arbitrary field~$K$ does not appear to have an immediate reduction to 
the case of algebraically closed fields. With respect to our intended 
application to maximal subgroups, the restriction to algebraically closed 
fields~$K$ is justified by property~(2) of the groups of class~$\mathcal{S}$,
as discussed at the beginning of this introduction. 

In Chapter~$2$ we collect some general results from representation theory which
are used to dispose of some subgroups as possible block stabilizers or to 
compare imprimitivity in different characteristics. We also comment on some
aspects of our notation here.

Chapter~\ref{SectionSporadicGroups} is devoted to the sporadic simple groups 
and their covering groups. Let~$G$ denote such a group. We determine all 
imprimitive irreducible $KG$-modules for fields~$K$ of 
arbitrary characteristic (Theorem~\ref{SporadicTheorem}). Thus this chapter 
contains the most complete results. 

Next let~$G$ be a quasisimple group such that $G/Z(G)$ is an alternating 
group. There is a complete classification of the imprimitive absolutely 
irreducible
$KG$-modules for fields~$K$ of characteristic~$0$. We cite the corresponding 
results from the literature in Chapter~\ref{alternating}. The results for the 
alternating groups are due to D.~\v{Z}.~Djokovi\'c and J.~Malzan \cite{DM1,DM2} 
and are about $35$ years old. More recent are the results for the covering groups 
due to D.~Nett and F.~Noeske \cite{NeNoe}. Here, the classification over fields
of positive characteristic remains open, although \cite{NeNoe} restricts the 
possible block stabilizers.

The remaining chapters deal with the quasisimple covering groups of the
finite groups of Lie type. Thus let~$G$ be a quasisimple group with $G/Z(G)$
a finite simple group of Lie type. The classification of the imprimitive
irreducible $KG$-modules in the defining characteristic case, i.e., where
the characteristic of~$K$ equals that of~$G/Z(G)$, is due to Gary Seitz. There are
only finitely many such instances. We state Seitz's result for completeness.
\begin{thm}\label{SameCharacteristic}{\rm (Seitz) \cite[Theorem~$2$]{seitz:2}}
Suppose that~$G$ is defined over a field of characteristic~$p$.
    If $M$ is an irreducible and imprimitive $KG$-module and
    $\mbox{\rm char}(K)=p$, then $G$ is $\SL_{2}(5)$, $\SL_{2}(7)$, 
    $\SL_{3}(2)$ or $\Sp_{4}(3)$ and $M$ is the Steinberg module for $G$.
\end{thm}
There is one further case not covered by the above theorem, since it involves
a group which has two defining characteristics, and is not a simple gorup of 
Lie type in the strict sense. Namely, let $G = {^2\!G}_2(3)' \cong \SL_2(8)$,
and $\mbox{\rm char}( K ) = 3$. Then the three $9$-dimensional irreducible 
$KG$-modules are imprimitive. In turn, the Steinberg module of ${^2\!G}_2(3)$
over~$K$ is induced from each of two $3$-dimensional $KH$-modules, where~$H$
is the normalizer of Sylow $2$-subgroup of ${^2\!G}_2(3)$.

From now on we will assume that the defining characteristic of $G/Z(G)$ is 
different from the characteristic of~$K$.
In all but finitely many cases,~$G$ is itself a finite group of Lie type. The 
exceptions occur if $G/Z(G)$ has an exceptional Schur multiplier. These cases 
are treated in Chapter~$5$, very much in the same way as the sporadic groups. 
Again, we obtain a complete classification for all characteristics with the 
only exception if $G/Z(G) = {^2\!E}_6(2)$ and $H/Z(G)$ is a parabolic subgroup
of $G/Z(G)$, in which case we obtain this classification only for
$\mbox{\rm char}(K) = 0$.

We are thus left with the cases that $G$ and $G/Z(G)$ are finite groups of
Lie type and that the characteristic of~$K$ is different from that
underlying~$G$. 
Our first main result here is Theorem~\ref{ReductionToParabolic}:
{\em If~$H$ is a maximal subgroup of~$G$ and a block stabilizer of an 
imprimitive irreducible $KG$-module, then $H$ is a parabolic subgroup of~$G$.} 
All of Chapter~\ref{CrossCharacteristic} is devoted to the proof of this 
theorem.

The remainder of our monograph is then devoted to the determination of
imprimitive irreducible $KG$-modules whose block stabilizer is a parabolic
subgroup~$H$ of~$G$. We first prove in this situation
that if $M$ is induced from a $KH$-module~$M_1$, then
the unipotent radical of~$H$ acts trivially in~$M_1$, in other words, $M$ is
Harish-Chandra induced (Proposition~\ref{ParabolicStabilizer}). Thus one can
apply the well developed machinery of Harish-Chandra theory. This is exploited
in Chapter~\ref{InductionFromParabolicSubgroups} to give sufficient conditions
on irreducible $KH$-modules to induce to irreducible $KG$-modules. These
conditions are given in terms of the Harish-Chandra parametrization 
(Theorem~\ref{HCresult}) and of the Lusztig parametrization 
(Theorem~\ref{FirstLusztigSeriesResult}) of the relevant irreducible 
$KH$-modules. These results are powerful enough to determine the asymptotic
proportion of the irreducible imprimitive $KG$-modules, when~$G$ runs through 
a series groups of fixed (twisted) Lie type (Theorem~\ref{Asymptotic}).
Moreover, for classical groups, these proportions tend to~$1$, if the Lie 
rank of the groups tends to infinity (Example~\ref{ExamplesForAsymptotics}).

Let us comment on a consequence of this result for the problem of classifying
maximal subgroups of a classical group~$X$. Suppose that~$L$ is a subgroup 
of~$X$ of class~$\mathcal{S}$ with $N_X( L ) = L$ and $G := F^*( L )'$ of Lie 
type of characteristic different from the characteristic of~$X$. Suppose that 
the natural module~$V$ of~$X$ is an imprimitive Harish-Chandra induced 
$KG$-module with block stabilizer~$P$. If the $G$-conjugacy class of~$P$ is 
$L$-invariant, the Frattini argument implies that $L = N_L(P)G$. 
Consequently,~$V$ is also an imprimitive $KL$-module. Indeed, let~$U$ denote
the unipotent radical of~$P$. As~$V$ is irreducible and Harish-Chandra induced
from~$P$, there is a decomposition $V = V_1 \oplus V_2 \oplus \cdots \oplus V_n$
such that $V_1 = C_V(U)$, and the $V_i$ are permuted by the action of~$G$. 
As~$U$ is fixed by~$N_L( P )$, the $KP$-module~$V_1$ is stabilized by~$N_L(P)$. 
It now follows from $L = N_L(P)G$, that~$L$ permutes the~$V_i$, and hence~$V$
is induced from the $KN_L(P)$-module~$V_1$. Thus~$L$ is a $\mathcal{C}_2$-type
subgroup of~$X$, and as such, not maximal in~$X$. Indeed, if $\tilde{L}$ is a 
maximal subgroup of~$X$ of $\mathcal{C}_2$-type such that $F^*(\tilde{L})'$ is 
quasisimple, then $F^*(\tilde{L})'$ does not act irreducibly on~$V$. Therefore 
our results imply the surprising fact that most (in the sense discussed at the
end of the previous paragraph) cross-characteristic irreducible 
modules for groups of Lie type do not lead to examples of maximal 
subgroups of classical groups.

Harish-Chandra theory involves the investigation of representations of certain
Iwahori-Hecke algebras. In case $\text{\rm char}(K) = 0$, this can be reduced
to the investigation of representations of Weyl groups. This programme is
carried out in Chapter~\ref{Characteristic0}, where we prove partial
converses to the results of Chapter~\ref{InductionFromParabolicSubgroups}.
In fact we show that the sufficient conditions derived there are 
necessary, if the groups~$G$ considered arise from algebraic groups with 
connected center (Theorem~\ref{HCresult:converse}), or, more generally, if the 
characters considered lie in Lusztig series corresponding to semisimple elements 
whose centralizer is connected (Theorem~\ref{FirstLusztigSeriesResultConverse}).
As a consequence, these results can not be applied to all characters of
quasisimple groups of Lie type.

The general results of Chapter~\ref{Characteristic0} are now specialized to
the classical groups in Chapter~\ref{ClassicalInCharacteristic0}. Again we
assume that $\text{\rm char}(K) = 0$, and that the groups~$G$ arise from
algebraic groups with connected center. We obtain a complete classification
of the Harish-Chandra imprimitive irreducible $KG$-modules in terms of
Harish-Chandra series (Propositions~\ref{PrimitivesInClassicalI} and
\ref{PrimitivesInClassicalII}) as well as in terms of Lusztig series
(Propositions~\ref{PrimitivesInGLn} and~\ref{PrimitivesInOtherClassicals}).

The final Chapter~\ref{ExceptionalGroups} is devoted to the exceptional
groups of Lie type. We obtain a complete, characteristic independent
and explicit classification of the imprimitive irreducible $KG$-modules
if~$G$ is an exceptional group of Lie type of small rank, i.e.,\ $G$
is one of ${^2\!B}_2(q)$, ${^2\!G}_2(q)$, ${^2\!F}_4(q)$,
$G_2(q)$, and ${^3\!D}_4(q)$. For the exceptional groups of Types~$E$
and~$F$ our results require the assumption $\text{\rm char}(K) = 0$.

\chapter{Generalities}

In this chapter we start with some comments on our notation, in particular
when it differs from that of our sources. Further notation is introduced as 
we go along.  We also collect 
a few general results of various types, needed in later chapters. This chapter 
can be skipped at a first reading as it mainly serves as a convenient reference.
Throughout,~$G$ denotes a group and~$K$ a field.

\section{Comments on the notation}

\subsection{Notation for groups}
Let~$G$ be a group. Our group actions are generally right actions. Thus for 
$x, y \in G$, we put
$$x^y := y^{-1}xy$$
and
$$[x,y] := x^{-1}y^{-1}xy = x^{-1}x^y.$$
If~$X$ and~$Y$ are subsets of~$G$ we write 
$$[X,Y] := \langle [x,y] \mid x \in X, y \in Y \rangle.$$
In particular, $[G,G]$ is the commutator subgroup of~$G$. If there is no danger
of confusion, we occasionally write~$G'$ for $[G,G]$. We write $Z(G)$ for the
center of~$G$, and $C_G( H )$ and $N_G( H )$ for the centralizer and normalizer
in~$G$ of a subgroup~$H$, respectively. Suppose that~$G$ is finite.
If~$H$ is a subgroup of~$G$, then $[G\colon\!H]$ denotes the index of~$H$ in~$G$. 
As usual, $F^*( G )$ denotes the generalized Fitting subgroup of~$G$. Also, 
if~$\ell$ is a prime, $O_\ell(G)$ is the largest normal $\ell$-subgroup of~$G$.

For group extensions we occasionally use the Atlas \cite{ATLAS} convention,
i.e., $A.B$ denotes a group~$G$ with a normal subgroup~$N$ isomorphic to~$A$ 
and $G/N$ isomorphic to~$B$. The symbol $A\colon\!B$ denotes a split extension,
$A \times B$ the direct product of~$A$ and~$B$, and $A \circ B$ a central 
product, i.e., a group~$G$ containing normal subgroups $A$ and~$B$ with
$[A,B] = 1$ and $AB = G$.
A cyclic group of order~$n$ is sometimes just denoted by the symbol~$n$.

The alternating and symmetric groups of degree~$n$ are denoted by~$A_n$ 
and~$S_n$, respectively. The quaternion group of order~$8$ is denoted by~$Q_8$.
Our notation for the classical groups is more traditional than that of the 
Atlas. We follow essentially the notation used by Wilson in \cite{WilsonBook}.
Thus $\GL_n( \mathbb{F} )$, $\SL_n( \mathbb{F} )$, and 
$\PSL_n( \mathbb{F} )$ denote, respectively, the general linear group, the 
special linear group and the projective special linear group of degree~$n$
over the field~$\mathbb{F}$. If~$\mathbb{F}$ is finite with~$q$ elements, we
write $\GL_n( q )$ etc. We write $\GU_n( q )$ for the unitary group which 
is a subgroup of $\GL_n( q^2 )$. Our notation for orthogonal groups is as
follows. If $n \geq 2$ is even, we write $\GO_n^+( q )$ and $\GO_n^-( q )$
for subgroup of $\GL_n( q )$ preserving a non-degenerate quadratic form on
$\mathbb{F}_q^n$ of Witt index $n/2$ and $n/2 - 1$, respectively.
If $n \geq 3$ is odd, $\GO_n( q )$ denotes the subgroup of $\GL_n( q )$ 
preserving a non-degenerate quadratic form on $\mathbb{F}_q^n$. If~$n$ is 
odd and~$q$ is even, by a non-degenerate quadratic form on $\mathbb{F}_q^n$
we mean a quadratic form whose polar form has a $1$-dimensional radical 
consisting of non-isotropic vectors. For odd~$n$ we write $\GO_n^0( q ) := 
\GO_n( q )$ so that $\GO^\epsilon_n( q )$ for $\epsilon \in \{ 0, + , - \}$ 
denotes our three types of orthogonal groups. Let $\epsilon \in 
\{ 0, + , - \}$. If~$q$ is odd, we put $\SO^\epsilon_n( q ) := 
\GO^\epsilon_n( q ) \cap \SL_n( q )$. We also
put $\Omega^\epsilon( q ) := [\GO^\epsilon_n( q ), \GO^\epsilon_n( q )]$ for 
the commutator subgroup of $\GO^\epsilon_n( q )$. Finally, 
$\Spin_n^\epsilon( q )$ denotes the corresponding spin groups, i.e., the 
(generic) Schur covering groups of the $\Omega^\epsilon_n( q )$.
Notice that $\Omega^\epsilon_n( q )$ is quasisimple if $n = 3$ and $q > 3$,
if $n = 4$ and $\epsilon \neq +$, if $n = 5$ and $q > 2$, or if $n \geq 6$.
For $n = 2$ we have that $\Omega^\epsilon_2( q )$ is cyclic of order 
$q - \epsilon \cdot 1$. If $\Omega^\epsilon_n( q )$ is quasisimple, the 
central factor group $P\Omega^\epsilon_n( q ) := 
\Omega^\epsilon_n( q )/Z( \Omega^\epsilon_n( q ) )$ is simple. As
$Z( \Omega^\epsilon_n( q ) )$ is trivial if $q$ is even or~$n$ is odd,
we just write $\Omega^\epsilon_n( q )$ for $P\Omega^\epsilon_n( q )$ in
these cases.

Of course, there is also a version of these symbols emphasizing the underlying 
vector space carrying a form. For example, if~$V$ is a vector space equipped 
with a non-degenerate quadratic form, $\GO( V )$ denotes the orthogonal group 
of the quadratic space~$V$, i.e., the automorphisms of~$V$ preserving the given 
form.

Suppose that~$V$ is a vector space equipped with a bilinear, sesquilinear or
quadratic form. A subspace~$W$ of~$V$ is called totally isotropic, if the 
form is constantly zero on $W \times W$ or $W$, respectively.

\subsection{Notation from linear algebra}
\label{NotationFromLinearAlgebra}
If~$n$ is a positive integer, we denote by~$I_n$ the $n \times n$ identity
matrix. We write~$J_n$ for the $n \times n$ matrix with ones along the
anti-diagonal and zeros elsewhere. Finally,~$\tilde{J}_{2n}$ denotes the
$2n \times 2n$ matrix defined by
$$\tilde{J}_{2n} = 
\left( \begin{array}{cc}
0_n & J_n \\
-J_n & 0_n
\end{array} \right),$$
where $0_n$ is the $n \times n$ zero matrix. These matrices are viewed as
matrices over the ring currently considered. If~$n$ is clear from the context,
we simply write~$I$,~$J$ and~$\tilde{J}$ for~$I_n$, $J_n$ and~$\tilde{J}_{2n}$,
respectively.

If~$A$ is a matrix, its transpose is denoted by $A^T$, and $A^{-T}$ is defined
by $A^{-T} := (A^{-1})^T$.

If $a$, $b$, $c$, \ldots are elements of a ring or square matrices over this 
ring (not necessarily of the same degree), we write $\mbox{\rm diag}( a, b, c, 
\ldots )$ for the (block) diagonal matrix with diagonal entries $a, b, c, \ldots$

Let~$V$ be a vector space and $G \leq \GL( V )$. We write~$I_V$ for the identity
element of $\GL( V )$. If~$U$ is a $G$-invariant
subspace of~$V$ and $x \in G$, we write $x_U$ for the restriction of~$x$
to~$U$, viewed as a linear map on~$U$. Moreover, we write
$$G_U := \{ x_U \mid x \in G \}.$$

\subsection{Notation from representation theory}

\label{NotationRepresentationTheory}
%Let~$G$ be a finite group and~$K$ a field. 
By a $KG$-mod\-ule we always mean
a finite dimensional right $KG$-module. It will be convenient occasionally to 
consider the characters of $KG$-modules. If~$K$ is a splitting field for
$KG$ and $\text{\rm char}(K) = \ell$, we write $\Irr(G)$ for the set of 
characters of the irreducible $KG$-modules if $\ell = 0$, and $\IBr_\ell(G)$ 
for the set of Brauer characters (with respect to a fixed $\ell$-modular system 
with residue class field~$K$) of the irreducible $KG$-modules, if $\ell > 0$. 
If~$M$ and~$N$ are $KG$-modules,
we put $[ M, N ]_{KG} = \mbox{\rm dim}_K \Hom_{KG}( M , N )$, omitting the
subscript~$KG$ if there is no danger of confusion.

If~$H$ is a subgroup of~$G$, we write $\Ind_H^G( - )$ and $\Res_H^G( - )$
for induction and restriction of $KH$-modules, respectively $KG$-modules.
The same symbols are used for induction and restriction of characters.

\section{Conditions for primitivity}

Since we are interested in imprimitive irreducible $KG$-modules, we collect
a few conditions guaranteeing that an induced $KG$-module is reducible. We
begin with the most general one.

\begin{lem}\label{SquaredSubgroup}
Suppose that~$K$ is a splitting field for~$KG$.
Let~$H$ be a subgroup of~$G$. If $|H|^2 < |G|$, then
$\Ind_H^G( M_1 )$ is reducible for all $KH$-modules~$M_1$.
\end{lem}
\begin{prf}
Let~$M_1$ be a $KH$-module. Then 
$$\mbox{\rm dim}_K(\Ind_H^G( M_1 ))^2 =
\mbox{\rm dim}_K(M_1)^2\,[G:H]^2 > \mbox{\rm dim}_K(M_1)^2\,|G| \geq |G|.$$ 
As $\mbox{\rm dim}_K(N)^2 \leq |G|$ for all irreducible $KG$-modules~$N$,
our assertion follows.
\end{prf}

We shall need the following easy consequence of Mackey's theorem, which
was also used implicitly in \cite{DM1,DM2}.
\begin{lem} \label{Mackey1}
Let~$H$ be a subgroup of~$G$. Suppose that there exists an element
$t \in G \setminus H$ such that~$t$ centralizes the intersection 
$H \cap H^t$. Then $\Ind_H^G( M_1 )$ is reducible for all 
$KH$-modules~$M_1$.
\end{lem}
\begin{prf}
With the notation introduced in Subsection~\ref{NotationRepresentationTheory}
and the Nakayama relations we have
\begin{equation}
\label{Mac1}
\mbox{\rm dim}_K( \End_{KG}( \Ind_H^G( M_1 ) ) =
\left[ M_1, \Res^G_H( \Ind^G_H(M_1) ) \right]_{KH}.
\end{equation}

By Mackey's subgroup theorem, we have
\begin{equation}
\label{Mac2}
\Res^G_H( \Ind^G_H(M_1) ) \cong 
M_1 \oplus \Ind_{H \cap H^t}^H( \Res_{H \cap H^t}^H( M_1^t ) ) \oplus M_2
\end{equation}
for some $KH$-module~$M_2$.
Inserting~(\ref{Mac2}) into~(\ref{Mac1}), we find
$$
\mbox{\rm dim}_K( \End_{KG}( \Ind_H^G( M_1 ) ) \geq 1 + 
\left[ M_1, \Ind_{H \cap H^t}^H( \Res_{H \cap H^t}^H( M_1^t ) ) \right]_{KH}.
$$
As~$t$ commutes with $H \cap H^t$, we have
$$\Res_{H \cap H^t}^H( M_1^t ) \cong \Res_{H \cap H^t}^H( M_1 ).$$
Using the Nakayama relations once more, this implies that
$$
\left[ M_1, \Ind_{H \cap H^t}^H( \Res_{H \cap H^t}^H( M_1^t ) ) \right]_{KH} \geq 1.
$$
Thus the dimension of the endomorphism ring of
$\Ind_H^G( M_1 )$ is at least two.
\end{prf}

The following lemma gives a sufficient condition for the existence of an 
element~$t$ as above.

\begin{lem} \label{2.3replace}
Suppose that~$H$ and $A$ are groups with $H \leq G \leq A$. Let $z \in A$ and
put $C = C_G(z)$. Suppose that $s \in G$ such that $H  \leq \langle C,s \rangle 
\leq N_A(C)$. Let $t \in G$ satisfy the following three conditions.
\begin{enumerate}
\item $t \in C_G(s)$;
\item $t \in \langle z, {z^t} \rangle$; 
\item ${C^t} \cap Cs^i = \emptyset$ for all $i \in \mathbb{Z}$ with $s^i \not\in C$.
\end{enumerate}
Then $t \in C_G( H \cap {H^t} )$.

(The assumption $t \in \langle z, {z^t} \rangle$ is in particular satisfied
if~$t$ has odd order and ${t^z} = t^{-1}$, since the latter implies
$[z,t] = (t^z)^{-1}t = t^2$. Assumption~(3) is trivially satisfied if $s = 1$.)
\end{lem}
\begin{prf}
As $t \in \langle z, {z^t} \rangle$ we have  $t \in C_G(C \cap {C^t} )$.
As $t$ commutes with $s$ we also have that~$t$ centralizes $\langle C 
\cap C^t, s \rangle$. We conclude by showing that $H \cap H^t \leq 
\langle C \cap C^t, s \rangle$. 

We have $H \leq \langle C, s \rangle = C\langle s \rangle = 
\cup_{i \in \mathbb{Z}} Cs^i$. Similarly, ${H^t} \leq
\cup_{i \in \mathbb{Z}} {C^t}s^i$. It follows that $H \cap H^t
\leq \cup_{i,j \in \mathbb{Z}} (Cs^i \cap {C^t}s^j)$. By our assumption,
$Cs^i \cap {C^t}s^j = \emptyset$, unless $s^{i-j} \in C$. 
Hence $H \cap H^t \leq  \cup_{i \in \mathbb{Z}} (Cs^i \cap {C^t}s^i) =
\langle C \cap C^t, s \rangle$.
\end{prf}

The case $s = 1$ is allowed in the above lemma. In fact this is the version of 
the lemma that we need to handle the $\mathcal{C}_1$-type subgroups. The case
$s \neq 1$ is relevant when considering subgroups of types $\mathcal{C}_{2}$,
$\mathcal{C}_{3}$ or $\mathcal{C}_{5}$.

Condition~(3) of Lemma~\ref{2.3replace} is sometimes hard to verify. The 
following variation describes some favorable circumstances allowing to
replace this verification. 

\begin{lem}
\label{2.3replace_var}
Assume the setup and the first two conditions of Lemma \ref{2.3replace}. Suppose 
also that~$s$ normalizes~$\langle z \rangle$, and that~$\langle t \rangle$
is a characteristic subgroup of $\langle z, z ^t \rangle$. Then~$t$ centralizes
$H \cap H^t$ if $N_G( \langle t \rangle ) = C_G( t )$.
\end{lem}
\begin{prf}
Since~$s$ normalizes~$\langle z \rangle$, the latter group is also normalized
by~$H$. Hence $H^t$ normalizes $\langle z^t \rangle$ and thus $H \cap H^t$
normalizes $\langle z, z^t \rangle$. It follows that $H \cap H^t \leq 
N_G( \langle t \rangle ) = C_G( t )$.
\end{prf}
%%%%%%%%%%%%%%%%%%%%%%%%%%%%%%%%
%%%%%%%%%%%%%%%%%%%%%%%%%%%%%%%%
%%
%% Hier die Lieferung von Kay vom 15.08.10 eingebaut, soweit es hierher passt.
%%

We will also need some generalizations of Lemma~\ref{2.3replace}.

\begin{lem} \label{2.3generalized2} Let~$S$ be a subgroup of~$G$. 
Let $t \in G$ satisfy the following three conditions.
\begin{enumerate}
\item $t \in L := \langle S, {S^t} \rangle$.
\item There is a subgroup $D \leq G$ with $[L,D] = 1$, and there is $s \in G$ 
normalizing~$L$ and~$D$ such that $N_G(L) = (L \circ D).\langle s \rangle$.
\item $t$ centralizes $N_{\langle L,s \rangle}(S) \cap 
N_{\langle L,s \rangle}(S)^t$.
\end{enumerate}
Then $t$ centralizes $N_G(S) \cap N_G(S)^t$. 
\end{lem}
\begin{prf}
Let $H := N_G(S)$ and $x \in H \cap {H^t}$. Then $x$ normalizes~$S$ 
and~${S^t}$, hence~$L$. As~$D$ centralizes~$L$, it centralizes~$S$ 
and~${S^t}$ and thus $D \leq H \cap {H^t} \leq N_G(L) = 
(L \circ D).\langle s \rangle$.
Thus $H \cap {H^t} \leq \langle N_{\langle L,s \rangle}(S) \cap 
N_{\langle L,s \rangle}(S)^t, D \rangle$. Also,~$t$ centralizes~$D$, 
as $t \in L$. The claim follows as~$t$ centralizes
$N_{\langle L,s \rangle}(S) \cap N_{\langle L,s \rangle}(S)^t$ by assumption.
\end{prf}

\begin{lem} \label{lift}
Let~$Z$ be a subgroup of $Z(G)$ and $t \in G$ of order coprime to $|Z|$.
Write $\,\bar{}: G \rightarrow G/Z(G)$ for the canonical epimorphism.
Then for all $x \in G$ such that $\bar{t} \in
C_{\bar{G}}(\bar{x})$ we have $[x,t] = 1$.
\end{lem}
\begin{prf}
This follows from \cite[8.2.2.(b)]{KurzStell}.
\end{prf}

This lemma has the following simple, but important, corollary.

\begin{cor} \label{spinlift}
Let $Z$ be a subgroup of $Z(G)$ and let $H \leq G$ with $Z \leq H$.
Suppose that $t \in G$ is such that $\bar{t} \in \bar{G} = G/Z$ has
order coprime to~$|Z|$. Then $t \in C_G( H \cap H^t )$ if
$\bar{t} \in C_{\bar{G}}( \bar{H} \cap {\bar{H}^{\bar{t}}} )$.
\end{cor}
\begin{prf}
As $\bar{t}$ has order coprime to~$|Z|$, there is an element~$t_1 \in G$
with $\bar{t}_1 = \bar{t}$ and $|t_1| = |\bar{t}|$. As $t_1 = tz$ for some
$z \in Z$, the claim follows from Lemma~\ref{lift}.
\end{prf}

\section{Some results on linear groups of small degree}

We collect a couple of results, needed later on, on linear groups of degrees 
two and three.

\begin{lem}
\label{ToriInSLSU}
Suppose that~$q$ is an odd prime power.

{\rm (a)} Let $G \in \{ \GL_2(q), \GU_2(q) \}$ and let $s \in G$
be a non-central involution. Then there are maximal tori $T_1$ and $T_2$ of
$\{ A \in G \mid \mbox{\rm det}(A) = 1 \} \cong \SL_2( q )$ with $|T_1| = 
q - 1$ and $|T_2| = q + 1$ such that~$s$ inverts $T_1$ and $T_2$.

{\rm (b)} Let $G = \PGL_2( q )$ and let $s \in G$ be an involution.
Then there are maximal tori $T_1$ and $T_2$ of $\PSL_2( q )$ with $|T_1| = 
(q - 1)/2$ and $|T_2| = (q + 1)/2$ such that~$s$ inverts $T_1$ and $T_2$.
\end{lem}
\begin{prf}
(a) Suppose first that~$G = \GL_2( q )$. As~$s$ is conjugate in~$G$ to the 
element
$$s_0 := \left(
\begin{array}{cc}
  0 & 1 \\
  1 & 0 \\
\end{array}
\right),$$
we may and will assume that $s = s_0$. Then
$$T_1 := \left\{ 
\left( \begin{array}{cc} x & 0 \\ 0 & x^{-1} \end{array} \right) \mid 
x \in \mathbb{F}_q^* \right\}$$ 
is a torus of $\SL_2( q )$ of order $q - 1$, inverted by~$s$.

Next, let $x \in \mathbb{F}_{q^2}^*$ be an element of order $q + 1$.
The minimal polynomial of~$x$ over~$\mathbb{F}_q$ has constant term equal  
to~$1$, as~$x$ and $x^q = x^{-1}$ are its zeroes. Define $a \in \mathbb{F}_q$ 
by $x^2 - ax + 1 = 0$. Then the matrix
\begin{equation}
\label{Z0Definition}
z_0 := \left(
\begin{array}{cc}
  0 & 1 \\
  -1 & a \\
\end{array}
\right)
\end{equation}
represents right multiplication with~$x$ on $\mathbb{F}_{q^2}$ with respect
to the basis $1, x$ of $\mathbb{F}_{q^2}$ over~$\mathbb{F}_q$. Thus 
$T_2 := \langle z_0 \rangle$ is a torus of~$\SL_2(q)$ of order $q + 1$, which 
is obviously inverted by~$s$.

Now suppose that $G = \GU_2( q ) \leq \GL_2( q^2 )$. We let the form 
defining~$G$ have Gram matrix~$I_2$. Then $s_0 \in G$. Moreover,
every involution of $G \setminus \SU_2( q )$ is conjugate to~$s_0$ in~$G$
(since any two elements of~$G$ are conjugate in~$G$ if and only if they are
conjugate in $\GL_2( q^2 )$).

The intersection of~$T_1$ (defined with respect to $\GL_2( q^2 )$) with 
$\SU_2( q )$ has order $q + 1$ and is inverted by~$s_0$. The general theory 
of twisting tori (see \cite[Section~$3.3$]{cart}) shows that there is $g \in 
\GL_2( \bar{\mathbb{F}}_{q^2} )$, where $\bar{\mathbb{F}}_{q^2}$ denotes an 
algebraic closure of~$\mathbb{F}_{q^2}$, such $T_1^g$ ($T_1$ defined with 
respect to $\GL_2( q )$) lies in~$G$, as well as $s_0^g$. As~$s^g_0$ is 
conjugate to~$s_0$ in~$G$,
the result follows.

(b) It is well known that $\PGL_2( q )$ has exactly two conjugacy classes 
of involutions. If $s \not\in \PSL_2( q )$, the claim follows from~(a).
If $s \in \PSL_2( q )$, the assertion follows from the fact that $\PSL_2( q )$
has dihedral subgroups of orders $q - 1$ and $q + 1$ (see 
\cite[S\"atze~II.$8.3$,~$8.4$]{HuppertI}).
\end{prf}

\begin{lem}
\label{PGL2}
Let $L = \PGL_2( q )$ or $L = \SL_2( q )$ with $q > 3$ and let~$S$ be a maximal 
torus of~$L$. If $q = 5$, assume that $|S| = 6$. Then there is $t \in L$ with 
$L = \langle S, {S^t} \rangle$ such 
that~$t$ centralizes $N_L( S ) \cap N_L( S )^t$. If~$q$ is odd, there is such 
a~$t$ of odd order. If also $L = \PGL_2( q )$, then $L' = \langle S \cap L', 
(S \cap L')^t \rangle$ with $L' := [L,L] \cong \PSL_2( q )$. 
\end{lem}
\begin{prf}
Suppose first that~$q$ is odd. To begin with, let $L = \PGL_2( q )$. Then~$S$ 
contains a unique involution~$z$ and $N_L( S ) = C_L ( z )$. By 
Lemma~\ref{ToriInSLSU}(b) there is an element $1 \neq t \in L$ of odd order 
coprime to~$q$ such that $z^{-1} t z = t^{-1}$, since one of $q - 1$ or 
$q + 1$ is not a $2$-power. Then $\langle z, {z^t} \rangle$ is a dihedral 
group of order $2|t|$ containing~$t$. It follows from Lemma~\ref{2.3replace}
that~$t$ centralizes $C_L( z ) \cap C_L( z )^t = N_L( S ) \cap 
N_L( S )^t$. 

Now write $L'$ for the subgroup of~$L$ of index~$2$, i.e., $L' = [L,L] \cong 
\PSL_2( q )$. Then $t \in L'$, as~$t$ has odd order. Put $S' := S \cap L'$. 
Then $|S'| = |S|/2$ and $N_L( S' ) = N_L( S )$ since $S = C_L( S' )$. Thus~$t$ 
centralizes $N_{L'}( S' ) \cap N_{L'}( S' )^t$. As $t$ does not commute with 
$z \in S$, it follows that $t \not\in S$, and thus $t \not\in N_{L'}( S' )$. 
If $q \geq 17$, Dickson's list of subgroups of~$L'$ 
(see \cite[II, Hauptsatz 8.27]{HuppertI}) now 
implies that $\langle S', {(S')}^t \rangle = L'$. Indeed, the assumption on~$q$
implies that~$S'$ cannot be isomorphic to a subgroup of~$A_4$, $S_4$ or~$A_5$,
nor to a subgroup of $\PGL_2(q_0)$ for a proper divisor~$q_0$ of~$q$. Any
two distinct conjugates of~$S'$ have trivial intersection (see \cite[S\"atze
II.$8.3$,~$8.4$]{HuppertI}), and thus $\langle S', {(S')}^t \rangle$ has at 
least $(q - 1)^2/4$ elements. In particular, $\langle S', {(S')}^t \rangle$ 
cannot be a cyclic nor a dihedral subgroup of~$L'$. Finally suppose that 
$\langle S', {(S')}^t \rangle$ is contained in a Borel subgroup~$B'$ of~$L'$.
Then $B' = \langle S', {(S')}^t \rangle$. Now~$B'$ is the stabilizer of a
point in the projective space of~$\mathbb{F}_q^2$, and~$S'$ is the subgroup
of~$B'$ fixing a further point. As~$S'$ and~$z$ are contained
in the same torus~$S$ of~$L$, the pair of points of this projective space
fixed by~$S'$ is also fixed by~$z$. Hence $B := \langle B', z \rangle$ is a 
Borel subgroup of~$L$. For the same 
reason we have $B = \langle B', z^t \rangle$. Thus $\langle z, z^t \rangle
\leq B$, a contradiction, since~$B$ does not contain any dihedral subgroup of 
order prime to~$q$. 
For $q \leq 13$, a computation
with GAP shows that there is $t \in L'$ of odd order with
$\langle S', {(S')}^t \rangle = L'$. Since $S \not\leq L'$, we 
obtain $\langle S, {S^t} \rangle = L$, proving all our assertions for~$q$ odd
and $L = \PGL_2( q )$.

Now let $L = \SL_2( q )$ (still assuming that~$q$ is odd), and write 
$\,\bar{} : \SL_2( q ) \rightarrow \PSL_2( q )$ for the canonical epimorphism.
By what we have already proved, there is an element $\bar{t} \in \bar{L}$ of
odd order such that $\bar{L} = \langle \bar{S}, \bar{S}^{\bar{t}} \rangle$
and such that~$\bar{t}$ centralizes $N_{\bar{L}}( \bar{S} ) \cap 
N_{\bar{L}}( \bar{S} )^{\bar{t}}$. Let $t \in L$ be a preimage of~$\bar{t}$ 
of odd order. As $L$ is a non-split central extension of $\bar{L}$, we have
$L = \langle S, S^t \rangle$. Moreover,~$t$ centralizes $N_L( S ) \cap 
N_L( S )^t$ by Corollary~\ref{spinlift}.

Now suppose that $q$ is even. Then $\PGL_2( q ) \cong \SL_2( q )$ and 
$N_L( S )$ is a dihedral group of twice odd order. Let~$z \in N_L( S )$ be an 
involution. Choose an
involution $t \in L$ centralizing~$z$ and not equal to~$z$. (This is 
possible since $q \neq 2$.) Then $t \not\in N_L( S )$ since no two
distinct involutions of $N_L( S )$ commute. As~$S$ is a TI subgroup of
$L$, it follows that $N_L( S ) \cap N_L( S )^t = \langle z \rangle$,
and thus $t$ commutes with this intersection. 

In order to show that $\langle S, {S^t} \rangle = L$, we distinguish
two cases. If~$S$ is a Coxeter torus of~$L$, then $\langle S, {S^t} 
\rangle = L$ since any two disjoint Coxeter tori of~$L$ generate~$L$.  
If~$S$ fixes a line
in the natural representation of $L$ on $\mathbb{F}_q^2$, i.e.,~$S$
lies in some Borel subgroup~$B$ of~$L$, then $\langle S, {S^t} \rangle 
= L$. Otherwise, $\langle S, {S^t} \rangle =B$. But~$z$ normalizes
$\langle S, {S^t} \rangle = B$, hence $z \in B$. This is impossible
as $N_B(S) = S$.
\end{prf}

\begin{lem} \label{gl3} Let $G = \GL_3(q)$, and $S \leq G$ a natural 
$\SL_2(q)$-subgroup. Then there exists a $t \in L := \SL_3( q ) \leq G$ such 
that~$t$ centralizes $N_G( S ) \cap N_G( S )^t$ and such that 
$L = \langle S, S^t \rangle$.
\end{lem}
\begin{prf}
Put $N := N_G( S )$. Then~$N$ fixes a unique $2$-space~$U$ of the natural 
module~$V$ of~$G$ and a unique $1$-space $\langle v \rangle$ lying outside 
of~$U$.  Conversely,~$N$ is the stabilizer of the pair of subspaces~$U$ 
and~$\langle v \rangle$. 

Now let $e_1,e_2,e_3$ be a basis of~$V$ and put $U_1 := 
\langle e_1, e_2 \rangle$, $U_2 := \langle e_2, e_3 \rangle$, 
$v_1 = e_3$,  $v_2 = e_1 + e_3$. Assume that~$N$ is the stabilizer of~$U_1$ 
and $\langle v_1 \rangle$. If $t \in G$ is the linear transformation defined 
by $v_1t = v_2$, $e_2t = e_2$ and $e_1t = e_3 = v_1$, then $U_1t = U_2$ and 
$v_1t = v_2$. Thus~$N^t$ is the stabilizer of~$U_2$ and $\langle v_2 \rangle$.
It follows that $N \cap {N^t}$ stabilizes $U_1 \cap U_2 = \langle e_2 \rangle$, 
$\langle v_1 \rangle$ and $\langle v_2 \rangle$. Let $x \in N \cap N^t$ and 
suppose that $v_1x = av_1$ and $v_2x = bv_2$ for some $a, b \in \mathbb{F}_q$. 
Then $a = b$, since~$x$ stabilizes~$U_1$. This implies that~$t$ commutes 
with~$x$, as the vectors $e_2, v_1, v_2$ form a basis of~$V$ and $v_2t = 
v_1 + v_2$. 

Finally we observe that $S = \SL(U_1)$ and the group generated by 
$\SL(U_1)$ and $\SL(U_2) = S^t$ contains a generating set of transvections for 
$L = \SL_3(q)$, and the lemma is proved. 
\end{prf}

%%
%% Ende der Lieferung vom 15.08.10
%%
%%%%%%%%%%%%%%%%%%%%%%%%%%%%%%%%%
%%%%%%%%%%%%%%%%%%%%%%%%%%%%%%%%%

\section{Reduction modulo~$\ell$ and imprimitivity}

So far our preliminary results have been independent of the characteristic
of~$K$. We will also need the following elementary observations relating 
imprimitivity of irreducible characters in different characteristics through 
reduction modulo~$\ell$. Here, we assume that~$K$ is large enough, i.e.,
that~$K$ is a splitting field for all subgroups of~$G$.

Let~$\ell$ be a prime and let~$\chi$ be an ordinary irreducible character 
of~$G$, i.e., $\chi \in \Irr(G)$. The restriction of $\chi$ to the 
$\ell$-regular conjugacy classes of~$G$ is called the {\em reduction of~$\chi$ 
modulo~$\ell$}. If~$K$ has characteristic~$\ell$, we call an irreducible 
$KG$-module {\em liftable}, if its Brauer character (with respect to a suitable 
$\ell$-modular system) is obtained by reduction modulo~$\ell$ of some ordinary 
irreducible character of~$G$. We refer the reader to \cite[Section~I.14]{laro} 
for more general definitions and properties of reduction modulo~$\ell$ and 
liftability of modules.

\begin{lem}
\label{Modp}
Let~$H$ be a proper subgroup of~$G$.
Let $\psi \in \Irr(H)$ and $\chi \in \Irr(G)$. Furthermore, let $\varphi \in
\IBr_\ell(H)$ and $\vartheta \in \IBr_\ell(G)$. Then the following hold.

{\rm (1)} Suppose that $\varphi$ is the reduction modulo~$\ell$ of $\psi$ and
that $\Ind_{H}^{G}(\psi)$ is reducible. Then $\Ind_{H}^{G}(\varphi)$ is also
reducible.

{\rm (2)} Suppose that $\varphi$ lifts to~$\psi$ and that
$\Ind_{H}^{G}(\varphi)$ is irreducible. Then $\Ind_{H}^{G}(\varphi)$ lifts
to the imprimitive irreducible ordinary character $\Ind_{H}^{G}(\psi)$.

{\rm (3)} Suppose that $\chi = \Ind_{H}^{G}(\psi)$, i.e., $\chi$ is
imprimitive. If~$\vartheta$ is the reduction modulo~$\ell$ of~$\chi$, then
$\vartheta$ is also imprimitive.

{\rm (4)} Suppose that $\psi$ is of $\ell$-defect zero and that $\chi =
\Ind_{H}^{G}(\psi)$. Then $\chi$ is of $\ell$-defect zero and the reduction
modulo~$\ell$ of~$\chi$ is imprimitive.
\end{lem}
\begin{prf}
Parts~(1) and~(3) are clear since induction of characters commutes with
reduction modulo~$\ell$, i.e., restriction to $\ell$-regular elements. Part~(2)
is just a restatement of Part~(1) from a different point of view. Finally,
Part~(4) is a special case of Part~(3).
\end{prf}

These observations allow us to eliminate many possible block stabilizers $H$
and characters $\psi$. In particular, the linear characters of $H$ remain
irreducible in any characteristic (although they may become trivial). Moreover
all modular linear characters arise as restrictions of some ordinary linear
character. Thus, knowledge of the behavior of the ordinary linear characters
of $H$ is sufficient to determine their modular behavior.

We also need a generalization of Lemma~\ref{Modp}.
\begin{lem}
\label{ModpGeneralization}
Let $\psi \in \Irr(H)$ and $\varphi \in \IBr_\ell(H)$ such that~$\varphi$ occurs
in the reduction modulo~$\ell$ of~$\psi$. If $\chi(1) < [G\colon\!H]\varphi(1)$ 
for all irreducible constituents~$\chi$ of $\Ind_H^G( \psi )$, then 
$\Ind_H^G( \varphi )$ is reducible.
\end{lem}
\begin{prf}
Again we use the fact that induction commutes with reduction modulo~$\ell$.
If $\Ind_H^G( \varphi )$ were irreducible, it would be a constituent of the 
reduction modulo~$\ell$ of some irreducible constituent $\chi$ of 
$\Ind_H^G( \psi )$, which contradicts the inequality on degrees.
\end{prf}

\section{A result on polynomials}
Finally, we record a lemma on the evaluation of polynomials which will be used
later on.

\begin{lem}
\label{LemmaOnPolynomials}
Let $f = \sum_{i=0}^n a_iX^i \in \mathbb{R}[X]$ be a non-zero polynomial of
degree~$n$ (i.e., $a_n \neq 0$). Put $B := \text{\rm max}\{ |a_i|/|a_n|
\mid 0 \leq i \leq n \}$ and let $b \in \mathbb{R}$ with $b \geq B + 1$.

Then $f( b ) \neq 0$. Moreover, $f( b ) < 0$ if and only if $a_n < 0$.
\end{lem}
\begin{prf}
We may assume that $|a_n| = 1$ and thus 
$$B = \text{\rm max}\{ |a_i| \mid 0 \leq i \leq n \}.$$ 
Then
$$\left| \sum_{i=0}^{n-1} a_i b^i \right| \leq \sum_{i=0}^{n-1} |a_i| b^i 
  \leq  B \frac{b^n - 1}{b - 1}  \leq  b^n - 1  <  b^n.$$
This implies both assertions.
\end{prf}

\chapter{Sporadic Groups and the Tits Group}
\label{SectionSporadicGroups}

In this chapter $G$ is a quasisimple group such that $G/Z(G)$ is a sporadic 
simple group, or $G = {^2\!F}_4(2)'$ is the Tits simple group. Also,~$H$
denotes a maximal subgroup of $G$ and~$K$ an algebraically closed field of 
characteristic $\ell \geq 0$. As our arguments in this chapter are based on 
characters rather than modules or representations, we work with the 
irreducible $K$-characters of $H$ and $G$ throughout. We determine all 
imprimitive irreducible $K$-characters of $G$ and the corresponding block 
stabilizers~$H$. Of course, it suffices to consider the cases $\ell = 0$ and
$\ell \mid |G|$. Let us begin with the Tits group.

\begin{prop}
\label{TitsGroup}
The Tits group ${^2\!F}_4(2)'$ does not have any imprimitive irreducible
$K$-characters.
\end{prop}
\begin{prf}
The (modular) character tables of the Tits group and of all of its maximal
subgroups are available in GAP~\cite{GAP4}. It is thus a routine task to check 
our assertion.
\end{prf}

From now on~$G$ is such that~$G/Z(G)$ is a sporadic simple group.
we mostly apply ad hoc methods, although we use the same general 
approach described below.

Step 1: Using a list of maximal subgroups of $G$ and the ordinary character 
tables of $G$ as found in the Atlas~\cite{ATLAS} or GAP~\cite{GAP4}, we generate 
a list of all maximal subgroups whose index is smaller than the maximal degree 
of an ordinary irreducible character for $G$. Note that every possible maximal 
block stabilizer, in any characteristic, must occur on this list. Moreover we 
produce an upper bound on the degree of a $K$-character of~$H$ which may induce 
to an irreducible $K$-character of~$G$. 

Step 2: If the degrees of the modular irreducible characters of $G$ are known, 
then we refine the list generated in Step~$1$ by checking whether the index of 
a given maximal subgroup divides the degree of some irreducible $K$-character 
of~$G$. We also determine the exact degree of a possible imprimitive character.

Step 3: If the character table for the maximal subgroup $H$ is also known,
we induce the characters of $H$ to $G$ using GAP and determine if any of the
induced characters are irreducible.

Step 4: The remaining groups on our list that can not be approached using
Step~$3$ must be dealt with separately. 

More details are given in the proof of Theorem~\ref{SporadicTheorem} below,
which is the main result of this chapter.
The instances of imprimitive irreducible characters are presented in
Table~\ref{tab:clredSporadic}. We are now going to describe how to read this 
table. More explanations on particular entries are given in 
Remark~\ref{SporadicExplanations} below.

The notation in the table follows, as far as possible, the Atlas~\cite{ATLAS},
except that we have chosen a more traditional notion for the classical groups.
The first two columns describe~$G$ and~$H$, respectively. The third column, 
headed~$\psi$, indicates the irreducible character of~$H$ which induces to the
irreducible character~$\chi$ of~$G$, identified in the sixth column.
In those cases where $H$ is given as a direct product $H_1 \times H_2$,
$\psi$ is written as an outer product of the form $\psi_1 \boxtimes \psi_2$, 
with $\psi_i$ an irreducible character of $H_i$, $i = 1, 2$. Whenever the
(compound) character tables of~$H$ or $H_i$ are contained in the Atlas, we 
use the Atlas notation to describe $\psi$ or $\psi_i$, 
respectively. Similarly, the character
$\chi$ of $G$ is identified. With $\zeta_i$ we denote a linear
$K$-character of order $i$, $i = 2, 3, \ldots$ Some other characters
are denoted by their degrees, usually with subscripts to distinguish
characters of the same degrees. There should be no problem to match
the characters with existing tables. The degrees of~$\psi$ and~$\chi$ 
are given in the columns labelled by~$\psi(1)$ and~$\chi(1)$, respectively.

Given~$G$, there is exactly one block of rows for the {\em faithful} irreducible
imprimitive characters of~$G$. These blocks are separated by two horizontal 
rules. As in the Atlas, a single row can describe more than one irreducible 
character of~$H$ or~$G$. This can be detected from the ``ind''-columns, giving 
the Frobenius-Schur indicators and, possibly, the number of the characters in a 
row. The characteristic is always odd in Table~\ref{tab:clredSporadic}, so that 
these indicators can easily be determined.

The column headed ``Primes'' lists those characteristics, for which the 
example described in the respective row exists for one of the reasons 
described in Lemma~\ref{Modp}.

\begin{rem}\label{SporadicExplanations}
{\rm
$M_{11}$: The character $1_3$ of the maximal subgroup $H =
3^{2}\colon\!Q_{8}.2$ is the non-trivial linear character which has the
elements of order $8$ of~$H$ in its kernel.

$M_{12}$ and $2.M_{12}$: The character $1_3$ of the maximal subgroup
$H = M_{10}\colon\!2$ of $M_{12}$ is the non-trivial linear character which
has the elements of order $10$ of~$H$ in its kernel.

There are two conjugacy classes of subgroups of $M_{12}$ isomorphic to $M_{11}$
and two conjugacy classes of subgroups isomorphic to $M_{10}\colon\!2$. In each
case, the two conjugacy classes are swapped by the outer automorphism. Since
the irreducible imprimitive characters of $G$ induced from these subgroups
are invariant under the outer automorphism, these characters are also induced
from the conjugate subgroups.

$M_{22}$: Note that the example in characteristic~$3$
arises from the reduction modulo~$3$ of the character $\chi_{40}$ of
$3.M_{22}$ (Atlas notation). Since this reduction is not faithful anymore, we
have introduced it in an extra row of Table~\ref{tab:clredSporadic}.

$3.M_{22}$: The characters $\chi_{13}$ and $\chi_{14}$ are the inflations to
$3.(2^4\colon\!A_6)$ of the characters of degree~$3$ of $3.A_6$.

$\Co_2$: The character $\psi_3'$ is (the inflation to $2^{10}\colon\!(M_{22}\colon\!2)$ of)
the extension to $M_{22}\colon\!2$ of $\psi_3$ with the property that its value
on class $2B$ equals $-3$.
}
\end{rem}

\begin{thm} 
\label{SporadicTheorem}
{\rm (1)}        The following groups have no irreducible and imprimitive
                characters in any characteristic: $J_1$, $2.J_2$, $M_{23}$, 
                $3.J_{3}$, $\He$, $2.\Ru$, $6.\Suz$, $\Co_{3}$, $\HN$, 
                $\Ly$, $\Th$, $\Fi_{23}$, $J_4$, $2.B$, and $M$.

{\rm (2)}       All examples of irreducible imprimitive characters for the
                groups $M_{11}$, $2.M_{12}$, $12.M_{22}$, $2.\HS$, $M_{24}$,
                $3.\McL$, $3.\ON$, $\Co_{2}$, $6.\Fi_{22}$, $2.\Co_{1}$ and 
                $3.\Fi_{24}'$ 
                are known and are listed in Table~\ref{tab:clredSporadic}.
\end{thm}
{\small
\begin{table}[h]
\caption{\label{tab:clredSporadic}Sporadic group examples}
\hspace*{-2.00cm}
\begin{tabular}{|l|c|c|c|c|c|c|c|l|}  \hline
$G$ & $H$   & $\psi$ & ind & $\psi(1)$ & $\chi$ & ind & $\chi(1)$ & Primes        \\ \hline\hline
$M_{11}$    & $M_{10}$                     & $\zeta_2$              & $+$ & 1 & $\chi_{5}$ & $+$ 
            & $11$ & 0, 5, 11   \\ \cline{2-8}
            & $3^{2}\colon\!Q_{8}.2$       & $1_3$                  & $+$ & 1 & $\chi_{10}$ & $+$ 
            & $55$ & 0, 5, 11   \\ \hline\hline
$M_{12}$    & $M_{11}$                     & $\psi_3$               & $\circ 2$ & 10 & $\chi_{13}$ & $+$ 
            & $120$ & 0, 5       \\ \cline{2-8}
            & $M_{10}\colon\!2$            & $1_3$                  & $+$ & 1 & $\chi_{11}$ & $+$
            & $66$  & 0, 5, 11   \\ \hline
$2.M_{12}$  & $2 \times M_{11}$            & $\zeta_2 \boxtimes \psi_1$ & $+$ & 1 & $\chi_{18}$ & $+$
            & $12$ & 0, 5, 11  \\ \cline{3-8}
            &                              & $\zeta_2 \boxtimes \psi_{2}$ & $+$ & 10 & $\chi_{24}$ & $+$
            & $120$ & 0, 5  \\ \cline{3-8}
            &                              & $\zeta_2 \boxtimes 9_1$  & $+$ & 9 & $108_1$ & $+$
            & $108$ & 11  \\ \hline\hline
$M_{22}$    & $2^4\colon\!A_6$                    & $3_1$                  & $+$ & 3 & $231_1$ & $+$
            & 231 & 3 \\ \cline{3-8}
            &                              & $3_2$                  & $+$ & 3 & $231_1$ & $+$
            & 231 & 3 \\ \hline
$3.M_{22}$  & $3.\PSL_{3}(4)$                 & $\psi_{33}$     & $\circ 2$ & 15 & $\chi_{41}$  & $\circ 2$
            & 330 & 0, 5, 7, 11        \\ \cline{3-8}
            &                              & $\psi_{34}$     & $\circ 2$ & 15 & $\chi_{41}$  & $\circ 2$
            & 330 & 0, 5, 7, 11        \\ \cline{2-8}
            & $3.(2^{4}\colon\!A_{6})$            & $\psi_{14}$     & $\circ 2$ & 3 & $\chi_{40}$ & $\circ 2$
            & 231  & 0, 5, 7, 11 \\ \cline{3-8}
            &                              & $\psi_{15}$     & $\circ 2$ & 3 & $\chi_{40}$ & $\circ 2$
            & 231  & 0, 5, 7, 11 \\ \cline{2-8}
            & $3 \times (2^{4}\colon\!S_{5})$     & $\zeta_6$       & $\circ 2$ & 1 & $\chi_{39}$ & $\circ 2$
            & 231 & 0, 5, 7, 11 \\ \hline
$4.M_{22}$  & $4_1.\PSL_{3}(4)$               & $\psi_{19}$ & $\circ 2$ & 8 & $\chi_{30}$ & $\circ 2 $
            & 176  & 0, 11     \\ \cline{3-8}
            &                              & $\psi_{20}$ & $\circ 2$ & 8 & $\chi_{30}$ & $\circ 2 $
            & 176  & 0, 11     \\ \hline
$6.M_{22}$  & $6 \times (2^{3}\colon\!\SL_{3}(2))$ & $\zeta_6$        & $\circ 2$ & 1& $\chi_{51}$ & $\circ 2$
            & 330 & 0, 5, 7, 11 \\ \hline\hline
$2.\HS$      & $2.M_{22}$                   & $\psi_{13}$            & $\circ 2 $ & 10 & $\chi_{32}$ & $-$
            & $1\,000$ & 0, 3, 5  \\ \cline{2-8}
            & $\PSU_{3}(5).(2\times2)^i$      & $\zeta_4$              & $\circ 2$ & 1 & $\chi_{26}$ & $\circ 2$
            & $176$ & 0, 3, 7, 11   \\ \hline\hline
$M_{24}$    & $2^{6}\colon\!3.S_{6}$              & $\zeta_2$              & $+$ & 1 & $\chi_{18}$ & $+$ 
            & $1\,771$ & 0, 5, 7, 11, 23 \\ \hline\hline
$\McL$      & $\PSU_{4}(3)$                   & $\psi_3$               & $+$ & 35 & $\chi_{20}$ & $+$ 
            & $9\,625$ & 0, 5, 7, 11  \\ \cline{3-8}
            &                              & $\psi_4$               & $+$ & 35 & $\chi_{20}$ & $+$ 
            & $9\,625$ & 0, 5, 7, 11  \\ \hline\hline
$3.\ON$     & $3 \times (\PSL_3(7)\colon\!2)$ & $\zeta_3 \boxtimes \zeta_2$ & $\circ 2$ & 1 & $ \chi_{45}$ & $\circ 2$
            & $122\,760$ & 0, 5, 11, 31 \\ \hline\hline
$\Co_{2}$   & $\PSU_{6}(2).2$                   & $\psi_{7}$           & $\circ 2$ & 560 & $\chi_{50}$ &$+$ 
            & $1\,288\,000$ & 0, 5, 7, 23      \\ \cline{2-8}
            & $2^{10}\colon\!(M_{22}\colon\!2)$            & $\psi_3'$             & $\circ 2$ & 45 & $\chi_{60}$ & $+$
            & $2\,095\,875$ & 0, 3, 5, 23      \\ \hline \hline
$6.\Fi_{22}$ & $6 \times \Omega_8^+(2)\colon\!S_3$ & $\zeta_6 \boxtimes 1$ & $\circ 2$ & 1 & $\chi_{167}$ & $\circ2$
            & $61\,776$ & 0, 5, 11, 13 \\ \cline{3-8}
            &                                & $\zeta_6 \boxtimes \zeta_2$ & $\circ 2$ & 1 & $\chi_{168}$ & $\circ2$
            & $61\,776$ & 0, 5, 11, 13 \\\cline{3-8}
            &                                & $\zeta_6 \boxtimes 2$ & $\circ 2$ & 2 & $\chi_{173}$ & $\circ2$
            & $123\,552$ & 0, 5, 11, 13 \\ \hline \hline

$2.\Co_{1}$  & $2 \times \Co_{2}$            & $\zeta_2 \boxtimes \psi_{12}$ & $\circ 2$ & 10\,395 & $\chi_{167}$ & $+$ 
            & $1\,021\,620\,600$ & 0, 7, 11, 13 \\ \hline \hline
$3.\Fi_{24}'$& $3 \times \Fi_{23}$           & $\zeta_3 \boxtimes \psi_{15}$ & $\circ 2$ & 837\,200 & $\chi_{175}$ & $\circ 2$
            & $256\,966\,819\,200$ & 0, 5, 7, 13, 23, 29         \\ \cline{3-8}
            &                               & $\zeta_3 \boxtimes \psi_{16}$ & $\circ 2$ & 837\,200 & $\chi_{175}$ & $\circ 2$
            & $256\,966\,819\,200$ & 0, 5, 7, 13, 23, 29         \\ \hline
\end{tabular}
\end{table}
}
\begin{prf} Our proof proceeds by completing Steps 1--4 above for all of the
sporadic simple groups and their covering groups. In this proof, by an 
{\em example} we always mean an imprimitive irreducible character of~$G$.
Characters can be ordinary or modular characters, and $\ell$ denotes the
characteristic of the underlying field~$K$, where $\ell = 0$ is allowed. If $H$
is a maximal subgroup of~$G$, a (hypothetical) character of~$H$ inducing to
an irreducible character of~$G$ is denoted by $\varphi$, if $\ell$ is not 
specified, and by $\psi$, if $\ell = 0$.

$G \in \{ M_{11}, 2.M_{12}, 2.J_{2}, M_{23}, 2.\HS, 3.\McL, \He \}$:
All the character tables are known in GAP for both $G$ and for all maximal
subgroups on the list generated in Step 2.
Here all examples can be determined using GAP.

$G = J_{1}$ or $G = 3.J_{3}$: The list of maximal subgroups for $G$ generated in 
Step~2 is empty, hence no examples exist.

$G = 12.M_{22}$: The modular character tables of $12.M_{22}$ modulo all primes are 
known as well as the ordinary character tables of all maximal subgroups of
$12.M_{22}$. It is easy to check using GAP and Lemma~\ref{Modp}, that the 
examples of imprimitive irreducible characters given in 
Table~\ref{tab:clredSporadic} for the various central quotient groups of 
$12.M_{22}$ are correct. 

There are a few candidates for imprimitive irreducible characters of~$G$ which
do not yield examples. However, the potential characters which could possibly 
induce to irreducible characters of~$G$ have small degrees, 
and these hypothetical cases are easily
ruled out by ad hoc arguments. We only give one such argument here. The group
$G = 2.M_{22}$ has two faithful irreducible characters of degree $154$ in
characteristic $\ell = 3$. Let~$\vartheta$ be one of them, and let 
$H = 2.(2^4:A_6)$ be the maximal subgroup of~$G$ of index~$77$. Suppose that 
$\varphi$ is an irreducible character of~$H$ of degree~$2$, inducing 
to~$\vartheta$. Then $Z(G)$ is not in the kernel 
of~$\varphi$. All ordinary irreducible characters of~$H$ with $Z(G)$ not 
in their kernels are faithful (see GAP). Now~$\varphi$ is a 
$\mathbb{Z}$-linear combination of such ordinary irreducible characters 
(restricted to the $3$-regular classes of~$H$). Since~$H$ does not have any 
non-trivial normal $3$-subgroup, this implies that $\varphi$ is faithful.
Since the character table of~$H$ has only quadratic irrationalities,  
this would yield an embedding of~$H$ into $\GL_2(9)$, which is absurd.

$G = M_{24}$: The ordinary character tables of the maximal subgroups of $M_{24}$ 
are available in GAP, as well as enough of the modular character tables, so
that the result can be derived with the help of Lemma~\ref{Modp}.

$G = 2.\Ru$: The only maximal subgroups that could lead to examples are 
$^{2}\!F_{4}(2)'$ and $2^{6}\colon\!\PSU_{3}(3)\colon\!2$. Moreover, all 
degrees of irreducible
$2.\Ru$ characters are known. We find that the only possibilities are for 
a~$13$ or~$26$-degree character of $^{2}\!F_{4}(2)'$ in characteristic~$0$ 
and~$7$. As~$7$ is coprime to the order of $^{2}\!F_{4}(2)'$, we only need
to check the ordinary characters using GAP. Here we find no examples.

$G = 6.\Suz$: The available character tables in GAP and Lemma~\ref{Modp} leave 
only one possibility for a potential imprimitive irreducible character of~$G$, 
a faithful character~$\vartheta$ of degree $68\,640$ in characteristic 
$\ell = 5$, with block stabilizer~$H$ of index $22\,880$, i.e.,~$H$ is the
second largest maximal subgroup of~$G$. Suppose that~$\varphi$ 
is an irreducible character of~$H$ inducing to~$\vartheta$. Then $\varphi(1) 
= 3$.  Since~$\vartheta$ is faithful,~$Z(G)$ is faithfully represented on a 
module with character~$\varphi$. But $Z(G) = Z(H) \leq H'$ (see GAP) implies
that the degree of~$\varphi$ is divisible by~$6$, a contradiction.

$G = 3.\ON$: The only subgroup with an index which is small enough is $3 \times 
(\PSL_{3}(7)\colon\!2)$. Moreover, the only character that may work is linear. 
Checking the ordinary characters, we see that exactly the two linear 
characters of $3 \times (\PSL_{3}(7)\colon\!2)$ of order~$6$ induce to 
irreducible characters of $3.\ON$. The corresponding ordinary irreducible
imprimitive characters of $3.\ON$ have degree $122\,760$. Since these two 
characters are reducible modulo~$\ell$ for $\ell \in \{ 2, 3, 7, 19\}$, and 
irreducible modulo all other primes, the results of
Table~\ref{tab:clredSporadic} follow from Lemma~\ref{Modp}.

$G = \Co_3$: The modular character tables of $G$ are known, as well as the 
ordinary character tables of all maximal subgroups of~$G$. Using this 
information, we find that only the largest maximal subgroup $\McL.2$ can 
possibly be a block stabilizer. But all modular 
character tables for $\McL.2$ are known as well, so that this case is easily 
ruled out.

$G = \Co_2$: The four possible maximal block stabilizers are $\PSU_{6}(2).2$, 
$2^{10}\colon\!M_{22}\colon\!2$, $\McL$ and $2^{1+8}\colon\!\Sp_{6}(2)$. The 
modular character tables of $\PSU_{6}(2).2$ are available in GAP, which makes it 
easy to complete the proof in this case. For 
$H \in \{2^{10}\colon\!M_{22}\colon\!2, \McL, 2^{1+8}\colon\!\Sp_{6}(2)\}$, we 
find that the only characters in any characteristic that could induce to 
irreducible characters are irreducible restrictions of characteristic~$0$ 
examples. Using Lemma~\ref{Modp} and GAP, we find the examples in these cases. 

$G = 6.\Fi_{22}$: All modular character tables for~$G$ are known (see 
\cite{HiLuII,Noeske}). Only the five largest maximal subgroups remain after Step~$1$,
and the smallest of these, $M_5$ say, of index $142\,155$, does not give any 
candidates inducing to faithful characters of~$G$. In fact the only candidates
for $M_5$ have degree~$6$, and this possibility only occurs for $\Fi_{22}$ and
$3.\Fi_{22}$ and in characteristics $\ell \neq 2$. Now $M_5$ has a normal subgroup~$N$
isomorphic to $2^{10}$, and any $\ell$-modular irreducible character of $M_5$
with $N$ not in its kernel has degree larger than~$6$ by Clifford theory.
Hence the candidates for $M_5$ have~$N$ in their kernel and are thus irreducible
$\ell$-modular characters of $M_{22}$ or $3.M_{22}$, respectively. But these groups
do not have irreducible $\ell$-modular characters of degree~$6$.

The ordinary character tables of the four largest maximal subgroups are 
also available in GAP, as well as their modular characters. This allows us
easily to deduce the result.

$G = \HN$: Here, the only subgroups with an index smaller than the largest 
irreducible character of~$G$ are~$A_{12}$ and~$2.\HS.2$. Moreover, the upper
bound for the degree is~$5$ for an $A_{12}$-character and~$3$ for a 
$2.\HS.2$-character. Thus, an example can arise only from a linear character 
of $2.\HS.2$. 
Checking the characteristic~$0$ case and using Lemma~\ref{Modp}, we find 
that there are no examples for this group.

$G = \Ly$: The only possibilities for $H$ are $G_{2}(5)$ and $3.\McL.2$. 
Moreover, the upper bound for the degree of an $H$-character is~$8$ in both 
cases. If $H=3.\McL.2$, then only linear characters are possible. Checking the 
characteristic~$0$ case and using Lemma~\ref{Modp}, we see that $H$ is not a 
block stabilizer. 

If $H=G_{2}(5)$, then a character~$\varphi$ inducing irreducibly to~$G$ must 
be the character of the minimal module for $G_{2}(5)$ in characteristic~$5$. 
Now there is a conjugate $H^t$ of $H$ such that $H^t \cap H = \PSU_{3}(3)$
(see \cite[p.~$174$]{ATLAS}). 
Moreover the restriction to this intersection of the minimal module for~$H$
is irreducible. Although $\PSU_{3}(3)$ has three $5$-modular characters of 
degree~$7$, the restriction of the minimal character of $H$ has value~$0$
on all $12$-elements of $\PSU_{3}(3)$, but no other $5$-modular character of 
degree~$7$ of $\PSU_{3}(3)$ has this property (irreducible or reducible). 
In particular, $\Res_{H^{t} \cap H}( \varphi ) = 
\Res_{H^{t} \cap H}(^{t}\!\varphi)$, which eliminates this case as well. 

$G = \Th$: The only subgroup that may provide an example is 
$^{3}\!D_{4}(2)\colon\!3$ with a linear character. Using the ordinary character 
tables and Lemma \ref{Modp}, we see that there are no examples for this group.

$G = \Fi_{23}$: According to the ``Improvements to the {\sc Atlas}'' in
\cite[Appendix~$2$]{modATL}, the list of maximal subgroups of $\Fi_{23}$
given in \cite[p.~$177$]{ATLAS} is complete, except that the last maximal
subgroup has to be replaced by $\mbox{\rm PSL}_2(23)$ (see 
\cite[p.~$304$]{modATL}). A potential block stabilizer is among the first
six maximal subgroups. We consider each of these possibilities in turn,
starting with the smallest one. In all cases, we denote by~$\varphi$ an
irreducible $K$-character inducing irreducibly to~$G$.

Suppose first that $H = 2^{11}.M_{23}$. Then $\varphi(1) \leq 2$ and, since~$H$ 
is perfect, we have $\varphi(1) = 2$. Since~$H$ does not have any irreducible
character of this degree, this subgroup does not provide any examples.

Next let $H = S_3 \times \Omega_7(3)$. Here, $\varphi(1) \leq 3$. Since the smallest
degree of a non-trivial irreducible character of~$\Omega_7(3)$ is~$7$, we must
have $\Omega_7(3)$ in the kernel of~$\varphi$, i.e., $\varphi$ is a character of~$S_3$.
In~$S_3$, every modular irreducible character lifts to an ordinary irreducible
character. By the first part of Lemma~\ref{Modp}, it suffices to check the
ordinary character tables. This does not give any example.

If $H = \Sp_8(2)$, we have $\varphi(1) \leq 7$, but a non-trivial irreducible
character of $H$ has degree at least~$8$. Thus we do not find an example 
with block stabilizer $\Sp_8(2)$.

Now let $H = 2^2.\PSU_6(2).2$. In this case $\varphi(1) \leq 10$. To show that only
characters of degree~$1$ satisfy this condition, it suffices to look at the
subgroup $2^2.\PSU_6(2)$ of~$H$. This group does not have any faithful irreducible
characters, and a non-trivial irreducible character of $2.\PSU_6(2)$ has degree 
at least~$20$. Hence $\varphi(1) = 1$. Checking the ordinary character table
and using Lemma~\ref{Modp}, we find no example for this subgroup either.

Next let $H = P\Omega_8(3)^+\colon\!S_3$. Put $F := \Fi_{24}$ and let $u, v, w$
be $3$-transpositions of~$F$ such that $uv=vu$, but $w$ does not commute with
$u$ and~$v$. Then $G = C_{F'}( u )$, and $H = C_G( uw ) = C_{F'}( u, w )$
(see \cite[(16.12), (15.14), (25.9)]{Asch3}). Put $L := C_{F'}( u, v, w )$. 
Then $L \cong \Omega_8^+(2)\colon\!S_3$ by \cite[(25.6)]{Asch3}. Put $s := uv$. 
Then $s \in \langle u, v, w \rangle$ and thus $s$ centralizes~$L$. It follows 
that $L = H \cap H^s$ since $L$ is a maximal subgroup of~$H$ (see 
\cite[p.~140]{ATLAS}). Hence $H$ is not a block stabilizer by 
Lemma~\ref{Mackey1}.

Finally, let~$H = 2.\Fi_{22}$, the largest maximal subgroup of~$G$. Let
$t \in G$ be an element of the $3A$-conjugacy class of~$G$. Then $N := 
N_G(\langle t \rangle) = \langle t, s \rangle \times L$ with $\langle t, 
s \rangle \cong S_3$ and $L \cong \Omega_7(3)$ (see \cite[p.~$177$]{ATLAS}).
Thus $s$ is an involution in class~$2A$, since no other involution
centralizer is divisible by~$|L|$.
We choose $H = C_G(s)$. Then $L \leq {H^t} \cap H$. Since $\langle s \rangle
\times L$ is a maximal subgroup of $H$ (see \cite[p.~$163$]{ATLAS}), and 
since ${s^t} \not\in H$, we have $L = {H^t} \cap H$. In particular, 
$t$~centralizes ${H^t} \cap H$. Hence $H$ is not a block stabilizer by 
Lemma~\ref{Mackey1}.

$G = 2.\Co_1$ or $G = \Co_1$: Suppose first that $\ell = 0$.
Step~$1$ leaves the seven largest maximal subgroups as candidates for block 
stabilizers, but only the largest and the fifth largest maximal subgroup 
remain after Step~$2$. For the latter subgroup the only candidates have 
degree~$1$. Now the subgroup $2^{1+8}_+.\Omega_8^+(2)$ of $\Co_1$ is perfect and 
thus cannot be a block stabilizer. Moreover, $2^{1+8}_+.\Omega_8^+(2)$ is not a 
subgroup of $2.\Co_1$, as can be shown with GAP. Hence $2.(2^{1+8}_+.\Omega_8^+(2))
\leq 2.\Co_1$  is perfect as well, and so it cannot be a block stabilizer 
either.

Next let $H = 2 \times \Co_2$ be the largest maximal subgroup of~$2.\Co_1$.
The only candidates have degrees $1\,080$, $5\,313$ and $10\,395$, and since 
$1\,080$ and $5\,313$ are not character degrees for $H$, only $10\,395$ remains. 
These candidates yield the characteristic~$0$ examples of Table 
\ref{tab:clredSporadic}. The modular character tables of $H$ and of $2.\Co_1$
are known for $\ell = 7, 11, 13$ and $23$ (see \cite{BBB} and 
\cite{ModAtlHomePage}). This allows us to conclude that the degree $10\,395$ 
characters also yield examples for $\ell = 7, 11$ and $13$ and that there are
no more examples with block stabilizer $H$ in these characteristics.

We are now ging to show that there are no more examples of irreducible
imprimitive characters. Assume from now on that $\ell \neq 0$. We consider the 
seven maximal subgroups of~$2.\Co_1$ in turn, starting with the smallest one. 

This is $H := 2.(A_4 \times G_2(4))\colon\!2)$. Here, $\varphi(1) = 1$ and thus 
$\varphi$ is liftable. But $2.\Co_1$ does not have an ordinary irreducible 
character of degree equal to the index of~$H$. Thus~$H$ does not give any 
example.

Next, let $H = 2.(\PSU_6(2)\colon\!S_3)$. Here, $\varphi(1) \leq 13$. Since the 
non-split extension $2.\PSU_6(2)$ is not a subgroup of~$2.\Co_1$ (use GAP to show
this), it follows that $H$ has a normal subgroup~$L$ of index~$12$, isomorphic 
to $\PSU_6(2)$. The smallest non-trivial character degree of~$\PSU_6(2)$ in any 
characteristic is larger than~$13$ (see \cite{ModAtlHomePage}). Hence~$\varphi$
has~$L$ in its kernel and 
is thus liftable. Moreover, $\varphi(1) \in \{ 1, 2, 3\}$. But~$2.\Co_1$ does 
not have an ordinary irreducible character of degree equal to $d \cdot 
[2.\Co_1\colon\!H]$ for $d = 1, 2, 3$, and so $H$ does not give any example.

Now let $G = \Co_1$ and $H = 2^{1+8}_+.\Omega_8^+(2)$. Then $H = C_G(x)$ for an 
element~$z$ in class $2A$ of~$G$. By computing class multiplication coefficients
we see that there is a conjugate~$y$ of~$z$ such that $zy$ lies in class $3B$.
Hence there is an element~$t$ of order $3$ such that $t \in \langle z, z^t \rangle$.
Lemma~\ref{2.3replace} and Corollary~\ref{spinlift} imply that neither~$H$ nor 
its inverse 
image in $2.G$ are block stabilizers of imprimitive irreducible representations 
of~$G$, respectively~$2.G$. 

Next let $G = 2.\Co_1$ and $H = 2 \times \Co_3$. All modular character tables of
$H$ are known. From these we conclude that the only candidates $\varphi$ which
are not liftable have degree~$22$, and occur for $\ell = 2$ and $3$. In this 
case,~$\varphi$ is a constituent of the reduction modulo~$\ell$ of one of the two
ordinary characters of $H$ of degree~$23$. We can use 
Lemma~\ref{ModpGeneralization} to show that~$H$ is not a block stabilizer.

Now let $G = 2.\Co_1$ and $H = 2.(2^{11}\colon\!M_{24})$. Here, $\varphi(1) 
\leq 123$. Using Clifford theory and the known modular character tables of 
$M_{24}$, we find that the only non-liftable Brauer characters $\varphi$ of 
degrees at most $123$ have degrees~$22$ if $\ell = 3$, and $11$ and $44$ if 
$\ell = 2$. Again, an application of Lemma~\ref{ModpGeneralization} rules out
these candidates.

Next, let $G = \Co_1$ and $H = 3.\Suz.2$. Here, $H = N_G (\langle z \rangle)$,
where~$z$ is an element of class~$3A$ of~$G$ (see \cite[p.~183]{ATLAS}). 
A computation of class multiplication coefficients shows that there is a 
conjuate~$z'$ of~$z$ such that $zz'$ is an involution in Class~$2B$. It
follows that $L := \langle z, z' \rangle$ is isomorphic to an alternating 
group~$A_4$, and $N_G( \langle z, z' \rangle ) \cong 
(A_4 \times G_2( 4 ))\colon\!2$ (again, see \cite[p.~183]{ATLAS}). Thus
$N_G(L) = (L \times D)\colon\!\langle s \rangle$ with $D \cong G_2(4)$,
and~$s$ an involution normalizing~$L$ and~$D$. Now $\langle L, s \rangle 
\cong S_4$. A computation in~$S_4$ shows that there exists an involution 
$t \in L$ such that $L = \langle z, z^t \rangle$ and such that~$t$ centralizes 
$N_{\langle L, s \rangle}( \langle z \rangle ) \cap 
N_{\langle L, s \rangle}( \langle z \rangle )^t$.
Lemma~\ref{2.3generalized2} implies that~$t$ centralizes $H \cap H^t$.

Now, take  $G = 2.\Co_1$ and $H = 6.\Suz.2$. It suffices to consider faithful 
characters of~$G$ and so we may assume that $\ell$ is odd. Ruling out the cases 
$\ell = 7, 11$ and $13$ with the known modular character tables for $G$, we may 
assume that $\ell = 3$ or~$5$. Here, $\varphi(1) \leq 660$.
The $\ell$-modular character tables and decomposition numbers for $H$ are known. 
If $\ell = 5$, there are just two non-liftable irreducible Brauer characters of
degrees at most $660$. Both of them have $Z(G)$ in their kernel, so they do not
induce to irreducible characters by what we have already proven.
Now let $\ell= 3$. Using Lemma~\ref{ModpGeneralization}, it is easy to rule out 
all but two potential examples. The remaining Brauer characters that could 
possibly induce to irreducible characters have degree~$12$. Let $\varphi \in
\IBr_3(H)$ with $\varphi(1) = 12$. Then $\bar{\varphi} \neq \varphi$ and 
$\bar{\varphi} + \varphi$ lifts to an ordinary irreducible character $\psi$ of
degree $24$. There is exactly one irreducible constituent~$\chi$ of 
$\Ind_H^G( \psi )$ with $\chi( 1 ) \geq 12[G\colon\!H]$, and $\chi$ occurs with
multiplicity~$1$ in $\Ind_H^G( \psi )$. Moreover, $\chi$ is real valued. (In 
fact, $\chi = \chi_{132}$ in Atlas notation.) Now if $\Ind_H^G( \varphi )$ were 
irreducible, so would be $\Ind_H^G( \bar{\varphi} )$, and both would be contained in 
the reduction modulo~$3$ of~$\chi$. However, $24[G\colon\!H] > \chi(1)$, so that
this possibility does not occur. 

Finally, let $G = 2.\Co_1$ and $H = 2 \times \Co_2$ if $\ell = 3$ or~$5$, and let
$G = \Co_1$ and $H = \Co_2$ if $\ell = 2$. The $\ell$-modular character tables and 
decomposition numbers for $H$ are known. Using Lemma~\ref{ModpGeneralization}
we can easily rule out $H$ as a block stabilizer. We only comment on the 
critical cases. Suppose first
that $\ell$ is odd. Then the only candidates arise from the liftable Brauer 
characters of~$H$ of degree $10\,395$. The lifts of these induce to
the irreducible ordinary character $\chi_{167}$ of $G$ of largest degree.
However, this character is neither irreducible modulo~$5$ (see 
\cite{ModAtlHomePage})
nor modulo~$3$. In each case, the restriction of $\chi_{167}$ to the $\ell$-regular
elements is a $\mathbb{Z}$-linear combination of restricted ordinary characters
of smaller degrees.
Now let $\ell = 2$. Then there exists $\varphi \in \IBr_2(H)$ with $\varphi( 1 ) 
= 748$. The complex conjugate character, $\bar{\varphi}$ is distinct from 
$\varphi$, and $\varphi$ and $\bar{\varphi}$ are both constituents of the
reduction modulo~$2$ of $\psi \in \Irr(H)$ with $\psi(1) = 1\,771$. There is 
exactly one irreducible constituent~$\chi$ of $\Ind_H^G( \psi )$ with $\chi( 1 )
 \geq 748 \cdot [G\colon\!H]$, and $\chi$ occurs with multiplicity~$1$ in 
$\Ind_H^G( \psi )$. Moreover, $\chi$ is real valued. (In fact, $\chi = 
\chi_{59}$ in Atlas notation.) Now if $\Ind_H^G( \varphi )$ were irreducible, 
so would be $\Ind_H^G( \bar{\varphi} )$, and both would be contained in the 
reduction modulo~$2$ of~$\chi$. However, $2 \cdot 748 \cdot [G\colon\!H] > 
\chi(1)$, so that this possibility does not occur.

This completes our proof for $2.\Co_1$.

$G = J_{4}$: Only the largest maximal subgroup $H = 2^{11}\colon\!M_{24}$ 
of $J_4$ arises in Step~$1$. 
Suppose first that~$\ell$ is odd or $\ell = 0$. Then the two non-trivial orbits 
of~$M_{24}$ on the set of irreducible characters (ordinary or $\ell$-modular) 
of the normal subgroup $2^{11}$ of~$H$ have lengths $759$ and $1\,288$, 
respectively (see the ordinary character table of~$H$ in~\cite{GAP4}). 
By~\cite{modATL}, a non-trivial character of~$M_{24}$ has dimension at 
least~$22$. But $22 \cdot [J_4\colon\!H]$ is larger than the largest degree 
of an ordinary irreducible character of~$J_4$.  Suppose then that $\ell = 2$. 
The $2$-modular character table of~$H$ is the same as the $2$-modular character 
table of~$M_{24}$. The only characters of~$H$ that can possibly induce to an
irreducible character of~$J_4$ are the two characters of degree~$11$. They 
arise as constitutents of the ordinary irreducible character~$\psi$ of~$M_{24}$
of degree~$23$. Inducing~$\psi$ from~$H$ to~$J_4$ using~\cite{GAP4} we
find 
$$\Ind_H^{J_4}(\psi) = \chi_{29} + \chi_{30} + \chi_{45} + \chi_{51}.$$
Each of these characters has degree less than $11 \cdot [J_4\colon\!H]$. 
We are done with Lemma~\ref{ModpGeneralization}.

$G = 3.\Fi_{24}'$: Step~$1$ leaves the six largest maximal subgroups as candidates
for block stabilizers. If $\ell = 0$, however, only $3 \times \Fi_{23}$ and 
$3.((3 \times P\Omega_8^+(3)\colon\!3)\colon\!2)$ remain after Step~$2$. We now
consider these six maximal subgroups in turn, starting with the largest one.

If $G = 3.\Fi_{24}'$ and $H = 3.(3^{1+10}\colon\!\PSU_5(2)\colon\!2)$, then 
$\varphi(1) = 1$. Hence,~$\varphi$ is liftable. By Lemma~\ref{Modp} there are 
no such examples.

Next, let $G = 3.\Fi_{24}'$ and $H = 3.(3^7.\Omega_7(3))$. Here, $\varphi(1) \leq 3$.
The smallest degree of a non-trivial projective character of $\Omega_7(3)$ is, in any
characteristic, larger than~$3$. It follows from Clifford theory, applied to
$O_3(H)$, that~$\varphi$ is liftable. Again, Lemma~\ref{Modp} rules out this
possibility.

In case $G = 3.\Fi_{24}'$ and $H = 3 \times \Omega_{10}^+(2)$, we must have 
$\varphi(1) \leq 8$, and a similar argument as above disposes of this 
possibility.

Next, let $G = \Fi_{24}'$ and $H = (3 \times P\Omega_8^+(3)\colon\!3)\colon\!2$.
In this case $H = N_G( \langle z \rangle )$ for an element~$z$ in the conjugacy 
class $3A$ of~$G$ (see \cite[p.~$207$]{ATLAS}). A computation of class 
multiplication coefficients shows that there is a conjuate~$z'$ of~$z$ such that 
$zz'$ is an involution in Class~$2A$. It follows that $L := 
\langle z, z' \rangle$ is isomorphic to an alternating group~$A_4$, and 
$N_G( \langle z, z' \rangle ) \cong 
(A_4 \times \Omega^+_8( 2 )\colon\!3)\colon\!2$ (again, see \cite[p.~207]{ATLAS}). 
Thus $N_G(L) = (L \times D)\colon\!\langle s \rangle$ with $D \cong 
\Omega^+_8( 2 )\colon\!3$, and~$s$ an involution normalizing~$L$ and~$D$. Now 
$\langle L, s \rangle \cong S_4$. A computation in~$S_4$ shows that there exists 
an involution $t \in L$ such that $L = \langle z, z^t \rangle$ and such that~$t$ 
centralizes $N_{\langle L, s \rangle}( \langle z \rangle ) \cap 
N_{\langle L, s \rangle}( \langle z \rangle )^t$.
Lemma~\ref{2.3generalized2} implies that~$t$ centralizes $H \cap H^t$ and
thus~$H$ is not a block stabilizer in $G$. Finally, Corollary~\ref{spinlift} 
excludes $3.H$ as a block stabilizer in $3.G$.

Next, let $G = \Fi_{24}'$ and $H = 2.\Fi_{22}.2$. In this case, $H = C_G( z )$ for
a $2A$-element~$z$. A class multiplication coefficient computation shows the 
existence of a conjugate~$y$ of~$z$ such that $zy$ lies in class $5A$. Hence 
there is an element~$t$ of order $5$ such that $t \in \langle z, z^t \rangle$.
Lemma~\ref{2.3replace} and Corollary~\ref{spinlift} imply that neither~$H$ nor 
its inverse
image in $3.G$ are block stabilizers of imprimitive irreducible representations
of~$G$, respectively~$3.G$. 

Let us finally consider the largest maximal subgroup. Suppose first that
$G = \Fi_{24}'$ and $H = \Fi_{23}$. Now $H = C_G( z )$ for an 
element $z \in \Fi_{24}$ lying in conjugacy class $2C$ 
(see \cite[p.~$207$]{ATLAS}). There is an element $z'$ conjugate to~$z$
such that $t := zz'$ has order~$3$. Hence $t \in \langle z, z^t \rangle$,
and there are no irreducible characters of~$G$ induced from~$H$ by 
Lemma~\ref{2.3replace}. Suppose then that $G = 3.\Fi_{24}'$ and $H = 3 \times 
\Fi_{23}$. Using GAP, we find exactly the characteristic~$0$ examples
given in Table~\ref{tab:clredSporadic}. The character $\zeta_3 \boxtimes
\chi_{15}$ of $H$ is of degree $837\,200$ and is thus of defect~$0$ for
$\ell = 5, 7, 13, 23, 29$. This yields the examples in these characteristics by
Lemma~\ref{Modp}. If $\ell = 11$ or $\ell = 17$, then $\zeta_3 \boxtimes 
\chi_{15}$ is irreducible modulo~$\ell$ (\cite[p.~$298$ and~$302$]{BBB}), but 
$\chi_{241}$ is not (\cite[p.~$366$ and~$378$]{BBB}). The reduction of $\zeta_3 
\boxtimes \chi_{15}$ modulo~$2$ is reducible (see \cite{HNN}), and modulo~$3$ 
it is not faithful. Suppose that $\ell \neq 2, 3$. All $\ell$-modular character 
tables of~$H$ are known. By applying Lemma~\ref{ModpGeneralization}, we find no 
other examples if~$\ell$ is odd. Since~$G$ has no faithful irreducible 
characters for $\ell = 3$, we may now asume that $\ell = 2$. In this case, the 
modular character table of $H$ is also known (see \cite{HNN}). Applying 
Lemma~\ref{ModpGeneralization} once more, only the irreducible Brauer characters 
of degrees $1\,494$ and $94\,588$ could possibly induce to irreducible 
characters of~$G$. These Brauer characters are faithful, lying in two dual 
blocks of~$H$. Each of these blocks contains one irreducible character of degree 
$1\,494$ and two irreducible characters of degree $94\,588$. To rule out these 
cases we need a refinement of the argument of Lemma~\ref{ModpGeneralization}. 
Let $\varphi_1$ denote one of the irreducible Brauer characters of $H$ of 
degree $1\,494$, and let $\varphi_2$ be the Brauer character of $H$ obtained by 
restricting the ordinary irreducible character of degree $274\,482$ (in the 
same block as~$\varphi_1$) to the $2$-regular elements. Then $\varphi_2$ 
contains each of the two irreducible Brauer characters of degree $94\,588$ 
(in the same block as~$\varphi_1$) as constituents (see \cite{HNN}). Thus 
$\varphi_2 = 94\,588_1 + 94\,588_2 + \varphi_3$ for some Brauer 
character~$\varphi_3$. Put $\Phi_i := \Ind_H^G( \varphi_i )$, $i = 1, 2$. We 
now express $\Phi_1$ and $\Phi_2$ in terms of suitable ordinary characters, 
restricted to the $2$-regular elements of~$G$. The result is given in 
Table~\ref{Fi24Mod2RuleOut1}. The first column of this table gives the degrees 
of the restricted ordinary characters, the second and third column give the 
multiplicities of these restrictions in $\Phi_1$ and $\Phi_2$, respectively.
\begin{table}[h]
\caption{\label{Fi24Mod2RuleOut1}Decomposition of some induced Brauer characters}
$$
\begin{array}{r|r|r} \hline
\multicolumn{1}{c|}{\mbox{\rm Degree}} & \multicolumn{1}{c|}{ \Phi_1 }
& \multicolumn{1}{c}{ \Phi_2 }
\rule[- 7pt]{0pt}{ 20pt} \\ \hline
          783 & -1 & -5 \rule[ 0pt]{0pt}{ 13pt}\\
      64\,584 & -1 &  0 \\
     306\,153 &  1 &  6 \\
  6\,724\,809 &  1 &  6 \\
 19\,034\,730 & -1 & -7 \\
 25\,356\,672 &  0 &  1 \\
 43\,779\,879 & -2 & -8 \\
195\,019\,461 &  0 &  7 \rule[- 7pt]{0pt}{ 5pt} \\ \hline
\end{array}\quad\quad
\begin{array}{r|r|r} \hline
\multicolumn{1}{c|}{\mbox{\rm Degree}} & \multicolumn{1}{c|}{ \Phi_1 }
& \multicolumn{1}{c}{ \Phi_2 }
\rule[- 7pt]{0pt}{ 20pt} \\ \hline
    203\,843\,871 & -1 & -6 \rule[ 0pt]{0pt}{ 13pt} \\
    216\,154\,575 &  1 &  4 \\
    216\,154\,575 &  1 &  2 \\
    330\,032\,934 &  1 &  8 \\
 1\,050\,717\,096 &  0 &  2 \\
 2\,801\,912\,256 &  0 &  6 \\
 8\,993\,909\,925 &  0 &  2 \\
21\,842\,179\,632 &  0 &  2
\rule[- 7pt]{0pt}{ 5pt} \\ \hline
\end{array}
$$
\end{table}
Thus the degree of any modular constituent of $\Phi_1$ is at most equal to 
$330\,032\,934 < 1\,494 \times 306\,936$, and the degree of any modular 
constituent of $\Phi_2$ is at most equal to
$21\,842\,179\,632 < 94\,588 \times 306\,936$. Since the second
of these factors is the index $[G\colon\!H]$, the characters
$\Ind_H^G( 1\,494 )$ and $\Ind_H^G( 94\,588_i )$, $i = 1, 2$ are reducible.
This completes the proof for $3.\Fi_{24}'$.

$G = 2.B$: There are only four maximal subgroups of~$2.B$ in the list generated 
in Step~$1$, namely $2^2.{^2\!E}_{6}(2)\colon\!2$, $2.(2^{1+22}\Co_{2})$, 
$2 \times \Fi_{23}$ and $2.(2^{9+16}\colon\!\Sp_{8}(2))$. If 
$H=2^2.{^2\!E}_{6}(2)\colon\!2$ or $2.(2^{1+22}\Co_{2})$, then the image of $H$ 
in~$B$ is the centralizer of a $2A$ or $2B$ element $z \in B$, respectively (see 
\cite[p.~$217$]{ATLAS}). In each case, there is a conjugate ${z^t}$ such 
that $\langle z, {z^t} \rangle \cong S_{3}$. (This is easily checked using
structure constants.) We may thus assume that~$t$ is
a $3$-element in $\langle z, {z^t} \rangle$. Moreover, as $t$ is 
a $3$-element, there is a $3$-element in the preimage of $t$ in $2.B$. Thus,~$H$ 
can not be a block stabilizer by Lemma~\ref{2.3replace} and 
Corollary~\ref{spinlift}.

If $H = 2 \times \Fi_{23}$, then $\varphi(1) \leq 29$. Therefore $\varphi$ must
be linear. Checking the characteristic~$0$ case, we see that the corresponding 
linear character induces to a reducible character. 

Finally suppose that $H = 2.(2^{9+16}.\Sp_{8}(2))$. Then $\varphi(1) < 7$ if
$Z(G)$ is in the kernel of~$\varphi$, and $\varphi( 1 ) < 12$, otherwise. 
The ordinary character table of the image~$\bar{H}$ of~$H$ in~$B$ 
is available in GAP. (It has been computed by Eamonn O'Brien.) We find
that $\bar{H}$ is perfect. Moreover, GAP shows that $\bar{H}$ is not a 
subgroup of $2.B$. Thus~$H$ is perfect as well. Suppose first that
$\ell = 2$. In this case every irreducible $2$-modular character of $H$ is 
a $2$-modular character of~$\Sp_8(2)$. However, the smallest degree of a 
non-trivial $2$-modular character of $\Sp_8(2)$ equals~$8$. Thus there is 
no imprimitive irreducible $2$-modular character of~$G$, since this would
have $Z(G)$ in its kernel. Suppose then that $\ell > 2$ or $\ell = 0$ and let 
$\varphi$ be an irreducible character of~$H$ of degree smaller than~$12$. 
Put $N := O_2(H)$ and let $\lambda$ be an irreducible constituent of
$\Res_N^{H}( \varphi )$. If~$\lambda$ were not invariant in $H$, the degree
of~$\varphi$ would be at least equal to the index of a proper subgroup of $\Sp_8(2)$, 
hence at least equal to~$120$. Thus~$\lambda$ is invariant in~$H$. Since the Schur 
multiplier of $\Sp_8(2)$ is trivial, and since a non-trivial character of $\Sp_8(2)$ 
has degree at least $35$, we conclude that~$\varphi$ is an extension of~$\lambda$.
Suppose now that $\ell \neq 0$. 
The restriction of~$\varphi$ to~$N$ equals~$\lambda$, which
lifts to characteristic~$0$, since $N$ is an $\ell'$-group. Thus $\varphi$ also
lifts to an ordinary character~$\psi$, say, namely to the extension to $H$ 
of the lift of $\Res^{H}_N( \varphi )$. But no ordinary character degree of
$2.B$ is divisible by the index of~$H$ in~$G$. In particular, 
$\Ind_{H}^{G}( \psi )$ is reducible. Then $\Ind_{H}^{G}( \varphi )$
is reducible as well by Lemma~\ref{Modp}(1). This completes the proof for
$2.B$.

$G = M$: It follows from the known information on the maximal subgroups of~$M$ 
(see \cite{WWWW}), that the only maximal subgroup which is large enough is $2.B$. 
Now $H = 2.B$ is the centralizer of a $2A$-element $z \in M$
(see \cite[p.~$234$]{ATLAS}). Using structure constants one checks that
there is a conjugate $z'$ of $z$ such that $t := zz'$ has order~$3$.
Hence $t \in \langle z, z^t \rangle$ and thus, by Lemma~\ref{2.3replace},
the subgroup $2.B$ does not give rise to imprimitive irreducible characters
of~$M$.
\end{prf}

\chapter{Alternating Groups}\label{alternating}

In this chapter $G$ is an alternating group or the double cover of
an alternating group (or possibly the triple or sixfold cover of
$A_{6}$ or $A_{7}$). The results presented here are known,
but included for the sake of completeness. Again, $K$ denotes an 
algebraically closed field of characteristic~$\ell \geq 0$, and the results
are formulated in terms of $K$-characters.

We first consider the groups $A_n$ for small values of~$n$. The examples are
given in Table~\ref{tab:clredSmallAlternating}, where we use similar conventions 
as for Table~\ref{tab:clredSporadic}.
\begin{thm}
\label{SmallAlternating}
Let $G$ be a covering group of $A_n$ with $5 \leq n \leq 9$. Then the irreducible
imprimitive $K$-characters of~$G$ are exactly those described in 
Table~\ref{tab:clredSmallAlternating}.
\end{thm}
\begin{prf}
The result is easily obtained with the help of GAP. These groups have also been 
treated in \cite{DM2} and \cite{NeNoe} in case $\text{\rm char}(K) = 0$.
\end{prf}

{\small
\begin{table}[h]
\caption{\label{tab:clredSmallAlternating}Small alternating group examples}
\hspace*{-2.00cm}
$$
\begin{array}{|l|c|c|c|c|c|c|c|l|}  \hline
G & H   & \psi & \text{\rm ind} & \psi(1) & \chi & \text{\rm ind} & \chi(1) 
& \text{\rm Primes} \\ \hline\hline

A_5 & A_4 &
\begin{array}{c} \zeta_3 \\ \bar{\zeta}_3 \end{array} & \circ & 1 &
\chi_{5} & + &  5 & 0, 5 \\ \hline

2.A_5 & 5\colon\!4 &
\begin{array}{c} \zeta_4 \\ \bar{\zeta}_4 \end{array} & \circ & 1 &
\chi_{9} & - &  6 & 0, 3 \\ \hline\hline

A_6 & 3^2\colon\!4 &
\begin{array}{c} \zeta_4 \\ \bar{\zeta}_4 \end{array} & \circ & 1 &
\chi_{7} & + & 10 & 0, 5  \\ \hline

2.A_6 & 3^2\colon\!8 &
\begin{array}{c} \zeta_8 \\ \bar{\zeta}_8 \end{array} & \circ & 1 &
\chi_{12} & - & 10 & 0, 5 \\ \cline{3-9}

& & \begin{array}{c} \zeta_8^{3} \\ \bar{\zeta}_8^{3} \end{array} & 
\circ & 1 & \chi_{13} & - & 10 & 0, 5 \\ \hline

3.A_6 & 3 \times A_5 &
\zeta_3 \boxtimes 1 & \circ2 & 1 & \chi_{16} & \circ2 & 6 & 0, 5 \\ \hline

6.A_6 & 3 \times 2.A_5 &
\zeta_3 \boxtimes 2_1 & \circ2 & 2 & \chi_{21} & \circ2 & 12 & 0  \\
\cline{3-9}
& &
\zeta_3 \boxtimes 2_2 & \circ2 & 2 & \chi_{22} & \circ2 & 12 & 0 \\
\hline\hline

3.A_7 & 3.A_6 & \begin{array}{c} \psi_{14} \\ \psi_{15} \end{array} &
\circ2 & 3 & \chi_{21} & \circ2 & 21 & 0, 7 \\ \cline{2-9}
& 3 \times \SL_3(2) & \zeta_3 \boxtimes 1 & \circ2 & 1 & \chi_{19} & 
\circ2 & 15 & 0, 5 \\ \cline{2-9}
& 3 \times S_5 & \zeta_3 \boxtimes \zeta_2 & \circ2 & 1 & \chi_{20} &
\circ2 & 21 & 0, 5, 7 \\ \hline\hline

A_8 & 2^3\colon\!\SL_3(2) &
\psi_2 & \circ2 & 3 & \chi_{10} & \circ2 & 45 & 0, 3, 5 \\ \cline{3-9}
& &
3 & + & 3 & \varphi_{10} & + & 45 & 7 \\ \cline{2-9}
       & 2^4\colon\!(S_3 \times S_3) &
\begin{array}{c} \zeta_2 \boxtimes 1 \\ 1 \boxtimes \zeta_2\end{array} & + & 1 &
\chi_{9} & + & 35 & 0, 3, 5, 7 \\ \cline{3-9}
& &
\begin{array}{c} \zeta_2 \boxtimes 2 \\ 2 \boxtimes \zeta_2\end{array} & + & 2 &
\chi_{14} & + & 70 & 0, 5, 7 \\ \hline\hline

A_9 & \SL_2(8)\colon\!3 & \begin{array}{c} \zeta_3 \\ \bar{\zeta}_3 \end{array} 
& \circ & 1 & \chi_{14} & + & 120 & 0, 5 \\ \hline
\end{array}
$$
\end{table}
}

If $n$ is large enough, the possible maximal block stabilizers are rather
restricted.

\begin{prop}
{\rm [Nett-Noeske~\cite{NeNoe}]}
Suppose that $n \geq 10$ and let $G = 2.A_n$ denote the twofold cover of $A_n$.
Let $H \leq G$ be the inverse image of a maximal, transitive subgroup of~$A_n$.
If $H$ is the block stabilizer of an irreducible imprimitive non-trivial
$K$-character, then $n = 2m$ is even and $H/Z(G) = (S_m \wr S_2) \cap A_n$.
\end{prop}

The next result gives all generic examples in case of $\text{\rm char}(K) = 0$.
\begin{thm}
\label{GenericAlternating}
{\rm [Djokovi\'c-Malzan, Nett-Noeske, Schur]}
Suppose that $\text{\rm char}(K) = 0$. Then the irreducible imprimitive
$K$-characters of $A_n$ or $2.A_n$, $n \geq 10$ are exactly those given
in Table~\ref{tab:clredGenericAlternating}. 

The group~$E$ in the middle column and second row of this table is a group of 
order~$4$; it is elementary abelian, if
$m$ is even, and cyclic, otherwise. The characters $1_3$ and $1_4$ in
the same column are characters of degree~$1$ which have $A_m \times A_m$ 
in their kernel, but not $(S_m \times S_m) \cap A_m$.

The Frobenius-Schur indicators in this table are computed as follows:
$$\iota_m = \left\{\begin{array}{cl}
+, & \text{\rm if\ } m(m-1)/2 \text{\rm\ is even,} \\
\circ, & \text{\rm if\ } m(m-1)/2 \text{\rm\ is odd;} \\
\end{array}\right.
$$

$$\kappa_m = \left\{\begin{array}{cl}
+, & \text{\rm if\ } m \text{\rm\ is even,} \\
\circ, & \text{\rm if\ } m \text{\rm\ is odd;} \\
\end{array}\right.
$$

$$\lambda_m = \left\{\begin{array}{cl}
+, & \text{\rm if\ } m(m-1)/2 \equiv  (0 \text{\rm\ mod\ } 8), \\
-, & \text{\rm if\ } m(m-1)/2 \equiv  (4 \text{\rm\ mod\ } 8), \\
\circ, & \text{\rm if\ } m(m-1)/2 \equiv  (2 \text{\rm\ mod\ } 4);
\end{array}\right.
$$

$$\nu_m = \left\{\begin{array}{cl}
+, & \text{\rm if\ } m(m-1)/2 \equiv  (0, 6 \text{\rm\ mod\ } 8), \\
-, & \text{\rm if\ } m(m-1)/2 \equiv  (2, 4 \text{\rm\ mod\ } 8).
\end{array}\right.
$$

\end{thm}
\begin{prf}
The result for the simple groups is contained in~\cite{DM2}, the result
for the covering groups in \cite{NeNoe}. The degrees of the characters involved
can be computed from the hook formulae (\cite[Theorem 2.3.31]{jake} and
\cite[Theorem 10.7]{HoHu}). The indicators of the characters are determined
from \cite[Theorem 2.5.13]{jake} and \cite[S\"atze IV, V]{SchurIndicators}.
\end{prf}

\begin{table}[h]
\caption{\label{tab:clredGenericAlternating}Generic alternating group examples
(characteristic~$0$)}
$$\begin{array}{|c||c|c|c|} \hline
G & A_{m^2+1}, m \geq 3 & A_{2m}, m \geq 5 & 2.A_n, n = \frac{m(m+1) + 2}{2},
\frac{m(m-1)}{2} \text{\rm\ even} \\ \hline
H & A_{m^2} & (A_m \times A_m)\colon\!E & 2.A_{n - 1} \\ \hline
\psi & \begin{array}{c} [m^m]_1 \\ \left[ m^m \right]_2  \end{array} &
\begin{array}{c} 1_3 \\ 1_4 \end{array} &
\begin{array}{c} < m, m-1, \ldots, 2, 1 >_1 \\ < m, m-1, \ldots, 2, 1 >_2 \end{array} \\ \hline
\text{\rm ind} & \iota_m & \kappa_m & \lambda_m \\ \hline
\psi(1) & \frac{(m^2)!}{2m^m\prod_{i=1}^{m-1}[i(2m-i)]^i} & 1 &
2^{m(m-1)/4 - 1}\cdot \frac{(n - 1)!}{m!} \prod_{i=1}^{m-1} \frac{(i+1)!}{(2i+1)!} \\ \hline
\chi & [m+1, m^{m-1}] & [m+1,1^{m-1}] & < m + 1, m-1, \ldots, 2, 1 > \\ \hline
\text{\rm ind} & + & + & \nu_m \\ \hline
\chi(1) & (m^2 + 1)\psi(1) & \left( \begin{array}{c} 2m-1 \\ m \end{array} \right) &
n\cdot \psi(1) \\ \hline
\end{array}
$$
\end{table}

\chapter{Exceptional Schur Multipliers and Exceptional Isomorphisms}
\markleft{5.\ EXCEPTIONAL SCHUR MULTIPLIERS}
\label{ExceptionalMultipliers}

In this chapter~$K$ is an algebraically closed field of characteristic 
$\ell \geq 0$  and we use, once more, the language of $K$-characters.
We determine the irreducible imprimitive $K$-characters for all 
quasisimple groups~$G$ of Lie type, such that $G/Z(G)$ has an exceptional 
Schur multiplier or such that $G/Z(G)$ is isomorphic to a simple group
of Lie type of a different characteristic, with the only restriction that
$\mbox{\rm char}(K) = 0$ if $G/Z(G) = {^2\!E}_6(2)$. 

The simple groups of Lie type with an exceptional Schur multiplier are listed 
in Table~\ref{ExcMGroups} (cf.\ \cite[Table 6.1.3]{GolYSo3}). 
There are two entries for the groups $\PSL_2(9)$ and $\Sp_4(2)'$, although these 
groups are isomorphic. The reason is that each of them has an exceptional Schur 
multiplier, one as a group of Lie type of characteristic~$3$, the other as a
group of Lie type of characteristic~$2$. On the other hand we have not listed 
the group $\PSp_4(3)$ isomorphic to $\SU_4(2)$, since $\PSp_4(3)$ does not have 
an exceptional multiplier as group of Lie type of characteristic~$3$. For the
same reason we have not listed $\PSL_2( 5 )$ nor $\PSL_2( 7 )$.

\begin{table}[h]
\centering
\caption{\label{ExcMGroups} The simple groups of Lie type with exceptional
Schur multiplier}
$$
\begin{array}{cccccc}  \hline\hline
\SL_{2}(4) & \PSL_{2}(9) & \SL_{3}(2) & \PSL_{3}(4)  & 
\SL_{4}(2) & \SU_{4}(2) \rule[ 0pt]{0pt}{ 15pt} \\ 
\PSU_{4}(3) & \PSU_{6}(2) & \Sp_{4}(2)' & \Sp_{6}(2) & \Omega_{7}(3) & 
\Omega_{8}^{+}(2) \\
^{2}\!B_{2}(8) & G_{2}(3) & G_{2}(4) & F_{4}(2) & 
{^2\!E}_{6}(2) \rule[- 7pt]{0pt}{  5pt} \\ \hline\hline
\end{array}
$$
\end{table}

\noindent
There are exactly six exceptional isomorphisms between simple groups of Lie 
type of different characteristics (see \cite[Theorem 2.2.10]{GolYSo3}). We give 
these isomorphisms in Table~\ref{ExcIsos}. Of these groups, only ${^2\!G}_2(3)' 
\cong \SL_2( 8 )$ and $G_2(2)' \cong \SU_3(3)$ do not occur in 
Table~\ref{ExcMGroups}.

\begin{table}[h]
\centering
\caption{\label{ExcIsos} Exceptional isomorphisms}
$$
\begin{array}{cc}  \hline\hline
\SL_2( 4 ) \cong \PSL_2( 5 ), & \SL_3(2) \cong \PSL_2( 7 ), \rule[ 0pt]{0pt}{ 15pt} \\
\PSL_2( 9 ) \cong \Sp_4(2)', & {^2\!G}_2(3)' \cong \SL_2( 8 ), \\
G_2(2)' \cong \SU_3(3), & \SU_4(2) \cong \PSp_4(3). \rule[- 7pt]{0pt}{  5pt} \\ \hline\hline
\end{array}
$$
\end{table}

The following theorem is the main result of this chapter.
\begin{thm}
\label{ExcMultThm}
Let~$G$ be a quasisimple group such that $G/Z(G)$ is one of the groups occurring 
in Tables~\ref{ExcMGroups} and~\ref{ExcIsos}, and let~$H$ be a maximal subgroup 
of~$G$.
In case~$G/Z(G) = {^2\!E}_6(2)$ and $H/Z(G)$ is a parabolic subgroup of~$G/Z(G)$, 
assume that $\mbox{\rm char}(K) = 0$. 

If~$\chi$ 
is an irreducible faithful $K$-character of~$G$ such that $\chi = 
\Ind_H^G(\psi)$ for some $K$-character~$\psi$ of~$H$, then~$\chi$ is as listed 
in one of the Tables~\ref{tab:clred1} -- \ref{tab:clred2d}. In the first two of 
these tables, $H/Z( G )$ is a non-parabolic subgroup of $G/Z(G)$, in the last 
five of these tables, $H/Z(G)$ is parabolic.
\end{thm}
The proof of Theorem~\ref{ExcMultThm} is given in 
Section~\ref{ProofExcMultThm} below.
Since we only consider faithful irreducible characters of~$G$,
the center $Z(G)$ is cyclic in all cases, even if the Schur multiplier of the 
simple group $G/Z(G)$ is non-cyclic. Also, if $\Ind_H^G(\psi)$ is faithful,
the restriction of~$\psi$ to $Z(G) \leq H$ is faithful as well.

\section{Description of the tables} We employ the same conventions as in
Table~\ref{tab:clredSporadic} for the sporadic groups. There are a couple of 
additions to the notations used there. If $H$ is of the form $L.2$ and $\psi$ 
is an extension of an irreducible character $\psi_i$ of $L$, where~$i$ is the 
number in the compound character table of the Atlas, we write $\psi = \psi_i'$ 
or $\psi_i''$. An example occurs in Table~\ref{tab:clred1}, in the entry for
$G = 2.\Sp_6(2)$ and $H = 2.\SU_4(2)\colon\!2$. Here, the group 
$L := 2.\SU_4(2)$ has a faithful irreducible character of degree~$20$ and 
Frobenius-Schur indicator~$-$, denoted $\chi_{23}$ in \cite{ATLAS}. The two
extensions of $\chi_{23}$ to $H$ are denoted by $\psi_{23}'$ and $\psi_{23}''$
in Table~\ref{tab:clred1}. A similar convention is used for characters which
are just denoted by their degrees. A further addition concerns the tables with
the parabolic examples. Suppose that~$G$ is of Lie type of characteristic~$\ell$
(i.e.,~$Z(G)$ is a quotient of the generic multiplier of $G/Z(G)$)
and that $H$ is a parabolic subgroup of~$G$. Then $H$ is of the form 
$U\colon\!L$, with $U = O_\ell(H)$, the unipotent radical of~$H$. In this
situation~$U$ is in the kernel of~$\psi$ by 
Proposition~\ref{ParabolicStabilizer}, and we give the Atlas name of $\psi$ as 
a character of~$L$. Consider, for example, the entry for $G = \SL_4(2)$ and 
$H = 2^3\colon\!\SL_3(2)$ in Table~\ref{tab:clred2}. Here, the 
character~$\psi_2$ denotes the irreducible character of the simple group 
$\SL_3(2)$ with Atlas name $\chi_2$.

Here, we also have an additional column headed ``Nts'' which contains further 
information of various kind, which we are now going to describe in more detail.

\begin{itemize}
\item[1.] Let~$p$ be the defining characteristic of~$G$. There is an isomorphism
of~$G$ with a group of Lie type~$G'$ of characteristic $p' \neq p$, which maps 
$H$ to a parabolic subgroup of~$G'$. ($H$ is a parabolic subgroup in the wrong 
characteristic.)

\item[2.] There are two $K$-characters $\psi_1 \neq \psi_2$ of~$H$ with
$\Ind_H^G(\psi_1) = \chi = \Ind_H^G(\psi_2)$.

\item[3.] This is an example with $\char(K) = \ell = p$, where $p$ is the defining
characteristic of~$G$. It also occurs in
the result of Seitz cited above (Theorem~\ref{SameCharacteristic}).

\item[4.] There is more than one maximal block stabilizer for $\chi$ (up to 
conjugation in $G$).

\item[5.] This is an example where~$\chi$ is not liftable to characteristic~$0$.

\item[6.] This is a non-canonical ``parabolic'' example, i.e.,~$H/Z(H)$
is a parabolic subgroup of~$G/Z(G)$, but~$G$ itself is not a group of
Lie type.

\item[$f_\ell$.] Let~$p$ be the defining characteristic of~$G$. There is an 
isomorphism of~$G$ with a group of Lie type~$G'$ of characteristic $\ell \neq p$.
(A defining, but wrong characteristic example). It also occurs in the result of 
Seitz cited above (Theorem~\ref{SameCharacteristic}).

\item[$\dagger$.] There is an outer automorphism of~$G$ of order~$2$ 
transporting the triple
$(H, \psi, \chi)$ to the triple $(H', \psi', \chi')$ such that $H$ and~$H'$ are 
not conjugate in~$G$ and $\chi \neq \chi'$. We have not included a row for 
$(H', \psi', \chi')$.

\item[$\dagger'$.] Similar to $\dagger$, but with an outer automorphism of 
order~$3$.
\end{itemize}
Additional explanations of various details are contained in 
Remark~\ref{ExceptionalExplanations} below.

\begin{rem}\label{ExceptionalExplanations}
{\rm
\textit{Table~\ref{tab:clred1}.} Entry for $2.\SU_4(2)$: The two characters 
$2_2$ and $2_3$ of $H = 3^3\colon\!\GL_2(3)$ denote the two faithful 
characters of $\GL_2(3)$ of degree~$2$, inflated to~$H$.

\textit{Table~\ref{tab:clred1a}.} (1) Entries for $3_1.\SU_4(3) = 
12_1.\PSU_4(3)$, $3_2.\SU_4(3) = 12_2.\PSU_4(3)$: The subgroup 
$4 \circ \Sp_4(3)$, the central product of $Z(\SU_4(3))$ with $\Sp_4(3)$ (see
\cite[Table~3.5 B]{kleidman:4}) is isoclinic to $2 \times 2.\SU_4(2)$. The 
characters $\psi_{21}'$ and $\psi_{22}'$ denote extensions of the two characters
of degree~$4$ of $\Sp_4(3)$ to~$H$.
There are two conjugacy classes of subgroups isomorphic to $4 \circ \Sp_4(3)$ in each
of $3_1.\SU_4(3)$ and $3_2.\SU_4(3)$. In $3_1.\SU_4(3)$, only one of these 
conjugacy classes contains block stabilizers. In $3_2.\SU_4(3)$ the two conjugacy
classes are swapped by an outer automorphism of $3_2.\SU_4(3)$ (see 
\cite[pp.~$52$--$53$]{ATLAS}), and so both conjugacy classes contain block 
stabilizers.

(2) At various places in Table~\ref{tab:clred1a} the following phenomenon occurs.
There is an orbit of length four of Galois conjugate irreducible characters
of~$H$, but only two of the four induce to irreducible characters of~$G$.
Consider, for example, the entry for $G = 3_2.\PSU_4(3)$ and $H = 3 \times
\SU_4(2)$. There are four irreducible characters of degree~$5$ which restrict
to faithful characters of~$Z(G)$, namely $\zeta_3 \boxtimes \psi_2$,
$\bar{\zeta}_3 \boxtimes \psi_2$, $\zeta_3 \boxtimes \psi_3$, and 
$\bar{\zeta}_3 \boxtimes \psi_3$, but only two of these, a pair of complex 
conjugates, induce to irreducible characters of~$G$. This happens exactly in 
the following cases, described in the form $(G,H,\psi)$, where~$\psi$ is one 
of the four characters in the Galois orbit. 
$$\begin{array}{c|c|c} \hline
G & H & \psi \\ \hline
3_2.\PSU_4(3) & 3 \times \SU_4(2) & \zeta_3 \boxtimes \psi_2 \\
6_2.\PSU_4(3) & 6 \times \SU_4(2) & \zeta_6 \boxtimes \psi_2 \\
3.\Omega_7(3) & 3 \times \PSL_4(3)\colon\!2_2 & \zeta_3 \boxtimes \zeta_2 \\
2.\SU_6(2) & 6 \times \SU_5(2) & \zeta_6 \boxtimes \psi_3 \\ \hline
\end{array}
$$

\textit{Table~\ref{tab:clred2}.} (1) Entries for $\PSL_3(4)$ and $\SL_3(4)$: The
character denoted by~$3$ is the inflation to~$H$ of the reduction modulo~$5$ 
of the two ordinary characters $\psi_2$ and $\psi_3$ of~$A_5$.

(2) Entry for $4_2.\SL_3(4) = 12_2.\PSL_3(4)$: The character denoted by $4_2$ 
is one of two faithful characters of $4.2^4\colon\!A_5$ of degree~$4$, inflated 
to~$H$.

\textit{Table~\ref{tab:clred2a}.} (1) Entries for $\PSU_4(3)$ and $G_2(3)$: The
two characters $2_2$ and $2_3$ of $H = 3_+^{1+4}\colon\!2.S_4$, respectively
$H = (3^{1+2}_+ \times 3^2)\colon\!2.S_4$ denote the two faithful characters 
of the Levi subgroup $L = 2.S_4$ of~$H$ of degree~$2$, inflated to~$H$.

(2) Entry for $2.\PSU_4(3)$: The two characters $2_4$ and $2_5$ of $H = 
3_+^{1+4}\colon\!2.(2.S_4)$ denote two irreducible characters of the Levi 
subgroup $L = 2.(2.S_4)$ of~$H$ of degree~$2$, inflated to~$H$. These two
characters are determined by the fact that $Z(H)$ is not in their kernel and 
that they are not rational.

(3) Entry for $\SU_4(3) = 4.\PSU_4(3)$: The two characters denoted by $1_3$ and
$1_4$ are linear characters of order~$8$. The characters denoted by $3_4$ 
and~$3_5$ denote two faithful characters of the Levi subgroup $L = 4.(2.S_4)$ 
of~$H$ of degree~$3$, inflated to~$H$.

\textit{Table~\ref{tab:clred2b}.} (1) Entries for $\Omega_7(3)$: The two 
characters $2_5$ and $2_6$ of $H = 3_+^{1+6}\colon\!(2.A_4 \times A_4).2$ 
denote the two non-rational irreducible characters of the Levi subgroup 
$L = (2.A_4 \times A_4).2$ of~$H$ of degree~$2$, inflated to~$H$.
Moreover, $6_3$ and $6_4$ are the two faithful irreducible characters of 
degree~$6$ of~$L$ inflated to~$H$.

(2) Entries for $2.\Omega_7(3)$: The two characters $2_7$ and $2_8$ of 
$H$ denote the two non-rational irreducible characters of the Levi 
subgroup $L = 2.(2.A_4 \times A_4).2$ 
of~$H$ of degree~$2$, which do not have $Z(H)$ in their kernel, inflated 
to~$H$. Moreover, $6_5$ and $6_6$ are the two non-rational irreducible 
characters of degree~$6$ of~$L$, which do not have $Z(H)$ in their kernel,
inflated to~$H$.
}
\end{rem}

\begin{table}[h]
\centering
\caption{\label{tab:clred1} Non-parabolic examples}
$$
\begin{array}{|c|c|c|c|c|c|c|c|l|l|}  \hline
G & H & \psi & \text{\rm ind} & \psi(1) & \chi & \text{\rm ind} &
\chi(1) & \text{\rm Primes} & \text{Nts} \\ \hline\hline

2.\SL_2(4) & 5\colon\!4 & 
\begin{array}{c} \zeta_4 \\ \bar{\zeta}_4 \end{array} & \circ & 1 & 
\chi_{9} & - &  6 & 0,3 & 1, 2 \\ \hline\hline

\SL_3(2)   & 7\colon\!3 & 
\begin{array}{c} \zeta_3 \\ \bar{\zeta}_3 \end{array} & \circ & 1 &
\chi_{6} & + & 8 & 0, 2 & 1, 2, 3 \\ \hline
2.\SL_3(2) & 7\colon\!6 & 
\begin{array}{c} \zeta_6 \\ \bar{\zeta}_6 \end{array} & \circ & 1 &
\chi_{11} & - & 8 & 0 & 1, 2 \\ \hline\hline

3.\PSL_2(9) & 3 \times A_5 & 
\zeta_3 \boxtimes 1 & \circ2 & 1 & \chi_{16} & \circ2 & 6 & 0, 5 & 4 \\ \hline
3.\SL_2(9) & 3 \times 2.A_5 & 
\zeta_3 \boxtimes \psi_6 & \circ2 & 2 & \chi_{21} & \circ2 & 12 & 0 & 4 \\ 
\cline{3-10}
& & 
\zeta_3 \boxtimes \psi_7 & \circ2 & 2 & \chi_{22} & \circ2 & 12 & 0 & 4 \\ 
\hline\hline

\Sp_4(2)' & 3^2\colon\!4 &
\begin{array}{c} \zeta_4 \\ \bar{\zeta}_4 \end{array} & \circ & 1 &
\chi_{7} & + & 10 & 0, 5 & 1, 2 \\ \hline
2.\Sp_4(2)' & 3^2\colon\!8 &
\begin{array}{c} \zeta_8 \\ \bar{\zeta}_8 \end{array} & \circ & 1 &
\chi_{12} & - & 10 & 0, 5 & 1, 2 \\ \cline{3-10}
& & 
\begin{array}{c} \zeta_8^{3} \\ \bar{\zeta}_8^{3} \end{array} & \circ & 1 &
\chi_{13} & - & 10 & 0, 5 & 1, 2 \\ \hline
3.\Sp_4(2)' & 3 \times A_5 & 
\zeta_3 \boxtimes 1 & \circ2 & 1 & \chi_{16} & \circ2 & 6 & 0, 5 & 4 \\ \hline
6.\Sp_4(2)' & 3 \times 2.A_5 & 
\zeta_3 \boxtimes \psi_6 & \circ2 & 2 & \chi_{21} & \circ2 & 12 & 0 & 4 \\ 
\cline{3-10}
& & 
\zeta_3 \boxtimes \psi_7 & \circ2 & 2 & \chi_{22} & \circ2 & 12 & 0 & 4 \\ 
\hline\hline

{^2\!G}_2(3)' & 2^3\colon\!7 &
\begin{array}{c} \zeta_7 \\ \bar{\zeta}_7 \end{array} & \circ & 1 &
\chi_{7} & + & 9 & 0, 3 & 1, 2 \\ \cline{3-10}
& & 
\begin{array}{c} \zeta_7^3 \\ \bar{\zeta}_7^3 \end{array} & \circ & 1 &
\chi_{8} & + & 9 & 0, 3 & 1, 2 \\ \cline{3-10}
& & 
\begin{array}{c} \zeta_7^2 \\ \bar{\zeta}_7^2 \end{array} & \circ & 1 &
\chi_{9} & + & 9 & 0, 3 & 1, 2 \\ \hline\hline

G_2(2)' & 3_+^{1+2}\colon\!8 &
\begin{array}{c} \zeta_8 \\ {\zeta}_8^5 \end{array} & \circ & 1 &
\chi_{11} & \circ & 28 & 0, 7 & 1, 2 \\ \cline{3-10}
& & 
\begin{array}{c} \bar{\zeta}_8 \\ \bar{\zeta}_8^5 \end{array} & \circ & 1 &
\chi_{12} & \circ & 28 & 0, 7 & 1, 2 \\ \hline\hline

2.\SL_4(2) & 2.A_7 & 
\begin{array}{c} 6_1 \\ 6_2 \end{array} & - & 6 &
\varphi_{15} & - & 48 & 3 & 2, 5 \\ \hline\hline

4_1.\PSL_3(4) & 4 \times A_6 & \zeta_4 \boxtimes 1 & \circ2 & 1  &
             \chi_{21} & \circ2 & 56 & 0, 5, 7 & 4 \\ \hline
4_1.\SL_3(4)& 12 \times \PSL_2(7) & \zeta_{12} \boxtimes 1 & \circ4 & 1 &
             \chi_{52} & \circ4 & 120 & 0, 5 & 4 \\ \hline\hline

2.\SU_4(2) & 3^3\colon\!\GL_2(3) & 
\begin{array}{c} 2_2 \\ 2_3 \end{array} & \circ & 2 & \chi_{34} & - & 80 
& 0, 5 & 1, 2 \\ \hline\hline

2.\Sp_{6}(2) & 2.\SU_4(2)\colon\!2 & 
\begin{array}{c} \psi_{23}' \\ \psi_{23}'' \end{array} & \circ & 20 & 
\chi_{42} & + & 560 & 0, 5, 7 & 2 \\ \cline{2-10}
           & 2 \times \SU_{3}(3)\colon\!2 & \zeta_2 \boxtimes \zeta_2 & + & 1 & 
\chi_{37}  & + & 120 & 0, 5, 7 & \\ \cline{3-10}
& & \begin{array}{c} \zeta_2 \boxtimes \psi_2' \\ \zeta_2 \boxtimes \psi_2'' \end{array} 
& \circ & 6 & \chi_{43} & + & 720 & 0, 5   & 2 \\ \hline
\end{array}
$$
\end{table}

\begin{table}[h]
\centering
\caption{\label{tab:clred1a} Non-parabolic examples (continued)}
$$
\begin{array}{|c|c|c|c|c|c|c|c|l|l|}  \hline
G & H & \psi & \text{\rm ind} & \psi(1) & \chi & \text{\rm ind} &
\chi(1) & \text{\rm Primes} & \text{Nts} \\ \hline\hline

3_1.\PSU_4(3)  & 3 \times \SU_4(2) & \begin{array}{c} 
                               \zeta_3 \boxtimes \psi_2 \\
			       \zeta_3 \boxtimes \psi_3 
			       \end{array} & \circ2 & 5 & 
             \chi_{68} & \circ2 & 630 & 0, 5, 7 & 2, 4 \\ \hline
3_2.\PSU_4(3) & 3 \times \SU_4(2) & \zeta_3 \boxtimes 1 & \circ2 & 1 &
             \chi_{102} & \circ2 & 126 & 0, 5, 7 & 4 \\ \cline{3-10}
	   &                 & \zeta_3 \boxtimes \psi_2 & \circ2 & 5 &
             \chi_{107} & \circ2 & 630 & 0, 5, 7 & 4 \\ \hline
6_1.\PSU_4(3) & 6 \times \SU_4(2) & \zeta_6 \boxtimes 1 & \circ2 & 1 & 
             \chi_{75} & \circ2 & 126 & 0, 5, 7 & 4 \\ \cline{3-10}
          &                & \begin{array}{c}
			     \zeta_6 \boxtimes \psi_2 \\
			     \zeta_6 \boxtimes \psi_3
			     \end{array} & \circ2 & 5 & 
             \chi_{82} & \circ2 & 630 & 0, 5, 7 & 2, 4 \\ \hline
6_2.\PSU_4(3) & 6 \times \SU_4(2) & \zeta_6 \boxtimes 1 & \circ2 & 1 &
             \chi_{113} & \circ2 & 126 & 0, 5, 7 & \dagger \\ \cline{3-10}
           &                 & \zeta_6 \boxtimes \psi_2 & \circ2 & 5 &
             \chi_{121} & \circ2 & 630 & 0, 5, 7 & 4 \\ \hline
3_1.\SU_4(3)& 3 \times (4 \circ \Sp_4(3)) & 
                                             \begin{array}{c}
                                             \zeta_3 \boxtimes \psi_{21}' \\ 
					     \zeta_3 \boxtimes \psi_{22}' \\ 
					     \end{array} & \circ4 & 4 &
             \chi_{95} & \circ4 & 504 & 0, 7 & 2 \\ \hline
3_2.\SU_4(3)& 3 \times (4 \circ \Sp_4(3)) & 
             \zeta_3 \boxtimes \psi_{21}' &\circ4 &4 &
             \chi_{128} & \circ4 & 504 & 0, 7 & 4 \\ \cline{3-10}
	  &                                & \zeta_3 \boxtimes \psi_{22}' &\circ4 &4 &
	     \chi_{129} & \circ4 & 504 & 0, 7 & 4 \\ \hline\hline

3.G_2(3)   & 3 \times (\SU_3(3)\colon\!2) & \zeta_3 \boxtimes 1 & \circ2 & 1 & 
             \chi_{28} & \circ2 & 351 & 0, 13 & 4 \\ \cline{3-10}
	   &  & \zeta_3 \boxtimes \zeta_2 & \circ2 & 1 & 
	     \chi_{29} & \circ2 & 351 & 0, 7, 13 & \dagger \\ \cline{2-10}
           & 3 \times (\SL_3(3)\colon\!2) & \zeta_3 \boxtimes \zeta_2 & \circ2 
	   & 1 & \chi_{31} & \circ2 & 378 & 0, 7, 13 & \dagger \\ \hline\hline

2.G_2(4)   & (2 \times \SU_3(4)).2 & 
             \begin{array}{c} \zeta_4 \\ \bar{\zeta}_4 \end{array} & \circ & 1 &
             \chi_{43} & - & 2016 & 0, 3, 7 & 2 \\ \hline\hline

3.\Omega_{7}(3) & 3 \times \PSL_4(3)\colon\!2_2 & \zeta_3 \boxtimes \zeta_2 & \circ2 & 1 &
               \chi_{92} & \circ2 & 378 & 0, 5, 7, 13 & \\ \hline
6.\Omega_{7}(3) & 3 \times \SL_4(3)\colon\!2_2 & 
              \begin{array}{c} \zeta_3 \boxtimes \psi_{30}' \\ 
	                       \zeta_3 \boxtimes \psi_{30}'' \end{array} & \circ2 & 40 &
              \chi_{127} & \circ2 & 15\,120 & 0, 5, 7 & 2 \\ \cline{2-10}
           & 2 \times 3.G_{2}(3) & \zeta_2 \boxtimes \psi_{25} & \circ2 & 27 & 
             \chi_{134} & \circ2 & 29160 & 0, 5 & 4 \\ \hline\hline

2.\SU_{6}(2) & 6 \times \SU_5(2) & \zeta_6 \boxtimes 1 & \circ2 & 1 &
             \chi_{122} & \circ2 & 672 & 0, 5, 7, 11 & \\ \cline{3-10}
           &                 & {\zeta}_6 \boxtimes \psi_3 & \circ2 & 11 &
	     \chi_{130} & \circ2 & 7392 & 0, 5, 7, 11 & \\ \hline
\end{array}
$$
\end{table}

\begin{table}[h]
\centering
\caption{\label{tab:clred2} Parabolic examples}
$$
\begin{array}{|c|c|c|c|c|c|c|c|l|c|}  \hline
G & H & \psi & \text{\rm ind} & \psi(1) & \chi & \text{\rm ind} &
\chi(1) & \text{\rm Primes} & \text{Nts} \\ \hline\hline

\SL_2(4) & 2^2\colon\!3 &
\begin{array}{c} \zeta_3 \\ \bar{\zeta}_3 \end{array} & \circ & 1 &
\chi_{5} & + &  5 & 0, 5 & 2, f_5 \\ \hline\hline

\SL_3(2) & 2^2\colon\!\SL_2(2) &
\zeta_2 & + & 1 & \chi_{5} & + & 7 & 0, 3, 7 & 4, f_7 \\ \hline\hline

\PSL_2(9) & 3^2\colon\!4 &
\begin{array}{c} \zeta_4 \\ \bar{\zeta}_4 \end{array} & \circ & 1 &
\chi_{7} & + & 10 & 0, 5 & 2 \\ \hline
\SL_2(9) & 3^2\colon\!8 &
\begin{array}{c} \zeta_8 \\ \bar{\zeta}_8 \end{array} & \circ & 1 &
\chi_{12} & - & 10 & 0, 5 & 2 \\ \cline{3-10}
& & 
\begin{array}{c} \zeta_8^{3} \\ \bar{\zeta}_8^{3} \end{array} & \circ & 1 &
\chi_{13} & - & 10 & 0, 5 & 2 \\ \hline\hline

\SL_4(2) & 2^3\colon\!\SL_3(2) &
\psi_2 & \circ2 & 3 & \chi_{10} & \circ2 & 45 & 0, 3, 5 & 4 \\ \cline{3-10}
& & 
3 & + & 3 & \varphi_{10} & + & 45 & 7 & 4 \\ \cline{2-10}
       & 2^4\colon\!(\SL_2(2) \times \SL_2(2)) & 
\begin{array}{c} \zeta_2 \boxtimes 1 \\ 1 \boxtimes \zeta_2\end{array} & + & 1 &
\chi_{9} & + & 35 & 0, 3, 5, 7 & 2 \\ \cline{3-10}
& &
\begin{array}{c} \zeta_2 \boxtimes 2 \\ 2 \boxtimes \zeta_2\end{array} & + & 2 &
\chi_{14} & + & 70 & 0, 5, 7 & 2 \\ \hline\hline
\PSL_3(4) & 2^4\colon\!A_5 & \psi_2 & + & 3 & \chi_8 & + & 63 & 0, 3, 7 
                                                           & 4 \\ \cline{3-10}
       &               & \psi_3 & + & 3 & \chi_9 & + & 63 & 0, 3, 7 & 4 \\ \cline{3-10}
       &               & 3 & + & 3 & \varphi_8 & + & 63 & 5 & 4 \\ \hline
\SL_3(4) & 2^4\colon\!( 3 \times A_5) & \zeta_3 \boxtimes 1 & \circ2 & 1 & 
                  \chi_{35} &  \circ2 & 21 & 0, 5, 7 & 4 \\ \cline{3-10}
         &              & \zeta_3 \boxtimes \psi_2 & \circ2 & 3 
		& \chi_{38} &  \circ2 & 63 & 0, 7 & 4 \\ \cline{3-10}
         &              & \zeta_3 \boxtimes \psi_3 & \circ2 & 3
                & \chi_{39} &  \circ2 & 63 & 0,  7 & 4 \\ \cline{3-10}
         &              & \zeta_3 \boxtimes 3 & \circ2 & 3
                & \varphi_{30} &  \circ2 & 63 & 5 & 4 \\ \cline{3-10}
         &              & \zeta_3 \boxtimes \psi_4 & \circ2 & 4 
                & \chi_{40} &  \circ2 & 84 & 0, 7 & 4 \\ \hline
4_2.\SL_3(4) & 3 \times (4.2^4\colon\!A_5) & \zeta_3 \boxtimes 4_2 &  \circ4 & 4 & 
                  \chi_{58} &  \circ4 & 84 & 0, 7 & 4, 6 \\ \hline
\end{array}
$$
\end{table}

{\small
\begin{table}[h]
\caption{\label{tab:clred2a} Parabolic examples (continued)}
$$
\begin{array}{|c|c|c|c|c|c|c|c|l|c|}  \hline
G & H & \psi & \text{\rm ind} & \psi(1) & \chi & \text{\rm ind} &
\chi(1) & \text{\rm Primes} & \text{Nts} \\ \hline\hline

\SU_4(2) & 2^4\colon\!A_5 &
\begin{array}{c} \psi_2 \\ \psi_3 \end{array} & + & 3 & \chi_{20} & + &
81 & 0, 3 & 2, f_3 \\ \hline\hline

{^2\!B}_2(8) & 2^{3+3}\colon\!7 &
\begin{array}{c} \zeta_7 \\ \bar{\zeta}_7 \end{array} & \circ & 1 & 
\chi_{8} & + & 65 & 0, 5, 13 & 2 \\ \cline{3-10}
& & \begin{array}{c} \zeta^3_7 \\ \bar{\zeta}^3_7 \end{array} & \circ & 1 & 
\chi_{9} & + & 65 & 0, 5, 13 & 2 \\ \cline{3-10}
& & \begin{array}{c} \zeta^2_7 \\ \bar{\zeta}^2_7 \end{array} & \circ & 1 & 
\chi_{10} & + & 65 & 0, 5, 13 & 2 \\ \hline\hline

\Sp_6(2) & 2^6\colon\!\SL_3(2) & 
\begin{array}{c} \psi_2 \\ \psi_3 \end{array} & \circ & 3 &
\chi_{28} & + & 405 & 0, 3, 5 & 2 \\ \hline\hline

\PSU_4(3) & 3^4\colon\!\PSL_2(9) & 
\begin{array}{c} \psi_4 \\ \psi_5 \end{array} & + & 8 &
\chi_{20} & + & 896 & 0, 2, 7 & 2 \\ \cline{2-10}
& 3_+^{1+4}\colon\!2.S_4 & 
\begin{array}{c} 2_2 \\ 2_3 \end{array} & - & 2 &
\chi_{16} & - & 560 & 0, 5, 7 & 2 \\ \hline
2.\PSU_4(3) & 3^4\colon\!(2 \times \PSL_2(9)) & 
\begin{array}{c} \zeta_2 \boxtimes \psi_4 \\ \zeta_2 \boxtimes \psi_5 \end{array} 
& + & 8 & \chi_{39} & + & 896 & 0, 7 & 2 \\ \cline{2-10}
& 3_+^{1+4}\colon\!2.(2.S_4) & 
\begin{array}{c} 2_4 \\ 2_5 \end{array} & \circ & 2 &
\chi_{35} & + & 560 & 0, 5, 7 & 2 \\ \hline
\SU_4(3) & 3^4\colon\!2.(2 \times \PSL_2(9)) & 
\begin{array}{c} \psi_{10}' \\ \psi_{11}' \end{array} & \circ2 & 8 &
\chi_{65} & \circ2 & 896 & 0, 7 & 2 \\ \cline{2-10}
& 3_+^{1+4}\colon\!4.(2.S_4) & 
\begin{array}{c} 1_3 \\ 1_4 \end{array} & \circ2 & 1 &
\chi_{47} & \circ2 & 280 & 0, 5, 7 & 2 \\  \cline{3-10}
& & \begin{array}{c} 3_4 \\ 3_5 \end{array} & \circ2 & 3 &
\chi_{54} & \circ2 & 840 & 0, 5, 7 & 2 \\ \hline\hline

G_2(3) & (3^{1+2}_+ \times 3^2)\colon\!2.S_4 &
\begin{array}{c} 2_2 \\ 2_3 \end{array} & \circ & 2 & 
\chi_{19} & + & 728 & 0, 7, 13 & 2, \dagger \\ \hline\hline

\Omega^+_8(2) & 2^6\colon\!A_8 & 
\begin{array}{c} \psi_6 \\ \psi_7 \end{array} & \circ & 21 & 
\chi_{47} & + & 2\,835 & 0, 7 & 2, \dagger' \\ \cline{3-10}
& & \begin{array}{c} \psi_{10} \\ \psi_{11} \end{array} & \circ & 45 & 
\chi_{53} & + & 6\,075 & 0, 3, 5 & 2, 4 \\ \hline
\end{array}
$$
\end{table}
}

{\small
\begin{table}[h]
\caption{\label{tab:clred2b} Parabolic examples (continued)}
$$
\begin{array}{|c|c|c|c|c|c|c|c|l|c|}  \hline
G & H & \psi & \text{\rm ind} & \psi(1) & \chi & \text{\rm ind} &
\chi(1) & \text{\rm Primes} & \text{Nts} \\ \hline\hline

G_2(4) & 2^{2+8}\colon\!(3 \times A_5) &
\begin{array}{c} \zeta_3 \boxtimes \psi_2 \\ 
\bar{\zeta}_3 \boxtimes \psi_2 \end{array} & \circ & 3 & 
\chi_{26} & + & 4\,095 & 0, 5, 7, 13 & 2 \\ \cline{3-10}
& & \begin{array}{c} \zeta_3 \boxtimes \psi_3 \\ 
\bar{\zeta}_3 \boxtimes \psi_3 \end{array} & \circ & 3 & 
\chi_{27} & + & 4\,095 & 0, 5, 7, 13 & 2 \\ \cline{2-10}
&  2^{4+6}\colon\!(A_5 \times 3) &
\begin{array}{c} 1 \boxtimes \zeta_3 \\ 
1 \boxtimes \bar{\zeta}_3 \end{array} & \circ & 1 & 
\chi_{16} & + & 1\,365 & 0, 5, 7, 13 & 2 \\ \cline{3-10}
& & \begin{array}{c} \psi_2 \boxtimes \zeta_3 \\ 
\psi_2 \boxtimes \bar{\zeta}_3 \end{array} & \circ & 3 & 
\chi_{24} & + & 4\,095 & 0, 5, 7, 13 & 2 \\ \cline{3-10}
& & \begin{array}{c} \psi_3 \boxtimes \zeta_3 \\ 
\psi_3 \boxtimes \bar{\zeta}_3 \end{array} & \circ & 3 & 
\chi_{25} & + & 4\,095 & 0, 5, 7, 13 & 2 \\ \cline{3-10}
& & \begin{array}{c} \psi_4 \boxtimes \zeta_3 \\ 
\psi_4 \boxtimes \bar{\zeta}_3 \end{array} & \circ & 4 & 
\chi_{32} & + & 5\,460 & 0, 7, 13 & 2 \\ \hline\hline

\Omega_7(3) & 3^{3+3}\colon\!\SL_3(3) &
\begin{array}{c} \psi_4 \\ \psi_5 \end{array} & \circ & 16 &
\chi_{55} & + & 17\,920 & 0, 2, 5, 7 & 2 \\ \cline{3-10}
& & \begin{array}{c} \psi_6 \\ \psi_7 \end{array} & \circ & 16 &
\chi_{54} & + & 17\,920 & 0, 2, 5, 7 & 2 \\ \cline{2-10}
& 3_+^{1+6}\colon\!(2.A_4 \times A_4).2 &
\begin{array}{c} 2_5 \\ 2_6 \end{array} & \circ & 2 &
\chi_{41} & + & 7\,280 & 0, 5, 7, 13 & 2 \\ \cline{3-10}
& & \begin{array}{c} 6_3 \\ 6_4 \end{array} & \circ & 6 &
\chi_{57} & + & 21\,840 & 0, 5, 7, 13 & 2 \\ \hline

2.\Omega_7(3) & 3^5\colon\!(2.\SU_4(2)\colon\!2) & 
\begin{array}{c} \psi_{23}' \\ \psi_{23}'' \end{array} & \circ & 20 &
\chi_{69} & + & 7\,280 & 0, 5, 7, 13 & 2 \\ \cline{3-10}
& & \begin{array}{c} \psi_{30}' \\ \psi_{30}'' \end{array} & \circ & 60 &
\chi_{87} & + & 21\,840 & 0, 5, 7, 13 & 2 \\ \cline{3-10}
& & \begin{array}{c} \psi_{33}' \\ \psi_{33}'' \end{array} & \circ & 64 &
\chi_{88} & + & 23\,296 & 0, 7, 13 & 2 \\ \cline{2-10}

& 3^{3+3}\colon\!(2 \times \SL_3(3)) &
\begin{array}{c} \zeta_2 \boxtimes \psi_4 \\ 
                 \zeta_2 \boxtimes \psi_5 \end{array} & \circ & 16 &
\chi_{85} & + & 17\,920 & 0, 2, 5, 7 & 2 \\ \cline{3-10}
& & \begin{array}{c} \zeta_2 \boxtimes \psi_6 \\ 
                 \zeta_2 \boxtimes \psi_7 \end{array} & \circ & 16 &
\chi_{84} & + & 17\,920 & 0, 2, 5, 7 & 2 \\ \cline{2-10}
& 3_+^{1+6}\colon2.(2.A_4 \times A_4).2 &
\begin{array}{c} 2_7 \\ 2_8 \end{array} & \circ & 2 &
\chi_{70} & + & 7\,280 & 0, 5, 7, 13 & 2 \\ \cline{3-10}
& & \begin{array}{c} 6_5 \\ 6_6 \end{array} & \circ & 6 &
\chi_{86} & + & 21\,840 & 0, 5, 7, 13 & 2 \\ \hline
\end{array}
$$
\end{table}
}

\begin{landscape}
{\small
\begin{table}[h]
\caption{\label{tab:clred2c} Parabolic examples (continued)}
$$
\begin{array}{|c|c|c|c|c|c|c|c|l|c|}  \hline
G & H & \psi & \text{\rm ind} & \psi(1) & \chi & \text{\rm ind} &
\chi(1) & \text{\rm Primes} & \text{Nts} \\ \hline\hline
\PSU_6(2) & 2^9\colon\!\PSL_3(4) & 
\begin{array}{c} \psi_6 \\ \psi_7 \end{array} & \circ & 45 &
\chi_{46} & + & 40\,095 & 0, 3, 5, 11 & 2 \\ \cline{2-10}
& 2^{4 + 8}\colon\!(\SL_2(2) \times \SL_2(4)) &
\begin{array}{c} 1 \boxtimes \psi_2 \\ 1 \boxtimes \psi_3 \end{array}
& + & 3 & \chi_{38} & + & 18\,711 & 0, 3, 7, 11 & 2 \\ \cline{3-10}
& & \begin{array}{c} \zeta_2 \boxtimes \psi_2 \\ \zeta_2 \boxtimes \psi_3 \end{array}
& + & 3 & \chi_{37} & + & 18\,711 & 0, 3, 7, 11 & 2 \\ \cline{3-10}
& & \begin{array}{c} 2 \boxtimes \psi_2 \\ 2 \boxtimes \psi_3 \end{array}
& + & 6 & \chi_{45} & + & 37\,422 & 0, 7, 11 & 2 \\ \hline
3.\PSU_6(2) & 2^9\colon\!\mbox{\rm SL}_3(4) & 
\begin{array}{c} \psi_{36} \\ \psi_{37} \end{array} & \circ 2 & 45 &
\chi_{120} & \circ 2 & 40\,095 & 0, 5, 11 & 2 \\ \cline{2-10}
& 2^{4 + 8}\colon\!(3 \times \SL_2(2) \times \SL_2(4)) & 
\begin{array}{c} \zeta_3 \boxtimes 1 \boxtimes \psi_2 \\ 
\zeta_3 \boxtimes 1 \boxtimes \psi_3 \end{array} & \circ2
 & 3 & \chi_{112} & \circ 2 & 18\,711 & 0, 7, 11 & 2 \\ \cline{3-10}
& & \begin{array}{c} \zeta_3 \boxtimes \zeta_2 \boxtimes \psi_2 \\ 
\zeta_3 \boxtimes \zeta_2 \boxtimes \psi_3 \end{array} & \circ2
 & 3 & \chi_{111} & \circ 2 & 18\,711 & 0, 7, 11 & 2 \\ \cline{3-10}
& & \begin{array}{c} \zeta_3 \boxtimes 2 \boxtimes \psi_2 \\ 
\zeta_3 \boxtimes 2 \boxtimes \psi_3 \end{array} & \circ2
 & 6 & \chi_{119} & \circ 2 & 37\,422 & 0, 7, 11 & 2 \\ \hline\hline

F_4(2) & [2^{20}]\colon\!(\SL_2(2) \times \SL_3(2)) &
\begin{array}{c} \zeta_2 \boxtimes \psi_2 \\ 
\zeta_2 \boxtimes \psi_3 \end{array} & \circ & 3 & 
\chi_{82} & + & 9\,398\,025 & 0, 3, 5, 13, 17 & 2 \\ \hline
\end{array}
$$
\end{table}
}
\end{landscape}

\begin{landscape}
{\small
\begin{table}[h]
\caption{\label{tab:clred2d} Parabolic examples (continued), $\mbox{\rm char}(K) = 0$}
$$
\begin{array}{|c|c|c|c|c|c|c|c|c|}  \hline
G & H & \psi & \text{\rm ind} & \psi(1) & \chi & \text{\rm ind} &
\chi(1) & \text{Nts} \\ \hline\hline
{^2\!E}_6(2) & 2^{8+16}\colon\!\Omega_8^-(2) & 
\begin{array}{c} \psi_{16} \\ \psi_{17} \end{array} & 
+ & 1\,071 & \chi_{94} & + & 24\,748\,759\,035 & 2 \\ \cline{3-9}
& &
\begin{array}{c} \psi_{23} \\ \psi_{24} \end{array} & 
+ & 2\,142 & \chi_{115} & + & 49\,497\,518\,070 & 2 \\ \cline{3-9}
& &
\begin{array}{c} \psi_{29} \\ \psi_{30} \end{array} & 
+ & 2\,835 & \chi_{121} & + & 65\,511\,420\,975 & 2 \\ \cline{3-9}
& & 
\begin{array}{c} \psi_{31} \\ \psi_{32} \end{array} & 
+ & 2\,835 & \chi_{122} & + & 65\,511\,420\,975 & 2 \\ \cline{2-9}
& [2^{29}]\colon\!(\SL_2(2) \times \PSL_3( 4 ) ) & 
\begin{array}{c} \zeta_2 \boxtimes \psi_{6} \\ 
\zeta_2 \boxtimes \psi_7 \end{array} & 
\circ & 45 & \chi_{116} & + & 53\,033\,055\,075 &  \\ \cline{2-9}
 & [2^{31}]\colon\!(\SL_2(4) \times \SL_3( 2 ) ) & 
\begin{array}{c} \psi_2 \boxtimes \psi_2 \\ 
\psi_3 \boxtimes \psi_3 \end{array} & 
\circ & 9 & \chi_{103} & \circ & 31\,819\,833\,045 & 2 \\ \cline{3-9}
& & 
\begin{array}{c} \psi_2 \boxtimes \psi_3 \\ 
\psi_3 \boxtimes \psi_2 \end{array} & 
\circ & 9 & \chi_{104} & \circ & 31\,819\,833\,045 & 2 \\ \hline
3.{^2\!E}_6(2) & 2^{8+16}\colon\!(3 \times \Omega_8^-(2)) & 
\begin{array}{c} \zeta_3 \boxtimes \psi_{16} \\ \zeta_3 \boxtimes \psi_{17} \end{array} & 
\circ 2 & 1\,071 & \chi_{287} & \circ 2 & 24\,748\,759\,035 & 2 \\ \cline{3-9}
& &
\begin{array}{c} \zeta_3 \boxtimes \psi_{23} \\ \zeta_3 \boxtimes \psi_{24} \end{array} & 
\circ 2 & 2\,142 & \chi_{323} & \circ 2 & 49\,497\,518\,070 & 2 \\ \cline{3-9}
& &
\begin{array}{c} \zeta_3 \boxtimes \psi_{29} \\ \zeta_3 \boxtimes \psi_{30} \end{array} & 
\circ 2 & 2\,835 & \chi_{335} & \circ 2 & 65\,511\,420\,975 & 2 \\ \cline{3-9}
& & 
\begin{array}{c} \zeta_3 \boxtimes \psi_{31} \\ \zeta_3 \boxtimes \psi_{32} \end{array} & 
\circ 2 & 2\,835 & \chi_{336} & \circ 2 & 65\,511\,420\,975 & 2 \\ \cline{2-9}
& [2^{29}]\colon\!(\SL_2(2) \times \SL_3( 4 ) ) & 
\begin{array}{c} \zeta_2 \boxtimes \psi_{36} \\ 
\zeta_2 \boxtimes \psi_{37} \end{array} & 
\circ 2 & 45 & \chi_{325} & \circ 2 & 53\,033\,055\,075 &  \\ \cline{2-9}
& [2^{31}]\colon\!(3 \times \SL_2(4) \times \SL_3( 2 ) ) & 
\begin{array}{c} \zeta_3 \boxtimes \psi_2 \boxtimes \psi_2 \\ 
\zeta_3 \boxtimes \psi_3 \boxtimes \psi_3 \end{array} & 
\circ 2 & 9 & \chi_{299} & \circ 2 & 31\,819\,833\,045 & 2 \\ \cline{3-9}
& & 
\begin{array}{c} \zeta_3 \boxtimes \psi_2 \boxtimes \psi_3 \\ 
\zeta_3 \boxtimes \psi_3 \boxtimes \psi_2 \end{array} & 
\circ 2 & 9 & \chi_{300} & \circ 2 & 31\,819\,833\,045 & 2 \\ \hline
\end{array}
$$
\end{table}
}
\end{landscape}

\section{The proofs}
\label{ProofExcMultThm}

The proof of Theorem~\ref{ExcMultThm} follows the same line as the one of 
Theorem~\ref{SporadicTheorem} for the sporadic groups. We will thus be rather 
sketchy here. Let~$G$ be a quasisimple group such that $G/Z(G)$ is one of the
simple groups of Tables~\ref{ExcMGroups} or~\ref{ExcIsos}. With a few 
exceptions, the 
ordinary and modular character tables of $G$ and all of its relevant
maximal subgroups are available in GAP. Many of the character tables
of the maximal subgroups have been computed by Sebastian Dany in his
Diploma thesis \cite{Dany}. For the construction of some of the 
subgroups involved, Rob Wilson's Atlas of Group Representations 
\cite{WWWW} was used. With all this information at hand, the proof
is tedious but straightforward. For the results in the case of
$\mbox{\rm char}(K) = \ell > 0$, Lemma~\ref{Modp} is invoked most of
the time. If~$G$ is such that $G/Z(G)$ is an alternating group, i.e.,
isomorphic to one of $A_5$,~$A_6$ or~$A_8$, some of our results are 
also contained in~\cite{DM2} and~\cite{NeNoe}, and, of course, in 
Chapter~\ref{alternating}.

In the following we comment on those groups~$G$ where some additional 
arguments were needed. This is the case if either some modular character 
table of $G$ is not (completely) available or some character tables of 
relevant subgroups of~$G$ are missing.

As before, if~$H$ is a maximal subgroup of~$G$, by $\varphi$ we denote a
$K$-character of~$H$ such that $\chi := \Ind_H^G( \varphi )$ is irrreducible.

\subsection{Induction from non-parabolic subgroups}
\label{ExceptionalMultipliersNonParabolic}
We begin with the cases where~$H/Z(G)$ is a non-parabolic subgroup 
of~$G/Z(G)$.

$G = 12_1.\PSU_4(3)$ or $12_2.\PSU_4(3)$: Let~$H$ be a maximal subgroup of~$G$
such that some $K$-character of~$H$ induces to an irreducible $K$-character 
of~$G$. Using
the known (modular) character tables of~$G$, one easily checks that this
implies $[G\colon\!H] = 126$, and thus $H = 12.\SU_4(2)$. (For a more precise
description of the structure of~$H$ see Remark~\ref{ExceptionalExplanations}.)
We also find that every modular character of~$H$ which could possibly induce
to an irreducible character of~$G$ must be liftable. In view of
Lemma~\ref{Modp} it suffices to look at the ordinary characters of~$H$.
Since these are available in GAP, the results given in Table~\ref{tab:clred1a}
are readily obtained.

$G = 2.G_2(4)$: There are three non-parabolic subgroups of $G/Z(G)$ whose
index is small enough to be possible block stabilizers. These are
$2.J_2$, $2.(\SU_3(4)\colon\!2)$, and $2.(3.\PSL_3(4)\colon\!2_3)$.
If $H = 2.(3.\PSL_3(4)\colon\!2_3)$, then $[G\colon\!H] = 2\,080$, 
and there is exactly one ordinary irreducible character of $G$, whose
degree is a multiple of this index. This degree equals
$2 \cdot 2\,080$, but $H$ does not have an irreducible character of
degree~$2$. The only candidate degree in characteristics $\ell > 0$
occurs for $\ell = 13$, and is the same as above. Since $13$ does not
divide $|H|$, this subgroup does not give rise to an imprimitive
irreducible character of~$G$.

Next, let $H = 2.(\SU_3(4)\colon\!2)$, where $[G\colon\!H] = 2016$.
The only multiple of $2\,016$ which is equal to the degree of 
an irreducible $K$-character of~$G$ is $2\,016$, and exactly one such 
character exists for $\ell = 0, 3, 7$. The subgroup $H$ is isoclinic
to but not isomorphic to $2 \times (\SU_3(4)\colon\!2)$. (GAP does
not find a possible fusion of conjugacy classes of the group
$2 \times (\SU_3(4)\colon\!2)$ into $G$.) Let~$K$ have one of the 
above characteristics. There are two $K$-characters of~$H$ of 
degree~$1$, arising from dual representations, inducing to the 
irreducible character of~$G$ of degree $2\,016$.

Finally, let $H = 2.J_2$. The only possible degree of a $K$-character 
inducing to an irreducible character of~$G$ is $14$. This can only
happen if $\ell = 0, 7, 13$. One checks with {GAP} that the two
characters of degree~$14$ induce to reducible characters of $G$ if
$\ell = 0$. This also rules out the cases $\ell = 7, 13$.

$G = 2.F_4(2)$: According to the ``Improvements to the {\sc Atlas}'' in
\cite[p.~297--327]{modATL}, the list of maximal subgroups of $2.F_4(2)$
given in \cite[p.~$170$]{ATLAS} is complete. There are only three 
non-parabolic maximal subgroups $H$ of $G$
with $|H|^2 \geq |G|$, namely $2 \times \Sp_8(2)$, 
$2.(\Omega^+_8(2)\colon\!S_3)$ 
and $2.({^3\!D}_4(2)\colon\!3)$. Let $H$ be one of these groups.

By the result of Seitz cited above (see Theorem~\ref{SameCharacteristic}), 
we may assume that  $\mbox{\rm char}(K) \neq 2$. Suppose in addition that 
$\mbox{\rm char}(K) \neq 3$. Using the ordinary character table and the 
known (modular) character tables of~$G$ (see \cite{HissF4}), we check that
no irreducible $K$-character of~$G$ has degree divisible by the index 
$[G\colon\!H]$, except possible if $H = 2.\Sp_8(2)$. In this case, however, 
the modular character tables of~$H$ are also known. It is then routine to 
rule out this possibility as well.

Finally, let $\char(K) = 3$. If $H = 2.({^3\!D}_4(2)\colon\!3)$ or
$2.(\Omega^+_8(2)\colon\!S_3)$, a $K$-character~$\varphi$ of~$H$ 
inducing irreducibly to~$G$ has degree at most~$5$. It follows (as in the 
proof of Theorem~\ref{SporadicTheorem} in the case $G = 2.\Co_1$ and 
$H = 2.(\PSU_6(2)\colon\!S_3)$) that~$\varphi$ is liftable. Since there are
no examples in characteristic~$0$, we are done for these two subgroups by
Lemma~\ref{Modp}. It remains to consider the subgroup $H = 2 \times \Sp_8(2)$.
Here, the only candidates for the character $\varphi$ have degrees $1$, $35$, $50$,
$118$, $135$, or $203$. We rule out these possibilities using the decomposition
matrix of $H$ modulo~$3$ and Lemma~\ref{ModpGeneralization}.

$G := 6.{^2\!E}_6(2)$: According to the ``Improvements to the {\sc Atlas}'' in
\cite[p.~297--327]{modATL}, the list of maximal subgroups of ${^2\!E}_6(2)$
given in \cite[p.~$191$]{ATLAS} is complete (up to a correction of the
structure of one of the maximal subgroups which does not affect its order).
If $H$ is a maximal subgroup of~$G$ such that $H/Z(G)$ is not a parabolic
subgroup of ${^2\!E}_6(2)$, then $H/Z(G)$ is one of $F_4(2)$,
$\mbox{\rm Fi}_{22}$ or $\Omega_{10}^-(2)$. By a result of Seitz 
(see Theorem~\ref{SameCharacteristic}), we may and will assume that
$\mbox{\rm char}(K) \neq 2$ in the following. We will also make use of the
ordinary character table of~$G$ computed by Frank L\"ubeck (see 
\cite{LuebeckChar}).

Suppose first that $H/Z(H) = F_4(2)$. Then $H$ is one of $6 \times F_4(2)$ or 
$3 \times 2.F_4(2)$ and $\varphi(1) \leq 4\,478$. All irreducible $K$-characters 
of these groups of degrees smaller than or equal to $4\,478$ are known (see 
\cite{WhiteF4,HissF4}). In case of $\mbox{\rm char}(K) \neq 0$, all of these 
characters are liftable, except for two characters of degree $1\,104$ for
$\mbox{\rm char}(K) = 7$, and one character of degree $1\,366$ for 
$\mbox{\rm char}(K) = 13$. Using GAP, we find that no ordinary irreducible
character of~$H$ induces to an irreducible character of~$G$. This also
rules out the cases of the liftable characters by Lemma~\ref{Modp}. The
non-liftable cases are ruled out using Lemma~\ref{ModpGeneralization}.
Suppose next that $H/Z(H) = \mbox{\rm Fi}_{22}$. Then 
$H = 6.\mbox{\rm Fi}_{22}$, the full covering group of $\mbox{\rm Fi}_{22}$
or $H = 2 \times 3.\mbox{\rm Fi}_{22}$
(GAP shows that there is no possible class fusion of $2.\mbox{\rm Fi}_{22}$
or $\mbox{\rm Fi}_{22}$ into $G$). Moreover, $\varphi(1) \leq 87$. 
This implies that $\varphi(1) \in \{ 77, 78 \}$ (see 
\cite{BBB,HiWhi,Noeske}), and thus that the elements of order~$3$ of $Z(G)$ 
lie in the kernel of~$\varphi$. This case is now easily ruled out using 
Lemmas~\ref{Modp} and~\ref{ModpGeneralization}.
Suppose finally that $H/Z(G) = \Omega_{10}^-(2)$. Then 
$\varphi(1) \leq 33$. Since $\Omega_{10}^-(2)$ has a trivial Schur 
multiplier, we have $H = 6 \times \Omega_{10}^-(2)$. By \cite{LaSe},
every irreducible $K$-character of~$H$ has degree larger than~$33$ 
if $\mbox{\rm char}(K) \neq 2$. 

\subsection{Induction from parabolic subgroups}
\label{ExceptionalMultipliersParabolic}
Here we consider the cases where $H/Z(H)$ is a parabolic subgroup of $G/Z(G)$.

$G = 2.(^{2}\!B_{2}(8))$: Every irreducible $K$-character of~$G$ whose degree 
is divisible by $65$ has $Z(G)$ in its kernel. Thus it suffices to consider 
the simple group $G/Z(G)$, and its maximal parabolic subgroup. The claimed
result follows from Harish-Chandra theory (or with GAP).

$G = 2.\Sp_6(2)$: The same arguments as above give the desired results.

$G = 12_1.\PSU_4(3)$ or $12_2.\PSU_4(3)$: There are two maximal subgroups~$H$ 
of~$G$ such that $H/Z(G)$ is a parabolic subgroup of~$G/Z(G)$. Let us first 
consider the subgroup $H$ of index $112$ in $G$, i.e., $H/Z(G) \cong 
3^4\colon\!\PSL_2(9)$. There are no faithful irreducible $K$-characters of~$G$ 
whose degree is divisible by~$112$. The analogous statement holds for the group 
$3_2.\PSU_4(3)$. The faithful irreducible characters of $3_1.\PSU_4(3)$ of 
degrees divisible by~$112$ have degree $3\cdot112 = 336$ and exist exactly in
characteristics $0, 5$ and $7$. It turns out that the inverse image of
$3^4\colon\!\PSL_2(9)$ in $3_1.\PSU_4(3)$ does not have irreducible 
characters of degree~$3$ in these characteristics. Now consider the central 
quotient $4.\PSU_4(3) = \SU_4(3)$ of~$G$. This is a finite 
group of Lie type, and thus we only have to consider Harish-Chandra induced 
representations from the parabolic subgroup $3^4\colon\!2.(2 \times \PSL_2(9))$
(see Proposition~\ref{ParabolicStabilizer}). It is not difficult to
construct the character table of this subgroup and check that only the examples
given in Table~\ref{tab:clred2a} exist.

Next, let~$H$ be the maximal subgroup of index $280$, i.e., $H/Z(G) \cong
3_+^{1+4}\colon\!2S_4$. Since $H$ is solvable, every modular irreducible
character of~$H$ is liftable. By Lemma~\ref{Modp},
it suffices to look at the ordinary characters of~$H$. The groups
$12_2.\PSU_4(3)$, $3_2.\PSU_4(3)$ and $3_1.\PSU_4(3)$ do not have any faithful 
irreducible characters of degrees divisible by $280$. For 
$G = 12_1.\PSU_4(3)$, there are such characters for $\ell \in \{0, 5, 7\}$. All of 
these have degree $3\cdot280 = 840$. There are faithful irreducible 
characters of $H$ of degree~$3$, but all of these induce to reducible
characters of~$G$. The same arguments as for the first parabolic subgroup 
of $G/Z(G)$ now give the desired results.

$G = 3.G_2(3)$: There are two conjugacy classes of maximal parabolic subgroups
of $G/Z(G) = G_2(3)$, conjugate by an outer automorphism of~$G/Z(G)$. The index
of these subgroups is $364$. There is no faithful irreducible $K$-character 
of~$G$ whose degree is divisible by $364$. Thus it suffices to consider the 
simple group $G/Z(G)$, and any one of its maximal parabolic subgroups. The only 
possible degree for an induced irreducible $K$-character of $G/Z(G)$ is $728$, 
and this can only occur for $\mbox{\rm char}(K) \in \{ 0, 7, 13 \}$. Using GAP, 
it is easy to find the example given in Table~\ref{tab:clred2a}.

$G = 2.\Omega^+_8(2)$: First, let $H \leq G$ be any maximal subgroup of 
index~$135$. The possible degrees for $\varphi$ are $21$, $24$ and $45$. If 
$\varphi(1) = 21$ or~$45$, 
then $\Ind_H^G( \varphi )(1) = 2\,835$ or $6\,075$. Any irreducible $K$-character
of~$G$ of one of these degrees has $Z(G)$ in its kernel, and we may thus look 
at the simple Chevalley group $\Omega^+_8(2)$. In this case we 
easily find the examples given in Table~\ref{tab:clred2a}. The possibility 
$\varphi(1) = 24$ only occurs if $\mbox{\rm char}(K) = 0$ or~$3$. However, there 
is no ordinary irreducible character of $H$ of degree $24$. Using MOC (see
\cite{MOC}), we can 
show that there is no irreducible $3$-modular character of degree $24$ of~$H$.
Next, let $H \leq G$ have index~$1\,575$. The only possible degree for
$\chi$ turns out to be~$1\,575$. Then~$\chi$ has $Z(G)$ in its kernel
and we may consider the parabolic subgroup $2^{1+8}_+\colon\!
(S_3 \times S_3 \times S_3)$ of $\Omega^+_8(2)$. Using GAP, one 
shows that no character of degree~$1$ of this parabolic subgroup induces to an 
irreducible character of~$G$.

$G = 2.G_2(4)$: There is no faithful irreducible $K$-character of~$G$ whose 
degree is divisible by the index of a maximal parabolic subgroup of $G/Z(G)$. 
Again, it suffices to consider the simple group $G/Z(G) = G_2(4)$. Using the 
known (modular) character tables and GAP, one easily obtains the result.

$G = 6.\Omega_7(3)$: 
Let~$H$ be a maximal subgroup of $6.\Omega_7(3)$ such that $H/Z(H)$ is a 
parabolic subgroup of $\Omega_7(3)$. The ordinary and modular character tables 
of $6.\Omega_7(3)$ are known and available in GAP. We use these to check that, 
with a few possible exceptions, no faithful irreducible character of 
$6.\Omega_7(3)$ or of $3.\Omega_7(3)$ has a degree divisible by the index of~$H$
in $6.\Omega_7(3)$. The exceptions occur 
for the largest parabolic subgroup $H = 6.(3^5\colon\!(\SU_4(2)\colon\!2))$ of 
index $364$ and only if $\mbox{\rm char}(K) \neq 2, 3$. Let $\mbox{\rm char}(K) 
\in 0, 5$; then a possible faithful character of $H$ inducing to an irreducible 
character of $6.\Omega_7(3)$ has degree $54$. But~$H$ does not have any 
faithful irreducible $K$-character of this degree. (The ordinary character 
table of~$H$ is available in GAP, and one checks that every ordinary
irreducible faithful character of~$H$ has degree divisible by~$108$.)
Next, let $H$ denote the maximal subgroup of index $364$ of $3.\Omega_7(3)$. In 
this case, a faithful character of $H$ inducing to an irreducible character 
of $3.\Omega_7(3)$ has degree $27$ and can only occur for 
$\mbox{\rm char}(K) = 0$. However, the ordinary irreducible characters of~$H$ 
of degree $27$ induce to reducible characters of $3.\Omega_7(3)$.

We are thus left
with the group of Lie type $G = 2.\Omega_7(3)$. The character tables of all
maximal subgroups of $G$ are contained in GAP. It is then routine to
obtain the results of Table~\ref{tab:clred2b}. Alternatively, one can use
Harish-Chandra theory.

$G = 6.\PSU_6(2)$: All ordinary and modular character tables of $G$ are
available in GAP. From this it follows that the central elements of 
order~$2$ of $G$ are in the kernel of every irreducible $K$-character of~$G$
which is induced from a character of a subgroup $H$ such that
$H/Z(G)$ is parabolic in $\PSU_6(2)$. Thus it suffcies to investigate the
group $3.\PSU_6(2) = \SU_6(2)$, which is a finite group of Lie type. By 
Proposition~\ref{ParabolicStabilizer} we may restrict to Harish-Chandra
induced characters. Using in addition the ordinary character tables of
the maximal parabolic subgroups of $\SU_6(2)$ as well as the modular 
character tables of their Levi subgroups, we obtain the results from 
Table~\ref{tab:clred2c}.

$G = 2.F_4(2)$: Let us first consider a maximal (parabolic) subgroup
$H$ of index $69\,615$. Then $H/Z(G)$ is the centralizer in $F_4(2)$ of
an involution $z \in F_4(2)$ in class $2A$ or $2B$. A computation of a class 
multiplication coefficient shows that there is a conjugate~$y$ of~$z$ such that 
$zy$ lies in class $3A$, respectively $3B$, of $F_4(2)$. Hence there is an 
element $t \in F_4(2)$ of order~$3$ such that $t \in \langle z, z^t \rangle$.
Lemma~\ref{2.3replace} and Corollary~\ref{spinlift} imply that neither $H/Z(G)$ 
nor~$H$ are block stabilizers of imprimitive irreducible modules of $G/Z(G)$, 
respectively~$G$. Next we consider a maximal subgroup $H$ with $H/O_2(H) \cong
\SL_2(2) \times \SL_3(2)$. There are no faithful irreducible characters of~$G$
whose degree is a multiple of $[G\colon\!H]$. (Although the $3$-modular
character table of $G$ is not completely known, the results of 
\cite[Theorem 2.2]{HissF4} are sufficient to prove the above claim even if 
$\mbox{\rm char}(K) = 3$.) Thus we may consider the Chevalley group $G/Z(G)$ 
and its parabolic subgroup $H/Z(H)$. By Theorem~\ref{SameCharacteristic}, we 
may assume that $\mbox{\rm char}(K) \neq 2$. No $7$-modular irreducible 
character of~$G/Z(G)$ has degree divisible by $[G\colon\!H]$. 
In the remaining characteristics, the only possible degree for $\varphi$ is~$3$.
(Once more, in characteristic~$3$ the incomplete results of \cite{HissF4}
suffice to see this.) The results of Table~\ref{tab:clred2c} now follow from 
Proposition~\ref{ParabolicStabilizer} and Harish-Chandra theory.

$G := 6.{^2\!E}_6(2)$: Here, we assume that $\mbox{\rm char}(K) = 0$.
The ordinary character table of~$G$ has been computed by Frank L\"ubeck
(see \cite{LuebeckChar}).
Let $H$ be a maximal subgroup of~$G$ such that $H/Z(G)$ 
is one of the four maximal parabolic subgroups of $G/Z(G) = {^2\!E}_6(2)$.
As a first step we show that the central
element of order~$2$ of $G$ is in the kernel of every $K$-character of~$H$
inducing to an irreducible $K$-character of~$G$. First, let $H$ be the largest 
of these four subgroups, so that $H/Z(H) = 2^{1+20}\colon\!\PSU_6(2)$. The 
image $\bar{H}$ of $H$ in the factor group $\bar{G} := 3.{^2\!E}_6(2)$ is equal
to the centralizer of an involution $z$. Using class multiplication coefficients
one can show that $z$ together with some conjugate generates a dihedral 
group of order $6$. Thus $\bar{H}$ and $H$ cannot be block stabilizers in 
$\bar{G}$, respectively~$G$ by Corollary~\ref{spinlift} and 
Lemma~\ref{2.3replace}. Now let $H$ be one of the remaining three maximal 
subgroups considered here. Then there is no faithful irreducible $K$-character 
of~$G$ whose degree is divisible by the index $[G\colon\!H]$. The analogous
statement holds for the factor group $\tilde{G} := 2.{^2\!E}_6(2)$ with one 
exception. This occurs for the second maximal subgroup $\tilde{H}$ of 
$\tilde{G}$, where $\tilde{H}/Z(\tilde{G}) = 
2^{8 + 16}\colon\!\Omega_8^-(2)$.
In this case, the only possible degree for $\varphi$ is $4\,480$. Using
Clifford theory it is not hard to show that $\tilde{H} = 2.(2^{8 + 16}
\colon\!\Omega_8^-(2)) = 
2.(2^{8 + 16})\colon\!\Omega_8^-(2)$ does not have an 
irreducible character of this degree. 

It suffices now to consider the group $\bar{G} = 3.{^2\!E}_6(2)$. Note that
$\bar{G}$ is a finite group of Lie type. Let $\bar{H}$ denote the image of~$H$ 
in $\bar{G}$. By Proposition~\ref{ParabolicStabilizer}, the unipotent radical
of~$\bar{H}$ is in the kernel of $\varphi$, and we may thus use Harish-Chandra
theory to obtain the results of Table~\ref{tab:clred2a}. We omit the details
(but see Chapters~\ref{InductionFromParabolicSubgroups}, \ref{Characteristic0},
\ref{ExceptionalGroups}).

\chapter{Groups of Lie type: Induction from non-parabolic subgroups}
\label{CrossCharacteristic}

\markleft{6.\ INDUCTION FROM NON-PARABOLIC SUBGROUPS}
In this chapter~$G$ is a quasisimple covering group of a simple group of 
Lie type. By the latter we understand any finite simple group which is
not isomorphic to a sporadic simple group nor to an alternating group~$A_7$ 
or~$A_n$, $n \geq 9$. Throughout this chapter we also assume that~$G/Z(G)$ 
does not have an exceptional multiplier, and that $G/Z(G)$ is not isomorphic 
to a simple group of Lie type of a different defining characteristic. The 
cases excluded here were treated in Chapter~\ref{ExceptionalMultipliers}. 
Moreover,~$K$ is an algebraically closed field of characteristic~$\ell \geq 0$,
but not equal to the defining characteristic of~$G$. In 
Proposition~\ref{TitsGroup} we have already observed that the Tits simple group 
${^2\!F}_4(2)'$ does not have any imprimitive irreducible module over~$K$. For
this reason we also exclude the Tits group from our considerations here. In 
particular,~$G$ is a finite group of Lie type. The main result of this chapter 
is the following reduction theorem.

\begin{thm} \label{ReductionToParabolic}
    Let~$G$ be a quasisimple covering group of a simple group of Lie type 
    subject to the restrictions formulated above, and let~$H$ be a maximal
    subgroup of~$G$ and a block stabilizer of an imprimitive irreducible 
    $KG$-module. Then~$H$ is a parabolic subgroup of~$G$.
\end{thm}

The proof of this theorem is given in Sections~\ref{ProofClassicalGroups} 
and~\ref{ProofExceptionalGroups}.

\section{Outline of the strategy} 
Our starting point is the observation that the index of a block stabilizer 
is at most equal to the dimension of an absolutely irreducible $KG$-module. 
Convenient upper bounds for the latter are due to Seitz.
\begin{lem}[Seitz]\label{Seitzlem}
Let~$M$ be an absolutely irreducible $KG$-module and~$B$ a Borel
subgroup of~$G$.
Then $\mbox{\rm dim}(M) \leq [G:B]c_G$ where either $c_G = 1$ or is as 
given below.
$$c_G = \left\{   \begin{array}{ll}
(q^{2}-1)/(q^{2}-q+1), & \text{\rm\ if\ } {G} \cong \SU_{n}(q),\quad
n \mbox{\rm\ odd}, \\
(q-1)/(q-\sqrt{2q}+1), & \text{\rm\ if\ } {G} \cong {^2\!B}_{2}(q),\quad
q = 2^{2m+1}, \\
(q-1)/(q-\sqrt{3q}+1), & \text{\rm\ if\ } {G} \cong {^2\!G}_{2}(q),\quad
q = 3^{2m+1}, \\
(q-1)^2/(q-\sqrt{2q}+1)^2, & \text{\rm\ if\ } {G} \cong {^2\!F}_{4}(q),\quad  
q = 2^{2m+1}.
\end{array} \right.$$
\end{lem}
\begin{prf}
This is due to Seitz \cite[Theorem 2.2]{seitz:1}.
\end{prf}

\noindent
In order to prove Theorem~\ref{ReductionToParabolic}, we may assume that~$H$ 
is a maximal subgroup of~$G$. In a first step we determine those non-parabolic 
maximal subgroups~$H$ of~$G$ for which $[G:H] \leq [G:B]c_G$, since no other~$H$ 
can be the block stabilizer of an imprimitive irreducible $KG$-module. When~$G$ 
is a classical group, we treat the various non-parabolic subgroups~$H$ according 
to their Aschbacher classes.
Given~$G$ and~$H$, we usually produce a lower bound~$b$ for the order of a Borel 
subgroup~$B$ of~$G$, and a strict upper bound~$h$ for the order of~$H$. Then 
$|H|c_G/|B| < hc_G/b$, and thus $hc_G/b \leq 1$ implies $[G:H] > [G:B]c_G$. 
In the Tables~\ref{c2}--\ref{c8}, the 
columns headed~$bc_G^{-1}$ and~$h$ give lower and upper bounds for $|B|c_G^{-1}$ 
and $|H|$, respectively. The column headed~$\tilde{H}$ describes groups 
closely related to~$H$, usually overgroups of~$H$ (in overgroups of~$G$); in
any case, $|H| \leq |\tilde{H}|$, and~$h$ is an upper bound for~$|\tilde{H}|$. 
These are derived from the following easily obtained upper bounds 
for the orders of certain classical groups:
\begin{eqnarray*}
|\GL_n( q )| & \leq & q^{n^2}; \\
|\GU_n( q )| & \leq & (q + 1)q^{n^2 - 1}; \\
|\Sp_n( q )| & \leq & q^{n(n+1)/2}; \\
|\GO_n^{\epsilon}( q )| & \leq & q^{n(n-1)/2},\quad\epsilon \in \{0, +, - \}.
\end{eqnarray*}
For those~$H$ surviving the first step, we attempt to construct an element 
$t \in G \setminus H$ such that~$t$ centralizes $H \cap {H^t}$, in which case 
Lemma~\ref{Mackey1} rules out~$H$ as possible block stabilizer. At this point 
very few possibilities for~$G$ and~$H$ remain. These are disposed of by 
character theoretic methods similar to those from 
Chapters~\ref{SectionSporadicGroups} and~\ref{ExceptionalMultipliers}. 
Here, Frank L\"ubeck's data base \cite{LuebeckWeb} containing the ordinary 
character degrees of groups of Lie type of small rank is particularly helpful.

\section{The classical groups of Lie type}
\label{ProofClassicalGroups}
For further reference, we collect the groups considered here in the 
following hypothesis.
\begin{hypo}
\label{ClassicalGroupCollection}
{\rm 
Throughout this section~$G$ denotes one of the following quasisimple 
classical groups.
\begin{itemize}
\item[(a)] $\SL_n(q)$ with $n = 2$ and $q \geq 11$ or $n \geq 3$ and 
$(n,q) \neq (3,2)$, $(3,4)$, $(4,2)$;
\item[(b)] $\SU_n(q)$ with $n \geq 3$, and $(n,q) \neq (3,2)$, $(3,3)$, 
$(4,2)$, $(4,3)$, $(6,2)$;
\item[(c)] $\Sp_{n}(q)$ with $n \geq 4$ even, and $(n,q) \neq 
(4,2)$, $(4,3)$, $(6,2)$;
\item[(d)] $\Omega_n( q )$ with~$q$ odd and $n \geq 7$ odd and $(n,q) 
\neq (7,3)$;
\item[(e)] $\Omega_n^+( q )$ with $n \geq 8$ even and $(n,q) \neq (8,2)$;
\item[(f)] $\Omega_n^-( q )$ with $n \geq 8$ even.
\end{itemize}
}
\end{hypo}

\noindent
The spin groups, which are excluded from our list, are treated as follows. 
Suppose that
$G = \Omega^\epsilon_n(q)$ is one of the orthogonal groups, and that~$H$
is a non-parabolic subgroup of~$G$. If the conclusion of 
Theorem~\ref{ReductionToParabolic} for~$G$ and~$H$ is due to the fact that
$[G:H] > [G:B]$ (notice that $c_G = 1$ in this case), then the inverse
image of~$H$ in $\Spin^\epsilon_n(q)$ is not a block stabilizer either.
If $\Spin^\epsilon_n(q) \neq \Omega^\epsilon_n(q)$ and we don't prove 
that $[G:H] > [G:B]$, we find an element $t \in G$ of odd order such 
that~$t$ centralizes $H \cap H^t$. Then the statement of 
Theorem~\ref{ReductionToParabolic} for the inverse image of~$H$ in 
$\Spin^\epsilon_n(q)$ follows from Corollary~\ref{spinlift}.

The natural module for~$G$ is denoted by~$V$. Thus, if~$G$ is
the special unitary group $\SU_n(q)$, then $V = \mathbb{F}_{q^2}^n$,
and $V = \mathbb{F}_q^n$ if~$G$ is one of the other groups.
Notice that by our restrictions on $(n,q)$, each of the classical groups
considered above has a unique natural module.

Recall from the introduction that the subgroups of~$G$ are grouped into 
Aschbacher classes. We use the definitions of these classes established 
in Kleidman-Liebeck's book~\cite{kleidman:4}, slightly adapted to our
situation, where we only consider maximal subgroups of~$G$.
In the sequel, we consider each Asch\-bacher class in turn. The respective
subgroups of~$G$ are roughly described in Tables~\ref{c2}--\ref{c8}
below. These tables are based on \cite[Tables~3.5A--3.5F]{kleidman:4},
except that we suppress most of the restrictions on the parameters given
there. Whenever such restrictions are relevant here, they will be discussed in
our proofs.

\subsection{The case $H$ is of type ${\mathcal C}_1$ but not parabolic}

In this subsection either~$G$ is one of the groups defined in
Hypothesis~\ref{ClassicalGroupCollection} or $G = \Omega_n( q )$ with 
$n \geq 5$ odd, where~$q$ can be any prime power, but $(n,q) \neq (5,2), (5,3),
(7,2), (7,3)$. Of course, the natural $n$-dimensional module for the groups
$\Omega_n(q)$ considered here is also denoted by~$V$. The reason for 
enlarging the collection of groups considered in this subsection will be 
explained below.

We call a subgroup~$H$ of~$G$ of ${\mathcal C}_1$-type, if~$H$ is the 
stabilizer in~$G$ of a subspace~$W$ of~$V$, for which we write $H = N_G(W)$.
The ${\mathcal C}_1$-type subgroup $H = N_G(W)$ is called parabolic, if~$G$ 
is linear and $W$ is any subspace of~$V$, if~$G$ is unitary or symplectic
and~$W$ is totally isotropic (i.e., $W \leq W^\perp$), or if~$G$ is orthogonal 
and $Q(v) = 0$ for every $v \in W$, where~$Q$ denotes the quadratic form 
defining~$G$.

We define a ${\mathcal C}_1^*$-type subgroup of~$G$ to be a non-parabolic 
subgroup of type~${\mathcal C}_1$. Clearly, if~$G$ is linear, then it does 
not have any ${\mathcal C}_1^*$-type subgroups. We have included the groups 
$\Omega_n(q)$ for $n \geq 5$ odd and~$q$ any prime power (subject to the
restrictions listed above), since we 
wish to deal with a certain ${\mathcal C}_7$-type subgroup of $G = \Sp_4(q)$ 
as a ${\mathcal C}_1^*$-type subgroup of~$\Omega_5( q )$ (see 
Lemma \ref{LemmaC7} and Proposition~\ref{C1Prop}), and with certain 
${\mathcal C}_8$-type subgroups of $G = \Sp_n(q)$, $q$ even, as 
${\mathcal C}_1^*$-type subgroups of $\Omega_{n+1}( q ) \cong G$ 
(see \ref{SectionC8}). 

We begin by describing the configurations we have to consider. Although the
result is well known, we sketch its proof for the convenience of the reader.

\begin{lem} \label{norm} Let~$H$ be a maximal subgroup of~$G$ of 
${\mathcal C}_1^*$-type.  Then one of the following occurs.

\mbox{\rm (1)} We have $G = \Omega^{\epsilon}_n(q)$ with $n \geq 8$ even and~$q$ 
even, and~$H$ is the stabilizer of a hyperplane~$Y$ of~$V$ with 
$\mbox{\rm dim}(Y \cap Y^\perp) = 1$. In particular, $H = \Omega_{n-1}(q)
\cong \Sp_{n-2}(q)$.

\mbox{\rm (2)} We have $G = \Omega_n(q)$ with $n \geq 5$ odd and~$q$ even, 
and either~$H$ is isomorphic to a maximal ${\mathcal C}_1^*$-type subgroup 
of $\Sp_{n-1}(q)$ through the isomorphism $G \cong \Sp_{n-1}(q)$, or~$H$ 
is the stabilizer of a non-degenerate hyperplane~$Y$ of~$V$. In the latter case, 
$H = \GO_{n-1}^\epsilon(q)$, where~$\epsilon$ depends on the
Witt index of~$Y$.

\mbox{\rm (3)} Otherwise~$H$ stabilizes a pair of complementary, mutually 
orthogonal, non-isometric, non-degenerate subspaces~$X$ and~$Y$ of~$V$. 
\end{lem}
\begin{prf}
If~$G$ is orthogonal, let~$Q$ denote the quadratic form on~$V$ defining~$G$.
Suppose that~$H$ stabilizes the proper subspace $0 \neq W \leq V$. Then it also 
stabilizes~$W^\perp$ and $\mbox{\rm rad}(W) := W \cap W^\perp$ (if~$G$ is 
orthogonal,
$\mbox{\rm rad}(W)$ is defined with respect to the polar form of~$Q$).
If $\mbox{\rm rad}(W) = 0$, i.e., if~$W$ is non-degenerate, we are in Case~(3).

Suppose first that~$q$ is odd or that~$G$ is not orthogonal. As~$H$ is maximal 
and not parabolic, it follows that~$W$ is non-degenerate. 

Now suppose that~$q$ is even, that~$G$ is orthogonal and that $\mbox{\rm rad}(W) 
\neq 0$. Put
$$U := \{ v \in \mbox{\rm rad}(W) \mid Q( v ) = 0 \}.$$
Then~$U$ is an $H$-invariant subspace of~$V$. As~$H$ is maximal and 
non-parabolic, it follows that $U = 0$. This implies that $\mbox{\rm rad}(W)$
is $1$-dimensional.

Now suppose that~$n$ is even. As~$H$ stabilizes $Y := \mbox{\rm rad}(W)^\perp$,
we are in Case~(1).

Finally suppose that~$n$ is odd. Then $\mbox{\rm rad}(V)$ is $1$-dimensional.
If $W + \mbox{\rm rad}(V) = V$, we are in the second alternative of Case~(2). 
Thus suppose that
$W + \mbox{\rm rad}(V) \neq V$. As $W + \mbox{\rm rad}(V)$ is $H$-invariant,
we find that $\mbox{\rm rad}(V) \leq W$ and $\mbox{\rm rad}(V) = 
\mbox{\rm rad}(W)$. It follows that $\mbox{\rm rad}(V) \neq W$. Now let
$$\,\bar{}: V \rightarrow \bar{V} := V/\mbox{\rm rad}(V)$$ 
denote the canonical epimorphism. This gives rise to the isomorphism 
$$\,\bar{}: \Omega(V) \rightarrow \Sp(\bar{V}).$$
As~$\bar{H}$ is a non-parabolic subgroup of~$\Sp(\bar{V})$, it follows
that $\bar{W}$ is a non-zero, proper, non-degenerate subspace of~$\bar{V}$.
This completes the proof.
\end{prf}

\noindent 
We first treat the case of odd defining characteristic.
\begin{lem}  \label{COneLin} 
Let~$H$ be a maximal subgroup of~$G$ of ${\mathcal C}_1^*$-type and suppose
that~$q$ is odd. Then there exists an element $t \in G \setminus H$ of odd 
order centralizing $H \cap H^t$.
\end{lem}
\begin{prf}
As~$q$ is odd, we are in Case~(3) of Lemma~\ref{norm}. Using the notation of 
that lemma, we choose $z \in \GL(V)$ so that it 
acts via $-I_X$ on~$X$ and via~$I_Y$ on~$Y$. Evidently~$z$ is an involution 
and $H \leq C_{G}(z)$.

We first assume that~$q \neq 3$ and
pick non-zero vectors $x \in X$ and $y \in Y$. Thus with respect to the basis 
$\{x, y \}$, the element~$z$ acts diagonally on $U := \langle x, y \rangle$ 
with eigenvalues $ \pm 1$.

If~$G$ is symplectic,~$x$ and~$y$ are necessarily isotropic and perpendicular, 
as $Y = X^{\perp}$. Consequently~$U$ is totally isotropic. Choose $x' \in X$
and $y' \in Y$ so that $(x,x')$ and $(y,y')$ are hyperbolic pairs, and put
$U' := \langle x', y' \rangle$. Then $U \oplus U'$ is a non-degenerate subspace
of~$V$. Put $P_U := N_G ( U )$. As~$U$ is totally isotropic,~$P_U$ is a parabolic 
subgroup of~$G$, and $L_U := N_{P_U}( U' )$ is a Levi subgroup of~$P_U$.
As $z \in G$ stabilizes~$U$ and~$U'$, we have $z \leq L_U$.
Now $L_U$ contains a $z_U$-invariant direct summand isomorphic to $\GL(U) \cong 
\GL_2( q )$. By Lemma~\ref{ToriInSLSU}, there are maximal tori~$T_1$ and~$T_2$ 
of $\SL(U)$ of orders $q - 1$ and $q + 1$, respectively, inverted by~$z_U$. 
As~$q > 3$, there is a non-trivial element $t_0 \in T_1 \cup T_2$ of odd order. 
Let~$t \in G$ be the element which acts as~$t_0$ 
on~$U$ and as the identity on $(U \oplus U')^\perp$. Then~$z$ inverts~$t$ and
hence $t \in C_G(H \cap {H^t})$ by Lemma~\ref{2.3replace}.

If $G$ is unitary or orthogonal we pick $x$ and $y$ so that they are non-isotropic 
with respect to the form defining $G$. Since~$X$ and~$Y$ are non-degenerate, this 
is always possible. Thus, as $x$ and $y$ are perpendicular,
$U$ is non-degenerate with respect to the form defining~$G$. 
Now if~$G$ is unitary, $\SU(U) \times \SU(U^{\perp}) \leq N_G(U)$ and $\SU(U)$ 
is $z_U$~invariant. Using Lemma~\ref{ToriInSLSU}, an argument analogous to the
one in the symplectic case shows that 
there exists a $t \in C_G(H \cap {H^t})$ with~$t$ of odd order.  

If~$G$ is orthogonal, $G = \Omega(V) \cong \Omega^\epsilon(q)$, we adjust the 
choice of~$y$ so that the order of $\SO(U)$ is twice an odd number. 
Since $q \neq 3$ this is always possible. 
Now $\SO(U) \times \SO(U^\perp) \leq \SO( V ) 
\cong \SO^\epsilon( q )$. Since $z_U \in \GO(U) \setminus \SO(U)$, 
we have $\langle z_U, \SO(U) \rangle = \GO(U)$; the latter is a dihedral group.
Let $t \in \SO(U)$ be a nontrivial element of odd order, viewed as an element
of~$\SO(V)$, by letting it act as the identity on $U^\perp$. As~$t$ has odd 
order, it lies in $\Omega(V) = G$. Since~$z$ inverts~$t$, our claim follows 
from Lemma~\ref{2.3replace}.

We now consider the case where $q=3$ and $G$ is symplectic. The group 
$\Sp_4( 3 )$ contains a unique class of non-central involutions, and a 
subgroup isomorphic to~$S_3$. We choose a non-degenerate subspace $U \leq V$ 
of dimension~$4$, whose intersections with~$X$ and~$Y$ are $2$-dimensional. 
Then~$z_U$ is a non-central involution in $\Sp( U ) \cong \Sp_4( 3 )$. Let 
$t \in G$ be an element of order~$3$ fixing~$U$ such that $t_U$ is inverted 
by~$z$ and such that~$t$ acts as the identity on~$U^\perp$. Then~$z$ inverts~$t$ 
and we are done by Lemma~\ref{2.3replace}.

We now consider the case where $q=3$ and $G = \SU_n( 3 )$, $n \geq 5$. We may 
assume that $\mbox{\rm dim}(X) \geq 2$ and choose a non-degenerate subspace 
$U \leq V$ of dimension~$3$, whose intersection with~$X$ is $2$-dimensional. 
Then~$z_U$ is an involution in $\SU( U ) \cong \SU_3( 3 )$. The latter group
contains a unique class of involutions, and a subgroup isomorphic to~$S_3$. 
Let $t \in G$ be an element of order~$3$ fixing~$U$ such that $t_U$ is 
inverted by~$z$ and such that $t$ acts as the identity on~$U^\perp$. Then~$z$ 
inverts~$t$ and we are done by Lemma~\ref{2.3replace}.

We now consider the case where $q=3$ and $G$ is orthogonal. The group 
$\Omega_5( 3 )$ has two classes of involutions, distinguished by their 
centralizer order. Let~$z_0$ be an involution in $\Omega_5( 3 )$ whose 
centralizer has order $2^6\cdot3^2$. (The centralizer of an involution in the 
other class has order $2^5\cdot 3$.) A computation of class multiplication
coefficients shows that there is a conjugate~$z_0'$ of~$z_0$ such that the
product $z_0z_0'$ has order~$3$. In particular, $z_0$ inverts an element
of order~$3$. These facts can easily be checked using the known character
table of $\Omega_5( 3 ) \cong \PSp_4( 3 )$ (see \cite{GAP4}). As we have
excluded the case $(n,q) = (5,3)$, we have $n \geq 7$, and may thus assume
that $\mbox{\rm dim}(X) \geq 4$. We now choose 
a non-degenerate subspace $U \leq V$ of dimension~$5$, whose intersection 
with~$X$ is $4$-dimensional. Now $z_U \in \Omega_5( 3 )$ since $-I_{U \cap X}$ 
is contained in $\Omega_4^+( 3 )$. Hence~$z_U$ is a non-central involution in 
$\Omega_5( 3 )$, and, by the above, inverts an element of order~$3$. Let 
$t \in G$ be an element of order~$3$ fixing~$U$ such that $t_U$ is inverted 
by~$z$ and such that~$t$ acts as the identity on~$U^\perp$. Then~$z$ inverts~$t$ 
and we are done by Lemma~\ref{2.3replace}.
\end{prf}

\noindent
Next we deal with the groups of even characteristic, where we distinguish 
various cases.
\begin{lem}
\label{Case1OfLemmanorm}
Let $G = \Omega^{\epsilon}_n(q)$ with $n \geq 8$ even and~$q$ even, and 
suppose that~$H$ is the stabilizer of a hyperplane~$Y$ of~$V$ with 
$\mbox{\rm dim}(Y \cap Y^\perp) = 1$. Then there is $t \in G \setminus H$ 
such that $t \in C_G(H \cap H^t)$.
\end{lem}
\begin{prf}
There is an involution $z \in \GO_6^\epsilon( 2 ) \leq \GO_n^\epsilon(q)$
with $H = C_G(z)$. Suppose first that $G = \Omega_n^+(q)$. Now $\GO_6^+(2) 
\cong S_8$ and there is $z' \in \GO_6^+(2)$ conjugate to~$z$ such that 
$t := zz' \in \Omega_6^+(2) \leq G$ has order~$3$. An analogous construction
works for $\Omega_n^-(q)$ since $\GO_6^-(2) \cong \SU_4( 2 )\colon\!2$.
The claim follows from Lemma~\ref{2.3replace}.
\end{prf}

\noindent
The following lemma deals with the generic case.
\begin{lem} 
\label{GenericC1qEvenCase}
Let~$G$ be one of the groups $\SU_n(q)$, $n \geq 3$,~$q$ even, $\Sp_n(q)$,
$n \geq 4$,~$q$ even, or $\Omega^\epsilon_n(q)$, $n \geq 8$ even and~$q$ even.
Suppose that~$H$ is a maximal subgroup of ${\mathcal C}_1^*$-type, fixing a
pair of non-degenerate, mutually orthogonal, complementary subspaces~$X$ 
and~$Y$.
Then there exists an element $t \in G \setminus H$ such that 
$t \in C_G(H \cap {H^t})$.
\end{lem}
\begin{prf}
For any subspace $W \leq V$, write $W' := W \otimes_{\mathbb{F}_q} 
{\mathbb{F}_{q^2}}$ and view $\GL(V)$ as a subgroup of $\GL(V')$. 
Let $\alpha \in \mathbb{F}^*_{q^2}$ be an element of order~$3$, and let 
$z \in \GL(V')$ be the element that acts as $\alpha I$ on~$X'$ and as 
$\alpha^2 I$ on~$Y'$. Then $H = C_G( z )$.

First we assume that both~$X$ and~$Y$ are at least $2$-dimensional, and pick 
$2$-dimensional subspaces~$X_0$ and~$Y_0$ of~$X$ and~$Y$, respectively such 
that $U := X_0 \oplus Y_0$ is a non-degenerate subspace of~$V$. If~$G$ is
orthogonal, we assume in addition that both~$X_0$ and~$Y_0$ are of plus-type.
(The case where this is not possible will be treated later on.)
If~$G$ is orthogonal, we choose bases of~$X_0$ and~$Y_0$ consisting of hyperbolic 
pairs. In the other cases we may choose bases of~$U$ so that the bilinear, 
respectively sesquilinear form of~$G$ restricted to~$U$ is represented 
by~$J_4$ (which, in the orthogonal case represents the polar form of the
quadratic form restricted to~$U$).
We may also assume that~$z_{U'}$ is represented by 
$$\left(
\begin{array}{cccc}
\alpha & 0 & 0 & 0  \\
0 & \alpha^2 & 0  & 0  \\
0 & 0 & \alpha^2 & 0 \\
0 & 0 & 0 & \alpha \\
\end{array}
\right).
$$ 
Now let 
$$t  = \left(
\begin{array}{cccc}
1 & 0 & 0 & 0  \\
1 & 1 & 0  & 0  \\
0 & 0 & 1 & 0 \\
0 & 0 & 1 & 1 \\
\end{array}
\right).
$$ 
Then $t \in G$ and $t = [z, (z^t)^{-1}]$ and thus $t \in \langle z,z^t \rangle$.
Our claim follows from Lemma~\ref{2.3replace}.

%%%%%%%%%%%%%%%%%%%%%%%%%%%%%%%%%%%%%%%%%%%%%%%%%%%%%%%%%%%%%%%%%%%%%%%%%%%%%%%
%%%%%%%%%%%%%%%%%%%%%%%%%%%%%%%%%%%%%%%%%%%%%%%%%%%%%%%%%%%%%%%%%%%%%%%%%%%%%%%
%%
%% Der folgende Absatz ist neu (17.07.2012)
%%
%\marginpar{The following section is new (17.07.2012)}
Next assume that~$G$ is unitary and that~$X$ is $1$-dimensional. As~$X$ is
non-degenerate and $\mbox{\rm dim}(Y) \geq 2$, we can find a non-degenerate 
subspace $Y_0$ of $Y$ of dimension~$2$, such that the Gram matrix of the form
defining~$G$, restricted to $U := X \oplus Y_0$ equals~$J_3$ (with respect to
a suitable basis of~$U$). With respect to the same basis, we define $t_0 \in 
\SU(U)$ by the matrix
$$t_0 := \left(
\begin{array}{ccc}
1 & 0 & 0  \\
1 & 1 & 0  \\
\alpha & 1 & 1 
\end{array}
\right).
$$
Let $t \in G$ be such that $t_U = t_0$ and $t_{U^{\perp}} = I_{U^{\perp}}$.
Again $t = [z,(z^t)^{-1}] \in \langle z,z^t \rangle$, so that our claim
follows from Lemma~\ref{2.3replace}.
%%
%% Ende des neuen Absatzes.
%%
%%%%%%%%%%%%%%%%%%%%%%%%%%%%%%%%%%%%%%%%%%%%%%%%%%%%%%%%%%%%%%%%%%%%%%%%%%%%%%%
%%%%%%%%%%%%%%%%%%%%%%%%%%%%%%%%%%%%%%%%%%%%%%%%%%%%%%%%%%%%%%%%%%%%%%%%%%%%%%%

Finally assume that~$G$ is orthogonal and that~$X$ is $2$-dimensional and 
totally anisotropic, i.e.\ of minus-type. As $\mbox{\rm dim}(Y) \geq 6$,
there is a non-degenerate $2$-dimensional subspace $Y_0 \leq Y$ of plus-type.
Put $U := X \oplus Y_0$. We identify $X$ with $\mathbb{F}_{q^2}$, equipped 
with the norm form. Then $\Omega(X)$ equals the set of right 
multiplications by elements of norm~$1$. Let $x \in \mathbb{F}^*_{q^2}$ be an 
element of order~$q+1$. Then multiplication by~$x$ is represented by the matrix
$$z_0 :=\left(
\begin{array}{cc}
  0 & 1 \\
  1 & a \\
\end{array}
\right)$$
with respect to the $\mathbb{F}_q$-basis $\{ 1, x \}$ of $\mathbb{F}_{q^2}$.
Choose a hyperbolic pair $\{y_1,y_2\}$ of~$Y_0$, and an ordered basis
of~$V$ of the form $\{ y_1, x_1, x_2, y_2, \ldots \}$, where the dots indicate
a basis of $U^\perp$, and where $x_1, x_2 \in X$
correspond to $1, x$, respectively, in our identification of $X$~with 
$\mathbb{F}_{q^2}$. Let $z \in G$ be the element represented by the matrix
$$z := \left( \begin{array}{cccc}
1 & 0 & 0 & 0\\
0 & z_0 & 0 & 0\\
0 & 0 & 1 & 0 \\
0 & 0 & 0 & I 
\end{array}\right)
$$
with respect to this basis.
Then $H = N_G( \langle z \rangle ) = \langle C_G( z ), s \rangle$, where 
$s \in G$ is an involution swapping the basis elements $x_1, x_2$ of~$X$ 
and leaving fixed~$y_1$ and~$y_2$. Now let $t$ be the element
$$t := \left( \begin{array}{cc}
t_0 & 0 \\
0 & I 
\end{array} 
\right),
$$
where~$t_0$ is the matrix
$$t_0  = \left(
\begin{array}{cccc}
1 & 0 & 0 & 0  \\
1 & 1 & 0 & 0  \\
1 & 0 & 1 & 0 \\
a & 1 & 1 & 1 \\
\end{array}
\right) 
$$
with respect to the basis $\{ y_1, x_1, x_2, y_2 \}$ of~$U$. Then~$t \in G$.

The fixed space of~$z$ and of $z^t$ equals $\langle y_1, Y_1 \rangle$, 
where~$Y_1$ denotes the orthogonal complement of~$Y_0$ in~$Y$.
Hence $H \cap H^t = N_G( \langle z \rangle ) \cap N_G( \langle z^t \rangle )$
stabilizes $\langle y_1, Y_1 \rangle$ and in turn also 
$\langle y_1, Y_1 \rangle^\perp = \langle y_1, X \rangle =: U_0$. Thus an element 
$g \in H \cap H^t$ has the form
\begin{equation}
\label{ElementOfHCapHt}
g  = \left(
\begin{array}{cc}
A & B \\
C & D
\end{array} \right),
\end{equation}
where $A$ is a $4 \times 4$-matrix,~$B$ a $4 \times (n - 4)$-matrix whose first
three rows are zero, and~$C$ is an $(n-4) \times 4$-matrix whose last three
columns are zero. 

Now 
$$g^t  = \left(
\begin{array}{cc}
t_0^{-1}At_0 & t_0^{-1}B \\
Ct_0 & D
\end{array} \right).
$$
The properties of~$C$ and~$D$ mentioned above imply $t_0^{-1}B = B$ and
$Ct_0 = C$.
Our claim follows once we can show that $t_0$ commutes with the matrices~$A$
arising in~(\ref{ElementOfHCapHt}). 

We have already observed above that $H \cap H^t$ stabilizes~$U_0$. It also
stabilizes $X \cap Xt = \langle x_1 + x_2 \rangle$, and hence 
$x_1 + x_2$. As $H \cap H^t$ also fixes $y_1$, it follows that
$(H \cap H^t)_{U_0} = \langle s_{U_0} \rangle$. As~$H$ stabilizes~$X$ 
and~$Y = X^\perp$ the matrices~$A$ arising in~(\ref{ElementOfHCapHt}) are
of the form
$$A  = \left(
\begin{array}{cccc}
1 & 0 & 0 & 0  \\
0 & 0 & 1 & 0  \\
0 & 1 & 0 & 0 \\
b & 0 & 0 & 1 \\
\end{array}
\right) 
$$
or
$$A  = \left(
\begin{array}{cccc}
1 & 0 & 0 & 0  \\
0 & 1 & 0 & 0 \\
0 & 0 & 1 & 0  \\
b & 0 & 0 & 1 \\
\end{array}
\right) 
$$
for suitable elements~$b \in \mathbb{F}_q$. Such elements obviously commute 
with~$t_0$. This completes our proof.
\end{prf}

\noindent
We finally consider the orthogonal groups of odd dimension and even characteristic.

\begin{lem}\label{SOnOddEven}
Let $G = \Omega_n(q)$ with $n \geq 5$ odd and~$q$ even, and suppose that~$H$ is the
stabilizer of a non-degenerate hyperplane~$Y$ of~$V$. Then there is 
$t \in G \setminus H$ such that~$t$ centralizes $H \cap H^t$.
\end{lem}
\begin{prf}
Notice that~$Y$ is a complement to $X := V^\perp$. Choose a non-degenerate
subspace~$Y_0$ of~$Y$ of dimension~$2$ of plus-type, i.e.,~$Y_0$ contains 
isotropic vectors. Let~$Y_1$ denote the orthogonal complement of~$Y_0$ 
in~$Y$, and put $U := X \oplus Y_0$.

Choose a hyperbolic pair $\{ y, y' \}$ of~$Y_0$ and a non-zero element 
$x \in X$. Let~$t_0$ be the element of~$\Omega(U)$ defined by the following matrix 
with respect to the basis $\{ x, y', y \}$ of~$U$:
$$\left(
\begin{array}{ccc}
1 & 0 & 0  \\
1 & 1 & 1  \\
0 & 0 & 1 
\end{array}
\right).
$$
Extend~$t_0$ to all of~$V$ by letting it act as the identity on~$Y_1$. In 
particular,~$t$ fixes~$y$ and~$x$. Now $H \cap H^t$ stabilizes $Y \cap Yt$.
As $y \in Y \cap Y^t$, we have $Y \cap Yt = Y_1 \oplus \langle y \rangle$,
with $\langle y \rangle = \mbox{\rm rad}(Y \cap Yt)$. In particular, the 
elements of $H \cap H^t$ fix~$y$ and~$x$.

Consider the basis $\{ x, y', y, y_3, \ldots , y_n \}$ of~$V$, where
$\{ y_3, \ldots , y_n \}$ is a basis of~$Y_1$. Let $g \in H \cap H^t$, 
written with respect to this basis as
$$g = \left(
\begin{array}{cc}
A & B  \\
C & D 
\end{array}
\right),
$$
where $A$ is a $3\times3$-matrix. Then
$$g^t = \left(
\begin{array}{cc}
t_0^{-1}At_0 & t_0^{-1}B  \\
Ct_0 & D 
\end{array}
\right).
$$
As~$g$ fixes~$y$ and~$x$, the first and last row of~$B$ consist of zeroes,
and thus $t_0^{-1}B = B$. As~$g$ stabilizes $\langle y \rangle \oplus Y_1$, 
the first two columns of~$C$ consist of zeroes, and thus $Ct_0 = C$.

Using the fact that~$g$ fixes~$x$ and~$y$, and applying the invariance of 
the polar form to the pair $(y,y')$, we see that~$A$ is of the form
$$\left(
\begin{array}{ccc}
1 & 0 & 0  \\
0 & 1 & \beta  \\
0 & 0 & 1 
\end{array}
\right),
$$
for some~$\beta \in \mathbb{F}_q$. It follows that~$D$ commutes with~$t_0$,
proving our assertion.
\end{prf}

\noindent 
We summarize the results of this subsection.

\begin{prop}
\label{C1Prop}
Let~$G$ be one of the following groups.

{\rm (a)} A group as in Hypothesis~\ref{ClassicalGroupCollection}.

{\rm (b)} A group $G = \Spin^\epsilon_n( q )$, $n \geq 5$, $q$ any prime power
$\epsilon \in \{ -1, 0, 1 \}$ and $(n,q) \neq (5,2), (5,3), (7,2 ), (7,3)$.

Let~$H$ be a non-parabolic maximal subgroup of~$G$ of ${\mathcal C}_1$-type.
(If $G = \Spin^\epsilon_n( q )$, this means that~$H$ is the inverse image of 
a maximal subgroup of $\Omega^\epsilon_n(q)$ of ${\mathcal C}_1$-type.)
Then~$H$ is not the stabilizer of an imprimitivity decomposition of an 
irreducible $KG$-module.
\end{prop}
\begin{prf}
The possibilities for~$H$ are described in Lemma~\ref{norm}. If~$q$ is odd 
and~$G$ is a spin group, our claim follows from Lemma~\ref{COneLin} in 
conjunction with Corollary~\ref{spinlift}. The other cases are treated in 
Lemmas~\ref{COneLin}, \ref{Case1OfLemmanorm}, \ref{GenericC1qEvenCase} and
\ref{SOnOddEven}.
\end{prf}

\subsection{The case $H$ is of type ${\mathcal C}_2$}

In this subsection we consider the case where~$H$ is a ${\mathcal C}_2$-type 
subgroup of~$G$. By definition, this means that there is a direct sum 
decomposition $V = V_1 \oplus \cdots \oplus V_k$ with $k \geq 2$ such
that the factors~$V_i$ are permuted by~$G$. The relevant bounds for~$|H|$ 
and $|B|c_G^{-1}$ are given in Table~\ref{c2}. The parameter~$k$ in 
Table~\ref{c2} is the same as above, and $n = mk$. We also assume that~$H$ 
is a maximal subgroup of~$G$.

\begin{lem} 
\label{C2Type1}
Let~$G$ be one of the groups of Hypothesis~\ref{ClassicalGroupCollection}.

{\rm (a)} Let $h$ and $bc_G^{-1}$ be as in Table \ref{c2}, and suppose that
$k = n$. Then $hc_G/b < 1$ in any of the following cases:

\begin{itemize}
\item[(i)] The group~$G$ is linear or unitary and $n \geq 5$,~$q \geq 3$
or $q = 2$ and $n \geq 7$.
\item[(ii)] The group~$G$ is symplectic or orthogonal and $n \geq 6$, 
$q \geq 5$, or $q = 3, 4$ and $n \geq 9$ or $q = 2$ and $n \geq 18$.
\end{itemize}

{\rm (b)}
If $H$ is a maximal subgroup of~$G$ of ${\mathcal C}_2$-type such that~$H$ 
stabilizes a decomposition of~$V$ into~$n$ blocks, then $H$ is not the
block stabilizer of an imprimitive irreducible $KG$-module.
\end{lem}

\begin{prf}
(a) In case $G = \SU_n(q)$ our estimate is  
$$\frac{hc_G}{b} \leq \frac{(q+1)^{n+1}n!}{(q^2-q+1)q^{n(n-1)/2}(q^2-1)^{(n - 3)/2}}
=: f(n,q).$$ 
For fixed $q \geq 2$, the function $f(n,q)$ is decreasing in~$n$. This is easily
verified by considering the quotient $f(n+1,q)/f(n,q)$ of two consecutive terms.
Now $f(5,q) < 1$ for all $q \geq 3$ and thus $f(n,q) < 1$ for all $q \geq 3$ 
and $n \geq 5$. Also, $f(7,2) < 1$ and thus $f(n,2) < 1$ for all $n \geq 7$.
The case $G = \SL_n( q )$ is treated similarly. We omit the details.

Suppose now that~$G$ is orthogonal. Here, $h = 2^nn!$, and, using the smallest
of the three lower bounds for the orders of the Borel subgroups, our estimate is
$$hc_G/b < \frac{2^{n+1}n!}{(q-1)^{(n-1)/2}q^{(n-1)^4/4}} =: f(n,q).$$
For fixed $q \geq 3$ the function $f(n,q)$ is decreasing in~$n$ for
$n \geq 6$ (in fact for $n \geq 4$). Since $f( 6, q ) < 1$ for all $q \geq 5$,
we have $f( n, q ) < 1$ for all $n \geq 6$ and $q \geq 5$. Also, $f( 9, 3 ),
f( 9, 4 )  < 1$ and and hence $f( n, q ) < 1$ for all $n \geq 9$ and $q = 3, 4$.
Finally, $f( 18, 2 ) < 1$ and our claim follows. The case $G = \Sp_n( q )$ is
treated similarly.

(b) We begin with the case $G = \SU_n( q )$. For~$n = 3, 4$ we use the exact 
values for~$|H|$, $|B|$ and~$c_G$ to prove that $|H|c_G/|B| < 1$ for all values 
of~$q$ not excluded in Hypothesis~\ref{ClassicalGroupCollection}. By~(a), it 
remains to consider the cases $q = 2$ and $3 \leq n \leq 6$, of which only the 
case $(n,q) = (5,2)$ is not excluded 
in Hypothesis~\ref{ClassicalGroupCollection}. The latter case is ruled out using 
the character table of $\SU_5( 2 )$ available in GAP~\cite{GAP4}. This concludes 
the case where~$G$ is unitary. The case of $G = \SL_n( q )$ is treated similarly. 
Note that the case $G = \Sp_n( q )$ does not occur, since~$m$ must be even in 
this case.

Suppose now that~$G$ is orthogonal. This case only occurs when~$q$ is odd 
(see \cite[Table 3.5]{kleidman:4}), which we assume from now on. By~(a)(ii),
we only have to consider the cases $q = 3$ and $n = 7, 8$, the first of
which is excluded in Hypothesis~\ref{ClassicalGroupCollection}. The cases 
$q = 3$ and $n = 8$ are treated by using
the exact values for the orders of~$B$ and~$H$ (notice that $G = \Omega^+_8(3)$
in this case; see \cite[Table 3.5.F]{kleidman:4}).
\end{prf}

\begin{lem} 
\label{C2LargeBlocks}
Let~$G$ be one of the groups of Hypothesis~\ref{ClassicalGroupCollection}.
If~$H$ is a maximal subgroup of~$G$ of ${\mathcal C}_2$-type such that~$H$ 
stabilizes a decomposition of~$V$ into~$k$ blocks with $k > 2$, then $H$ is 
not the block stabilizer of an imprimitive irreducible $KG$-module.
\end{lem}
\begin{prf}
We fix $k \geq 3$. First assume that~$G$ is one of the groups listed in
Lemma~\ref{C2Type1}(a). We show $hc_G/b < 1$ by induction on~$m$. By this
lemma, this assertion is true for $m = 1$. Now assume that the assertion 
holds for all $1 \leq m \leq j$ for some $j \geq 1$. Setting $e = 1$ if~$G$ 
is linear or unitary and $e = 2$, 
otherwise, we see that increasing~$m$ by~$1$ increases $h$ by at most a factor 
of 
$$\Delta_h := q^{(2km+2k)/{e}}.$$ 
On the other hand, the estimate for $b$ increases by at least a factor of 
$$\Delta_b:=(q-1)^{k/2}q^{k(2km + k - 2)/{2e}}.$$ 
One easily checks that $\Delta_h/\Delta_b \leq 1$ for all $k \geq 3$ and 
$m \geq 2$, implying our claim.

To finish the proof, we have to consider the groups satisfying our hypothesis
($k \geq 3$), that are excluded in Lemma~\ref{C2Type1}(a). First assume
that $G = \Sp_n( q )$, Since $k \geq 3$, we have $n \geq 6$. If~$n = 6$, only
$k = 3$, $m = 2$ is possible. In this case, $h/b < 1$ for all $q \geq 3$.
So we are left with the case $q = 2$, $6 \leq n \leq 18$. Here,~$h$ is 
smallest for $m = 2$. In this case, 
$$h/b \leq \frac{ (n/2)!}{2^{n(n-6)/4}} < 1$$
for all $n \geq 8$. As $\Sp_6(2)$ is excluded by 
Hypothesis~\ref{ClassicalGroupCollection}, we are left with the case 
$G = \Sp_8(2)$ and $m = 2$, which is ruled out using the exact value 
for~$|H|$.

Now assume that~$G$ is not symplectic. By Lemma~\ref{C2Type1}(b), we may 
assume that $m \geq 2$, and thus $n \geq 6$. As $\SU_6( 2 )$ is excluded by
Hypothesis~\ref{ClassicalGroupCollection}, we are left with $G = \SL_6(2)$ 
or~$G$ orthogonal. The former case yields $h/b < 1$. 
Now suppose that $G$ is orthogonal with $n \geq 8$ even. Here, the bound for~$h$ 
is largest if $m = 2$. Using the smallest valued for~$b$, we obtain
$$\frac{h}{b} \leq \frac{2(n/2)!}{q^{n(n - 6)/4}} < 1$$
for all $n \geq 10$. Now $\Omega^+_8(2)$ is excluded by
Hypothesis~\ref{ClassicalGroupCollection}, and $\Omega^-_8(2)$ is ruled
out using its ordinary character table.
If~$G$ is orthogonal and~$n$ odd, we
only have to consider the case $G = \Omega_9( 3 )$, and $m = k = 3$, which 
again yields $h/b < 1$. The proof is complete.
\end{prf}

\begin{table}[h] 
 \caption{Bounds for ${\mathcal C}_2$-type groups \label{c2}}
    \begin{center}
    \renewcommand{\arraystretch}{1.24}
     \begin{tabular}{||c|c|c|c||}
      \hline
      $G$ & $bc_G^{-1}$ & $\tilde{H}$ & $h$ \\
      \hline
      \hline
      $\SL_n(q)$ & $q^{n(n-1)/2}(q-1)^{n-1}$ & $\GL_m(q) \wr S_k$ & $q^{km^2}k!$ \\
      \hline
      $\SU_n(q)$ & $\frac{q^2-q+1}{q^2-1}q^{n(n-1)/2}(q^2-1)^{[n/2]}(q+1)^{-1}$ 
                 & $\GU_m(q) \wr S_k$ 
      & $(q+1)^kq^{k(m^2-1)}k!$ \\
       $n \geq 3$ & & $\GL_{n/2}(q^2).2$ & $2q^{n^2/2}$ \\
      \hline
      $\Sp_n(q)$ & $q^{n^2/4}(q-1)^{n/2}$ & $\Sp_m(q) \wr S_k$  
                 & $q^{km(m + 1)/{2}}k!$  \\
      $n \geq 4$ & & $\GL_{n/2}(q).2$ &  $2q^{n^2/4}$\\
      \hline
      $\Omega_n(q)$ & $\frac{1}{2}q^{(n-1)^2/4}(q-1)^{(n-1)/2}$ &
      $\GO_m(q) \wr S_k$ & $2^kq^{km(m - 1)/{2}}k!$ \\
       $2 \nmid qn $ & & $\GO_1(p) \wr S_n$ & $2^nn!$ \\
      \hline
      $\Omega_n^+(q)$ & $\frac{1}{2}q^{n(n-2)/4}(q-1)^{n/2}$ &
      $\GO_m^{\epsilon}(q) \wr S_k$ & $2^kq^{km(m - 1)/{2}}k!$ \\
      $ n \geq 8$ & & $\GO_1(q) \wr S_n$ & $2^nn!$ \\
      & & $\GL_{n/2}(q).2$  & $2q^{n^2/4}$ \\
      & & $\GO_{n/2}(q)^2$ &  $q^{{n(n-2)}/{4}}$\\
      \hline
      $\Omega_n^-(q)$ & $\frac{1}{2}q^{n(n-2)/4}(q-1)^{(n-2)/2}(q+1)$ &
      $\GO_m^{\epsilon}(q) \wr S_k$ & $2^kq^{km(m-1)/{2}}k!$ \\
      $ n \geq 8$ & & $\GO_1(q) \wr S_n$ & $2^nn!$ \\
      & & $\GO_{n/2}(q)^2$ & $q^{{n(n-2)}/{4}}$ \\
      \hline
      \end{tabular}
      \end{center}
      \end{table}

\begin{lem} 
\label{C2tExistence}
Let~$G$ be one of the groups of Hypothesis~\ref{ClassicalGroupCollection}.
If $H$ is a maximal subgroup of $G$ of ${\mathcal C}_2$-type such that $H$ 
stabilizes a decomposition of the natural module of $G$ into two blocks, then 
there exists an element $t \in G \setminus H$ such that $t \in C_G(H \cap H^t)$. 
Moreover, if~$q$ is odd, there exists such a~$t$ of odd order.
\end{lem}

\begin{prf}
We first observe that $H$ stabilizes a pair of complementary and isometric 
subspaces $X$ and $Y$ of $V$ and moreover there exists an involution 
$s \in \GL(V)$ which normalizes $H$ and interchanges the two subspaces $X$ 
and $Y$. To prove our lemma we achieve the hypotheses of 
Lemma \ref{2.3replace}. 

With respect to a suitable basis $\{x_1,\dots,x_m\}$ for $X$ and 
$\{y_1,\dots,y_m\}$ for $Y$ we choose 
$$s=  \left(
\begin{array}{cc}
  0 & I \\
  I & 0 \\
\end{array}
\right),$$ 
unless $q$ is odd, $G$ is symplectic, and $X$ and $Y$ are totally isotropic,
in which case we choose 
$$s_1 =  \left(
\begin{array}{cc}
  0 & I \\
  -I & 0 \\
\end{array}
\right).$$

If $q = 2, 3$, we put $\mathbb{F}' := \mathbb{F}_{q^2}$. Otherwise, let
$\mathbb{F}' := \mathbb{F}_q$. For any subspace $W \leq V$, write 
$W' := W \otimes_{\mathbb{F}_q} \mathbb{F}'$ and view $\GL(V)$ as a 
subgroup of $\GL(V')$. 
For $q \neq 2, 3$ we choose $\alpha$ to be a primitive element in the prime
field of~$\mathbb{F}_q$. For~$q = 2, 3$, we choose $\alpha \in \mathbb{F}'$
of order~$3$ and~$8$, respectively. Then $\alpha \neq \alpha^{-1}$.

Let $z \in \GL(V')$ be the element that acts as $\alpha I_{X'}$ on~$X'$ and 
as $\alpha^{-1} I_{Y'}$ on~$Y'$ i.e., with respect to any pair of bases of~$X'$ 
and~$Y'$ we have
$$z=  \left(
\begin{array}{cc}
  \alpha I & 0 \\
  0 & \alpha^{-1}I  \\
\end{array}
\right).$$
Evidently $H \leq \langle C_G(z),s \rangle $, respectively  
$H \leq \langle C_G(z),s_1 \rangle $ and so we have achieved part of 
the hypothesis of Lemma \ref{2.3replace}. To find the required element 
$t \in C_G(s)$, respectively $t \in C_G(s_1)$, we need to further analyze the 
centralizers of $s$, respectively $s_1$ in $G$, for our possible groups $G$.

Our strategy is to fix $X$, $Y$ and $s$, respectively $s_1$ and to adjust the 
form defining $G$ in such a way as to guarantee that our $s$ respectively $s_1$ 
and $t$ lie in $G$ and $X,Y$ are non-degenerate or totally isotropic. 

We observe that in the $x_i,y_i$ basis above, the elements of $C_{\GL(V')}(s)$ 
have the form 
$$\left(
\begin{array}{cc}
  A & B \\
  B & A \\
\end{array}
\right),$$ 
and the elements of $C_{\GL(V')}(s_1)$ have the form 
$$\left(
\begin{array}{cc}
  A & B \\
  -B & A \\
\end{array}
\right).$$ 

Next we exhibit elements $t \in C_G(s)$, respectively $t \in C_G(s_1)$ such 
that $t \in \langle z, z^t \rangle$ for the various situations described above. 

If $G$ is linear let
$$t =  \left(
\begin{array}{cc}
  I & N_{m,1} \\
  N_{m,1} & I \\
\end{array}
\right),$$ 
where $N_{m,1}$ is the $(m \times m)$-matrix with entry~$1$ at position $(m,1)$,
and zeroes, elsewhere.
As $N_{m,1}^2 = 0$ we have 
$$t^k = \left(
\begin{array}{cc}
  I & kN_{m,1} \\
  kN_{m,1} & I \\
\end{array}
\right),$$ 
for all integers~$k$ and hence~$t$ is an element of prime order, equal to the 
characteristic of~$G$. We also find that $[z^{-1}, z^t] = t^k$ with 
$k = 2 - \alpha^2 - \alpha^{-2}$. As $2 - \alpha^2 - \alpha^{-2}$ is not equal 
to~$0$ and lies in the prime field of $\mathbb{F}'$, we conclude that 
$t \in \langle z, z^t \rangle$. Clearly,~$t$ commutes with~$s$, and thus the 
first two hypotheses of Lemma~\ref{2.3replace} are satisfied.

To deal with the unitary case, we assume that $n = 2m \ge 4$, and define the 
$(m \times m)$-matrix
$$N(a,b) =  \left(
\begin{array}{ccc}
  0_{m-2} & 0 & 0\\
  0 & a & b \\
 0 & -b^q & a^q \\
\end{array}
\right),$$ 
where $a,b \in \mathbb{F}_{q^2}$ such that $a + a^q = 0$ and 
$a^{q+1} + b^{q+1} = 0$, and where $0_{m-2}$ denotes the
$(m -2) \times (m - 2)$ zero matrix. One checks that $N(a,b)N(a^q,b^q)^T = 0 = 
N(a,b) + N(a^q,b^q)^T$. In particular 
this implies that $N(a,b)^2 = -N(a,b)N(a^q,b^q)^T= 0$.
We put 
$$t = \left(
\begin{array}{cc}
  I & N(a,b) \\
  N(a,b) & I \\
\end{array}
\right)$$
and find $t \in \langle z, z^t \rangle$, as in the case $G = \SL_n(q)$.

Finally we  pick the defining unitary form of $G = \SU_n(q)$ to be represented 
by the identity matrix respectively the matrix~$s$ depending on whether~$X$ 
and~$Y$ are to be non-degenerate respectively totally isotropic subspaces 
of~$V$, and observe that in both cases $s, t \in G$. 

If $G = \Sp_n( q )$ is symplectic, $n  = 2m \ge 4$ and~$X$ and~$Y$ are 
non-degenerate, then we define~$t$ as in the linear case, that is 
$$t = \left(
\begin{array}{cc}
  I & N_{m,1} \\
  N_{m,1} & I \\
\end{array}
\right).$$   

Here we pick the defining symplectic form of~$G$ to be represented by the 
matrix
$$ \left(
\begin{array}{cc}
  \tilde{J} & 0 \\
   0 & \tilde{J} \\
\end{array}
\right).$$ 
To see that $t \in G$ we observe that 
\begin{eqnarray*}
\left(
\begin{array}{cc}
  I & N_{m,1} \\
  N_{m,1} & I \\
\end{array}
\right) \left(
\begin{array}{cc}
  \tilde{J} & 0 \\
   0 & \tilde{J} \\
\end{array}
\right)\left(
\begin{array}{cc}
  I & N_{m,1} \\
  N_{m,1} & I \\
\end{array}
\right)^T & = & \\
\left(
\begin{array}{cc}
  \tilde{J} & N_{m,1} \tilde{J}  \\
  N_{m,1}\tilde{J} & \tilde{J} \\
\end{array}
\right)\left(
\begin{array}{cc}
  I & N_{m,1}^T \\
  N_{m,1}^T & I \\
\end{array}
\right) & = & \\
\left(
\begin{array}{cc}
  \tilde{J} + N_{m,1}\tilde{J}N_{m,1}^T & \tilde{J}N_{m,1}^T + N_{m,1}\tilde{J} \\
  N_{m,1}\tilde{J} + \tilde{J}N_{m,1}^T &  N_{m,1}\tilde{J}N_{m,1}^T  +    \tilde{J}\\
\end{array}
\right). 
\end{eqnarray*}
As $N_{m,1}\tilde{J}N_{m,1}^T = 0 = \tilde{J}N_{m,1}^T + N_{m,1}\tilde{J}$ 
we have our desired conclusion. 

If $G = \Sp_n(q )$ is symplectic, $n = 2m \ge 4$ and~$X$ and~$Y$ are totally 
isotropic, then we define~$t$ as 
$$t = \left(
\begin{array}{cc}
  I & N_{m,1} \\
  -N_{m,1} & I \\
\end{array}
\right).$$
The commutator $[z^{-1}, z^t]$ shows once more that $t \in 
\langle z, z^t \rangle$. Again we check that $s_1, t \in G$, where the form 
defining~$G$ is represented by the matrix $\tilde{J}_{2m}$.

If $G$ is orthogonal and $n = 2m \ge 8$, then we define~$t$ as 
$$t = \left(
\begin{array}{cc}
  I & N \\
  N & I \\
\end{array}
\right),$$ 
where $N = N_{m-1,1} - N_{m,2}$. Now $N^2 = 0$ as $m \geq 4$, and we
conclude as before that $t \in \langle z, z^t \rangle$.

To deal with the case that~$X$ and~$Y$ are totally isotropic, we define the 
quadratic form $Q$ on~$V$ via the formula $Q(v) = \sum_{i = 1}^m a_ib_{m-i+1}$ 
for $v \in V$ with $v = \sum_{i=1}^m a_i x_i + b_i y_i$. To deal with the case
that~$X$ and~$Y$ are non-degenerate, we 
define the quadratic form $Q$ of plus-type on $V$ via the formula 
$Q(v) = \sum_{i = 1}^m a_ia_{m-i+1} + b_ib_{m-i+1}$
for $v \in V$ with $v = \sum_{i=1}^m (a_ix_i + b_i y_i)$. (As $X, Y$ are
both totally isotropic or both non-degenerate, the form~$Q$ is necessarily 
of plus-type.)

Then $s \in G$. Moreover for our $t$ and $v \in V$ we have $Q(vt) = Q(v)$. 
Indeed,
$vt = (a_1 + b_{m-1})x_1 + (a_2 - b_m)x_2 + a_3x_3 + \cdots + a_mx_m + 
(b_1 + a_{m-1})y_1 + (b_2 - a_m)y_2 + b_3y_3 + \cdots + b_my_m$.
Now if $X, Y$ are non-degenerate, we have
$Q(vt) = (a_1 + b_{m-1})a_m + (a_2 - b_m)a_{m-1} + \sum_{i = 3}^m a_ia_{m-i+1}
+ (b_1 + a_{m-1})b_m + (b_2 - a_m)b_{m-1} + \sum_{i = 3}^m b_ib_{m-i+1}  = 
\sum_{i = 1}^m a_ia_{m-i+1} + b_ib_{m-i+1} = Q(v)$. 
The calculation for the case $X$, $Y$ totally isotropic is very similar. 
Our claim is now proved. 

Put $C := C_G( z )$. We finally show that $C^t \cap Cs = \emptyset$ in 
each case. The elements of~$C$ are of the form
$$\left(
\begin{array}{cc}
A & 0 \\
0 & B
\end{array}
\right),
$$
and thus the elements of $Cs$ are of the form
$$\left(
\begin{array}{cc}
0 & A \\
B & 0
\end{array}
\right).
$$
Now 
$$t =
\left(
\begin{array}{cc}
I & N \\
N & I
\end{array}
\right),
$$
for a matrix $N$ with $N^2 = 0$, and thus the elements of $C^t$ are of the 
form
$$\left(
\begin{array}{cc}
A - NBN & AN - NB \\
BN -AN  & B - NAN 
\end{array}
\right).
$$
If such an element also lies in $Cs$, then $A - NBN = 0$ and $B - NAN = 0$.
This implies $AN = NA = 0$ and $BN = NB = 0$, a contradiction. A similar
calculation shows $C^t \cap Cs_1 = \emptyset$ for the~$t$ we have chosen in
this case. Thus all the 
hypotheses of Lemma~\ref{2.3replace} are satisfied, which, together with
Corollary~\ref{spinlift} proves our claim.

In all cases we have chosen~$t$ to be a unipotent element of~$G$, and thus~$t$
has odd order if~$q$ is odd.
\end{prf}

We summarize the results of this subsection.

\begin{prop}
\label{C2Prop}
Let~$G$ be one of the following groups.

{\rm (a)} A group as in Hypothesis~\ref{ClassicalGroupCollection}.

{\rm (b)} A group $G = \Spin^\epsilon_n( q )$, $n \geq 7$, with~$q$ odd,
$\epsilon \in \{ -1, 0, 1 \}$ and $(n,q) \neq (7,3)$.

Let~$H$ be a maximal subgroup of~$G$ of ${\mathcal C}_2$-type.
Then~$H$ is not the stabilizer of an imprimitivity decomposition of an
irreducible $KG$-module.
\end{prop}
\begin{prf}
If~$G$ is not a spin group, the result follows from Lemmas~\ref{C2LargeBlocks}
and \ref{C2tExistence}. If~$G$ is a spin group, our claim follows from these
lemmas together with Corollary~\ref{spinlift}.
\end{prf}

\subsection{The case $H$ is of type ${\mathcal C}_3$}

Here we assume that~$H$ is a maximal subgroup of~$G$ of type ${\mathcal C}_3$.
This means that~$H$ stabilizes an $\mathbb{F}'$-structure of~$V$, where
$\mathbb{F}'$ is an extension field of degree~$r$ of the natural field of~$G$,
with $r$ a prime and $n = mr$ for some integer~$m$. A more precise description
is given below in the case of $r = 2$. In contrast to Aschbacher's original
definition \cite{Asch1}, we impose the somewhat stricter condition of 
\cite[Definition in~$\S 4.3$]{kleidman:4}. Table \ref{c3} shows the possible 
pairs $(G,H)$ that need to be considered. In this table, the parameter~$r$ 
has the same meaning as above. We begin by reducing to the case $r = 2$.

\begin{lem}
\label{CaseC3r2}
Let~$G$ be one of the groups of Hypothesis~\ref{ClassicalGroupCollection},
and let~$H$ be a maximal subgroup of $G$ of ${\mathcal C}_3$-type.

If $r \geq 3$ or if $r = 2$ and both~$G$ and~$H$ are orthogonal, then the 
maximal dimension of an irreducible $KG$-module is smaller than $[G\colon\!H]$.
\end{lem}
\begin{prf}
The assertion for $r \geq 3$ is easily established by showing 
$|H|c_G/|B| \leq 1$, using the bounds given in Table~\ref{c3}.

Suppose now that $r = 2$ and that~$G$ and~$H$ are orthogonal. Then the bounds
in Table~\ref{c3} give $h/b \leq 4 / |T|$ where~$T$ is a maximal split torus 
of~$G$. As $n \geq 7$, we again obtain $|H|c_G/|B| \leq 1$ if~$q > 2$.
Now suppose that $2 = q$. In this case~$n$ is divisible by~$4$ (see
\cite[Tables $3.5$]{kleidman:4}), and $H \leq \GO^\epsilon_{n/2}(4)\colon\!2$, 
if $G = \Omega^\epsilon_n(2)$, $\epsilon \in \{ 1, -1 \}$. 
One easily proves that $|H|^2 < |G|$ if $n \geq 12$, and thus we are done with 
Lemma~\ref{SquaredSubgroup}. The result for $n = 8$ (and $G = \Omega^-_8(2)$) 
is derived form the character table given in the Atlas \cite[p.~$88$]{ATLAS}.
\end{prf}

\begin{table}[h] 

 \caption{\label{c3} Bounds for ${\mathcal C}_3$-type groups; $n = mr$, $r \geq 2$}
    \begin{center}
    \renewcommand{\arraystretch}{1.24}
     \begin{tabular}{||c|c|c|c||}
      \hline
      $G$ & $bc_G^{-1}$ & $\tilde{H}$ & $h$ \\
      \hline
      \hline
      $\SL_n(q)$ & $q^{n(n-1)/2}(q-1)^{n-1}$ & $\GL_m(q^r).r$ & $rq^{rm^2}$ \\
      \hline
      $\SU_n(q)$ & $\frac{q^2-q+1}{q^2-1}q^{n(n-1)/2}(q^2-1)^{[n/2]}(q+1)^{-1}$ & 
                 $\GU_m(q^r).r $ & $r(q+1)q^{rm^2-1}$  \\
      $ r \geq 3$ & & &  \\
      \hline
      $\Sp_n(q)$ & $q^{n^2/4}(q-1)^{n/2}$ & 
                   $\Sp_m(q^r).r$ & $rq^{rm(m+1)/2}$ \\
      $ n \geq 4$ & & $\GU_{n/2}(q).2$ &  $2(q+1)q^{n^2/4 -1}$ \\
      \hline
      $\Omega_n(q)$ & $\frac{1}{2}q^{(n-1)^2/4}(q-1)^{(n-1)/2}$ &
      $\GO_{m}(q^r).r$ & $2rq^{rm(m-1)/2}$ \\
       $2 \nmid qn $ & & &  \\
      \hline
      $\Omega_n^+(q)$ & $\frac{1}{2}q^{n(n-2)/4}(q-1)^{n/2}$ &
      $\GU_{n/2}(q)$& $(q+1)q^{n^2/4 -1}$ 
      \\
      $ n \geq 8$ & & $\GO_{m}^+(q^r).r$  & $2rq^{{rm(m-1)}/{2}}$  \\
      & & $\GO_{n/2}(q^2)$ & $2q^{n(n-2)/{4}}$ \\
      \hline
      $\Omega_n^-(q)$ & $\frac{1}{2}q^{n(n-2)/4}(q-1)^{(n-2)/2}(q+1)$ &
      $\GU_{n/2}(q)$& $(q+1)q^{n^2/4 -1}$ 
      \\
      $ n \geq 8$ & & $\GO_{m}^-(q^r).r$ & $2rq^{{rm(m-1)}/{2}}$ \\
      & & $\GO_{n/2}(q^2)$ & $2q^{n(n -2 )/{4}}$   \\
      \hline
      \end{tabular}
      \end{center}
      \end{table}

We now make a few general remarks about $\mathcal{C}_3$-type subgroups where 
$r=2$. Then $G = \SU_n(q)$ is excluded by Table~\ref{c3}, so that we may
assume that our ground field is~$\mathbb{F}_q$. We begin with an $m$-dimensional 
$\mathbb{F}_{q^2}$-vector space~$V$ and basis $v_1,\dots , v_m$. Let $x \in 
\mathbb{F}_{q^2} \setminus \mathbb{F}_q$.
Then ${1,x}$ is an $\mathbb{F}_q$ basis of $\mathbb{F}_{q^2}$ and $v_1, xv_1, 
v_2, xv_2, \dots, v_m, xv_m$ is an $\mathbb{F}_q$-basis of~$V$. In the 
following, elements of $\Aut_{\mathbb{F}_q}( V ) = \GL(V)$ are written with 
respect to this basis, unless explicitly stated otherwise. Now assume that 
$x^{q + 1} = 1$ and define $a \in \mathbb{F}_q$ by $x^2 = ax - 1$.
Then the polynomial $X^2 - a X + 1$ is irreducible over $\mathbb{F}_q$
(with roots~$x$ and $x^q$). 

Let 
$$z_0 :=\left(
\begin{array}{cc}
  0 & 1 \\
  -1 & a \\
\end{array}
\right)$$ 
and observe that the minimal polynomial of $z_0$ equals  $X^2 - a X + 1$, and that 
$z_0$ is the matrix representing right muliplication with~$x$ in 
$\mathbb{F}_{q^2}$ with respect to the basis ${1,x}$. 

Define $z$ to be the $(n \times n)$-block diagonal matrix all of whose diagonal 
blocks are equal to~$z_0$. We have the following well known fact.

\begin{lem} 
\label{Fq2Linearity}
If $t \in \GL(V)$ and $M$ is the matrix representing $t$ with respect to the 
basis $v_1,xv_1,v_2,xv_2,\dots, v_m,xv_m$, then $t$ is $\mathbb{F}_{q^2}$-linear 
if and only if $zM = Mz$. 
\end{lem}

Thus~$z$ defines an $\mathbb{F}_{q^2}$-structure on~$V$, and a subgroup~$H$ of 
$\GL(V)$ is of ${\mathcal C}_3$-type (with respect to this structure), if 
$H \leq N_{\GL(V)}( \langle z \rangle )$.

\begin{lem}
\label{C3tExistence}
Let~$G$ be one of the groups of Hypothesis~\ref{ClassicalGroupCollection}.
Suppose that~$H$ is a maximal subgroup of~$G$ of ${\mathcal C}_3$-type defined 
over a field extension of~$\mathbb{F}_q$ of degree~$2$, and that~$H$ is not 
orthogonal if~$G$ is orthogonal. Then there exists an element $t \in G 
\setminus H$ such that $ t \in C_G(H \cap H^t)$. Moreover, if~$G$ is orthogonal
and~$q$ is odd, there exists such a~$t$ of odd order.
\end{lem}
\begin{prf}
Notice that Table~\ref{c3} excludes the case that~$G$ is unitary. By
Lemma~\ref{Fq2Linearity} we may assume $H = N_G( \langle z \rangle )$.
Moreover, if $s \in \GL_n(q)$ with $z^s = z^{-1}$ and $G^s = G$, then~$s$ 
induces a field automorphism of order two on $C_G( z )$ and 
$H = N_G( \langle z \rangle ) = \tilde{H} \cap G$ with $\tilde{H} :=  
\langle \tilde{C}, s \rangle$ where $\tilde{C} := C_{\GL_n(q)}( z )$. 

Put 
$$s_0 := \left(
\begin{array}{cc}
  0 & 1 \\
  1 & 0 \\
\end{array}
\right)$$
and let~$s$ be the $(n \times n)$-block diagonal matrix all of whose diagonal 
blocks are equal to~$s_0$. Then $z^s = z^{-1}$.
As in the previous subsection we see that any matrix that commutes 
with~$s$ is a block matrix all of whose blocks must have the form 
$$
\left(
\begin{array}{cc}
  \alpha & \beta \\
  \beta & \alpha \\
\end{array}
\right)$$
with $\alpha, \beta \in \mathbb{F}_q$.

Suppose first that~$G$ is linear, i.e., $G = \SL_n(q)$. Then $G^s = G$
and thus $H = \tilde{H} \cap G$. Define~$t$ to be a block diagonal matrix 
with diagonal blocks~$t_0$ and~$I_{n-4}$ where 
$$t_0 = \left(
\begin{array}{cc}
  I_2 & 0  \\
  s_0 & I_2 \\
\end{array}
\right).$$ 
Then $t \in G \setminus H$ and~$t$ commutes with~$s$. A straightforward 
calculation shows that the upper left $(4 \times 4)$-block of $[ z^{-1}, z^t ]$
equals
$$
\left(
\begin{array}{cc}
  I_2 & 0 \\
  s_0(2I_2 - z_0^{-2} - z_0^{2})  & I_2 \\
\end{array}
\right).
$$
Now 
$2I_2 - z_0^{-2} - z_0^{2} = (4 - a^2)I_2$ and thus $[ z^{-1}, z^t ] = 
t^{4 -a^2}$. As $a \neq \pm 2$ and~$t$ is of prime order, it follows that 
$t \in \langle z, z^t \rangle$. So the first two hypotheses of 
Lemma~\ref{2.3replace} are satisfied. We now verify its third 
hypothesis. We view the elements of~$\tilde{C}$ as block matrices 
$(A_{ij})_{1 \leq i,j \leq m}$, where each block $A_{ij}$ is a 
$(2 \times 2)$-matrix commuting with~$z_0$. Let $(A_{ij})$ be an element 
of~$\tilde{C}$. Then the first row of $(2 \times 2)$-blocks of 
$(A_{ij})^t$ consists of the matrices $A_{1,1} - A_{1,2}s_0, A_{1,2}, 
\ldots , A_{1,m}$. On the other hand, the first row of $(A_{ij})s$ consists of 
the matrices $A_{1,1}s_0, A_{1,2}s_0, \ldots , A_{1,m}s_0$. One easily checks 
that a $(2 \times 2)$-matrix~$M$ commuting with~$z_0$ and satisfying $M = Ms_0$ 
is the zero matrix or $a = -2$. Since $a \neq \pm 2$, the latter 
case cannot occur. Thus ${\tilde{C}^t} \cap \tilde{C}s = \emptyset$.
It follows from Lemma~\ref{2.3replace} that~$t$ commutes with $\tilde{H} \cap
\tilde{H}^t$, hence~$t$ commutes with $H \cap H^t$.

Now let~$G$ be symplectic, i.e., $G = \Sp_n(q)$. In this case $H$ is either 
symplectic or unitary. To deal with the case that~$H$ is unitary, let~$F$ be 
the block diagonal matrix with blocks
$$\tilde{J}_2:=
\left(
\begin{array}{cc}
  0 & 1 \\
  -1 & 0 \\
\end{array}
\right),$$
and assume that~$F$ represents the form defining~$G$. Then $z \in G$ and thus
$C_G(z)$ is unitary, and we may take $H = N_G( \langle z \rangle )$.

In the following, we choose a non-degenerate, $2k$-dimensional, $z$-in\-variant 
subspace~$X$ of~$V$, e.g., $X := \langle v_1, xv_1, \ldots , v_k, xv_k \rangle$ 
in the notation introduced before Lemma~\ref{Fq2Linearity}, and put $S := 
\{ \alpha \in \Stab_G( X ) \mid \alpha_{X^\perp} = I_{X^\perp} \}$. Then 
$S \cong \Sp(X) \cong \Sp_{2k}( q )$. We view $z_X$ as an element of~$S$.

First assume that $q > 3$. Here, we let $k = 1$, and thus $S \cong \SL_2(q)$
and $\langle z_0 \rangle$ is a torus of order $q + 1$ of~$S$.
Put $N := N_S( \langle z_0 \rangle )$. By Lemma~\ref{PGL2}, there is an element 
$t \in S$ such that~$t$ centralizes $N \cap N^{t} = Z( S )$ and 
$\langle z_0, z_0^{t} \rangle = S$. As~$S$ is perfect, we find that~$S$ is the 
derived subgroup 
of $\langle z, z^t \rangle$. Now $H \cap H^t \leq N_G( \langle z, z^t \rangle )$, 
and thus $H \cap H^t$ normalizes~$S$. It follows that $H \cap H^t \leq N_G( S )
= \Sp( X ) \times \Sp( X^\perp )$. By construction, $t \in C_G( H \cap H^t )$.

When $q=3$, we let $k = 2$.
Then $S \cong \Sp(X) \cong \Sp_4( 3 )$. We note that~$z_X$ 
corresponds to an element of order~$4$ in $\Sp_4(3)$ which squares to~$-I_4$.
Using GAP, we find an
element~$t$ of order~$5$ in~$S$ which is inverted by~$z_X$, and such that
$t \in C_S( N_S( \langle z_X \rangle ) \cap N_S( \langle z_X \rangle )^t)$. The
first of these properties implies that $t \in \langle z, z^t \rangle$, and 
that~$t$ generates a normal Sylow $5$-subgroup of $\langle z, z^t \rangle$.
As an element of order~$5$ of~$S$ clearly does not fix any non-trivial vector 
of~$X$, the fixed space of $\langle t \rangle$ equals~$X^\perp$. Thus 
$H \cap H^t \leq N_G( \langle z, z^t \rangle ) \leq N_G( \langle t \rangle )$ 
fixes~$X^\perp$. The second property of~$t$ now implies that 
$t \in C_G( H \cap H^t )$.

Finally suppose that $q = 2$. This time we let $k = 3$, which is possible as 
$n \geq 8$ by Hypothesis~\ref{ClassicalGroupCollection}. Then
$S \cong \Sp(X) \cong \Sp_6( 2 )$. Using GAP, we find an element $t \in S$ 
such that $\langle z_X, z_X^t \rangle$ is isomorphic to the simple group 
$\SL_8( 2 )$. Moreover,~$t$ can be chosen such that
$t \in C_S( N_S( \langle z_X \rangle ) \cap N_S( \langle z_X \rangle )^t)$.
(In fact~$t$ can be chosen to be an involution whose centralizer in~$S$ has
order~$384$.)
These choices imply that the derived subgroup of $\langle z, z^t \rangle$ is
equal to $\langle z_X, z_X^t \rangle \cong \SL_8(2)$. The maximal subgroup
$\SL_2(8)\colon\!3$ of $\Sp_6(2)$ acts irreducibly on the natural 
$6$-dimensional module of $\Sp_6(2)$, and the restriction of this module to 
$\SL_2(8)$ is also irreducible. It follows that the fixed space of 
$\langle z_X, z_X^t \rangle$ on~$V$ equals~$X^\perp$. Hence
$H \cap H^t \leq N_G( \langle z, z^t \rangle ) \leq 
N_G( \langle z_X, z_X^t \rangle )$ stabilizes~$X^\perp$. In turn,
this implies that $t \in C_G( H \cap H^t )$. 

We now consider the case where $H$ is symplectic. In this case $n = 2m = 4k$ 
and~$H$ is of type $\Sp_{2k}(q^2).2$. We choose the $\mathbb{F}_{q^2}$-bilinear
symplectic form on~$V$ stabilized by~$H$ in such a way that its Gram matrix 
with respect to our $\mathbb{F}_{q^2}$-basis $\{ v_1, \ldots , v_m \}$ is the 
block diagonal matrix whose diagonal blocks are 
$$\tilde{J}_2:=
\left(
\begin{array}{cc}
  0 & 1 \\
  -1 & 0 \\
\end{array}
\right).$$
Composing this $\mathbb{F}_{q^2}$-bilinear form with the trace 
$\mathbb{F}_{q^2} \rightarrow \mathbb{F}_q$, we obtain an 
$\mathbb{F}_q$-bilinear symplectic form on the $2m$-dimensional 
$\mathbb{F}_q$-vector space~$V$. We take~$G$ to be the group stabilizing this
latter from. With respect to the $\mathbb{F}_{q}$-basis $\{ v_1, xv_1, \ldots ,
v_m, xv_m \}$, the Gram matrix of this latter form is a block diagonal matrix 
consisting of~$k$ diagonal blocks of the ($4 \times 4$)-matrix~$E_0$ with 
$$E_0:=
\left(
\begin{array}{cccc}
0 & 0 & 2 & a \\
0 & 0 & a & a^2-2 \\
-2 & -a & 0 & 0 \\
-a & -a^2 + 2 & 0 & 0 \\
\end{array}
\right).$$ 
To simplify our calculations we change our $\mathbb{F}_q$-basis of 
$\mathbb{F}_{q^2}$, resulting in a change of the $\mathbb{F}_q$ basis of $V$.

Recall that the minimal polynomial of $z_0$ is $X^2 - a X + 1$. Whenever 
possible, we choose~$z_0$ such that $a \not \in \{0,1\}$. Such a choice is 
possible if $q \neq 2,3,5$. When $q = 2,5$ we choose $a = 1$ and when 
$q = 3$ we choose $a = 0$. If $q \neq 2,3,5$ the base change matrix we use is 
$$B_0:=
\left(
\begin{array}{cc}
  1 & u \\
  -u & u \\
\end{array}
\right)$$ 
where  $u = (1-a)^{-1} = -(a-1)^{-1}$.
Note that $B_0z_0 = z_0B_0$ and that the new Gram matrix has block diagonal 
form with each diagonal block equal to 
$$F_0:= (1+u)
\left(
\begin{array}{cccc}
0 & 0 & 1-u & -2u \\
0 & 0 & -2u & 1-u \\
u-1 & 2u & 0 & 0 \\
2u & u-1 & 0 & 0 \\
\end{array}
\right).$$ 
If $q = 2,5$ the base change matrix we use is 
$$B_0:=
\left(
\begin{array}{cc}
  0 & 1 \\
  -1 & 1 \\
\end{array}
\right) = z_0,$$ 
and the new Gram matrix has block diagonal 
form where each diagonal block is equal to 
$$F_0:= 
\left(
\begin{array}{cccc}
0 & 0 & -1 & -2 \\
0 & 0 & -2 & -1 \\
1 & 2 & 0 & 0 \\
2 & 1 & 0 & 0 \\
\end{array}
\right).$$ 
If $q = 3$ the base change matrix we use is 
$$B_0:=
\left(
\begin{array}{cc}
  1 & -1 \\
 1 & 1 \\
\end{array}
\right).$$ 
Again we note that $B_0z_0 = z_0B_0$ and that the new Gram matrix has block 
diagonal form where each diagonal block is equal to 
$$F_0:= 
\left(
\begin{array}{cccc}
0 & 0 & 0 & 1 \\
0 & 0 & 1 & 0 \\
0 & -1 & 0 & 0 \\
-1 & 0 & 0 & 0 \\
\end{array}
\right).$$ 
We now view~$s$ as an element of $\GL(V)$, written with respect to our new 
$\mathbb{F}_q$-basis of~$V$. Then $s \in G$ and $H = 
N_G( \langle z \rangle ) = \langle C_G( z ), s \rangle$. Let $t \in \GL(V)$ 
be the element represented, with respect to the new 
basis, by the block diagonal matrix $\mbox{\rm diag}(t_0,I_{n-4})$ with
$$t_0 = \left(
\begin{array}{cc}
  I_2 & 0 \\
  s_0  & I_2 \\
\end{array}
\right).$$ 
Then $t \in G \setminus H$. As in the linear case we find that the 
hypotheses of Lemma~\ref{2.3replace} are satisfied.

We now consider the cases where~$G$ is orthogonal. Here~$H$ is either 
orthogonal or unitary. The first of these cases is excluded by the 
hypothesis of the lemma, so let us assume that~$H$ is unitary. 

We may assume that the quadratic space~$V$ is the orthogonal direct sum of 
the $2$-dimensional, $z$-invariant, anisotropic spaces $V_1, \ldots, V_{m}$, 
with $V_i := \langle v_i, xv_i \rangle$, where the form on $V_i$ is given by 
$(\alpha, \beta ) \mapsto \alpha^2 + a\alpha\beta + \beta^2$ with respect to 
the basis $\{ v_i, xv_i \}$. Let $V = X \oplus Y$ with $X := 
V_1 \oplus \cdots \oplus V_4$ and $Y := V_5 \oplus \cdots \oplus V_m$. 
Then~$X$ is a non-degenerate $8$-dimensional 
subspace of plus-type, invariant under~$z$, and $Y = X^\perp$. Our aim 
is to produce an element $t \in G$ acting trivially on~$Y$, such that
$t \in C_G( H \cap H^t )$.

We choose a new basis $w_1, \ldots , w_8$ of~$X$ such that $(w_i, w_{9-i})$ is 
a hyperbolic pair for $1 \leq i \leq 4$, such that the totally isotropic subspaces
$X_0 := \langle w_1, \ldots , w_4 \rangle$ and 
$X_0' := \langle w_5, \ldots , w_8 \rangle$ of~$X$ are $z$-invariant, and such 
that the matrix of~$z_{X}$ with respect to $\{ w_1, \ldots , w_8 \}$ equals 
$\mbox{\rm diag}(z_0,z_0,z_0^*,z_0^*)$ with
$$z_0^* = J_2^{-1}z_0^{-T}J_2 = 
\left( \begin{array}{cc} 0 & -1 \\ 1 & a \end{array} \right).$$
To see that such a basis of~$X$ exists, let $w'_1 := (1,0,\alpha,\beta)$ and 
$w'_2 := (0,1, -\beta, \alpha + a\beta)$ be elements of $V_1 \oplus V_2$, 
written with respect to the basis $\{ v_1, xv_1, v_2, xv_2 \}$ 
of~$V_1 \oplus V_2$. Here, $\alpha, \beta \in \mathbb{F}_q$ satisfy 
$\alpha^2 + a\alpha\beta + \beta^2 = -1$. Then $\langle w'_1, w'_2 \rangle$ is 
a totally isotropic subspace of $V_1 \oplus V_2$, and $w'_1 z = w'_2$ and 
$w'_2z = -w'_1 + aw'_2$. Next, choose an isotropic vector 
$w'_3 \in V_1 \oplus V_2$, orthogonal to~$w'_1$ and with
scalar product~$1$ with~$w'_2$. In other words, $(w'_2,w'_3)$ is a hyperbolic
pair of $V_1 \oplus V_2$, orthogonal to~$w'_1$. By subtracting a suitable
scalar multiple of~$w_1'$, we may assume that $w_3'$ is of the form $w_3' = 
(0, \gamma, \delta, \varepsilon)$. Now put $w'_4 := -w'_3z$.
Then $(w'_1, w'_4)$ is a hyperbolic pair orthogonal to $(w'_2,w'_3)$, and
$w'_4z = w'_3 + aw'_4$ (notice that $z^2 = az -I_V$). Producing an analogous 
basis for $V_3 \oplus V_4$ and rearranging the basis vectors, we obtain the 
desired basis $\{ w_1, \ldots , w_8 \}$ of~$X$.

By replacing the first eight elements of $\{ v_1, xv_1, \ldots , v_m, xv_m \}$ 
by $\{ w_1, \ldots , w_8 \}$, we obtain a new basis of~$V$. Elements of $\GL(V)$ 
are now written with respect to this new basis. Again, we view~$s$ as an element 
of $\GL(V)$, written with respect to the new basis. Then $s \in G$ and 
$H = N_G( \langle z \rangle ) = \langle C_G( z ), s \rangle$.

Now assume that $q \neq 3$. Then $a \neq 0$.
Define $t_0 \in \Omega( X )$ by the matrix
$$
\left(
\begin{array}{cccc}
I_2 & 0 & 0 & 0 \\
0 & I_2 & 0 & 0 \\
R & 0 & I_2 & 0 \\
0 & -R & 0 & I_2 \\
\end{array}
\right),$$
with 
$$R := \left(
\begin{array}{cc} -2a^{-1} & 1 \\ 1 & -2a^{-1} \end{array}
\right).$$
The cyclic group $Z := \langle \mbox{\rm diag}(z_0,z_0,z_0^*,z_0^*) \rangle$ 
of order $q + 1$ acts linearly on the $4$-dimensional space of matrices
$$
U := \left\{ 
\left(
\begin{array}{cccc}
I_2 & 0 & 0 & 0 \\
0 & I_2 & 0 & 0 \\
S & 0 & I_2 & 0 \\
0 & S' & 0 & I_2 \\
\end{array}
\right) 
\mid S \in \mathbb{F}_q^{2 \times 2}, S' = -J_2S^TJ_2 \right \}$$ 
by conjugation. As an $\mathbb{F}_qZ$-module,~$U$ is the direct sum 
$C_U( Z ) \oplus [U,Z]$. The space $C_U( Z )$ is $2$-dimensional, 
hence so is~$[U,Z]$. It is easy to see that $t_0 \in [U,Z]$, in fact
$t_0 = [u,z_X]$ with 
$$u := \left(
\begin{array}{cccc}
I_2 & 0 & 0 & 0 \\
0 & I_2 & 0 & 0 \\
S & 0 & I_2 & 0 \\
0 & S' & 0 & I_2 
\end{array}
\right),
$$
where 
$$S := \left(
\begin{array}{cc} 2a^{-1} & 0 \\ -1 & 0 \end{array}
\right).$$
Moreover,~$t_0$ and $[t_0,z_X]$ are linearly independent in $[U,Z]$. 
Thus $[U,Z]$ is an irreducible $\mathbb{F}_qZ$-module. It follows that 
$\langle z_X, z_X^{t_0} \rangle = \langle z_X, [t_0,z_X] \rangle = 
\langle z_X, t_0 \rangle = Z[U,Z]$. In particular, 
$t_0 \in \langle z_X, z_X^{t_0} \rangle$. 
Also, the derived subgroup of $\langle z_X, z_X^{t_0} \rangle$ equals 
$[U,Z]$, whose fixed space on~$X$ equals~$X_0$.
Now let $t \in G$ be the element which acts as~$t_0$ on~$X$ and as~$I_Y$ on~$Y$. 
Then the fixed space of the derived subgroup of $\langle z, z^t \rangle$ 
equals $X_0 \oplus Y$. As $H \cap H^t \leq N_G( \langle z, z^t \rangle )$,
it follows that $H \cap H^t$ fixes $X_0 \oplus Y$ and 
$X_0 = (X_0 \oplus Y)^\perp$.
Now let 
$$g = \left(
\begin{array}{cc}
A & B  \\
C & D 
\end{array}
\right)
$$
be an element of $H \cap H^t$, written with respect to our new basis of~$V$,
where~$A$ is an $(8 \times 8)$-matrix. Then
$$g^t = \left(
\begin{array}{cc}
t_0^{-1}At_0 & t_0^{-1}B  \\
Ct_0 & D 
\end{array}
\right).
$$
Since~$g$ stabilizes~$X_0$, the first four rows of~$B$ are zero, and thus
$t_0^{-1}B = B$. Similarly, as~$g$ stabilizes $X_0 + Y$, the last four columns 
of~$C$ are zero, and hence $Ct_0 = C$. It also follows that~$A$ is invertible,
and thus $A \in N_{\GL_8(q)}( \langle z_X \rangle ) \cap 
N_{\GL_8(q)}( \langle z_X^{t_0} \rangle )$. Now
$N_{\GL_8(q)}( \langle z_X \rangle ) = 
\langle C_{\GL_8(q)}( z_X ), s_X \rangle$, and~$s_X$ commutes with~$t_0$. 
Thus, by Lemma~\ref{2.3replace}, in order to prove that~$A$ commutes with~$t_0$,
we are left to show that $C_{\GL_8(q)}( z_X )^{t_0} \cap C_{\GL_8(q)}( z_X )s_X 
= \emptyset$. This is done in a similar way as in the proof for the linear case.

Now let $q = 3$. Using GAP, we find that the conjugacy class of~$z_X$ in 
$\GO(X) \cong \GO_8^+( 3 )$ is characterized by the property that $|z_X| = 4$ 
and $z_X^2 = -I_X$. Using GAP, one can show that there is an element 
$t_0 \in \GO(X)$ of order~$5$, centralizing $N_{\GO(X)}( \langle z_X \rangle ) 
\cap N_{\GO(X)}( \langle z_X^{t_0} \rangle )$, such that 
$\langle z_X, z_X^t \rangle$ has order~$20$ and such that its derived
subgroup is spanned by~$t_0$ 
(in particular, $t_0 \in \langle z_X, z_X^{t_0} \rangle$). Moreover,~$t_0$
does not have any non-trivial fixed vector on~$X$. Now let $t \in G$ 
be defined by letting it act as $t_0$ on~$X$ and as $I_Y$ on~$Y$. Then the
fixed space of the derived subgroup of $\langle z, z^t \rangle$ equals~$Y$.
As $H \cap H^t \leq N_G( \langle z, z^t \rangle )$, it follows that 
$H \cap H^t$ stabilizes~$Y$ and hence also~$X$. 
As $(H \cap H^t)_X$ equals $N_{\GO(X)}( \langle z_X \rangle ) 
\cap N_{\GO(X)}( \langle z_X^{t_0} \rangle )$, it follows that~$t$ 
centralizes $H \cap H^t$.

In any case, the element~$t$ constructed for~$G$ orthogonal and~$q$ odd
has odd order.
\end{prf}

\begin{rem}
{\rm
As $\Sp_4(2) \cong S_6$ is not simple, this group has been excluded in
Lemma~\ref{C3tExistence}. Nevertheless, there is an overfield subgroup
$H = \Sp_2( 4 ).2 \cong S_5$, and our argument in the proof of this lemma
shows that $S_6$ has an element $t$ such that $t \in C_{S_6}(H \cap {H^t})$. 
However $A_6$ contains no such element.
}
\end{rem}

We summarize the results of this subsection.

\begin{prop}
\label{C3Prop}
Let~$G$ be one of the following groups.

{\rm (a)} A group as in Hypothesis~\ref{ClassicalGroupCollection}.

{\rm (b)} A group $G = \Spin^\epsilon_n( q )$, $n \geq 7$, with~$q$ odd,
$\epsilon \in \{ -1, 0, 1 \}$ and $(n,q) \neq (7,3)$.

Let~$H$ be a maximal subgroup of~$G$ of ${\mathcal C}_3$-type.
Then~$H$ is not the stabilizer of an imprimitivity decomposition of an
irreducible $KG$-module.
\end{prop}
\begin{prf}
If~$G$ is not a spin group, the result follows from Lemmas~\ref{CaseC3r2}
and \ref{C3tExistence}. If~$G$ is a spin group, our claim follows from these
lemmas together with Corollary~\ref{spinlift}.
\end{prf}

\subsection{The case $H$ is of type ${\mathcal C}_4$ or ${\mathcal C}_7$}

Here we assume that~$H$ is a maximal subgroup of~$G$ of type ${\mathcal C}_4$
or~${\mathcal C}_7$. This means that~$H$ stabilizes certain tensor product
decompositions of~$V$. For the precise description see \cite[Definitions in
$\S 4.4$, $\S 4.7$]{kleidman:4}. Table~\ref{c4} gives the relevant estimates 
for~$bc_G^{-1}$ and~$h$. 

\begin{lem}
\label{C4Lemma}
Let~$G$ be one of the groups of Hypothesis~\ref{ClassicalGroupCollection}.
If~$H$ is a maximal subgroup of~$G$ of ${\mathcal C}_4$-type, then $hc_G/b < 1$.
\end{lem}
\begin{prf}
We may assume that $m$ and $n/m$ are distinct, since otherwise~$H$ is contained 
in a $\mathcal{C}_7$-type subgroup. Hence $n \geq 6$. Fix~$n$ and hence~$b$
(as a function of~$q$). 
We first identify the values of~$m$ for which~$h$ and consequently $hc_G/b$ 
is maximal. 

If $G \neq \Sp_n(q)$, then we observe that the exponent of~$q$ in~$h$ is 
of the form $t(x) + t(n/x)$ where $t$ is a quadratic polynomial in~$x$ which is 
positive and increasing on the interval $[2,n-2]$. So by symmetry we see that 
on the interval $[2,n-2]$ the function $t(x) + t(n/x)$ attains a minimum at 
the point $x = n/x$ and maxima at the end points. Evidently the endpoint maxima 
are equal. If $G = \Sp_n(q)$, the exponent of~$q$ in~$h$ is of the form 
$t(x) + s(n/x)$, where~$t(x)$ is as before and $s(x) = t(x)-x$. Here we see 
that the absolute maximum occurs at the right end point, i.e., when 
$x = n/2$. 

For example, if $G = \SU_n(q)$ and $m = n/2$, the exponent of~$q$ in~$h$ is 
$n^2/4 + 2$, whereas the exponent of~$q$ in~$b$ is $n(n-1)/2$. Thus, as 
$n \geq 6$, we find
$$ \frac{hc_G}{b} \leq 
\frac{(q^2 - 1)(q+1)^3q^{n^2/4 + 2 }}{(q^2 - q + 1)(q^2-1)^{[n/2]}q^{(n(n-1)/2)}} 
\leq \frac{1}{q^{(n^2-2n-8)/{4}}} < 1.$$
Similar calculations for the other cases prove our claim.
\end{prf}

\begin{table}[h] 
      \caption{\label{c4} Bounds for ${\mathcal C}_4$-type groups; $n = ms$, $m \neq s$}
      \begin{center}
      \renewcommand{\arraystretch}{1.24}
      \begin{tabular}{||c|c|c|c||}
      \hline
      $G$ & $bc_G^{-1}$ & $\tilde{H}$ & $h$ \\
      \hline
      \hline
      $\SL_n(q)$ & $q^{n(n-1)/2}(q-1)^{n-1}$ & $\GL_{m}(q) \otimes GL_{s}(q)$ &
       $q^{m^2 + s^2}$ \\
      \hline
      $\SU_n(q)$ & $\frac{q^2-q+1}{q^2-1}q^{n(n-1)/2}(q^2-1)^{[n/2]}(q+1)^{-1}$ &
      $\GU_{m}(q) \otimes \GU_{s}(q) $&
        $(q+1)^2q^{m^2 - 1 + s^2 - 1}$\\
      $ n \geq 3$ & & & \\
      \hline
      $\Sp_n(q)$ & $q^{n^2/4}(q-1)^{n/2}$ &
      $\Sp_m(q)\otimes \GO^{\epsilon}_{s}(q) $ &
      $ 2q^{(m(m+1) + s(s-1))/2}$ \\
      $ n \geq 4$ & & & \\
      \hline
      $\Omega_n(q)$ & $\frac{1}{2}q^{(n-1)^2/4}(q-1)^{(n-1)/2}$ &
      $\GO_m(q) \otimes \GO_{s}(q)$ 
      & $4q^{(m(m-1) + s(s-1))/2}$\\
      $ 2 \nmid (qn)$ & & & \\
      \hline
      $\Omega_n^+(q)$ & $\frac{1}{2}q^{n(n-2)/4}(q-1)^{n/2}$ &
      $\Sp_{m}(q) \otimes \Sp_{s}(q) $ & $q^{(m(m + 1) + s(s+1))/2}$
      \\
      $ n \geq 8$ & & $\GO^{\epsilon_1}_{m}(q) \otimes \GO^{\epsilon_2}_{s}(q) $
      &  $4q^{(m(m - 1) + s(s-1))/2}$\\
      \hline
      $\Omega_n^-(q)$ & $\frac{1}{2}q^{n(n-2)/4}(q-1)^{(n-2)/2}(q+1)$ &
      $\GO_m(q) \otimes \GO^-_{s}(q)$ 
      & $4q^{(m(m-1) + s(s-1))/2}$
      \\
      $ n \geq 8$ & & & \\
      \hline
      \end{tabular}
      \end{center}
      \end{table}

The case of ${\mathcal C}_7$-type subgroups is easily reduced to the previous
case.
\begin{lem} 
\label{LemmaC7}
Let~$G$ be one of the groups of Hypothesis~\ref{ClassicalGroupCollection}.
If~$H$ is a maximal subgroup of~$G$ of ${\mathcal C}_7$-type, 
then $|H|c_G/|B| < 1$ unless $G = \Sp_4( q )$. In the latter case~$H$ is a 
${\mathcal C}_1$-type subgroup (which has been ruled out as a block stabilizer
of an imprimitive irreducible $KG$-module in Proposition~\ref{C1Prop}). 
\end{lem}
\begin{prf}
Here,~$H$ is the stabilizer of a tensor product decomposition 
$$V = V_1 \otimes V_2 \otimes \cdots \otimes V_t$$
with $m$-dimensional subspaces $V_1, \ldots , V_t$, where $t, m \geq 2$. Now 
consider the maximal ${\mathcal C}_4$-type subgroup~$H_4$ of~$G$ stabilizing 
the tensor decomposition $V = V_1 \otimes W$ with $W = V_2 \otimes \cdots 
\otimes V_t$. 

Let~$h_4$ and~$h$ denote the upper bounds for~$|H_4|$ and~$|H|$ given in 
Tables~\ref{c4} and~\ref{c7}, respectively. If $t \geq 3$, then $|H| < h_4$. 
This is proved by showing $h < h_4$, except in the case $G = \Sp_8( q )$ with
$t = 3, m = 2$. In the latter case, we use the exact value of $|\tilde{H}|$ as
an upper bound for~$|H|$ and the fact that~$q$ is odd (see 
\cite[Table 3.5.C]{kleidman:4}). Now $|H|c_G/|B| \leq |H|c_G/b < h_4c_G/b$,
and thus Lemma~\ref{C4Lemma} establishes our claim in this case. If $t = 2$ and 
$m \geq 3$, then $hc_G/b < 1$ and we are done. 

If $t = 2 = m$, then $n=4$ and so~$G$ is not orthogonal. If~$G$ is symplectic, 
then~$G$ is isomorphic to a $5$-dimensional spin group and~$H$ corresponds to
${\mathcal C}_1$-type of this spin group. If~$G$ is linear or unitary, we use 
the exact values of $|\tilde{H}|$ (see Table~\ref{c7}) as upper bounds 
for~$|H|$. These give $|H|c_G/|B| < 1$.
\end{prf}

\begin{table}[h] 
      \caption{\label{c7} Bounds for ${\mathcal C}_7$-type groups; $n=m^t$}
      \begin{center}
      \renewcommand{\arraystretch}{1.24}
      \begin{tabular}{||c|c|c|c||}
      \hline
      $G$ & $bc_G^{-1}$ & $\tilde{H}$ & $h$ \\
      \hline
      \hline
      $\SL_n(q)$ & $q^{n(n-1)/2}(q-1)^{n-1}$ & $\GL_m(q) \wr S_{t}$ & $ q^{tm^2} \ t! $\\
      \hline
      $\SU_n(q)$ & $\frac{q^2-q+1}{q^2-1}q^{n(n-1)/2}(q^2-1)^{[n/2]}(q+1)^{-1}$ 
      & $\GU_m(q) \wr S_{t} $& $ (q+1)^tq^{t(m^2-1)}\ t!$ \\
      $ n \geq 3$ & & & \\
      \hline
      $\Sp_n(q)$ & $q^{n^2/4}(q-1)^{n/2}$ & $\Sp_m(q) \wr S_{t}$ &  
      $ q^{tm(m + 1)/2} \ t! $\  \\
      $ n \geq 4$ & & & \\
      \hline
      $\Omega_n(q)$ & $\frac{1}{2}q^{(n-1)^2/4}(q-1)^{(n-1)/2}$ &
      $\GO_{m}(q) \wr S_{t}$ &  $ 2^t q^{tm(m - 1)/2} \ t! $\ \\
      $ 2| (qn)$ & & & \\
      \hline
      $\Omega_n^+(q)$ & $\frac{1}{2}q^{n(n-2)/4}(q-1)^{n/2}$ & $\Sp_{m}(q) \wr S_{t}$&
         $ q^{tm(m + 1)/2} \ t! $ \\
      $ n \geq 8$ & & $\GO_{m}^{\epsilon}(q) \wr S_{t}$  &  $ 2^t q^{tm(m-1)/2} \ t! $ \\
      \hline
      \end{tabular}
      \end{center}
      \end{table}

We summarize the results of this subsection.

\begin{prop}
\label{C4u7Prop}
Let~$G$ be one of the following groups.

{\rm (a)} A group as in Hypothesis~\ref{ClassicalGroupCollection}.

{\rm (b)} A group $G = \Spin^\epsilon_n( q )$, $n \geq 7$, with~$q$ odd,
$\epsilon \in \{ -1, 0, 1 \}$ and $(n,q) \neq (7,3)$.

Let~$H$ be a maximal subgroup of~$G$ of ${\mathcal C}_4$-type
or ${\mathcal C}_7$-type.
Then~$H$ is not the stabilizer of an imprimitivity decomposition of an
irreducible $KG$-module.
\end{prop}
\begin{prf}
The result follows from Lemmas~\ref{C4Lemma} and \ref{LemmaC7}. 
\end{prf}

\subsection{The case $H$ is of type ${\mathcal C}_5$}

Suppose now that~$H$ is a maximal subgroup of~$G$ of type~${\mathcal C}_5$.
This means that~$H$ is the normalizer in~$G$ of a subgroup that can be 
defined over a subfield $\mathbb{F}_{q_0}$ (respectively 
$\mathbb{F}_{q_0^2}$) of the natural field for~$G$. For the precise 
description see \cite[Definition in $\S 4.5$]{kleidman:4}.
Table~\ref{c5} gives the relevant estimates for~$b$ and~$h$.

\begin{lem} 
\label{C5Lemma}
Let~$G$ be one of the groups of Hypothesis~\ref{ClassicalGroupCollection}.
If~$H$ is a maximal subgroup of~$G$ of ${\mathcal C}_5$-type, then 
$|H|c_G/|B| < 1$ unless $G = \SU_n( q )$ and~$H$ is symplectic. 
\end{lem}
\begin{prf}
Suppose that $G = \SL_n( q )$. We first consider the case $n = 2$. Then
$|H| \leq 2|\SL_2(q_0)|$. Hence
$$\frac{|H|}{|B|} \leq \frac{2q_0(q_0^2 - 1)}{q_0^k(q_0^k-1)} \leq 
\frac{2}{q_0^{k-1}} < 1$$
if $q_0 \neq 2$ or $q_0 = 2$ and $k \geq 3$. As $\SL_2( 4 )$ is excluded
by Hypothesis~\ref{ClassicalGroupCollection}, our claim follows in this case.
Now let $n = 3$. Then $|H| \leq 3|\SL_3(q_0)|$. As $3|\SL_3(q_0)|$ is 
largest for $q_0 = \sqrt{q}$, we find
$$\frac{|H|}{|B|} \leq \frac{3q_0^3(q_0^2 - 1)(q_0^3 - 1)}{q_0^6(q_0^2-1)^2} 
= \frac{3(q_0^2 + q_0 + 1)}{q_0^3(q_0 + 1)} < 1$$
for all $q_0 \geq 2$.
Next, let $n = 4$. Here, we use the obvious bound $|H| \leq 
|\GL_4( q_0 )|(q - 1)$. Again, we may assume that $q = q_0^2$, and so
$$\frac{|H|}{|B|} \leq \frac{(q_0 - 1)(q_0^3 - 1)(q_0^2 + 1)}{q_0^6} < 1$$
for all $q_0 \geq 2$.
Suppose now that $n \geq 5$. Then~$h$ is maximal if $q = q_0^2$. In this case 
$$\frac{h}{b} = \frac{q_0^n}{(q_0^2-1)^{n-2}} < 1$$ 
if $q_0 \geq 3$ or $q_0 = 2$ and $n \geq 6$. 
The case $q_0 = 2$ and $n = 5$ can be ruled out using the bound $|H| \leq
|\GL_5(2)|3$.

Let $G = \SU_n( q )$ and~$H$ unitary. Then $q = q_0^k$ with $k \geq 3$
(see \cite[Table 3.5.B]{kleidman:4}). Clearly,~$h$ is maximal if $k = 3$.
Suppose first that $n = 3$. Here we use the sharper bound $|H| \leq
3|\SU_3( q_0 )|$ and the exact value $|B| = q^3(q^2 - 1)$ to obtain
\begin{eqnarray*}
\frac{|H|c_G}{|B|} & \leq & 
\frac{3q_0^3(q_0^2-1)(q_0^3+1)(q^2 - 1)}{q^3(q^2-1)(q^2-q+1)} \\
& \leq & \frac{3(q_0^2-1)(q_0^3+1)}{q_0^6(q_0^6-q_0^3+1)} \\
& \leq & \frac{3 \cdot 3 \cdot 9}{64 \cdot 57} < 1
\end{eqnarray*}
for all $q_0 \geq 2$. For $n \geq 4$ we have
\begin{eqnarray*}
\frac{hc_G}{b} & = &
\frac{(q+1)^3}{q_0^{(n^2-3n+2)/2}(q^2 - 1)^{[n/2] - 1}(q^2 - q + 1)} \\
& \leq & \frac{(q_0^3 + 1)^2}{q_0^{(n^2-3n+2)/2}(q_0^6 - 1)^{[n/2] - 1}} \\
& \leq & \frac{(q_0^3 + 1)^2}{q_0^3(q_0^6 - 1)} \\
& \leq & \frac{9^2}{8 \cdot 7 \cdot 9} < 1
\end{eqnarray*}
for all $q_0 \geq 2$.

If $G = \SU_n( q )$ and~$H$ is neither unitary nor symplectic then
$$\frac{hc_G}{b} \leq \frac{2(q+1)}{(q^2-1)^{[n/2] - 1}} < 1$$
for all $n \geq 4$.
For $n=3$ our claim can be checked by using the exact value for~$|B|$.

If $G = \Sp_n( q )$, then~$h$ is maximal if $q = q_0^2$ and~$H$ is symplectic.            
In this case 
$$\frac{h}{b} = \frac{2q_0^{n/2}}{(q_0^2-1)^{n/2}} 
\leq \frac{2^{(n+2)/2}}{3^{n/2}} < 1$$ 
for all $n \geq 4$.

If~$G$ is orthogonal, then~$h$ is maximal if $q = q_0^2$. In this case 
$$\frac{h}{b} = \frac{4q_0^{[n/2]}}{(q_0^2-1)^{[n/2]}} < 1$$
for $n \geq 8$ and $n = 7$, $q_0 \geq 3$,
which implies our claim. Our proof is now complete.
\end{prf}

\begin{lem} 
\label{C5tExistence}
If $G = \SU_n(q)$ and~$H$ is a maximal subgroup of~$G$ of symplectic 
${\mathcal C}_5$-type, then there exists an element $t \in G \setminus H$ with 
$t \in C_G(H\cap {H^t})$.
\end{lem}
\begin{prf}
Here $n=2m$ and we may define $G$ as 
$$G := \{ A \in \SL_n( \mathbb{F}_{q^2} ) \mid A^T \hat{J} \rho(A) = \hat{J} \},$$
where $\rho$ is the automorphism of $\GL_n( \mathbb{F}_{q^2} )$ raising
every matrix entry to its $q$th power, and 
$$\hat{J} := \left( \begin{array}{cc}
0&I_m\\
-I_m&0\\
\end{array} \right).$$ 
Notice that~$\hat{J}$ defines a skew Hermitian form on $\mathbb{F}_{q^2}^n$
if~$q$ is odd. In this case we choose $\alpha \in \mathbb{F}_{q^2}^*$
with $\alpha^q = - \alpha$. Then $\alpha \hat{J}$ defines a non-degenerate 
Hermitian form with isometry group~$G$.

In this setting, we put
$$C := \{ A \in G \mid \rho( A ) = A \},$$
i.e.\ $C = C_G( \rho )$.
Then $H = N_G( C )$. In fact, $H = \langle C, s \rangle \cap G$
with $s = \beta s_0$, where
$$s_0 := \left( \begin{array}{cc}
0&I_m\\
I_m&0\\
\end{array} \right),$$
and $\beta \in \mathbb{F}_{q^2}^*$ of order $2(q+1)$ if~$q$ is odd, and of
order~$q+1$ if~$q$ is even. In the latter case, $H = C \times Z(G)$.

We choose
$$t =  \left( \begin{array}{cc}
A& 0\\
0 & \rho( A )^{-T}\\
\end{array} \right),$$ 
with 
$$A = \left( \begin{array}{cc}
\lambda & 0\\
0 & I_{m-1}\\
\end{array} \right),$$ 
if~$q$ is even, where $\lambda \in \mathbb{F}_{q^2}$ is a non-trivial element with
$\lambda^{q+1} = 1$,
and 
$$A = \left( \begin{array}{ccc}
1 &\lambda & 0\\
0 & 1 & 0 \\
0 & 0 & I_{m-2}\\
\end{array} \right),$$ 
if~$q$ is odd,
where $\lambda \in \mathbb{F}_{q^2}$ is a non-zero
element with $\lambda + \lambda^q = 0$. Then, in any case, $\rho( t ) = t^{-1}$.
Thus, $\langle t, \rho \rangle$ is a dihedral subgroup of $\Aut( G )$ of twice 
odd order. In particular,~$t$ is a power of ${\rho^t}\rho = t^2$, and hence 
$t \in \langle \rho, {\rho^t} \rangle$. Our claim follows from 
Lemma~\ref{2.3replace} in case~$q$ is even. If~$q$ is odd, we let~$\bar{\ }$
denote the natural homomorphism $G \rightarrow G/Z(G)$, and $\bar{\rho}$ the
automorphism of~$\bar{G}$ induced by~$\rho$. Then $\bar{H} = 
C_{\bar{G}}( \bar{\rho} )$. By Lemma~\ref{2.3replace}, we may conclude that
$\bar{t} \in C_{\bar{G} }( \bar{H} \cap \bar{H}^{\bar{t}} )$. As $\bar{t}$ has 
order prime to~$q + 1$, our assertion follows from Corollary~\ref{spinlift}.
\end{prf}

\begin{table}[h]
      \caption{\label{c5} Bounds for ${\mathcal C}_5$-type groups; $q = q_0^k$}
      \begin{center}
      \renewcommand{\arraystretch}{1.24}
      \begin{tabular}{||c|c|c|c||}
      \hline
      $G$ & $bc_G^{-1}$ & $\tilde{H}$ & $h$ \\
      \hline
      \hline
      $\SL_n(q)$ & $q^{n(n-1)/2}(q-1)^{n-1}$ & $\GL_n(q_0).(q-1)$ &  
      $q_0^{n^2}(q-1)$ \\
      \hline
      $\SU_n(q)$ & $\frac{q^2-q+1}{q^2-1}q^{n(n-1)/2}(q^2-1)^{[n/2]}(q+1)^{-1}$ & 
      $N_{\GU_n(q)}( \SU_n(q_0) ) $ & $q_0^{n^2-1}(q+1)^2$ \\
      $n \geq 3$ & & $\GO_n^{\epsilon}(q).(q+1) $ &  $2q^{n(n - 1)/{2}}(q+1)$  \\
       & & $\Sp_n(q).(q+1)$ &  $q^{n(n + 1)/{2}}(q+1)$  \\
      \hline
      $\Sp_n(q)$ & $q^{n^2/4}(q-1)^{n/2}$ & $\Sp_n(q_0).2$
      & $2q_0^{n(n+1)/2}$\\
      $ n \geq 4$ & &  & \\
      \hline
      $\Omega_n(q)$ & $\frac{1}{2}q^{(n-1)^2/4}(q-1)^{(n-1)/2}$ &
      $\GO_{n}(q_0)$ & $2q_0^{n(n-1)/{2}}$\\
      $ 2 \nmid qn$ & & & \\
      \hline
      $\Omega_n^+(q)$ & $\frac{1}{2}q^{n(n-2)/4}(q-1)^{n/2}$ & $\GO^-_n(q^{1/2})$&
        $2q^{n(n - 1)/{4}}$ \\
      $ n \geq 8$ & & $\GO_{n}^+(q_0)$ & $2q_0^{n(n - 1)/{2}}$\\
      \hline
      $\Omega_n^-(q)$ & $\frac{1}{2}q^{n(n-2)/4}(q-1)^{(n-2)/2}(q+1)$ &
      $\GO_{n}^-(q_0)$ &  $2q_0^{n(n - 1)/{2}}$ \\
      $ n \geq 8$ & &  & \\
      \hline
      \end{tabular}
      \end{center}
      \end{table}

We summarize the results of this subsection.

\begin{prop}
\label{C5Prop}
Let~$G$ be one of the following groups.

{\rm (a)} A group as in Hypothesis~\ref{ClassicalGroupCollection}.

{\rm (b)} A group $G = \Spin^\epsilon_n( q )$, $n \geq 7$, with~$q$ odd,
$\epsilon \in \{ -1, 0, 1 \}$ and $(n,q) \neq (7,3)$.

Let~$H$ be a maximal subgroup of~$G$ of ${\mathcal C}_5$-type.
Then~$H$ is not the stabilizer of an imprimitivity decomposition of an
irreducible $KG$-module.
\end{prop}
\begin{prf}
The result follows from Lemmas~\ref{C5Lemma} and \ref{C5tExistence}.
\end{prf}

\subsection{The case $H$ is of type ${\mathcal C}_6$}
Here we consider maximal subgroups~$H$ of~$G$ of ${\mathcal C}_6$-type. This 
means that $H = N_G( R )$, 
where~$R$ is a group of extraspecial type (we use the notation and
terminology of \cite[\S 4.6]{kleidman:4}). In particular,~$R$ is an
$r$-group for some prime~$r$, and~$R$ is the largest normal $r$-subgroup 
of~$H$.

\begin{lem} 
\label{C6Lemma}
Let~$G$ be one of the groups of Hypothesis~\ref{ClassicalGroupCollection}.
If~$H$ is a maximal subgroup of~$G$ of ${\mathcal C}_6$-type, then 
$|H|c_G/|B| < 1$.
\end{lem}
\begin{prf}
Let~$R$ be an $r$-group of extraspecial type such that $H = N_G( R )$. 

If either~$r$ is odd, or $Z(R)$ is of order~$4$, then $H/RZ(G) \cong
\Sp_{2m}(r)$ and $G = \SL_n( q )$ or $\SU_n( q )$ with $n = r^m$
(see \cite[Table 4.6.B, p.~150]{kleidman:4}). If $r = 2$, then $r \leq q$ and
$c_G = 1$ as~$n$ is even. If~$r$ is odd, then $r \mid q - 1$ or $r \mid 
q + 1$ and $c_G \leq (q^2 - 1)/(q^2 - q + 1) \leq q$. Thus, unless $r = q + 1$
is a Fermat prime, we have $|H|c_G \leq q^{2+2m}q^{2m^2+m}|Z(G)|$. Also,
$|B| \geq q^{n(n-1)/2}|Z(G)|$. Considering exponents we see that 
$|H|c_G/|B| < 1$ whenever
$$2m^2+3m+2 < r^m(r^m-1)/2.$$ 
The latter is true unless $m=1$ and $r = 2,3$, or $m=2=r$ or $m = 3$ and 
$r = 2$. Using sharper estimates for~$|B|$ and~$|H|$, we find that 
$|H|c_G/|B| < 1$ in these cases as well.
The case that $r = q + 1$ is a Fermat prime is ruled out using similar 
considerations.

If $Z(R)$ is of order~$2$, then $r=2$ and $H/R$ is isomorphic to a
subgroup of $\GO^{\pm}_{2m}(2)$ with $m \geq 1$, 
and $G = \Sp_n( q )$ or $\Omega^+_n( q )$ with $n = 2^m$ and~$q$ odd
(see \cite[Table 4.6.B, p.~150]{kleidman:4}).
If $m = 1$, then $G = \Sp_2( q )$ with~$q$ odd, and~$H$ is the normalizer of a
quaternion group of order~$8$. As $|H| \leq 48$ and $q \geq 11$ by assumption,
we have $|H| < q(q - 1) = |B|$.
If $m = 2$, then $G = \Sp_4( q )$ with~$q$ odd (the case $G = \Omega^+_4( q )$ 
is excluded by our assumptions on~$G$), and $H \leq 2^{1+4}\GO_4^-(2)$. 
Thus $|B| = q^4(q-1)^2$ and $|H| \leq 32 \cdot 72$. As $q > 3$, it follows that
$|H|/|B| < 1$.
If $m \geq 3$, then $G = \Sp_{2^m}(q)$ and $H \leq 2^{1+2m} \GO_{2m}^-(2)$ or 
$G = \Omega^+_{2^m}(q)$ and $H  \leq 2^{1+2m}\GO_{2m}^+(2)$. 
Using our standard bounds for~$|H|$ and~$|B|$, we find
$$\frac{|H|}{|B|} \leq \frac{2^{2m^2+m+1}}{(q-1)^{2^{m-1}}q^{2^{m-1}(2^{m-1} - 1)}}.$$ 
The latter is clearly smaller than~$1$ if $q \geq 5$ or if $q = 3$ and $m \geq 4$.
The case $q = 3$ and $m = 3$ can be settled using better estimates for~$|H|$.
\end{prf}

We summarize the results of this subsection.

\begin{prop}
\label{C6Prop}
Let~$G$ be one of the following groups.

{\rm (a)} A group as in Hypothesis~\ref{ClassicalGroupCollection}.

{\rm (b)} A group $G = \Spin^\epsilon_n( q )$, $n \geq 7$, with~$q$ odd,
$\epsilon \in \{ -1, 0, 1 \}$ and $(n,q) \neq (7,3)$, or $G = 
\Spin^+_6( q )$.

Let~$H$ be a maximal subgroup of~$G$ of ${\mathcal C}_6$-type.
Then~$H$ is not the stabilizer of an imprimitivity decomposition of an
irreducible $KG$-module.
\end{prop}
\begin{prf}
The result follows from Lemma~\ref{C6Lemma}.
\end{prf}

\subsection{The case $H$ is of type ${\mathcal C}_8$}
\label{SectionC8}
Suppose that~$H$ is a maximal subgroup of~$G$ of type ${\mathcal C}_8$.
This means that~$H$ is the normalizer in~$G$ of a classical subgroup on~$V$.
Table \ref{c8} displays the possible pairs $(G,H)$ that need to be considered,
but there are further conditions (see \cite[Table~3.5.A]{kleidman:4}) on~$n$ 
and~$q$ if~$G$ is linear. The parameter~$d$ in Table~\ref{c8} is defined only
for~$n$ even and equals $\mbox{\rm gcd}( n/2, q - 1 )$.
Notice that the case $G = \Sp_n(q)$ and~$q$ even has already been treated in 
the subsection on ${\mathcal C}_1$-type subgroups (cf.\ Lemma~\ref{SOnOddEven}).

\begin{lem} 
\label{C8Lemma}
Let $G = \SL_n(q)$ with $n \geq 2$, where~$n$ and~$q$ are subject to the 
restrictions formulated in Hypothesis~\ref{ClassicalGroupCollection}(a).  
Suppose that~$H$ is a maximal subgroup of~$G$ of orthogonal or unitary 
${\mathcal C}_8$-type; if~$H$ is orthogonal, assume that $q$ is odd, and 
if $H$ is unitary, assume that~$q$ is a square. Then $h/b < 1$, unless 
$n = 2$ and~$H$ is unitary.
\end{lem}
\begin{prf}
This is clear from the estimates in Table \ref{c8}. 
\end{prf}

\begin{lem} 
\label{C8tExistence}
If $G = \SL_n(q)$ with $n = 2m \geq 4$ even and if~$H$ is 
a maximal subgroup of~$G$ of symplectic ${\mathcal C}_8$-type, then there 
exists an element $t \in G \setminus H$ with $t \in C_G(H\cap {H^t})$.
\end{lem}
\begin{prf}
By $\sigma$ we denote the 
automorphism of~$G$ sending $A$ to $\tilde{J}^{-1} A^{-T} \tilde{J}$.
Then $C := C_G( \sigma ) = \Sp_{2m}( q )$. We have $N_{\GL_n(q)}( C ) =
\langle C, s_\lambda \rangle$ 
with 
$$
s_\lambda := \left( \begin{array}{cc}
\lambda I & 0 \\
0 & I \\
\end{array} \right),$$
where $\lambda$ is a generator of $\mathbb{F}_q^*$. Let $d := 
\mbox{\rm gcd}( m, q - 1 )$ and put $s := s_\lambda^{(q-1)/d}$. Then~$s$ has
order~$d$ and determinant~$1$, and $H = N_G( C ) = \langle C, s \rangle$.

If~$q$ is even and $m = 2, 4$, or if~$q = 2$, we let
$$t =  \left( \begin{array}{cc}
A& 0\\
0 & J A^T J \\
\end{array} \right),$$
with
$$A =  \left( \begin{array}{cc}
A_0 & 0\\
0 & I_{m-2} \\
\end{array} \right),$$
where $A_0 \in \SL_2( q )$ has order $q + 1$.
In all other cases we put
$$t =  \left( \begin{array}{cc}
A& 0\\
0 & A\\
\end{array} \right),$$
where~$A$ is chosen as follows.
If~$q$ is even (and hence $m \neq 2, 4$ and $q > 2$), we let~$A$ be a diagonal 
matrix of 
determinant~$1$ and order~$q-1$ such that the eigenspaces of~$A$ for different 
eigenvalues have distinct dimension, and such that $J^{-1}A^{-T}J = A^{-1}$.
If $q$ is odd, we put
$$A = \left( \begin{array}{ccc}
1 & 0 & 1\\
0 & I_{m-2} & 0 \\
0 & 0 & 1 \\
\end{array} \right).$$ 
Then $\sigma( t ) = t^{-1}$ and thus
$\langle t, \sigma \rangle$ is a dihedral subgroup of $\Aut(G)$ of twice odd 
order. In particular,~$t$ is a power of ${\sigma^t}\sigma = t^2$, and hence
$t \in \langle \sigma, {\sigma^t} \rangle$. As~$t$ commutes with~$s$, 
the first two hypotheses of Lemma~\ref{2.3replace} are satisfied.

Clearly,~$\langle t \rangle$ is a characteristic subgroup of 
$\langle z, z^t \rangle$. In the even~$q$ case we have $N_G( \langle t \rangle )
= C_G( t )$ by our conditions on the eigenvalues. Thus~$t$ centralizes 
$H \cap H^t$ by Lemma~\ref{2.3replace_var}.

Suppose now that~$q$ is odd. In this case we verify Condition~(3) of 
Lemma~\ref{2.3replace}.
Let $X \in C$ such that 
$t^{-1}Xt = Ys^i$ for some $Y \in C$ and some $0 < i < d$. Thus
$$t^{-1} X t s^{-i} \in C.$$
Now $\sigma( t ) = t^{-1}$ and $\sigma( s ) = \nu^{-1} s$ for some $\nu \in 
\mathbb{F}_q^*$ of order~$d$, and so
$$t^{-1} X t s^{-i} = \sigma( t^{-1} X t s^{-i} ) = \nu^{i} t X t^{-1} s^{-i}.$$
It follows that 
$$t^2X \nu^i = Xt^2.$$
As $\nu^i \neq 1$, this implies that $X = 0$, a contradiction.
\end{prf}

\begin{table}[h]
      \caption{\label{c8} Bounds for ${\mathcal C}_8$-type groups}
      \begin{center}
      \renewcommand{\arraystretch}{1.24}
      \begin{tabular}{||c|c|c|c||}
      \hline
      $G$ & $b$ & $\tilde{H}$ & $h$ \\
      \hline
      \hline
      $\SL_n(q)$ & $q^{n(n-1)/2}(q-1)^{n-1}$ & $\Sp_n(q).d$ & $dq^{n(n+1)/2}$  \\
       & & $\GO_n(q)$ & $2q^{n(n-1)/2}$\\
      & & $\GO_n^{\pm}(q)$ & $2q^{n(n-1)/2}$\\
      & & $\GU_n(q^{1/2})$ & $(q^{1/2}+1)q^{(n^2-1)/2}$ \\
      \hline
      $\Sp_n(q)$ & $q^{n^2/4}(q-1)^{n/2}$ & $\GO^{\pm}_n(q)$ & $q^{n(n-1)/2}$\\
      $n \geq 4$, $\mbox{ $q$ even}$ & & & \\
      \hline
      \end{tabular}
      \end{center}
      \end{table}

We summarize the results of this subsection.

\begin{prop}
\label{C8Prop}
Let~$G$ be one of the groups of Hypothesis~\ref{ClassicalGroupCollection}(a)
or~(c).

Let~$H$ be a maximal subgroup of~$G$ of ${\mathcal C}_8$-type, where the
restrictions of \cite[Tables 3.5A, 3.5C]{kleidman:4} apply.
Then~$H$ is not the stabilizer of an imprimitivity decomposition of an
irreducible $KG$-module.
\end{prop}
\begin{prf}
The result follows from Lemmas~\ref{C8Lemma}, \ref{C8tExistence} and 
\ref{SOnOddEven}, except for $G = \SL_2( q )$ and~$H$ unitary. The latter case
is treated in Lemma~\ref{C3tExistence}.
\end{prf}

\subsection{The case $H$ is of type ${\mathcal S}$}
In this subsection we let~$G$ be one of the groups of 
Hypothesis~\ref{ClassicalGroupCollection}. 
In \cite[p.~$3$]{kleidman:4}, Kleidman and Liebeck give a definition of
the class $\mathcal{S}$ of subgroups of the simple group $\bar{G} := 
G/Z(G)$. If $\hat{G}$ is a quasisimple central extension of~$\bar{G}$,
we define the $\mathcal{S}$-type subgroups of $\hat{G}$ to be the inverse
images of the $\mathcal{S}$-type subgroups of~$\bar{G}$. The definition
of Kleidman and Liebeck implies in particular, that subgroups of~$G$
of type~$\mathcal{S}$ are not of type~${\mathcal C}_8$.

Suppose that~$H$ is a maximal subgroup of~$G$ from class~${\mathcal S}$. 
Then $Z( H ) = Z( G )$ and we write $\bar{H} := H/Z(G)$.
Also, $\bar{H}$ is almost simple and thus $F^*( \bar{H} )$ is a non-abelian
finite simple group. Finally, $F^*( H )'$ is quasisimple and acts absolutely
irreducibly on~$V$. 

Recall that~$V$ is the natural vector space of~$G$ and that~$B$ denotes a Borel 
subgroup of~$G$. In this situation the main theorem of Liebeck 
\cite[Theorem 4.1]{liebeck:1} states that either $|\bar{H}| \leq |V|^3$ 
or $F^*( \bar{H} )$ is an alternating group and $V$ is the reduced permutation 
module for $F^*( \bar{H} )$. We consider each possibility in turn. 

In the first case we easily obtain the following estimate.

\begin{lem} 
If $|\bar{H}| \leq |V|^3$, then $|H|c_G/|B| < 1$ unless $G$ is linear and 
$\mbox{\rm dim}(V) \leq 7$, or $G$ is unitary or symplectic and 
$\mbox{\rm dim}(V) \leq 12$, or $G$ is orthogonal and 
$\mbox{\rm dim}(V) \leq 14$.
\end{lem}
\begin{prf}
This is proved using $|\bar{H}| \leq q^{3n}$, respectively $|\bar{H}| \leq 
q^{6n}$ if~$G$ is unitary, together with the lower bounds for~$|B|$ given 
in Table~\ref{c2}, for example.
\end{prf}

Now the tables in Hiss-Malle \cite{HiMa} and L{\"u}beck \cite{LuebeckDefining} 
leave only the following list of possibilities.

\begin{lem}
\label{TypeSFirstCollection}
If~$H$ is a maximal subgroup of~$G$ of ${\mathcal S}$-type and either~$G$ is 
linear and $\mbox{\rm dim}(V) \leq 7$, or~$G$ is unitary or symplectic and 
$\mbox{\rm dim}(V) \leq 12$, or $G$~is orthogonal 
and $\mbox{\rm dim}(V) \leq 14 $, then one of the following is true.

{\rm (1)} If $F^*( \bar{H} )$ is a simple group of Lie type, cross 
characteristically embedded into~$G$, or if $F^*( \bar{H} )$ is a sporadic 
group, then $F^*( \bar{H} )$ is one of the following groups: $\PSL_2(q)$ with 
$q \leq 29$, $\PSL_3(3)$, $\PSL_3(4)$, $\PSp_4(3)$, $\PSp_4(5)$, $\PSp_6(2)$, 
$\PSp_6(3)$, $\PSU_3(3)$, $\PSU_3(4)$, $\PSU_4(3)$, $\PSU_5(2)$, 
$\Omega^+_8(2)$, $G_2(3)$, $G_2(4)$, ${^2\!B}_2(8)$, $M_{i}$ where 
$i \in \{11, 12, 22, 23, 24 \}$, $J_1$, $J_2$, $J_3$, $\Suz$.

{\rm (2)} If $F^*( \bar{H} )$ is a simple group of Lie type of characteristic 
equal to that of~$G$, then 
$F^*( \bar{H} )$ is of type $A_1$, of type~$A_2$ or 
${^2\!A}_2$ and $\mbox{\rm dim}(V) \in \{ 6, 7,8,10 \}$, of type~$A_3$ or 
${^2\!A}_3$ and $\mbox{\rm dim}(V) \in \{ 10,14 \}$, of type~$A_4$ or 
${^2\!A}_4$ and $\mbox{\rm dim}(V) = 10$, of type $B_2$ and 
$\mbox{\rm dim}(V) \in \{ 10, 12, 13, 14 \}$, of type ${^2\!B}_2$ and 
$\mbox{\rm dim}(V) = 4$, of type $B_3$ and $\mbox{\rm dim}(V) \in \{ 8, 14 \}$, 
of type $C_3$ and $\mbox{\rm dim}(V) \in \{ 8, 13, 14 \}$, or of type $G_2$ 
or ${^2\!G}_2$ and $\mbox{\rm dim}(V) \in \{ 6, 7, 14 \}$.

{\rm (3)}
If $F^*(\bar{H})$ is alternating and $V$ is not the reduced permutation module, 
then $F^*(H)' = A_n$ or $2.A_n$ with $n \leq 8$, or $F^*(H)' = 2.A_9$ or 
$2.A_{10}$ and $\mbox{\rm dim}( V ) = 8$.
\end{lem}

To deal with the case where $G$ and $F^*(H)'$ are of Lie type of equal 
characteristic, we need the following lemma, whose formulation and proof 
is due to Frank L\"ubeck. For the notation used see 
\cite[Section~$2$]{LuebeckDefining}.
\begin{lem}
\label{FranksLemma}
Let~$p$ be a prime and $\mathbb{F}$ an algebraic closure of the finite
field~$\mathbb{F}_p$ with~$p$ elements.

Let~$\mathbf{H}$ be a simple, simply connected algebraic group 
over~$\mathbb{F}$, defined over~$\mathbb{F}_p$, and let~$F$ denote the 
corresponding standard Frobenius map. Furthermore, let 
$$\Psi: \mathbf{H} \rightarrow \GL_n( \mathbb{F} )$$
be an irreducible rational representation of~$\mathbf{H}$ with
highest weight~$\lambda$. 

Let~$r$ be a positive integer, $q := p^r$, and put 
$H := \mathbf{H}^{F^r}$. Assume that $\lambda$ is $q$-restricted such that
$\Res_H^{\mathbf{H}}( \Psi )$ is irreducible by Steinberg's theorems
\cite[Theorems~$7.4$,~$9.2$]{Steinberg}. Let $\lambda = \lambda_0 + 
p \lambda_1 + \cdots + p^{r - 1} \lambda_{r - 1}$ be the $p$-adic 
decomposition of~$\lambda$ (i.e., the $\lambda_i$ are $p$-restricted,
$0 \leq i \leq r - 1$).

Let~$s$ be a positive integer and write~$i'$ for the smallest non-negative
remainder 
of $s + i$ modulo~$r$ for all $0 \leq i \leq r - 1$. Suppose that $\Psi( H ) 
\leq \GL_n( p^s )$. Then $\lambda_i = \lambda_{i'}$ for all $0 \leq i 
\leq r - 1$

In particular, $s \geq r$, if $\lambda = p^i\lambda_i$ for some
$0 \leq i \leq r - 1$. 
\end{lem}
\begin{prf}
By Steinberg's tensor product theorem \cite[Theorem~$1.1$]{Steinberg}, 
the representation~$\Psi$ 
is equivalent to a representation on the twisted tensor product 
$$L(\lambda) \cong L(\lambda_0) \otimes L(\lambda_1)^{(1)} \otimes \cdots 
\otimes L(\lambda_{r - 1})^{(r-1)},$$
where the tensor products are taken over~$\mathbb{F}$. The fact that 
$\Psi( H ) \leq \GL_n( p^s )$ implies that $\Res_H^{\mathbf{H}}( L( \lambda ) )$ 
is isomorphic to $\Res_H^{\mathbf{H}}( L( \lambda )^{(s)} )$. Now 
$L( \lambda )^{(s)} \cong L(\lambda_0)^{(s)} \otimes L(\lambda_1)^{(s + 1)} 
\otimes \cdots \otimes L(\lambda_{r - 1})^{(s + r-1)}$. The restriction
of the latter module to~$H$ is isomorphic to the restriction to~$H$ of
$L(\lambda_0)^{(0')} \otimes L(\lambda_1)^{(1')} 
\otimes \cdots \otimes L(\lambda_{r - 1})^{((r-1)')}$. 
The assertion follows from the uniqueness result of Steinberg's theorems 
\cite[Theorems~$7.4$,~$9.2$]{Steinberg}.
\end{prf}

\noindent
The 
examples that still need closer inspection are collected in the next lemma.

\begin{lem}
\label{Alternative1}
Let~$H$ be a maximal subgroup of~$G$ of ${\mathcal S}$-type. 
If $|\bar{H}| \leq |V|^3$ and $\mbox{\rm dim}(V) \leq 14$, 
then $|H|c_G/|B| < 1$ unless one of the following holds.

{\rm (1)}
If $F^*( \bar{H} )$ is of Lie type, cross characteristically embedded 
into~$G$, or if~$F^*( \bar{H} )$ is a sporadic group, then
$(G,H) = (\Omega_8^+(3),2.\Omega_8^+(2))$.

{\rm (2)} If $F^*( \bar{H} )$ is of Lie type of characteristic equal 
to that of~$G$, then $(G,H)$ is one of $(\Sp_6(q),G_2(q))$, 
$(\Omega_8^+(q),\Sp_6(q))$ with~$q$ even, or one of 
$(\Omega_8^+(q),\Spin_7(q))$, $(\Omega_7(q),G_2(q))$ with $q$ odd.

{\rm (3)} If $F^*( \bar{H} )$ is alternating and~$V$ is not the reduced 
permutation module, then $(G,H) \in 
\{ (\SU_3(5),3.A_7), (\Omega_8^+(2), A_9) \}$.
\end{lem}
\begin{prf}
For each simple group~$S$ from the lists in Lemma~\ref{TypeSFirstCollection}, 
we consider the quasisimple classical groups~$G$ possibly containing a 
maximal subgroup~$H$ of type~$\mathcal{S}$ with $F^*( \bar{H} ) = S$. We then 
use the precise values for the orders of $\Aut( S )$ and sharp 
bounds for the orders of the Borel subgroups of~$G$ to exclude further 
possibilities.

The first and third part are straightforward computations.
For the second part we use Lemma~\ref{FranksLemma} to bound the order of the
underlying field of~$G$ from below.
\end{prf}

The second alternative of Liebeck's result is treated in the following lemma.

\begin{lem}
\label{Alternative2}
If $F^*( \bar{H} ) = A_m$ and $V$ is the reduced permutation module then 
$|H|c_G/|B| < 1$ unless $q = 2$ and $n \leq 14$ or $q=3$ and $n \leq 7$.
\end{lem}
\begin{prf}
By definition, $n \in \{ m - 1, m - 2 \}$, and~$V$ carries a non-degenerate
quadratic or symplectic form. In particular, $c_G = 1$. Writing 
$\bar{B} := B/Z(G)$ we thus have $|\bar{B}| \geq q^{n(n-2)/4}$. Also, 
$|\bar{H}| \leq m! \leq (n+2)!$ unless $m = 6$. Notice that 
$|\bar{H}|/|\bar{B}| = |H|/|B|$.

Put $$f(n,q) := \frac{(n+2)!}{q^{n(n-2)/4}}.$$
Then 
$$f(n+1,q)/f(n,q) = \frac{n+3}{q^{(2n-1)/4}} < 1$$
for all $q \geq 2$ and all $n \geq 8$. Thus $f(n,q)$ is increasing for fixed 
$q \geq 2$ and $n \geq 8$.

We have $f(8,4) < 1$, hence $f(n,q) < 1$ for all $n \geq 8$ and $q \geq 4$.
Using better bounds for~$|\bar{H}|$ and~$|\bar{B}|$ we show that $|H| < |B|$
for all $q \geq 4$ and all $3 \leq n \leq 7$.

We have $f(10,3) < 1$, and hence $f(n,3) < 1$ for all $n \geq 10$. The case
$q = 3$ and $n = 8$ cannot occur, and for $q = 3$ and $n = 9$ we use sharper 
bounds for $|\bar{H}|$ and~$|\bar{B}|$ to obtain our result.

We finally have $f(16,2) < 1$ and thus $f(n,2) < 1$ for all $n \geq 16$. The
case $q = 2$ and $n = 15$ cannot occur.
\end{prf}

\begin{lem} 
Suppose that $({G},{H})$ occurs in the conclusion of Lem\-ma~\ref{Alternative1},
but $$(G,H) \not\in \{(\SL_2(11),2.A_5), (\SU_3(5),3.A_7), 
(\Sp_6(q), G_2(q)), \ q \mbox{\ even}\}.$$
If $G = \Omega_7( q )$ or $\Omega_8^+( q )$ with~$q$ odd, we write $\hat{G}$
for the corresponding spin group and $\hat{H}$ for the inverse image of~$H$ 
in~$\hat{G}$.

Then there is $t \in G \setminus H$ respectively $\hat{t} \in \hat{G} \setminus 
\hat{H}$ such that $t \in C_G(H \cap H^t)$, respectively $\hat{t} \in 
C_{\hat{G}}(\hat{H} \cap \hat{H}^{\hat{t}})$.
\end{lem}
\begin{prf}
Suppose that $(G, H) = (\Omega_8^+(3),2.\Omega_8^+(2))$. The following results
are easily obtained by a computation in GAP in the permutation representation 
of~$\bar{G}$ on the cosets of $\bar{H}$. The number $|H \backslash G /H|$ of 
double cosets of~$H$ in~$G$ equals~$6$. The largest double coset is of size 
$22\,400$. The corresponding point stabilizer $H \cap H^x$ 
contains a Sylow $3$-subgroup of~$H$.
No other double coset point stabilizer contains a Sylow $3$-subgroup of~$H$. 
The centralizer $C_G (H \cap H^x)$ is cyclic of order~$6$ and contains an 
element~$t$ of order~$3$ which does not lie in~$H$. Thus 
$H \neq H \cap H^t \geq H \cap H^x$. So in fact $H \cap H^t = H \cap H^x$, 
as $HxH$ is the unique double coset whose point stabilizer contains a 
full Sylow $3$-subgroup of~$H$. Thus $t \in C_G(H \cap H^t)$ which is our claim. 
The corresponding result for $\Spin_8^+( 3 )$ follows from the above in
conjunction with Corollary~\ref{spinlift}

If $({G}, {H})$ occurs in conclusion~(2) and if $H \neq G_2(q)$, then 
$G = \Omega_8^+(q)$ and $H = \Spin_7( q )$ ($\cong \Sp_6( q )$ if~$q$ is even). 
There is a triality automorphism~$\tau$ of $\bar{G} = P\Omega_8^+( q )$ sending 
$\bar{H}$ to a ${\mathcal C}_1^*$-type subgroup $\bar{H}^\tau$ of~$\bar{G}$ 
(see \cite[Proposition~$2.2.4$]{kleidman:1}). Write~$H_1$ for the inverse image
of~$\bar{H}^\tau$ in~$G$. By Lemmas~\ref{COneLin} and~\ref{Case1OfLemmanorm}, 
there is an element $t_1 \in G \setminus H_1$ such that~$t_1$ centralizes 
$H_1 \cap H_1^{t_1}$. Moreover, if~$q$ is odd,~$t_1$ has odd order. Write 
$\bar{t}_1$ for the image of~$t_1$ in~$\bar{G}$. Then $\bar{t}_1 \not\in 
\bar{H}^\tau$ and $\bar{t}_1$ centralizes $\bar{H}^\tau \cap 
(\bar{H}^\tau)^{\bar{t}_1}$.

Put $\bar{t} := \bar{t}_1^{\tau^{-1}}$. Then $\bar{t} \in \bar{G} \setminus 
\bar{H}$ and $\bar{t} \in C_{\bar{G}}( \bar{H} \cap \bar{H}^{\bar{t}})$.
If~$q$ is even, $\bar{G} = G$, and we are done. If~$q$ is odd, $\bar{t}$ has 
odd order, and the result for $(G,H)$ and for $(\hat{G}, \hat{H})$ follows 
from Corollary~\ref{spinlift}.

Suppose now that $(G,H) = (\Omega_7(q), G_2(q))$ and~$q$ is odd. 
The group $H \cong G_2(q)$ contains maximal subgroups $H^\epsilon \cong 
\SL_3^\epsilon(q)$, $\epsilon \in \{ 1, -1 \}$, such that $H^\epsilon$ 
stabilizes a vector $v_{\epsilon} \in V$ whose stabilizer in $G \cong 
\Omega_7(q)$ equals $\SO^{\epsilon}_6(q)$ (see \cite{Asch2}).
Choose $\epsilon \in \{ 1, -1 \}$ such that $q \equiv \epsilon 
\ \mbox{mod} \ 4$. 
Then, by \cite[p.~$168$]{Taylor},~$G$ contains an involution~$t$ which 
fixes $v_{\epsilon}$ and acts as $-I$ on 
$\langle v_{\epsilon} \rangle^{\perp}$. Now~$t$ is not in~$H$, as~$t$ 
centralizes $H^\epsilon = C_H(v_{\epsilon})$. Thus $C_H(v_{\epsilon})$ is 
contained in $H \cap H^t \neq H$.
So by the maximality of $C_H(v_{\epsilon})$ in~$H$ we obtain equality and so~$t$ 
centralizes $H \cap H^t$. 

The element~$t$ lifts to an element~$\hat{t}$ of order~$4$ in $\Spin_7(q)$.
In fact,~$\hat{t}$ generates the center of $\Spin_6^\epsilon(q) \cong 
\SL^\epsilon_4(q)$, hence commutes with the inverse image of $H^\epsilon$.
Thus our claim also holds for $\hat{G} = \Spin_7(q)$.
This completes our proof.
\end{prf}

\begin{lem}
Suppose that $({G},{H})$ occurs in the conclusion of 
Lem\-ma~\ref{Alternative2}. Then $G$ does not have an 
irreducible imprimitive cross characteristic representation with block 
stabilizer~$H$.
\end{lem}
\begin{prf}
If $({G},{H})$ occurs in the conclusion of Lemma~\ref{Alternative2}, 
then $F^*(\bar{H})$ is an alternating group, and~$V$ is the reduced permutation 
module. If $q = 3$, then $n \leq 7$. If $n = 7$, then $\bar{G} \cong 
\Omega_7(3)$, a case excluded by Hypothesis~\ref{ClassicalGroupCollection}.
If $n = 6$, then $\bar{G} \cong P\Omega_6^-(3) \cong \PSU_4(3)$, again excluded
by Hypothesis~\ref{ClassicalGroupCollection}. The case $n = 5$ does not occur.
If $n = 4$, then $\bar{H} \cong A_6$ and $\bar{G} \cong P\Omega_4^-( 3 ) 
\cong A_6$.

If $q=2$ then~$n$ is even and $4 \leq n \leq 14$. Let~$V$ denote the reduced
permutation module for $A_{n+2}$. Then~$V$ restricts to the reduced 
permutation module for~$A_{n+1}$. Also,~$V$ carries a non-degenerate, 
$A_{n+2}$-invariant quadratic form, if and only if $n = \mbox{\rm dim}(V)$ 
is not divisible by~$4$ (see \cite[p.~187]{kleidman:4} or \cite[p.~637]{GoWi}).
In this case, this form is also
$S_{n+2}$-invariant, but the elements of $S_{n+2} \setminus A_{n+2}$ do
not lie in $\Omega(V)$. If $4 \mid n$, then we get an embedding of 
$A_{n+1}$ as a maximal subgroup of $\Omega_n^\epsilon(2)$ 
for $\epsilon \in \{ +, - \}$, and an embedding of $S_{n+2}$ as a maximal 
subgroup of $\Sp_n(2) = \PSp_n(2)$.

For $n = 4$ we obtain
embeddings of $A_5$ into $\Omega_4^-( 2 ) \cong \PSL_2( 4 ) \cong A_5$
and of $A_6$ into $\PSp_4( 2 )$, so this case need not be considered.
The case $n = 6$ leads to $\bar{G} \cong \Omega_6^+( 2 ) \cong \PSL_4( 2 )
\cong A_8$, so again we are done. If $n = 8$, we get an embedding of~$A_9$
into $\Omega_8^+( 2 )$ and of~$S_{10}$ into $\PSp_{8}( 2 )$. The former
case is excluded by Hypothesis~\ref{ClassicalGroupCollection}. In the latter
case we use the known ordinary and modular character tables of $\bar{G} = 
\PSp_{8}( 2 )$ to show that no irreducible $KH$-module induces to an 
irreducible $KG$-module. If $n = 10$ we obtain an embedding of $A_{12}$
into $\Omega_{10}^-( 2 )$. All character tables of $\Omega_{10}^-( 2 )$
are available in GAP (except for the $2$-modular table, which is not relevant
in this case by the global assumption at the beginning of this chapter). 
It is then easy to check that $[\bar{G}\colon\!\bar{H}]$
does not divide the degree of any irreducible cross characteristic
representation of $\Omega_{10}^-( 2 )$. Next, let $n = 12$. Suppose that
$H = A_{13}$, embedded in the $12$-dimensional orthogonal group $G = 
\Omega^\epsilon_{12}(2)$. As 
$[ 13!/2^{31} ] = 2$, it follows that only $1$-dimensional $KH$-modules can 
possibly induce to irreducible $KG$-mod\-ules. As $H = A_{13}$ is perfect,
the only $1$-dimensional $KH$-module induces to a reducible $KG$-module.
We also have an embedding of $H = S_{14}$ as a maximal subgroup of 
$\Sp_{12}(2)$. The ordinary character table of $\Sp_{12}(2)$ is available
in GAP, so that it is trivial to check that the largest dimension of an
irreducible $KG$-module is smaller than the index of $S_{14}$ in $\Sp_{12}(2)$.
Finally, suppose that $n = 14$. Then $A_{16}$ embeds as a maximal subgroup
into $\Omega_{14}^\epsilon(2)$. As $[ 16!/2^{43} ] = 2$, we are done as in
the case $n = 12$. Our proof is complete now.
\end{prf}

\begin{lem}
If $q > 2$ is even and $(G,H) = (\Sp_6(q), G_2(q))$, then~$G$ does not have an 
irreducible imprimitive $KG$-module with block stabilizer~$H$.
\end{lem}
\begin{prf}
According to our assumption at the beginning of this chapter, 
$\mbox{\rm char}( K )$ is odd, as~$q$ is even. Then the smallest dimension 
of a non-trivial irreducible $KH$-module equals $q^3 - 1$ 
(see \cite[Section~$9.2$, p.~$126$]{HiHa}).
We have
\begin{eqnarray*}
\frac{[G:H](q^3-1)}{[G:B]} & = & 
\frac{q^3(q^4-1)(q^3-1)}{(q+1)^3(q^2+1)(q^4+q^2+1)} \\
& = & \frac{1}{q^4(q-3) - q - 1} > 1
\end{eqnarray*}
for $q > 4$.
Hence, by Lemma~\ref{Seitzlem}, no nontrivial $KH$-module induces to an 
irreducible $KG$-module. Clearly, the trivial $KH$-module does not induce
to an irreducible module, and our claim follows.
\end{prf}

We finally rule out the remaining possibilities.

\begin{lem}
If $(G,H) \in \{(\SL_2(11),2.A_5), (\SU_3(5),3.A_7)\}$, then~$G$ does not 
have an irreducible imprimitive $KG$-module with block stabilizer~$H$.
\end{lem} 
\begin{prf}
If $(G,H) = (\SU_3(5),3.A_7)$, then $[\SU_3(5):3.A_7)] = 50$ but $G$ does not 
have irreducible ordinary or $\ell$-modular characters whose degree is divisible 
by $50$ (see \cite{modATL}), hence the conclusion. If 
$(G,H) = (\SL_2(11),2.A_5)$, then $[G:H] = 11$ and~$G$ has a unique irreducible
character of degree~$11$ for $\ell \neq 2, 3$, and no irreducible character of 
degree divisible by~$11$ for $\ell = 2, 3$ (again see \cite{modATL}). However~$H$ 
is perfect and again we get the desired conclusion. 
\end{prf}

We summarize the results of this subsection.

\begin{prop}
\label{SProp}
Let~$G$ be one of the following groups.

{\rm (a)} A group as in Hypothesis~\ref{ClassicalGroupCollection}.

{\rm (b)} A group $G = \Spin^\epsilon_n( q )$, $n \geq 7$, with~$q$ odd,
$\epsilon \in \{ -1, 0, 1 \}$ and $(n,q) \neq (7,3)$.

Let~$H$ be a maximal subgroup of~$G$ of ${\mathcal S}$-type.
Then~$H$ is not the stabilizer of an imprimitivity decomposition of an
irreducible $KG$-module.
\end{prop}
\begin{prf}
This follows from the above lemmas.
\end{prf}

\section{The exceptional groups of Lie type}
\label{ProofExceptionalGroups}

Here, we prove Theorem~\ref{ReductionToParabolic} for quasisimple groups~$G$ 
of exceptional Lie type, where we exclude the groups already considered in 
Chapter~\ref{ExceptionalMultipliers}, i.e., 
$$G \not\in \{ {^2\!B}_2( 8 ), {^2\!G}_2( 3 )', G_2( 2 )', G_2( 3 ), G_2( 4 ),
{^2\!F}_4( 2 )', F_4( 2 ), {^2\!E}_6( 2 ) \}.$$ 
The Suzuki group ${^2\!B}_2( 2 )$ is excluded as it is not quasisimple.
Throughout,~$H$ denotes a non-parabolic maximal subgroup of~$G$.

\subsection{The Case $G_2(q)$.}\label{CaseG2}

Let $G = G_2(q)$, $q > 4$. The maximal subgroups
of~$G$ have been classified, see Aschbacher~\cite{Asch2}, 
Kleidman~\cite{kleidman:2}, Cooperstein~\cite{CoopersteinG2} and 
Migliore~\cite{Migliore}. The ordinary character degrees of~$G$ are also known 
(see the web page \cite{LuebeckWeb}
by Frank L\"ubeck). By inspection we see that the only possibilities for our 
block stabilizer~$H$ are $\SL_3(q).2$ and $\SU_3(q).2$. Moreover, inspection 
of the character degrees of~$G$ and $H$ shows that if~$M_1$ is an irreducible 
$KH$-module such that $\Ind_H^G(M_1)$ is irreducible, then $M_1$ is 
$1$-dimensional, and hence $\mbox{\rm dim}_K(\Ind_H^G(M_1)) = [G:H]$. If 
$\mbox{\rm char}(K) = 0$, we use the table of character degrees of~$G$ to 
show that there is no irreducible character $\chi$ of~$G$ with $\chi(1) = 
\frac{1}{2}\,q^3\,(q + 1)\,(q^2 - q + 1)$
or $\frac{1}{2}\,q^3\,(q - 1)\,(q^2 + q + 1)$, i.e., such an~$M_1$
does not exist. By Lemma~\ref{Modp}, this implies that no such~$M_1$ exists 
for~$K$ of any characteristic. 

\subsection{The Case $F_4(q)$.}\label{CaseF4}
Let $G = F_4(q)$, $q > 2$.
The non-parabolic subgroups $H$ whose index is smaller than the largest 
character degree of~$G$ are $\Spin_9(q)$, $N_G(\Spin_8^+(q))$, and 
${^3\!D}_4(q).3$ (see Liebeck-Saxl \cite{liebeck:2} in conjunction with 
Lemma~\ref{Seitzlem}). 

We first consider the case $G = F_4(q)$, $q$ odd, and $H = \Spin_9(q)$.
Here $H = C_G(a)$ for an involution $a \in G$ whose trace on the 
$26$-dimensional ${\mathbb F}_qG$-module~$N$ is $-6$. It is well known 
that $G$ contains exactly two classes of involutions which are distinguished 
by their traces on~$N$. It is also known that~$G$ contains a subgroup 
isomorphic to $\SO_3(q) \times G_2(q)$ whose action on~$N$ is $(N_3 \otimes 
N_7) \oplus N_5$, where~$N_3$ and~$N_7$ are the natural modules for $\SO_3(q)$ 
and $G_2(q)$, respectively, and~$N_5$ is the $\SO_3(q)$-module of bivariate 
homogeneous polynomials of degree~$4$ (and on which $G_2(q)$ acts trivially). 
For these facts see \cite[p.~$60$]{Asch2b}. Thus the trace of an involution 
from $\Omega_3(q)$ on~$N$ is $-6$, and hence it is conjugate to~$a$.
Thus every involution in $\Omega_3(q)$ is conjugate to~$a$. Now $\Omega_3(q)
\cong \PSL_2(q)$ contains two involutions whose product has odd order $\neq 1$.
By Lemma~\ref{2.3replace} (with $s = 1$), an irreducible $KH$-module will 
never induce to an irreducible $KG$-module.

In case~$q$ even and $H = \Spin_9(q) = \Sp_8(q)$, there does not exist an 
element~$t$ such that $t \in C_G({H^t} \cap H)$; thus we argue as follows. 
The index of $H$ in $G$ equals $q^8(q^8 + q^4 + 1)$.

Suppose that $\text{\rm char}(K)$ is~$0$ or odd, and let $M_1$ be an irreducible
$KH$-module with $M := \Ind_H^G( M_1 )$ irreducible. Since $H$ is perfect,
the dimension of~$M_1$ is larger than~$1$. Also, $\text{\rm dim}(M_1)
< (q^4 - 1)(q^3 - 1)q^2$, since otherwise $\text{\rm dim}( M ) > [G\colon\!B]$,
where~$B$ denotes a Borel subgroup of~$G$, contradicting the upper bound of
Seitz given in Lemma~\ref{Seitzlem}. Hence, by the result 
\cite[Theorem~$1.1$]{GuraTiep} of Guralnick and Tiep, $M_1$ is a Weil module.

Suppose now that $\text{\rm char}( K ) = 0$. The degrees of the ordinary Weil 
characters of~$H$ are given in \cite[Table~$1$]{GuraTiep} as polynomials (over 
the rationals) in~$q$.  The degrees of the ordinary irreducible characters 
of~$G$ can be computed from Lusztig's Jordan decomposition of characters. These 
degrees, again as polynomials in~$q$, are explicitly given on the web site 
\cite{LuebeckWeb} of Frank L\"ubeck. Using these data together with 
Lemma~\ref{LemmaOnPolynomials}, one can easily show that no ordinary irreducible 
character~$\chi$ of~$G$ satisfies $\chi(1) = [G\colon\!H]\psi(1)$ for some 
ordinary Weil character~$\psi$ of~$H$.

Hence $\text{\rm char}( K ) = \ell \neq 0$. Let~$\varphi$ denote the 
Brauer character of~$M_1$. By \cite[Table~$1$]{GuraTiep}, there is an ordinary 
Weil character~$\psi$ of~$H$ such that $\varphi$ occurs in the reduction
modulo~$\ell$ of~$\psi$ and such that $\varphi( 1 ) = \psi( 1 )$ or $\varphi( 1 )
= \psi( 1 ) - 1$. If $\varphi( 1 ) = \psi( 1 )$, the character $\varphi$ is
liftable to $\psi$, and thus $\Ind_H^G( \varphi )$ is reducible by 
Lemma~\ref{Modp} and by the above result in the $\text{\rm char}( K ) = 0$ case.
Thus $\varphi( 1 ) = \psi( 1 ) - 1$.

Since induction commutes with reduction modulo~$\ell$, there is an irreducible
constituent $\chi$ of $\Ind_H^G( \psi )$ such that the irreducible Brauer
character $\Ind_H^G( \varphi )$ of~$M$ is a constituent in the reduction
modulo~$\ell$ of~$\chi$. Hence
\begin{equation}
\label{SqueezedInDegree}
[G\colon\!H]( \psi( 1 ) - 1 ) \leq \Ind_H^G( \varphi )( 1 ) \leq \chi( 1 )
\leq [G\colon\!H]\psi( 1 ).
\end{equation}
Using Lemma~\ref{LemmaOnPolynomials} once more, one can show that no ordinary 
irreducible character~$\chi$ of~$G$ satisfies~(\ref{SqueezedInDegree}). This 
contradiction proves that~$H$ is not a block stabilizer of some imprimitive 
irreducible $KG$-module.

Next, let $H = N_G(\Spin_8^+(q)) = \Spin_8^+(q).S_3$. Using the result of
Seitz given in Lemma~\ref{Seitzlem}, together with the bounds on minimal
degrees by Landazuri and Seitz \cite{LaSe} for $\Spin_8^+(q)$, we find that
if~$M_1$ is an irreducible $KH$-module such that $\Ind^G_H( M_1 )$ is
irreducible, then $\Spin_8^+(q)$ is in the kernel of~$M_1$, i.e.,~$M_1$ may
be viewed as a $KS_3$-module and $\mbox{\rm dim}_K( M_ 1 ) \leq 2$. Since~$S_3$
is solvable we may assume that $\mbox{\rm char}(K) = 0$ in order to rule out~$H$
as a block stabilizer (cf.\ Lemma~\ref{Modp}). If~$q$ is not divisible by~$2$
or~$3$, then~$q^{12}$ divides $[G\colon\!H]$. In any case~$q^{11}$ divides
this index. Using Frank L\"ubeck's explicit list of character degrees of~$G$
(see \cite{LuebeckWeb}), we find that no degree of an ordinary irreducible
character of~$G$ is divisible by $[G\colon\!H]$.

Finally, suppose that $H = {^3\!D}_4(q).3$. We choose a subgroup~$Z$ of~$H$
isomorphic to~$G_2(q)$. We note that~$Z$ is unique up to conjugacy in~$H$ 
and that $C_H(Z)Z$ is the unique maximal overgroup of~$Z$ in~$H$ (see 
\cite{kleidman:3}). Also, $C_H(Z) = \langle z \rangle$ is cyclic of order~$3$.
We claim that $C_G(Z) \cong \SO_3(q)$.
To see this note that the action of~$H$ on the $27$-dimensional module~$V$ 
of $E := \mbox{\rm gcd}(3, q - 1).E_6(q)$ is described in 
\cite[Lemma $7.3$]{Asch2b}. In particular we see 
that $\mbox{\rm dim}(C_V(H')) = 3$, every nonzero vector in $C_V(H')$ is 
centralized by a subgroup of~$E$ isomorphic to~$G$, and $[H',V]$ is an 
irreducible $H'$-module of dimension~$24$, but not absolutely irreducible. 
In fact it is the $8$-dimensional module of~$H'$ over the field of oder~$q^3$. 
Thus $\mbox{\rm dim}(C_{[H',V]}(Z)) = 3$ and hence $\mbox{\rm dim}(C_V(Z)) = 6$.

Now \cite[Lemma $7.1$]{Asch2b} implies that $C_V(Z)$ can not consist entirely 
of dark points. Cohen-Cooperstein \cite[p.~$467$]{CohenCooper} 
describe the $E$-orbits on the one-dimensional subspaces of~$V$. In their 
terminology the 
dark points are the black points. The singular points of Aschbacher are the 
white points, whereas Aschbacher's brilliant non-singular points are the gray 
points. The gray and white points are distiguished by the $E$-invariant map 
$\#: V \rightarrow V$. A simple calculation using the representative $e_1 + e_2$
shows that $(e_1 + e_2)^{\#} = e_3$, a white point. Translating this back into 
Aschbacher's language we see that if~$Z$ fixes a point in~$V$ which is
not dark, then~$Z$ fixes a singular point. Thus $C_V(Z)$ contains singular 
points. The stabilizer of a singular point $v \in C_V(Z)$ is a parabolic 
subgroup~$P$ of~$E$ of type~$D_5$. Now $q^2\colon\!\SL_2(q) \leq C_P(Z)$
and thus, as $C_E(Z) \geq C_E(H')$, which is cyclic of order $q^2+q+1$, we 
see that $C_E(Z)Z$ is the maximal subgroup of~$E$ listed in [82, Table III].

The action of $C_E(Z)Z$ on~$V$ is described explicitly in \cite[5.7]{Asch2b}. 
In particular, we see in \cite[4.5.7]{Asch2b}, that $C_V(Z)$ is an irreducible 
$C_E(Z)$-module when~$q$ is odd.
A six dimensional $C_E(Z) = \SL_3(q)$-module is necessarily a Galois twist of 
the symmetric square of the natural $\SL_3(q)$-module. Using the notation in 
\cite[4.5]{Asch2b}, the vector $d = v_1 + v_2 + v_3$ has the property that 
$G = C_E(d)$. Now~$d$ represents a non-degenerate symmetric bilinear form of 
the natural $C_E(Z)$-module and hence $C_Z(E) \cap G = \SO_3(q)$.

When~$q$ is even, then we see, using \cite[4.5.7]{Asch2b}, that $Z \leq 
C_G(v_1) \leq G$, and $C_G(v_1) \cong \Sp_8(q) \cong \Spin_9(q)$. Thus 
$C_{C_G(v_1)}(Z) \cong \SL_2(q) \cong \SO_3(q)$ as~$q$ is even. The centralizer 
of~$Z$ in~$G$ does not grow. To see this, using the notation of 
\cite[4.5.7]{Asch2b}, we note that $G = C_E(d)$, with $d = v_1 + v_2 + v_3$ and 
that $Q(d,-)$ (for the definiton of~$Q(x,y)$ see \cite[Section~$2$]{Asch2a}) 
restricted to the $C_E(Z)$-module $\langle w_1,w_2,w_3 \rangle,$ see 
\cite[5.7.8]{Asch2b}, is a nondegenerate quadratic form. 

Thus our claim is established, i.e., $C_G(Z) \cong \SO_3(q)$.
Since~$q > 2$, we can choose $1 \neq t \in C_G(Z)$ such that $C_H(Z) \cap 
C_H(Z)^t = 1$.
We claim that $H \cap {H^t} = Z$. If not, $H \cap {H^t}
= Z \times \langle z \rangle$, and thus $z \in {H^t}$. It follows
that $z \in C_{H^t}( Z ) = \langle z \rangle^t$, a contradiction.
We conclude that~$t$ centralizes $H \cap {H^t}$ and thus, as $t \not\in H$, 
no irreducible $KG$-module is induced from a $KH$-module.

\subsection{The Case $E_6(q)$ or ${^2\!E}_6(q)$}\label{CaseE6}
Let $G = d.E_6(q)$ with $d = \mbox{\rm gcd}(3,q-1)$ or $G = d.{^2\!E}_6(q)$ 
with $d = \mbox{\rm gcd}(3,q+1)$ and $q \neq 2$. Again, by 
Lemma~\ref{Seitzlem} and the results of Liebeck and Saxl \cite{liebeck:2}, 
the non-parabolic subgroups~$H$ of~$G$, whose index is smaller than the 
largest character degree of~$G$ are $Z(G)F_4(q)$ and, in case 
$G = d.{^2\!E}_6(q)$, also $N_G(\Spin_{10}^-(q))$. (To rule out the
particular case $G = E_6(2)$ and $H = \SL_2( 2 ) \times \SL_6( 2 )$ we use
the ordinary character degrees of~$G$ available in GAP \cite{GAP4}.)

Let $H = Z(G)Z$ with $Z \cong F_4( q )$. We may assume that $Z = C_G(\sigma)$ 
for the graph automorphism~$\sigma$ of~$G$ (see 
\cite[Table~4.5.2 and Proposition 4.9.2(b)(4)]{GolYSo3}). There is a 
$\sigma$-invariant Levi subgroup~$L$ of~$G$ of type $A_5(q)$ or ${^2\!A}_5(q)$, 
respectively, and~$\sigma$ acts as a graph automorphism on~$L$. Next, there is 
a non-trivial element $t \in L$ of order prime to~$6$, which is inverted 
by~$\sigma$; for example we can choose~$t$ to be a non-trivial element of order 
prime to~$6$ of a $\sigma$-invariant Singer cycle of~$L$. Now write
$\,\bar{}: G \rightarrow G/Z(G)$ for the canonical epimorphism. 
As $\bar{t} \not\in \bar{H}$, Lemma~\ref{2.3replace} (with $s = 1$) shows that 
no irreducible $K\bar{H}$-module can induce to an irreducible $K\bar{G}$-module.
By Corollary~\ref{spinlift}, the same conclusion holds for~$H$.

Now let $G = d.{^2\!E}_6(q)$ and $H = N_G(\Spin^-_{10}(q))$. If~$q$ is odd, 
$H = C_G(a)$, 
where~$a$ is an involution of the subgroup $F_4(q) \leq G$ with centralizer 
$\Spin_9(q)$. In Subsection~\ref{CaseF4} we have constructed two conjugates 
of~$a$ in $F_4(q)$ whose product has odd order~$\neq 1$. Thus we are done with
Lemma~\ref{2.3replace} (with $s = 1$).

If~$q$ is even, $H = \Spin^-_{10}( q ) \times (q + 1 )$ (see
\cite[Table~$5.1$]{LiebeckSaxlSeitz}). Using Lemma~\ref{Seitzlem}, together
with the Landazuri and Seitz bounds \cite{LaSe} for $\Spin_{10}^-(q)$,  
we find that if~$M_1$ is an irreducible $KH$-module
inducing to an irreducible $KG$-module, then $\mbox{\rm dim}_K( M_1 ) = 1$.
Now use the table of the irreducible character degrees of~$G$ provided on
Frank L\"ubeck's web page (see \cite{LuebeckWeb}), to prove that $[G\colon\!H]$ 
does not divide any such degree if $\mbox{\rm char}( K ) = 0$. The same 
conclusion then holds for any characteristic of~$K$ by Corollary~\ref{spinlift}.

\subsection{The Case $E_8(q)$}\label{CaseE8}
Let $G =E_8(q)$. The non-parabolic subgroups~$H$ whose index is smaller than 
the largest character degree of~$G$ are $H = (\SL_2(q) \circ E_7(q)).2$ if~$q$ 
is odd, and $H = \SL_2(q) \times E_7(q)$ if~$q$ is even (see Liebeck-Saxl 
\cite{liebeck:2} in conjunction with Lemma~\ref{Seitzlem}). The first factor~$S$ 
of $F^*(H)$, isomorphic to $\SL_2(q)$, is generated by long root subgroups. Now 
let~$L$ be a subgroup of~$E_8(q)$ isomorphic to $\SL_3(q)$,
generated by long root subgroups and  
containing~$S$. It is known (see e.g., \cite[Table 4.7.3A]{GolYSo3})
that $N_G(L) = L \times E_6(q)$ if $q \not\equiv 1 (\mbox{\rm mod\ }3)$ and 
$N_G(L) = (L \circ E_6(q))\colon\!\langle s \rangle$, otherwise, where~$s$ 
induces a diagonal automorphism of order~$3$ on each factor. Clearly~$L$ 
respectively $\langle L,s \rangle$ is isomorphic to a subgroup of $\GL_3(q)$, 
so Lemma~\ref{gl3} guarantees the existence of an element $t \in L$ such 
that~$t$ centralizes $N_L(S) \cap N_L(S)^t$ respectively 
$N_{\langle L,s \rangle}(S) \cap N_{\langle L,s \rangle}(S)^t$. Thus the 
hypotheses of 
Lemma~\ref{2.3generalized2} are satisfied (with $s = 1$ if 
$q \not\equiv 1 (\mbox{\rm mod\ }3)$), 
showing that~$t$ centralizes 
$H \cap {H^t}$. Thus~$H$ is not the block stabilizer of an imprimitivity 
decomposition of any irreducible $KG$-module.

\subsection{The Case $E_7(q)$}
Let $G = d.E_7(q)$ with $d = \mbox{\rm gcd}(2, q - 1)$. To simplify the 
exposition, we put $E^+_6(q) := \mbox{\rm gcd}(3,q-1).E_6( q )$, and 
$E_6^-( q ) := \mbox{\rm gcd}(3,q+1).{^2\!E}_6( q )$. Using once more 
the result of Liebeck and Saxl \cite{liebeck:2} together with 
Lemma~\ref{Seitzlem}, we find that the non-parabolic subgroups~$H$ whose 
index is smaller than the largest character degree of~$G$ are $H = 
N_G( {E}^\epsilon_6(q) )$, $\epsilon \in \{ +, - \}$, and $H = N_G(S)$, 
where $S \cong \SL_2(q)$ is generated by long root subgroups of~$G$. 
If $H = N_G(S)$, then the 
argument from Subsection~\ref{CaseE8} goes through verbatim with one minor 
difference. Namely that now $C_G(L)$ is $A_5(q)$ rather than $E_6(q)$. 

Let $H = N_G( {E}^\epsilon_6(q) )$. Then $H = 
(S \circ {E}^\epsilon_6(q)).\langle \sigma \rangle = N_G( S )$, where~$S$ is a 
cyclic torus of order $q - \epsilon \cdot 1$, and $(S \circ {E}^\epsilon_6(q))$ 
is a (not necessarily split) Levi subgroup of~$G$ (see 
\cite[Table~$1$]{liebeck:2}; notice that we have incorporated the~$e_{\pm 1}$ 
of \cite[Table~$1$]{liebeck:2} in our notation for ${E}^\epsilon_6(q)$). 
Moreoever,~$\sigma$ induces the graph automorphism on ${E}^\epsilon_6(q)$. This
follows, e.g., from the fact that $\langle \sigma \rangle$ is the relative Weyl 
group of the Levi subgroup $S \circ {E}^\epsilon_6(q)$, and the action of
$\langle \sigma \rangle$ on the latter can be determined by a computation in
the Weyl group~$W$ of~$G$. Indeed, the stabilizer in~$W$ of the set of six 
fundamental roots forming a subsytem of type~$E_6$, is of order~$2$ and acts 
as the graph automorphism on this set of six roots.
As already observed in Subsection 
\ref{CaseE6}, we have $Z := C_{{E}^\epsilon_6(q)}( \sigma ) \cong F_4(q)$.
Let $\,\bar{}: G \rightarrow G/Z(G)$ denote the canonical epimorphism. We 
identify~$Z$ with its image in~$\bar{G}$. Then $N_{\bar{G}}( Z ) = L \times Z$
with $L \cong \PSL_2( q )$.
For $q > 3$ this follows from 
\cite[Theorem 2(e)]{LiebeckSeitz}. For $q = 2, 3$, we may apply the same 
theorem to the groups $E_7( q^2 )$ and consider the fixed points of the 
field automorphism. As~$\bar{S}$ centralizes~$Z$, we have $\bar{S} \leq L$. 
Now suppose that $q > 3$. If $q = 5$, assume that $\epsilon = -$.
Then Lemma~\ref{PGL2} guarantees the existence of ${t} \in {L}$ with 
${L} = \langle \bar{S}, {\bar{S}^{t}} \rangle$ and ${t}$ centralizing 
$N_{{L}}( \bar{S} ) \cap N_{{L}}( \bar{S} )^{t}$. It follows from 
Lemma~\ref{2.3generalized2} that $t$ centralizes $\bar{H} \cap {\bar{H}^{t}}$, 
and  so no irreducible $K\bar{H}$-module induces to an irreducible 
$K\bar{G}$-module. As we may choose such a $t$ of odd order if~$q$ is odd, the 
same conclusion holds for~$H$ and~$G$ by Corollary~\ref{spinlift}.

Now suppose that $q \in \{ 3, 5 \}$ and $\epsilon = +$. (The case $q = 3$ and 
$\epsilon = -$ will be treated below.) If $q = 3$, we have 
$\bar{H} = E_6( 3 ).2$, and if $q = 5$ we have $\bar{H} = (E_6(5) \times 2).2$.
We thus have 
$\bar{H} = (E_6(q) \times \langle x \rangle ).\langle \bar{\sigma} \rangle$, 
where $q = 3$ or~$5$, and $x \in C_G( E_6( q ) )$ is trivial if $q = 3$ 
and and has order~$2$, otherwise. 
Now $N_{\bar{G}}( Z ) = L \times Z$ with $\langle x, \bar{\sigma} \rangle 
\leq C_{\bar{G}}( Z ) = L \cong \PSL_2( q )$. Clearly, there is an element
$t \in L$ of order~$3$ such that $\langle x, \bar{\sigma} \rangle \cap 
\langle x, \bar{\sigma} \rangle^t$ is trivial. As $t \not\in \bar{H} = 
N_{\bar{G}}( \bar{H}' )$, the intersection $\bar{H} \cap \bar{H}^t$ does not 
contain $\bar{H}' \cong E_6( q )$. 
As $C_{\bar{H}}(Z) = \langle x, \bar{\sigma} \rangle$ and
$Z \times C_{\bar{H}}(Z)$ is the unique maximal subgroup of~$\bar{H}$
containing~$Z$ but not $\bar{H}'$, we have

$$\bar{H} \cap \bar{H}^t \leq Z \times C_{\bar{H}}(Z) \leq Z \times L = N_{\bar{G}}(Z)$$
and
$$\bar{H} \cap N_{\bar{G}}(Z) = Z \times \langle x, \bar{\sigma} \rangle.$$
This implies
\begin{eqnarray*}
\bar{H} \cap \bar{H}^t  & = & \bar{H} \cap \bar{H}^t \cap N_{\bar{G}}(Z) \\
            & = & [\bar{H}  \cap N_{\bar{G}}(Z)] \cap [ \bar{H}^t \cap N_{\bar{G}}(Z)] \\
            & = & Z \times 
                  [ \langle x, \sigma \rangle \cap \langle x, \sigma \rangle^t ] \\
            & = & Z 
\end{eqnarray*}
by our choice of~$t$.
As~$t$ centralizes~$Z$, it follows that $t \in C_{\bar{G}}( \bar{H} \cap \bar{H}^t )$.
Since~$t$ has order~$3$, Lemma~\ref{Mackey1} and Corollary~\ref{spinlift} show 
that~$H$ is not the block stabilizer of an imprimitive irreducible $KG$-module.

Next, assume that $q = 3$ and $\epsilon = -$. Then 
$H = ({^2\!E}_6(3) \times 4).2$ and $\bar{H} = C_{\bar{G}}( x )$ for an 
involution $x \in \bar{G}$. There is a subgroup $Y \leq Z \leq {^2\!E}_6(3)$ 
with $Y \cong G_2( 3 )$, such that $C_{\bar{G}}( Y ) \cong \PSp_6( 3 )$ 
(this follows from \cite[Theorem 1(b), Theorem 2(e)]{LiebeckSeitz}).
A computation of class multiplication coefficients shows that there is a 
conjugate~$x'$ of~$x$ in $C_{\bar{G}}( Y )$ such that the product $xx'$ has 
order~$3$. Putting $t := xx'$, Lemma~\ref{2.3replace} and Corollary~\ref{spinlift}
prove that~$H$ is not a block stabilizer in~$G$.

Finally, suppose that $q = 2$. Then $G = E_7( 2 )$ and $H = E_6( 2 )\colon\!2$ 
or $H = (3.{^2\!E}_6(2))\colon\!S_3$ (see \cite[Table~$1$]{liebeck:2} or
\cite[p.~$219$]{ATLAS}). The case $H = E_6( 2 )\colon\!2$ is settled as 
follows. We have $N_G( Z ) = Z \times L$ with 
$L \cong \GL_2( 2 )$, and $N_H( Z ) = Z \times \langle \sigma \rangle$. 
Let~$t \in L$ be an element of order~$3$. Then $t \not\in H$, and thus
$H \cap H^t = Z$. As $t \in C_G( Z )$, we are done.
Now let $H = (3.{^2\!E}_6(2))\colon\!S_3$. Then $H = N_G( \langle z \rangle )$,
where~$z$ is an element of order~$3$ in the center of $F^*( H ) = 
3.{^2\!E}_6(2)$. By the tables of conjugacy classes of elements of order~$3$
given on Frank L\"ubeck's web page (see \cite{LuebeckWeb}),~$z$ is the unique 
element of order~$3$, up to conjugacy, whose centralizer contains an element of 
order~$13$. Now~$G$ contains a subgroup~$D$ isomorphic to ${^3\!D}_4( 2 )$ 
centralizing a subgroup~$L$ isomorphic to~$\SL_2( 8 )$ (see 
\cite[Table~$5.1$]{LiebeckSaxlSeitz} or \cite[p.~$219$]{ATLAS}). As~$D$ 
contains an element of order~$13$, we may assume that $z \in L$.
We have $N_G(L) = (L \times D).\langle s \rangle$ with~$s$ an element of 
order~$3$ normalizing~$L$ and~$D$ and acting non-trivially on each factor. 
In particular, $\langle L, s \rangle \cong \Aut( \SL_2( 8 ) ) \cong 
{^2\!G}_2( 3 )$. A computation inside ${^2\!G}_2( 3 )$ shows that there 
exists an element $t \in L$ such that $L = \langle z, z^t \rangle$ and such 
that~$t$ centralizes $N_{\langle L, s \rangle}( \langle z \rangle ) \cap
N_{\langle L, s \rangle}( \langle z \rangle )^t$. Lemma~\ref{2.3generalized2} 
implies that~$t$ centralizes $H \cap H^t$ and thus~$H$ is not a block 
stabilizer in~$G$. 

\subsection{The remaining cases}\label{Remaining}
Finally we consider the Steinberg triality groups ${^3\!D}_4(q)$,
the Suzuki groups ${^2\!B}_2( 2^{2m + 1} )$, $m \geq 1$, 
and the Ree groups ${^2\!F}_4( 2^{2m + 1} )$, $m \geq 1$,
and ${^2\!G}_2( 3^{2m + 1} )$, $m \geq 1$.
(Recall that the Tits group ${^2\!F}_4( 2 )'$ has already been considered 
in Proposition~\ref{TitsGroup}.)
The maximal subgroups of the Steinberg triality groups were 
determined by Kleidman \cite{kleidman:3},
those of the large Ree groups by Malle \cite{mallemax},
those of the small Ree groups by Kleidman \cite{kleidman:2} and
those of the Suzuki groups by Suzuki \cite{Suzuki}.
In every case the smallest index of a non-parabolic maximal subgroup 
exceeds the bound of Seitz given in Lemma~\ref{Seitzlem}, eliminating
these cases.

\chapter{Groups of Lie type: Induction from parabolic subgroups}
\label{InductionFromParabolicSubgroups}

\markleft{7.\ INDUCTION FROM PARABOLIC SUBGROUPS}
Here, we begin the investigation of imprimitive irreducible modules
whose block stabilizers are parabolic subgroups. First we show that
such an irreducible module is Harish-Chandra induced.
Then we give two sufficient conditions for the irreducibility of a
Harish-Chandra induced module. Our first condition is in terms 
of Harish-Chandra theory, and the second one in terms of 
Deligne-Lusztig theory. The two conditions are sufficient to show 
that the bulk of all cross-characteristic irreducible modules
of groups of Lie type are in fact Harish-Chandra induced, hence
imprimitive.
 
\section{Harish-Chandra series}
In this section we let~$G$ be a finite group with a split $BN$-pair
of characteristic~$p$, satisfying the commutator relations (see e.g., 
\cite[\S\S65,69]{cr1} for the definition and principal results arising 
from this set of axioms). By a parabolic subgroup, respectively Levi
subgroup of~$G$ we mean some $N$-conjugate of a standard parabolic
subgroup, respectively standard Levi subgroup of~$G$. As above, we 
let~$K$ denote an algebraically closed field of characteristic 
$\ell \geq 0$, with $\ell \neq p$.

\begin{prop}
\label{ParabolicStabilizer}
Let~$M$ be an irreducible $KG$-module which is imprimitive with block
stabilizer~$P$, where~$P$ is a parabolic subgroup of~$G$. Then 
$M = R_L^G( M_1 )$ for some $KL$-module~$M_1$, where~$L$ is a Levi 
complement in~$P$. In other words,~$M$ is Harish-Chandra induced.
\end{prop}
\begin{prf}
In this proof we write $[X,Y] := \mbox{\rm dim}_K \Hom_{KH}( X, Y )$ for 
subgroups $H \leq G$ and $KH$-modules~$X$ and~$Y$.

Let~$M_2$ be an irreducible $KP$-module such that $\Ind_P^G(M_2) \cong M$.
Let $Q$ be the parabolic subgroup opposite to~$P$.
Thus $P \cap Q = L$ is a Levi complement of~$P$ and of~$Q$.
If~$M_1$ is an irreducible constituent of the head of $\Res^P_L(M_2)$,
then $[ \Res_L^P( M_2 ), M_1 ] > 0$. Denote the inflation of~$M_1$
to a $KQ$-module by~$\tilde{M}_1$. Let~$D$ be a set of double coset 
representatives, containing~$1$, for the $P$-$Q$-double cosets of~$G$.
 By Mackey's theorem,
\begin{eqnarray*}
[ \Ind_P^G( M_2 ), \Ind_Q^G ( \tilde{M}_1 ) ] & = &
[ \Res_Q^G \Ind_P^G ( M_2 ), \tilde{M}_1 ] \\
 & = & \sum_{x \in D} [ \Ind_{{P^x}\cap Q}^Q \Res_{{P^x}\cap Q}^{{P^x}}
(M^x_2 ), \tilde{M}_1 ] \\
 & = & \sum_{x \in D} [ \Res_{{P^x}\cap Q}^{{P^x}} (M^x_2 ),
\Res_{{P^x}\cap Q}^{Q} ( \tilde{M}_1 ) ].
\end{eqnarray*}
Our choice of~$M_1$ implies that the 
summand corresponding to $x = 1$ in the last line of the formula above is 
non-zero. It follows that $[ \Ind_P^G( M_2 ), \Ind_Q^G ( \tilde{M}_1 ) ] > 0$.

The latter implies that $\Ind_P^G(M_2)$ is isomorphic to a submodule of 
$\Ind_Q^G(\tilde{M}_1)$, as $\Ind_P^G(M_2)$ is irreducible. On the other hand, 
our choice of~$M_1$ and the
fact that~$P$ and~$Q$ have the same order, implies that
$\mbox{\rm dim}_K ( \Ind_Q^G(\tilde{M}_1) ) \leq \mbox{\rm dim}_K ( \Ind_P^G(M_2) )$.
Thus $\Res_L^P ( M_2 ) = M_1$ and $\Ind_P^G(M_2) \cong \Ind_Q^G(\tilde{M}_1)$. 
Since $R_L^G( M_1 ) \cong \Ind_Q^G(\tilde{M}_1)$, the claim follows.
\end{prf}

An ordinary irreducible $KG$-module is called {\em Harish-Chandra
imprimitive}, if it is Harish-Chandra induced from a proper Levi 
subgroup of~$G$. Otherwise it is called {\em Harish-Chandra primitive}. 
By the above proposition, an irreducible $KG$-module is imprimitive
with parabolic block stabilizer, if and only if it is Harish-Chandra
imprimitive.

We now apply Harish-Chandra theory to produce Harish-Chandra imprimitive 
irreducible modules. If~$L$ is a Levi subgroup of~$G$ we 
denote Harish-Chandra induction from $L$ to~$G$ by~$R_L^G$.
Thus if $M$ is a $KL$-module, then $R_L^G(M) = \Ind_P^G( \tilde{M} )$,
where~$P$ is a parabolic subgroup of~$G$ with Levi complement~$L$,
and $\tilde{M}$ is the inflation of~$M$ to~$P$.
It is known that, up to isomorphism, $R_L^G(M)$ is independent of
the choice of the parabolic subgroup~$P$ with Levi complement~$L$
(see \cite{DD,HL3}).
If~$M$ is irreducible and cuspidal, $W_G(L,M) := (N_G(L,M) \cap N)L/L$ denotes 
the ramification group of $R_L^G(M)$ (here,~$N$ is the subgroup of~$G$
specified by the $BN$-pair; see \cite[Section~3]{GeHiMa2}).

\begin{thm}
\label{HCresult}
Let~$L_0$ be a Levi subgroup of~$G$ and let~$M_0$ be an 
irreducible cuspidal $KL_0$-module. Suppose that~$L$ is a Levi 
subgroup of~$G$ with $L_0 \leq L$ and $W_G(L_0,M_0) =
W_L( L_0, M_0 )$. 
Then $R_L^G(M_1)$ is irreducible for every irreducible quotient or 
submodule $M_1$ of $R_{L_0}^L(M_0)$.
\end{thm}
\begin{prf}
The dimension of $\End_{KG}( R_{L_0}^G(M_0) )$ equals the order of 
$W_G(L_0,M_0)$ (see e.g., \cite[Theorem 2.9]{DiFl}).
Moreover, $\End_{KL}( R_{L_0}^L(M_0) )$ is embedded into
$\End_{KG}( R_{L_0}^G(M_0) )$ as a unital subalgebra (see
\cite[2.5]{GruHi}). Since the two algebras have the same
dimension, they are isomorphic.

The $(L_0,M_0)$ Harish-Chandra series of~$L$ consists of the set of 
irreducible submodules of $R_{L_0}^L(M_0)$ (up to isomorphism),
and this set is equal to the set of irreducible quotients of
$R_{L_0}^L(M_0)$ (see \cite[2.2]{GeHi} for more explanations and
references). The result now follows from \cite[Proposition 2.7]{GeHiMa2}.
\end{prf}

\noindent Notice that the assumption of the theorem is 
satisfied for every intermediate Levi subgroup $L_0 \leq L \leq G$ 
if $W_G( L_0, M_0 ) = \{ 1 \}$.

\section{Lusztig series}\label{LusztigSeries}
Now we give a sufficient condition for Harish-Chandra imprimitivity
of an irreducible module of a finite group of Lie type in terms of
Deligne-Lusztig theory, more precisely, in terms of Lusztig series 
of characters. We introduce the relevant notation. Our basic references
for the terminology and notation in this section are \cite{cart,DiMi2}.
Let~${\mathbf G}$ be a connected reductive algebraic group over the 
algebraic closure of the prime field~${\mathbb F}_p$ and let~$F$ be a Frobenius 
morphism of~${\mathbf G}$. Let $G = {\mathbf G}^F$ denote the corresponding 
finite group of Lie type.

Let $\mathbf{G}^*$ denote the group dual to $\mathbf{G}$ in
the sense of Deligne and Lusztig (see \cite[Chapter~$4$]{cart}). By
abuse of notation, the Frobenius morphism of $\mathbf{G}^*$ induced 
by~$F$ is also denoted by~$F$. The finite group in duality with~$G$ 
is $G^* := {\mathbf{G}^*}^F$. The set of ordinary irreducible characters 
of~$G$ is partitioned into rational Lusztig series $\mathcal{E}( G, [s])$, 
where $[s]$ ranges over the ${G^*}^F$-conjugacy classes of semisimple 
elements of~${G^*}^F$ (see \cite[Proposition 14.41]{DiMi2}).

An $F$-stable Levi subgroup of $\mathbf{G}^*$ is called {\em split},
if it is the Levi complement of an $F$-stable
parabolic subgroup of~$\mathbf{G}^*$. The $F$-fixed points of split
Levi subgroups of $\mathbf{G}^*$ are Levi subgroups of
$G^*$ in the sense of groups with split $BN$-pairs.

Let~$s$ be a semisimple 
$\ell$-regular element in $G^*$. It was shown by Brou\'e and Michel 
(\cite{BM}, Th{\'e}or{\`e}me 2.2), that the set
\[{\mathcal E}_{\ell}(G,[s]):=\bigcup\limits_{t\in C_{G^*}(s)_{\ell}}
{\mathcal E}(G,[ts]),\]
where $C_{G^*}(s)_{\ell}$ is the set of $\ell$-elements of
$C_{G^*}(s)$, is the union of the ordinary irreducible characters in a 
set of $\ell$-blocks of $G$.

\begin{thm}
\label{FirstLusztigSeriesResult}
Let $s \in G^*$ be semisimple such that $C_{\mathbf{G}^*}(s)$ is 
contained in a proper split Levi subgroup $\mathbf{L}^*$ of 
$\mathbf{G}^*$. Let~$\mathbf{L}$ be a split Levi subgroup of
$\mathbf{G}$ dual to $\mathbf{L}^*$. 

Then every ordinary irreducible character of~$G$ contained in 
$\mathcal{E}(G,[s])$ is Harish-Chandra induced from a character
of $\mathcal{E}(L,[s])$.

If~$s$ is $\ell$-regular for some prime~$\ell$ not dividing~$q$,
then every irreducible $\ell$-modular character of~$G$ contained
in ${\mathcal E}_{\ell}(G,[s])$ is Harish-Chandra induced from
a Brauer character lying in ${\mathcal E}_{\ell}(L,[s])$.
\end{thm}
\begin{prf}
Under our assumption on~$s$,
the Lusztig map $\varepsilon_{\mathbf{G}}\varepsilon_{\mathbf{L}}
R_{L}^{G}$ induces a bijection between $\mathcal{E}( L, [s] )$ and
$\mathcal{E}( G, [s] )$ (see \cite[Theorem 13.25(ii)]{DiMi2} and 
\cite[Th{\'e}or{\`e}me 11.10]{CeBo2}).
Since $\mathbf{L}$ is split, the signs $\varepsilon_{\mathbf{G}}$
and $\varepsilon_{\mathbf{L}}$ are equal (see \cite[Corollary 6.5.7]{cart})
and the map $R_{L}^{G}$ is just Harish-Chandra
induction (see \cite[p.~81]{DiMi2}).

This proves the first part of the theorem. By a result of Brou{\'e}
(\cite[p.~62, {\em Remarque}]{brouMo}), Harish-Chandra induction 
from~$L$ to~$G$ induces a Morita equivalence between the unions of 
$\ell$-blocks $\mathcal{E}_{\ell}( L, [s] )$ and 
${\mathcal E}_{\ell}(G,[s])$. In particular, irreducible Brauer 
characters are mapped to irreducible Brauer characters.
\end{prf}

\section{Asymptotics}
The above theorem implies that a large proportion of the irreducible
characters of a finite group of Lie type in cross-characteristic are 
Harish-Chandra induced. To make this statement more precise, we have
to talk about infinite families of groups, in our case of series of 
finite groups of Lie type. That is, we fix the type of the group and 
vary the underlying field. The formal setting for considerations of 
this type is that of a {\em generic finite reductive group} as introduced
by Brou{\'e} and Malle in~\cite{BrouMa}; see also \cite[Section 1.A]{BMM}.

A generic finite reductive group $\mathbb{G} = (\Gamma,F_0)$ consists of a 
root datum $\Gamma = (X,R,Y,R^\vee)$ and an automorphism~$F_0$ of~$\Gamma$ 
of finite order. The components~$X$ and~$Y$ of~$\Gamma$ are free abelian
groups of the same finite rank, the {\em rank of~$\mathbb{G}$}. They are
called ``character group'' and ``cocharacter group'' of $\mathbb{G}$,
respectively. According to the usage in~\cite{BrouMa}, we choose notation
such that~$F_0$ acts on~$X$ (contrary to the usage in~\cite[Section 1.A]{BMM},
where~$F_0$ acts on~$Y$). The Weyl group~$W$ of~$\Gamma$ also acts
linearly on~$X$. Extending scalars we obtain actions of~$W$ and of~$F_0$ 
on the $\mathbb{R}$-vector space $X_{\mathbb{R}} := X \otimes_{\mathbb{Z}} 
\mathbb{R}$.

Assume that $\mathbb{G}$ is not equal to one of the particular types 
described in \cite[(drc.t2), (drc.t3)]{BrouMa} giving rise to the Suzuki 
and Ree groups. Then any prime number~$p$ determines a connected reductive 
group~$\mathbf{G}$ over the algebraic closure of the prime field of 
characteristic~$p$, together with a maximal torus~$\mathbf{T}$, such that 
$\Gamma$ is the root datum associated to the pair $(\mathbf{G}, \mathbf{T})$. 
The Weyl group~$W$ of~$\Gamma$ is isomorphic to the Weyl group 
$N_{\mathbf{G}}(\mathbf{T})/\mathbf{T}$ of~$\mathbf{G}$, and the first 
component~$X$ of~$\Gamma$ is isomorphic to the character group of 
$\mathbf{T}$, the isomorphism being compatible with the actions of 
$N_{\mathbf{G}}(\mathbf{T})/\mathbf{T}$ and~$W$. Any power $q > 1$ of~$p$ 
determines a Frobenius morphism~$F$ of~$\mathbf{G}$ such that~$\mathbf{T}$ 
is $F$-stable. In particular,~$F$ acts on the character group of~$\mathbf{T}$, 
inducing the action $qF_0$ on~$X$. We shall write $\mathbb{G}(q)$ for the 
finite reductive group $\mathbf{G}^F$ determined in this way by~$q$. The set 
of groups $\{ \mathbb{G}(q) \mid q > 1 \text{\rm\ a prime power}\}$ is also 
called a {\em series of finite groups of Lie type}. For example, $\{ \SL_n(q) 
\mid q > 1 \text{\rm\ a prime power}\}$ (for fixed~$n$) or $\{ {^2\!E}_6( q )
\mid q > 1 \text{\rm\ a prime power}\}$ are series of finite groups of 
Lie type.

The Suzuki and Ree groups, excluded in the theorem below,
will be dealt with in a later chapter.

\begin{lem}
\label{NonRegularElements}
Let $\mathbb{G} = (\Gamma,F_0)$ be a generic finite reductive
group of rank~$m$. Assume that $\mathbb{G}$ is not equal to one of the
particular types described in \cite[(drc.t2), (drc.t3)]{BrouMa}.
We let~$W$ denote the Weyl group of $\mathbb{G}$ and write~$X$ for
the character group of~$\Gamma$, i.e., the first component of
the root datum~$\Gamma$.

Then there is a bound~$M$ (depending only on $\mathbb{G}$) such that for 
all prime powers~$q$ the following holds: If $(\mathbf{G}, F)$ is the pair 
of a connected reductive group~$\mathbf{G}$ and a Frobenius morphism~$F$
of~$\mathbf{G}$ determined by~$\mathbb{G}$ and~$q$, and if $\mathbf{T}$ is 
an $F$-stable maximal torus of~$\mathbf{G}$, then
$$\left| \{  s \in \mathbf{T}^F \mid C_{\mathbf{G}}( s )^F \neq 
\mathbf{T}^F \} \right| \leq M (q+1)^{m-1}.$$
\end{lem}
\begin{prf}
Let~$q$ be a prime power, and let $(\mathbf{G},F)$ and~$\mathbf{T}$
be as in the statement of the lemma.

An element $s \in T := \mathbf{T}^F$ is {\em regular}, if its connected 
centralizer $C_{\mathbf{G}}(s)^\circ$ equals~$\mathbf{T}$. By Part~(3) 
of the proof of \cite[Theorem 2.1]{Luebeck}, the number of non-regular 
elements in~$T$ is at most equal to $m^2 \cdot 2^{m} \cdot (q + 1)^{m-1}$. 

Suppose that $s \in T$ is regular, but that $C_{\mathbf{G}}( s )^F$ strictly 
contains~$T$. Then, by \cite[Theorem 3.5.3]{cart}, there is a non-trivial
$F$-stable $w \in N_{\mathbf{G}}(\mathbf{T})/\mathbf{T} = W$ centralizing~$s$. 
Thus~$s$ is contained in the closed $F$-stable subgroup $\mathbf{T}' = 
\{ t \in \mathbf{T} \mid \dot{w}^{-1} t \dot{w} = t \}$ of~$\mathbf{T}$, 
where $\dot{w}$ is an element of $N_{\mathbf{G}}(\mathbf{T})$ mapping to~$w$ 
under the natural epimorphism.

To estimate the order of ${\mathbf{T}'}^F$, we introduce the following 
notation. If $A \leq X$ is a subgroup of~$X$, we write $\hat{A}$ for the
smallest pure subgroup of~$X$ containing~$A$. Following 
\cite[Section 1.12]{cart}, we put $A^\perp := \{ t \in \mathbf{T} \mid
\chi(t) = 1 \text{\ for all\ } \chi \in A\}$. This is a closed subgroup
of~$\mathbf{T}$ and it is a sub-torus if $A = \hat{A}$. Moreover, for 
subgroups $A \leq B \leq X$ such that $B/A$ is finite, we have 
$A^\perp / B^\perp \cong (B/A)_{p'}$, where $p$ is the prime dividing~$q$ 
(see \cite[Chapter III]{Borel} for proofs of these facts). 

For $v \in W$, let $[v,X]$ denote the image of $v - I_X$ in~$X$.
Then $\mathbf{T}' = [w,X]^\perp$. Since~$w$ is $F$-stable, $[w,X]$ and 
$\widehat{[w,X]}$ are $F$-stable as well and $\widehat{[w,X]}^\perp$ is an 
$F$-stable sub-torus of $\mathbf{T}'$. Since $w \in W$ is non-trivial, the 
rank of~$[w,X]$, and hence that of~$\widehat{[w,X]}$ is at least equal to~$1$. 
By \cite[Propositon 2.4]{Luebeck}, the set of $F$-stable elements of the 
torus $\widehat{[w,X]}^\perp$ has at most $(q + 1)^{m-1}$ elements.

Let $M_0$ be the largest order of a group of the form $\widehat{[v,X]}/[v,X]$, 
where~$v$ runs through the elements of~$W$. Then the order of $\mathbf{T}'/ 
\widehat{[w,X]}^\perp$ is at most equal to~$M_0$. Hence $|{\mathbf{T}'}^F| 
\leq M_0 (q + 1)^{m-1}$ and the claim follows by summing up over all $w \in W$.
\end{prf}

\begin{thm}
\label{Asymptotic}
Let $\mathbb{G} = (\Gamma,F_0)$,~$W$ and~$X$ be as in 
Lemma~{\rm \ref{NonRegularElements}}. Denote by~$X'$ the root lattice
of~$\mathbb{G}$, i.e., the sublattice of~$X$ spanned by the roots~$R$
in~$X$. (Note that $F_0$ and $W$ act on~$X'$ and hence on $X'_{\mathbb{R}}
:= X' \otimes_{\mathbb{Z}} {\mathbb{R}}$.)

{\rm (a)} Put $W_{\text{\rm split}} := \{ w \in W \mid F_0 \circ w 
\text{\ has an eigenvalue~$1$ on\ } X'_{\mathbb{R}} \}$.
Then the number 
of absolutely irreducible ordinary characters of a group $\mathbb{G}(q)$
in the series $\{ \mathbb{G}(q) \mid q > 1 \text{\ a prime power} \}$ equals 
$q^m + O(q^{m-1})$, and
$$\frac{|W_{\text{\rm split}}|}{|W|}\,q^m + O(q^{m-1})$$
of these characters are Harish-Chandra imprimitive. 

{\rm (b)} Let~$\ell$ be a prime number and let~$e$ be a positive integer
dividing $\ell - 1$. If~$q$ is a prime power, we write $e(q,\ell) = e$, if
$\ell$ does not divide~$q$ and if~$e$ equals the order of $q \mod \ell$ in
$\mathbb{F}_\ell^*$. Let $W_{(e,\ell)'}$ denote the set of $w \in W$ for
which $F_0 \circ w$, acting on $X'_{\mathbb{R}}$, has no eigenvalue of
order~$\ell^ie$ for some non-negative integer~$i$. Put
$W_{\text{\rm split}, (e,\ell)'} := W_{\text{\rm split}} \cap W_{(e,\ell)'}$.
Then the number of absolutely irreducible Harish-Chandra imprimitive
$\ell$-modular Brauer characters of a group in the set
$\{ \mathbb{G}(q) \mid q > 1 \text{\ a prime power with\ } e(q,\ell) = e \}$
is at least 
$$\frac{|W_{\text{\rm split},(e,\ell)'}|}{|W|}\,q^m + O(q^{m-1}).$$
\end{thm}
\begin{prf}
Choose a prime power~$q$, hence a group~$\mathbf{G}$ and a Frobenius 
morphism~$F$ and write $G := \mathbb{G}(q) = \mathbf{G}^F$. 
Let $\mathbf{G}^*$, $F$ be the dual pair. A {\em geometric conjugacy
class of $G^* := {\mathbf{G}^*}^F$} is the intersection of a conjugacy
class of $\mathbf{G}^*$ with $G^*$, if this intersection is not empty.
We have a partition of the set of irreducible characters 
$\text{\rm Irr}(G)$ of~$G$ into Lusztig series,
$$\text{\rm Irr}(G) = \bigcup_{(s)} \mathcal{E}( G, (s) ),$$
where~$(s)$ runs through the geometric conjugacy classes of semisimple
elements of $G^*$ (see \cite[Proposition 13.17]{DiMi2}). For each
semisimple element $s \in G^*$, there is a bijection between 
$\mathcal{E}( G, (s) )$ and $\mathcal{E}( C_{G^*}(s), (1) )$, the 
{\em unipotent characters} of $C_{G^*}(s)$ (see \cite[Theorem 13.23]{DiMi2}).
Moreover, the latter are parameterized independently of~$q$. These 
results are due to Lusztig; for groups $\mathbf{G}$ with a connected center
they are contained as a special case in \cite[4.23 Main Theorem]{luszbuch},
and for groups $\mathbf{G}$ with non-connected center, they are proved
in~\cite{Lu}.

Our strategy is to count the geometric conjugacy classes of semisimple
elements $s \in G^*$ for which $C_{G^*}(s)$ is a maximal torus of~$G^*$. 
For such an~$s$, the set $\mathcal{E}( C_{G^*}(s), (1) )$ just consists 
of the trivial character, and hence $\mathcal{E}( G^*, (s) )$ contains a 
unique element. Asymptotically, the elements lying in such Lusztig series
account for nearly all irreducible characters of~$G$.

Up to conjugation in~$G^*$, the number of maximal tori in $G^*$ equals the 
number of $F$-conjugacy classes of~$W$ (see \cite[Proposition 3.3.3]{cart}). 
It follows from Lemma~\ref{NonRegularElements}, applied to the dual generic
reductive group $\mathbb{G}^*$, that the number of geometric conjugacy 
classes of semisimple elements~$s$ of $G^*$ such that $C_{G^*}(s)$ is not a 
maximal torus of~$G^*$ is at most a constant (independent of~$q$) times 
$(q + 1)^{m-1}$. Since the total number of geometric conjugacy classes of 
semisimple elements of~$G^*$ is of the form $q^m + O(q^{m-1})$ (see 
\cite[Theorem 3.7.6(i)]{cart}), we have proved the first assertion of 
Part~(a) of the theorem. 

To prove the second, let $\mathbf{T_0}$ be a maximally split torus of 
$\mathbf{G}$.  Let~$\mathbf{T_0}^*$ be the dual torus of $\mathbf{G}^*$. We 
choose a set of representatives for the $F$-conjugacy classes of~$W$, and 
for each representative~$w$ an $F$-stable maximal torus~$\mathbf{T}^*_w$ 
of~$\mathbf{G}^*$ obtained from $\mathbf{T_0}^*$ by twisting with~$w$ 
(see \cite[p.~85]{cart}). By \cite[Theorem 3.5.3]{cart},
an $F$-stable element~$s$ in $\mathbf{T}^*_w$ is regular and satisfies
$C_{G^*}(s) = {\mathbf{T}^*_w}^F$ if and only if its orbit under 
$N_{\mathbf{G}^*}(\mathbf{T}^*_w)^F$ has maximal possible length, namely
$|N_{\mathbf{G}^*}(\mathbf{T}^*_w)^F/{\mathbf{T}^*_w}^F|$. The latter number 
equals $|C_{W,F}(w)|$, the order of the $F$-centralizer $C_{W,F}(w)$ of~$w$
(see \cite[Proposition 3.3.6]{cart}). Thus~$\mathbf{T}^*_w$ gives rise to
$([W\colon\!C_{W,F}(w)] / |W|)\,q^m + O(q^{m-1})$ conjugacy classes
of regular elements of $G^*$ whose centralizer in~$G^*$ equals
${\mathbf{T}^*_w}^F$. Note that the number in the numerator
of the coefficient at $q^m$ in the above formula equals the
number of elements in the $F$-conjugacy class of~$W$ containing~$w$.

To conclude, we have to observe that $\mathbf{T}^*_w$ is contained in a 
proper split Levi subgroup of $\mathbf{G}^*$ if and only if $F_0 
\circ w$ has an eigenvalue~$1$ on~$X'_{\mathbb{R}}$. Indeed, $\mathbf{T}^*_w$ 
is contained in a proper split Levi subgroup of $\mathbf{G}^*$ if and only 
if it contains a non-trivial and non-central $\mathbb{F}_q$-split torus
of $\mathbf{G}^*$ (see \cite[Section 3.D]{BrouMa}). This is the case if
and only if $\mathbf{T}^*_w/Z^\circ(\mathbf{G}^*)$ contains a non-trivial 
$\mathbb{F}_q$-split torus of $\mathbf{G}^*/Z^\circ(\mathbf{G}^*)$, where
$Z^\circ(\mathbf{G}^*)$ denotes the connected component of the
center of~$\mathbf{G}^*$. Using the considerations in \cite[Section 1.A]{BMM},
we find that $\mathbf{T}^*_w/Z^\circ(\mathbf{G}^*)$ contains a non-trivial
$\mathbb{F}_q$-split torus, if and only if $F_0 \circ w$ has an eigenvalue~$1$ 
on $X'_{\mathbb{R}}$.

To prove Part~(b), assume that~$q$ satisfies $e(q,\ell) = e$. Let 
$s \in G^*$ be a semisimple $\ell'$-element such that $C_{G*}(s)$
is a maximal torus. By the main result of Brou{\'e} and Michel \cite{BM},
the (unique) ordinary character in $\mathcal{E}( G^*, (s) )$ remains
irreducible on reduction modulo~$\ell$. Moreover, if~$s'$ is another
such element not conjugate to~$s$ in~$G^*$, then the two irreducible
$\ell'$-modular characters arising from $\mathcal{E}( G^*, (s) )$
and $\mathcal{E}( G^*, (s') )$ are distinct. We thus have to count
the number of $G^*$-conjugacy classes of semisimple $\ell'$-elements
of $G^*$ whose centralizer is a maximal torus.

Let $\mathbf{T}^*_w$ be a maximal torus of $\mathbf{G}^*$ as above.
All non-central elements in ${\mathbf{T}^*_w}^F$ are $\ell'$-elements
if $\ell$ does not divide the order of 
$({\mathbf{T}^*_w}/Z^\circ(\mathbf{G}^*))^F$. This order
is given by evaluating the characteristic polynomial of $F_0 \circ w$,
acting on $X'_{\mathbb{R}}$, at~$q$ (see \cite[Proposition 3.3.5]{cart}).
The characteristic polynomial is a product of cyclotomic polynomials
$\Phi_j$. Now $\Phi_j(q)$ is divisible by~$\ell$ if and only if
$j = \ell^ie$ for some non-negative integer~$i$ (see. e.g., 
\cite[Lemma IX.8.1]{HuBII}). This completes the proof.
\end{prf}

Let us consider some examples. 

\begin{exl}
\label{ExamplesForAsymptotics}
{\rm
Let $q$ be a power of the prime~$p$, and let~$\bar{\mathbb{F}}_p$ denote an 
algebraic closure of the prime field~$\mathbb{F}_p$.

(1) Let $\mathbf{G} = \SL_n( \bar{\mathbb{F}}_p )$. 
The Weyl group~$W$ of~$\mathbf{G}$ (defined with respect to the maximal 
torus~$\mathbf{T}_0$ of~$\mathbf{G}$ consisting of the diagonal matrices) 
is isomorphic to the symmetric group on~$n$ letters and $X'_{\mathbb{R}} = 
X_{\mathbb{R}}$ is the reduced permutation module for~$W$. The rank of 
$\mathbf{G}$ equals $n - 1$. 

(a) Suppose first that~$F$ is the standard Frobenius morphism of~$\mathbf{G}$
(raising every entry of a matrix to its $q$th power), so that $G = \SL_n(q)$. 
Then~$F_0 = 1$. Now $w \in W$ has eigenvalue~$1$ on~$X_{\mathbb{R}}$, if 
and only if it has more than one cycle (in its action as a permutation on~$n$ 
points). Hence $|W_{\text{\rm split}}| = n! - (n-1)!$. It follows that for 
fixed~$n$ the proportion of (Harish-Chandra) imprimitive characters of 
$\SL_n(q)$ tends to $1 - 1/n$ as~$q$ tends to infinity. Moreover, this
asymptotic value tends to~$1$, as~$n$ tends to infinity.

(b) Now suppose that~$F$ is a Frobenius morphism giving rise to the unitary 
group, e.g., $F( (a_{ij}) ) = J_n { (a_{ij}^q)^{-T} } J_n$, for a matrix 
$(a_{ij}) \in \mathbf{G}$. (For the notation see 
Subsection~\ref{NotationFromLinearAlgebra}.) This~$F$ yields $G = \mathbf{G}^F 
= \SU_n(q)$. Let~$w_0$ denote the longest element in~$W$, and let~$w \in W$. 
One checks that $F_0 \circ w$ has eigenvalue~$1$ on $X_{\mathbb{R}}$ if and 
only if $-w_0 w$ has eigenvalue~$1$ on $X_{\mathbb{R}}$. A straightforward 
computation shows that $-w_0 w$ has eigenvalue~$1$ on $X_{\mathbb{R}}$, if 
and only if $w_0 w$ has at least one cycle of even length (in its action as a 
permutation on~$n$ points) i.e., $w_0w$ has even order. Thus 
$|W_{\text{\rm split}}| = 
| \{ w \in W \mid w \text{\rm\ has even order}\}|$.

From \cite[Lemma~4.2(a)]{LuebeckEtAl}, which is based upon 
\cite[Theorem 2.3(b)]{BealsEtAl}, we get
$$1 - \frac{3}{\sqrt{n}} < \frac{|W_{\text{\rm split}}|}{|W|} < 
1 - \frac{1}{4 \sqrt{n}}.$$
Again, the asymptotic value (for $q \rightarrow \infty$) of the proportion
of (Harish-Chandra) imprimitive characters of $\SU_n(q)$ tends to~$1$, 
as~$n$ tends to infinity.

(2) Let $\mathbf{G} = \Sp_{2m}( \bar{\mathbb{F}}_p )$ or $\mathbf{G} = 
\Spin_{2m + 1}( \bar{\mathbb{F}}_p )$, $m \geq 2$, and let~$F$ denote the 
standard Frobenius morphism of~$\mathbf{G}$. Then $F_0 = 1$ and 
$G = \Sp_{2m} ( q )$ or $G = \Spin_{2m+1}( q )$.

As~$\mathbf{G}$ is semisimple, we again have $X'_{\mathbb{R}} = X_{\mathbb{R}}$.
The Weyl group~$W$ of~$\mathbf{G}$ is a Coxeter group of type~$B_m$, and~$W$
permutes the $2m$ elements $\pm e_1, \ldots , \pm e_m$ for some basis 
$e_1, \ldots , e_m$ of~$X_{\mathbb{R}}$. A cycle of an element $w \in W$ on
the set $\{ \pm e_1, \ldots , \pm e_m \}$ is called negative, if it contains
$e_j$ and $-e_j$ for some $1 \leq j \leq m$. Otherwise, the cycle is called 
positive. It is easy to see that $w \in W$ has eigenvalue~$1$ 
on~$X_{\mathbb{R}}$, if and only if at least one cycle of~$w$ is positive.
Put $f(m) := |W \setminus W_{\text{\rm split}}|/|W|$, the proportion of
elements of~$W$ with only negative cycles.
It has been shown by Frank L\"ubeck that
\begin{equation}
\label{WBnAsymptotics}
f(m) = \frac{1\cdot 3 \cdot 5 \cdots (2m - 1)}{|W|},
\end{equation}
and that this proportion tends to~$0$, as~$m$ tends to infinity
(see \cite[Proposition 2(a), Remark~$3$]{LuebeckNegativeCycles}).

(3) Let $\mathbf{G} = \Spin_{2m}( \bar{\mathbb{F}}_p )$, $m \geq 4$, and 
let~$F^\epsilon$ be a Frobenius morphism such that $G^\epsilon := 
\mathbf{G}^{F^\epsilon} = \Spin_{2m}^\epsilon( q )$ for 
$\epsilon \in \{ +, - \}$.

Again, $X'_{\mathbb{R}} = X_{\mathbb{R}}$, and the latter may be identified
with the group $X_{\mathbb{R}}$ arising in Part~(2) above. Moreover, the Weyl 
group~$\tilde{W}$ of~$\mathbf{G}$ is a Coxeter group of type~$D_m$. We may 
embed~$\tilde{W}$ into the Weyl group~$W$ of type~$B_m$ as a subgroup of 
index~$2$ such that the action of~$\tilde{W}$ on $X_{\mathbb{R}}$ is the 
restriction of the action of~$W$ (with the above identification of the vector 
spaces $X_{\mathbb{R}}$). In fact,~$\tilde{W}$ consists of the elements of~$W$ 
which have an even number of negative cycles.

If $\epsilon = +$, we may assume that $F_0 = 1$. In this case, 
$$
\frac{|\tilde{W} \setminus {\tilde{W}}_{\text{\rm split}}|}{|\tilde{W}|} = 
f(m)\frac{2m - 2}{2m-1}
$$
by \cite[Proposition 2(b)]{LuebeckNegativeCycles}, where~$f(m)$ is as in 
Equation~(\ref{WBnAsymptotics}). 

If $\epsilon = -$, we may assume that~$F_0$ acts on~$X_{\mathbb{R}}$ as an
element of~$W \setminus \tilde{W}$. In this case,
$$
\frac{|\tilde{W} \setminus {\tilde{W}}_{\text{\rm split}}|}{|\tilde{W}|} = 
f(m)\frac{2m}{2m-1}
$$
by \cite[Proposition 2(c)]{LuebeckNegativeCycles},  where~$f(m)$ is as in 
Equation~(\ref{WBnAsymptotics}). Clearly, these proportions tend to~$0$ as~$m$
tends to infinity.
}
\end{exl}

\chapter{Groups of Lie type: $\text{\rm char}(K) = 0$}
\label{Characteristic0}

\markleft{GROUPS OF LIE TYPE: {$\text{\rm char}(K) = 0$}}

In this chapter we assume that~$G$ is a finite group of Lie type (not 
necessarily quasisimple). We also assume that~$K$ is an algebraically closed 
field of characteristic~$0$. If, in addition,~$G$ arises from an algebraic 
group with connected center, we can prove a converse to Theorem~\ref{HCresult},
and thus obtain a classification of the irreducible imprimitive $KG$-modules
with parabolic block stabilizers. Under less restrictive assumptions, we
obtain a corresponding result about Lusztig series 
(Theorem~\ref{FirstLusztigSeriesResultConverse}). As a consequence of the 
latter, all unipotent characters of~$G$ are Harish-Candra primitive.

\section{Some results on Weyl groups}

In the following,~$W_m$ denotes a Weyl group of type~$B_m$. The 
maximal parabolic subgroups of~$W_m$ are of the form $W_{m-k}
\times S_k$ for $1 \leq k \leq m$. Here,~$S_k$ denotes a symmetric
group on~$k$ letters. The ordinary irreducible characters of~$W_m$ 
are parameterized by bi-partitions of~$m$, and we write $\chi^\alpha$ 
for the irreducible character of~$W_m$ labelled by the bi-partition~$\alpha$. 
Similarly, $\zeta^\beta$ denotes the irreducible character of~$S_k$ 
labelled by the partition~$\beta$ of~$k$. For partitions $\alpha$,
$\beta$, and $\gamma$ of $m-k$,~$k$, and $m$, respectively, we write
$g_{\alpha,\beta}^\gamma$ for the multiplicity of $\zeta^\gamma$
in the induced character $\Ind_{S_{m-k} \times S_k}^{S_m}( \zeta^\alpha
\times \zeta^\beta )$.

\begin{lem}
\label{WeylBm}
Let $m \geq 2$ and $1 \leq k \leq m$ be positive integers.
Moreover, let~$\gamma$ be a partition of~$k$, and let
$\beta = (\beta^0, \beta^1)$ be a bi-partition of $m - k$.
Then 
\begin{equation}
\label{IndBm}
\Ind_{W_{m-k} \times S_k}^{W_m}( \chi^\beta \times \zeta^\gamma )
\end{equation}
has at least three distinct irreducible constituents. If $k > 1$ and 
$(m,k) \neq (2,2), (3,3)$, then~{\rm (\ref{IndBm})} has at least five
distinct irreducible constituents.
\end{lem}
\begin{prf}
Suppose first that $k = 1$. The irreducible constituents of 
$\Ind_{W_{m-1}}^{W_m}
( \chi^\beta )$ are of the form $\chi^\alpha$ with bi-partitions
$\alpha = ( \alpha^0, \alpha^1 )$ of~$m$ satisfying $\alpha^0 =
\beta^0$ or $\alpha^1 = \beta^1$. In the first case $\alpha^1$
is obtained from adding a node to $\beta^1$, and in the second case
$\alpha^0$ is obtained from $\beta^0$ in this way. (This is a special
case of the general formula~(\ref{branchingB}) below.) Since $m \geq 2$,
at least one of $\beta^0$ or $\beta^1$ is not the empty partition.
If~$\beta^0$ is not the empty partition, there are at least two
distinct partitions which can be obtained from $\beta^0$ by adding
a node. The result follows in this case.

Now suppose that $k > 1$. If $(m,k) = (2,2)$ or $(3,3)$ we use the
CHEVIE~\cite{chevie} share package of GAP~\cite{GAP4} to show 
that~(\ref{IndBm}) has at least three distinct irreducible constituents
(in fact in these cases it has at least four). Thus suppose in addition 
that $(m,k) \neq (2,2), (3,3)$. For a bi-partition $(\alpha^0,\alpha^1)$
of~$m$, put $j := |\alpha^0| - |\beta^0|$. Then the multiplicity of
$\chi^\beta \times \zeta^\gamma$ in the restriction of $\chi^\alpha$
to $W_{m-k} \times S_k$ is equal to~$0$ unless $0 \leq j \leq k$.
In the latter case this multiplicity equals
\begin{equation}
\label{branchingB}
\sum_{\delta^0} \sum_{\delta^1}
                          g^{\alpha^0}_{\beta^0 \delta^0}
                          g^{\alpha^1}_{\beta^1 \delta^1}
                          g^{\gamma}_{\delta^0 \delta^1},
\end{equation}
where~$\delta^0$ and~$\delta^1$ run through the partitions of~$j$ and 
$k - j$, respectively. Formula~(\ref{branchingB}) follows from 
\cite[Theorem~6]{Osima} and an application of the Littlewood-Richardson
rule (see \cite[2.8.14]{jake}) for symmetric groups. For a detailed
deduction of the formula see \cite[Lemma 3.5.2]{GoetzDiss}).

Let $\delta = (\delta^0,\delta^1)$ be one of the following bi-partitions 
of~$k$:
\begin{equation}
\label{Deltas}
( \text{-}, \gamma ), ( \gamma, \text{-} ), ( (1), \mu ), ( \mu, (1) ),
( \nu, \kappa ),
\end{equation}
where $\mu$, $\nu$, and $\kappa$ are partitions of $k - 1$, $2$, and
$k - 2$, respectively, such that $g^\gamma_{ (1), \mu } \neq 0$ and 
$g^\gamma_{ \nu, \kappa } \neq 0$.

For each $(\delta^0,\delta^1)$ as above, there is a bi-partition $\alpha = 
(\alpha^0,\alpha^1)$ of~$m$ such that $g^{\alpha^0}_{\beta^0,\delta^0}
\neq 0$ and $g^{\alpha^1}_{\beta^1,\delta^1} \neq 0$. By~(\ref{branchingB})
this implies that $\chi^\alpha$ occurs in~(\ref{IndBm}) as a constituent.
If $k \geq 4$, then the five $\delta$'s of~(\ref{Deltas}) give rise to 
five distinct $\alpha$'s. Indeed, letting $\beta^0$ be a partition of~$m_0$,
the first component $\alpha^0$ of $\alpha$ is a partition of $m_0$, $m_0 + k$,
$m_0 + 1$, $m_0 + k - 1$, and $m_0 + 2$, respectively. If $k = 2$ or~$3$,
then $m - k > 0$ by assumption. Thus at least one of $\beta^0$ or $\beta^1$
is not the empty partition. Assume without loss of generality that $\beta^0$
is not empty. Then, by the Littlewood-Richardson rule, there are at least
two distinct partitions $\alpha^0$ of $m_0 + k$ with 
$g^{\alpha^0}_{\beta^0, \gamma } \neq 0$ and two distinct partitions
$\alpha^0$ of $m_0 + 1$ with $g^{\alpha^0}_{\beta^0, (1) } \neq 0$.
For $\alpha^0 = \beta^0$ we also have $g^{\alpha^0}_{\beta^0, (\text{-}) } 
\neq 0$, and we are done.
\end{prf}

\begin{lem}
\label{WeylInd}
Let~$W$ be a finite Weyl group and let~$W_0$ be a proper parabolic
subgroup of~$W$ (for this notation see e.g. \cite[Section 2.3]{cart}).
Then $\Ind_{W_0}^W(\psi)$ is reducible for every ordinary irreducible 
character~$\psi$ of~$W_0$. 
%In fact, if~$W$ is of type~$D_m$, $m \geq 4$,
%or of type~$E_7$, then $\Ind_{W_0}^W(\psi)$ has at least three distinct
%irreducible constituents. If~$W$ is of type~$E_6$, then $\Ind_{W_0}^W(\psi)$
%has at least four distinct irreducible constituents.
\end{lem}
\begin{prf}
We may assume that~$W$ is irreducible and that~$W_0$ is a maximal
parabolic subgroup of~$W$. Suppose first that~$W$ is not of type $A_m$, 
$D_m$, or $E_6$. Then~$w_0$, the longest element in~$W$ 
(\cite[Proposition 2.2.11]{cart}) is contained in the center 
of~$W$ (see \cite[Planches I--IX]{Bourbaki}). Since~$w_0$ is not contained 
in~$W_0$, the result follows from Lemma~\ref{Mackey1}. This argument is 
due to Meinolf Geck.  (Of course, the desired result for $W$ of type~$B_m$
also follows from Lemma~\ref{WeylBm}.)

It remains to consider the cases that $W$ is of type~$A_m$, $m \geq 2$,
$D_m$, $m \geq 4$, and $E_6$. The latter case is easily checked with 
CHEVIE\@. If $W$ is of type~$A_m$, the Littlewood-Richardson rule (see
\cite[2.8.14]{jake}) shows that no irreducible character of~$W_0$ induces 
to an irreducible character of~$W$.

In case~$W$ is of type~$D_m$, $m \geq 4$, we proceed as follows. 
We write $\tilde{W}_m$ for the Weyl group of type~$D_m$, and embed 
$\tilde{W}_m$ into a Weyl group $W_m$ of type~$B_m$ of index~$2$ 
(exactly as in the first of the two embeddings described in 
\cite[4.6]{luszbuch}).
A maximal parabolic subgroup of $\tilde{W}_m$ is of the form
$\tilde{W}_{m-k} \times S_k$ for some $1 \leq k \leq m$.
This is embedded into $W_{m-k} \times S_k$ by embedding 
$\tilde{W}_{m-k}$ into $W_{m-k}$.

Now let $\psi$ be an irreducible character of $\tilde{W}_{m-k} \times S_k$
inducing to an irreducible character~$\chi$ of $\tilde{W}_m$. Then
$\Ind_{\tilde{W}_m}^{W_m}(\chi)$ has at most two irreducible constituents.
This implies that the irreducible constituents of
$\Ind_{\tilde{W}_{m-k} \times S_k}^{W_{m-k} \times S_k}(\psi)$
induce to characters of $W_m$ with at most two irreducible
constituents. This contradicts Lemma~\ref{WeylBm}.
\end{prf}

\section{Harish-Chandra series}\label{Converse}

Our next aim is to prove a converse of Theorem~\ref{HCresult} for groups of Lie 
type arising from an algebraic group with connected center. More specifically, 
let~${\mathbf G}$ be a connected reductive algebraic group over the 
algebraic closure of~${\mathbb F}_q$ and let~$F$ be a Frobenius morphism 
of~${\mathbf G}$. Let $G = {\mathbf G}^F$ be the corresponding finite group of 
Lie type. Assume that the center of~${\mathbf G}$ is connected. We exclude the
pairs $({\mathbf G}, F)$ that give rise to the Suzuki and Ree groups, i.e., we
do not allow~$G$ to be one of ${^2\!B}_2(q)$, ${^2\!G}_2(q)$, ${^2\!F}_4(q)$.
The latter groups will be treated in Section~\ref{ExplicitExceptional}.

\begin{thm}
\label{HCresult:converse}
For~$G$ as above, let~$M$ be an irreducible imprimitive $KG$-module whose 
block stabilizer is a parabolic subgroup of~$G$ with Levi complement~$L$.
Then there exists a Levi subgroup~$L_0$ of~$L$, 
an irreducible cuspidal $KL_0$-module $M_0$ with $W_G(L_0,M_0) 
= W_L( L_0, M_0 )$, and an irreducible constituent~$M_1$ of 
$R_{L_0}^L( M_0 )$ such that $M \cong R_L^G(M_1)$.
\end{thm}
\begin{prf}
By Proposition~\ref{ParabolicStabilizer}, there is an irreducible
$KL$-module~$M_1$ such that $M \cong R_L^G(M_1)$.
Let $L_0$ be a Levi subgroup of~$L$ and~$M_0$ an
irreducible cuspidal $KL_0$-module such that $M_1$ is a
composition factor of $R_{L_0}^L(M_0)$. 

By the results of Lusztig in \cite[8.3--8.5]{luszbuch}), the ramification 
group $W_G( L_0, M_0 )$ is a Weyl group and $W_L( L_0, M_0 )$
is a parabolic subgroup of $W_G( L_0, M_0 )$. 
By the comparison theorem of Howlett 
and Lehrer \cite[Theorem 5.9]{HowLeh2}, there are irreducible 
characters~$\psi$ 
of $W_L( L_0, M_0 )$ and $\chi$ of $W_G( L_0, M_0 )$ corresponding 
to~$M_1$ and $M$ respectively, such that the multiplicity of $M$ in 
$R_L^G(M_1)$ equals the multiplicity of $\chi$ in the induced character 
$\Ind_{W_L( L_0, M_0 )}^{W_G( L_0, M_0 )} (\psi)$.
The first of these multiplicities being~$1$, this is true for the
second as well. By Lemma~\ref{WeylInd} this implies that
$W_L( L_0, M_0 ) = W_G( L_0, M_0 )$ and we are done.
\end{prf}

\section{Lusztig series}\label{Section73}

Here we give a characterization of the irreducible Harish-Chandra 
imprimitive $KG$-modules in terms of Lusztig series of characters, 
thus proving a converse to Theorem~\ref{FirstLusztigSeriesResult}
in characteristic~$0$. Let ${\mathbf G}$,~$F$ and $G = {\mathbf G}^F$
be as in Section~\ref{Converse} above. We do no longer, however, assume 
that the center of ${\mathbf G}$ is connected. As in 
Section~\ref{LusztigSeries}, we let $(\mathbf{G}^*,F)$ be the pair
dual to $(\mathbf{G},F)$, and we write $G^* = {\mathbf{G}^*}^F$.

\begin{thm}
\label{FirstLusztigSeriesResultConverse}
Let $s \in G^*$ be semisimple such that $C_{\mathbf{G}^*}( s )$ is connected.
If $\mathcal{E}(G,[s])$ contains a Harish-Chandra imprimitive element, then
$C_{\mathbf{G}^*}(s)$ is contained in a proper split $F$-stable Levi subgroup
of $\mathbf{G}^*$, and thus every element of $\mathcal{E}(G,[s])$ is
Harish-Chandra imprimitive.
\end{thm}
\begin{prf}
First assume that the center of ${\mathbf G}$ is connected.
Let $\chi \in \mathcal{E}(G,[s])$ be Harish-Chandra induced from the
Levi subgroup $L = \mathbf{L}^F$ of~$G$. Suppose that~$\chi$ lies
in the $(L_0,\psi)$-Harish-Chandra series for some Levi subgroup 
$\mathbf{L}_0^F = L_0 \leq L$ and some cuspidal irreducible character
$\psi$ of~$L_0$. We let $\mathbf{L}^*_0$ and $\mathbf{L}^*$ denote split
$F$-stable Levi subgroups of $\mathbf{G}^*$ in duality with $\mathbf{L}_0$
and $\mathbf{L}$ respectively. As Harish-Chandra induction preserves
Lusztig series (see \cite[Proposition~$15.7$]{CaEn}) we may assume that 
$s \in \mathbf{L}^*_0 \leq \mathbf{L}^*$.

We now use \cite[(8.5.7) and (8.5.8)]{luszbuch}.
This implies that $W_G(L_0,\psi)$ is a Coxeter group with canonical
generators in bijection to the $F$-orbits on $S_{C_{\mathbf{G}^*}(s)}
\setminus S_{C_{\mathbf{L}^*_0}(s)}$. Here, $S_{C_{\mathbf{G}^*}(s)}$
denotes the set of fundamental reflections of the Weyl group of 
$C_{\mathbf{G}^*}(s)$. Now $W_G(L_0,\psi) = W_L(L_0,\psi)$ by 
Theorem~\ref{HCresult:converse}, and hence the number of $F$-orbits on
$S_{C_{\mathbf{G}^*}(s)} \setminus S_{C_{\mathbf{L}^*}(s)}$ is zero.
Since $C_{\mathbf{L}^*}(s)$ is a split Levi subgroup of
$C_{\mathbf{G}^*}(s)$, it follows that $C_{\mathbf{G}^*}(s) =
C_{\mathbf{L}^*}(s)$, i.e., that $C_{\mathbf{G}^*}(s) \leq \mathbf{L}^*$.

%The last assertion follows from Theorem~\ref{FirstLusztigSeriesResult}.
%\end{prf}
%
%The above result can be generalized to the case where the center of 
%${\mathbf G}$ is not necessarily connected. 
%
%\begin{thm}
%\label{SeconLusztigSeriesResultConverse}
%Let $s \in G^*$ be semisimple such that $C_{\mathbf{G}}^*( s )$ is connected.
%If $\mathcal{E}(G,[s])$ contains a Harish-Chandra imprimitive element, then 
%$C_{\mathbf{G}^*}(s)$ is contained in a proper split $F$-stable Levi subgroup 
%of $\mathbf{G}^*$, and thus every element of $\mathcal{E}(G,[s])$ is 
%Harish-Chandra imprimitive.
%\end{thm}
%\begin{prf}

Now we deal with the general case.
Let~$\tilde{\mathbf{G}}$ denote a connected reductive algebraic group over the 
algebraic closure of $\mathbb{F}_q$, equipped with a Frobenius endomorphism~$F$,
such that the following conditions are satisfied. The group~$\mathbf{G}$ is a
closed, $F$-stable subgroup of~$\tilde{\mathbf{G}}$, and the Frobenius map 
on~$\mathbf{G}$ is the restriction of~$F$ to~$\mathbf{G}$. Moreover, the
derived subgroups of~$\mathbf{G}$ and of~$\tilde{\mathbf{G}}$ coincide
and~$\tilde{\mathbf{G}}$ has connected center. Such a group always
exists (see \cite[Section~$15.1$]{CaEn}).

By duality, there is a surjective, $F$-equivariant homomorphism of algebraic groups, 
\begin{equation}
\label{DualRegularEmbedding}
\tilde{\mathbf{G}}^* \rightarrow {\mathbf{G}}^*, 
\end{equation}
whose kernel is connected and contained in the center of~$\tilde{\mathbf{G}}^*$.

Let $\tilde{s} \in \tilde{\mathbf{G}}^*$ denote an $F$-stable element mapping
to~$s$ under~(\ref{DualRegularEmbedding}). (The existence of such an element
follows from \cite[Corollaire~$2.7$]{CeBo2}.)
As $C_{\mathbf{G}}^*( s )$ is connected, we obtain a bijection
\begin{equation}
\label{RestrictionBijection}
\mathcal{E}(\tilde{G},[\tilde{s}]) \rightarrow \mathcal{E}(G,[s]),\quad 
\tilde{\chi} \mapsto \Res^{\tilde{G}}_G( \tilde{\chi} ).
\end{equation}
This follows from the main result of Lusztigs paper~\cite{Lu}.
Now suppose that $\chi \in \mathcal{E}(G,[s])$ is Harish-Chandra imprimitive.
Thus there is an $F$-stable split proper Levi subgroup~$\mathbf{L}$ 
of~$\mathbf{G}$, and an irreducible character $\psi$ of~$L$ such that
$\chi = R_L^G( \psi )$. Let~$\mathbf{L}^*$ denote an $F$-stable split Levi
subgroup of~$\mathbf{G}^*$ dual to~$\mathbf{L}$. As above we may assume that 
$s \in L^*$, and thus $\tilde{s} \in \tilde{L}^*$, where $\tilde{\mathbf{L}}^*$ is the 
inverse image of~$\mathbf{L}^*$ under the 
homomorphism~(\ref{DualRegularEmbedding}). Let~$\tilde{\mathbf{L}}$ denote
an $F$-stable split Levi subgroup of~$\tilde{\mathbf{G}}$ dual to 
$\tilde{\mathbf{L}}^*$ and
with $\mathbf{L} \leq \tilde{\mathbf{L}}$. As the center of $\tilde{\mathbf{G}}$ 
is connected, so is $C_{\tilde{\mathbf{G}}^*}( \tilde{s} )$ (see
\cite[Theorem 4.5.9]{cart}) and thus 
$C_{\tilde{\mathbf{L}}^*}( \tilde{s} ) = \tilde{\mathbf{L}}^* \cap 
C_{\tilde{\mathbf{G}}^*}( \tilde{s} )$ is connected as well. 
By~(\ref{RestrictionBijection}), applied to 
$\mathbf{L} \leq \tilde{\mathbf{L}}$, there is an irreducible character
$\tilde{\psi} \in \mathcal{E}(\tilde{L},[\tilde{s}])$ with 
$\Res^{\tilde{L}}_L( \tilde{\psi} ) = \psi$.

We have
$$\Res^{\tilde{G}}_G( R_{\tilde{L}}^{\tilde{G}}( \tilde{\psi} ) ) =
R_L^G( \Res^{\tilde{L}}_L( \tilde{\psi} ) ) = \chi,$$
by \cite[Equation~$(15.5)$]{CaEn}. As every irreducible constituent of
$R_{\tilde{L}}^{\tilde{G}}( \tilde{\psi} )$ is contained in 
$\mathcal{E}(\tilde{G},[\tilde{s}])$, we conclude that 
$R_{\tilde{L}}^{\tilde{G}}( \tilde{\psi} )$ is irreducible from the 
bijection~(\ref{RestrictionBijection}). As $\tilde{L} \neq \tilde{G}$ 
is a split Levi subgroup of~$\tilde{G}$, it follows from the first part 
of the proof that
$C_{\tilde{\mathbf{G}}^*}( \tilde{s} )$ contained in a proper split
Levi subgroup of $\tilde{\mathbf{G}}^*$. Now $C_{\mathbf{G}^*}( s )$
is the image of $C_{\tilde{\mathbf{G}}^*}( \tilde{s} )$ under the 
homomorphism~(\ref{DualRegularEmbedding}) (see \cite[Proposition~$2.3$]{DiMi2}). 
The result follows, as the
split Levi subgroups of~$\tilde{\mathbf{G}}^*$ and of $\mathbf{G}^*$ are
in bijection via the map~(\ref{DualRegularEmbedding}).

The last assertion follows from Theorem~\ref{FirstLusztigSeriesResult}.
\end{prf}

Notice that $C_{\mathbf{G}^*}( s )$ is connected for every semisimple element
$s \in G^*$, if the center of~$\mathbf{G}$ is connected (see 
\cite[Theorem 4.5.9]{cart}).
\begin{cor}
\label{UnipotentsAreImprimitive}
The unipotent characters of~$G$ are Harish-Chandra primitive.
\end{cor}
\begin{prf}
The unipotent characters are the elements of $\mathcal{E}(G,[1])$,
and $C_{\mathbf{G}}^*( s ) = \mathbf{G}^*$ is connected.
\end{prf}

\chapter{Classical groups: $\text{char}(K) = 0$}
\label{ClassicalInCharacteristic0}

\markleft{CLASSICAL GROUPS: {$\text{\rm char}(K) = 0$}}
Here we apply the results of the preceding chapters to certain 
classical groups for which we obtain explicit descriptions of 
their Harish-Chandra induced irreducible characters. Since we
are only dealing with representations over algebraically closed
fields of characteristic~$0$, it seems appropriate to use the
language of characters rather than modules. A character thus
means a $K$-character, where~$K$ is an algebraically closed
field of characteristic~$0$.

We begin with a list of groups satisfying the hypotheses of 
Theorem~\ref{HCresult:converse}.

\section{The groups}\label{TheClassicalGroups} For the remainder 
of this section we let~$G$ be one of the following finite classical 
groups:
\begin{itemize}
\item[(a)] a general linear group $\GL_n(q)$ with $n \geq 2$,
\item[(b)] a general unitary group $\GU_n(q)$ with $n \geq 3$,
\item[($\mbox{\rm c}'$)] a symplectic group $\Sp_{n}(q)$ with~$q$ even 
and $n \geq 4$ even,
\item[(c)] a conformal symplectic group $\CSp_{n}(q)$ 
with~$q$ odd and $n \geq 4$ even,
\item[(d)] a Clifford group $D_{n+1}^0(q)$ with~$q$ odd and $n \geq 4$ even,
\item[($\mbox{\rm e}'$)] an orthogonal group $\Omega^\pm_{n}(q)$ with 
$q$ even and $n \geq 8$ even,
\item[(e)] a conformal special orthogonal group 
$\CSO^\pm_{n}(q)$ with~$q$ odd and $n \geq 8$ even,
\item[(f)] a Clifford group $D_{n}^{\pm}(q)$ with~$q$ odd and $n \geq 8$
even.
\end{itemize}
For the definition of these groups we refer to \cite[Sections 1.19]{cart} 
and \cite[Chapter IX, \S\ 9, $\mbox{\rm n}^\circ$5]{Bourbaki2}. 
The groups in (c) and (d) are dual groups in the sense of Deligne and 
Lusztig, as well as the groups of types (e) and (f) (see \cite[8.1]{lusz2}).
As before, we denote the natural $n$-dimensional (projective) module 
for~$G$ by~$V$.

\section{Harish-Chandra series}
The aim in this section is a characterization of the ordinary 
irreducible Harish-Chandra imprimitive characters of the above 
classical groups in terms of Harish-Chandra theory. 
%These results will 
%be used in describing the descent to the commutator subgroups in the
%second part of our paper.

Let~$W$ denote the Weyl group of~$G$ (as finite group with $BN$-pair).
Thus, for example, the Weyl group of $\GU_{2m}(q)$ is the Weyl group of
type~$B_m$. The set of fundamental reflections of~$W$ is denoted by~$S$.
The Weyl groups and root systems of the groups in~\ref{TheClassicalGroups}
are irreducible of types $A$,~$B$,~$C$ and~$D$.

For every subset $J \subseteq S$ there is a standard Levi subgroup~$L_J$
of~$G$, whose Weyl group is the parabolic subgroup $W_J \leq W$, generated
by the reflections in~$J$. As before, a Levi subgroup of~$G$ is an 
$N$-conjugate of some standard Levi subgroup, where~$N$ denotes the
monomial subgroup of the $BN$-pair of~$G$. The relative Weyl group $W_G(L)$
of the Levi subgroup $L$ is defined as $W_G(L) := (N_G(L) \cap N)L/L$.

\subsection{General linear groups}
We first investigate the easy case of the general linear group
$G = \GL_n(q)$. In this case~$W$ may be identified with the subgroup
of permutation matrices in~$G$, and~$S$ with the subset~$\{ s_1,
\ldots , s_{n-1} \}$, where~$s_i$ corresponds to the transposition 
$(i,i+1)$, $1 \leq i \leq n-1$.

By a composition of~$n$ we understand a sequence $ (\lambda_1, \ldots , 
\lambda_r)$ of positive integers summing up to~$n$. Let~$J$ be a proper 
subset of~$S$. By putting $\lambda_1 := \mbox{\rm min}\{k \mid 1 \leq 
k \leq n - 1, s_k \not\in J\}$, and $\lambda_{j+1} := 
\mbox{\rm min}\{k \mid 1 + \sum_{i=1}^j \lambda_i \leq k \leq n - 1, 
s_k \not\in J\} - \sum_{i=1}^j \lambda_i$, for $j \geq 1$ and 
$\sum_{i=1}^j \lambda_i < n$, we inductively associate a composition
$(\lambda_1, \lambda_2, \ldots , \lambda_r)$ of~$n$ to~$J$. The standard 
Levi subgroup $L_J$ of~$G$ is of the form
\begin{equation}
\label{L_JinGL}
L_J \cong \GL_{\lambda_1}(q) \times \GL_{\lambda_2}(q) \times \cdots 
\times \GL_{\lambda_r}(q),
\end{equation}
with the natural diagonal embedding of the right hand side of 
(\ref{L_JinGL}) into~$G$. 

An irreducible cuspidal character~$\psi$ of~$L_J$ as above is of the form
\begin{equation}
\label{psiinGL}
\psi = \psi_1 \otimes \psi_2 \otimes \cdots \otimes \psi_r,
\end{equation}
where $\psi_i$ is an irreducible cuspidal character of
$\GL_{\lambda_i}(q)$ for $1 \leq i \leq r$.

Let us shortly describe $W_G(L_J)$ and $W_G(L_J,\psi)$ for~$L_J$ 
and~$\psi$ as above. For every pair $(i,j)$, $1 \leq i < j \leq r$ 
such that $\lambda_i = \lambda_j$, there is an element $\tau_{ij} 
\in W_G(L_J)$ swapping $\GL_{\lambda_i}(q)$ with $\GL_{\lambda_j}(q)$ 
and centralizing the other factors of~$L_J$. In fact, with respect to a 
suitable basis of~$V$, the natural (projective) module for~$G$, the 
element~$\tau_{ij}$ is represented by the matrix
\begin{equation}
\label{ti}
\left( \begin{array}{ccccc}
{I}_{a} & 0 & 0 & 0 & 0 \\
0 & 0 & 0 & I_{\lambda_i} & 0 \\
0 & 0 & I_b & 0 & 0 \\
0 & I_{\lambda_i} & 0 & 0 & 0 \\
0 & 0 & 0 & 0 & {I}_{c}
\end{array}\right),
\end{equation}
with $a = \sum_{k=1}^{i-1} \lambda_k$, $b = \sum_{k=i+1}^{j-1} \lambda_k$,
and $c = \sum_{k = j+1}^r \lambda_k$ (recall that $I_d$ denotes the $d \times d$ 
identity matrix).
By \cite{How} we have $W_G(L_J) = \langle \tau_{ij} \mid 1 \leq i < j \leq r, 
\lambda_i = \lambda_j \rangle$. Moreover, $W_G(L_J,\psi) = \langle \tau_{ij} 
\mid 1 \leq i < j \leq r, \lambda_i = \lambda_j, \psi_i = \psi_j \rangle$. 
This follows from one of the main results of Chapter~$8$ of Lusztig's book
(namely \cite[(8.5.13)]{luszbuch}).

\begin{prop}
\label{PrimitivesInClassicalI}
Let~$G = \GL_n(q)$, $n \geq 2$, let $L_0 = L_J$ be a proper Levi subgroup 
of~$G$ as in~{\rm (\ref{L_JinGL})}, and let~$\psi$ be a cuspidal irreducible 
character of~$L_0$ as in~{\rm (\ref{psiinGL})}.

If $\lambda_1 = \lambda_2 = \cdots = \lambda_r$ and $\psi_1 = \psi_2 =
\cdots = \psi_r$, then every irreducible constituent of $R_{L_0}^G(\psi)$ 
is Harish-Chandra primitive. Otherwise every such constituent
is Harish-Chandra imprimitive.
\end{prop}
\begin{prf}
Suppose first that $\lambda_1 = \lambda_2 = \cdots = \lambda_r$ and 
$\psi_1 = \psi_2 = \cdots = \psi_r$. Then $W_G(L_0,\psi) = W_G(L_0)$
is isomorphic to the symmetric group on~$r$ letters. Let~$L$ be a proper
Levi subgroup of~$G$ such that~$L_0$ is a Levi subgroup of~$L$. Then some 
$W_G(L_0)$-conjugate~$L_1$ of~$L$ is a standard Levi subgroup of~$G$.
Obviously,~$W_{L_1}(L_0)$ is a proper subgroup of~$W_G(L_0)$. In 
particular, $W_G(L_0,\psi) \neq W_{L_1}(L_0,\psi)$ and so $W_G(L_0,\psi)
\neq W_L(L_0,\psi)$. It follows from Theorem~\ref{HCresult:converse}
and Harish-Chandra theory (see \cite[Theorem (70.15A)]{cr1}), that no 
irreducible constituent of $R_{L_0}^G(\psi)$ is Harish-Chandra imprimitive.

Suppose now that $\lambda_i \neq \lambda_j$ for some $1 \leq i,j \leq r$,
or that $\lambda_1 = \lambda_2 = \cdots = \lambda_r$, but $\psi_i \neq
\psi_j$ for some $1 \leq i,j \leq r$. By conjugating with a suitable
element of~$W$ we may assume that there is some $i < r$ such that, in
the first case, $\lambda_1 = \cdots = \lambda_i$, but $\lambda_j \neq 
\lambda_1$ for all $j > i$, and, in the second case, $\psi_1 = \cdots = 
\psi_i$, but $\psi_j \neq \psi_1$ for all $j > i$. (Such a conjugation
changes~$L_0$ and~$\psi$, but it does not change $R_{L_0}^G(\psi)$ by
\cite[(70.11)]{cr1}.) Put $k = \sum_{j = 1}^i \lambda_j$ and 
$m = n - k$. Let $L \cong \GL_k(q) \times \GL_m(q)$ be the standard 
Levi subgroup corresponding to the composition $(k,m)$ of~$n$. Then
$L_0 \leq L$ and $W_G(L_0,\psi) = W_L(L_0,\psi)$ by the remarks preceding
the proposition.  It follows from Theorem~\ref{HCresult} that every
constituent of $R_{L_0}^G(\psi)$ is Harish-Chandra imprimitive.
\end{prf}

\subsection{Other classical groups}
We now aim at a similar result for the other classical groups,
excluding the groups of 9.1(f).
If~$G$ is one of the groups of \ref{TheClassicalGroups}(b)--(d) 
or one of the orthogonal groups or Clifford groups of 
\ref{TheClassicalGroups}($\text{\rm e}'$),(e) of $+$~type, 
we put $m := \lfloor n/2 \rfloor$. In the remaining cases we put $m := 
n/2 - 1$. Then $m \geq 1$. The Weyl group $W$ of~$G$ is of type $B_m$, 
$m \geq 1$, of type $C_m$, $m \geq 2$, or of type~$D_m$, $m \geq 4$, and
we write $S = \{ s_1, s_2, \ldots , s_m \}$ for the set of fundamental 
reflections of~$W$. If~$W$ is of type~$B_m$ or $C_m$ with $m \geq 3$,
we choose notation so that 
$s_2, \ldots , s_m$ are conjugate, but $s_1$ is not conjugate to~$s_2$.
If~$W$ is of type~$B_2$ or~$C_2$, we let $s_1$ correspond to a short root
or a long root, respectively.
If~$W$ is of type~$D$ we let $s_1$, $s_2$, and $s_m$ be the end nodes
of the Dynkin diagram with~$s_1$ and $s_2$ of distance~$2$ (if~$m$ is
larger than~$4$, this uniquely determines the set $\{ s_1, s_2 \}$).

Let~$J$ be a proper subset of~$S$.
Let $k := m + 1  - \mbox{\rm min}\{ 1 \leq i \leq m \mid s_i \not\in J \}$. 
We call the subset $\{ s_1, \ldots , s_{m-k} \}$ of~$J$ its component of 
type~$B$,~$C$ or~$D$, respectively (in the latter case only if $s_1$ and 
$s_2$ are contained in~$J$). The subset 
$\{s_{m-k+2}, \ldots , s_m\} \cap J$ of $\{s_{m-k+2}, \ldots , s_m\}$ 
determines a composition $(\lambda_1, \ldots , \lambda_r)$ of~$k$ as above 
(in case of the general linear group). The standard Levi subgroup $L_J$ 
of~$G$ is of the form
\begin{equation}
\label{L_J}
L_J \cong G_{n-2k}(q) \times \GL_{\lambda_1}(q^\delta) \times
\GL_{\lambda_2}(q^\delta) \times \cdots \times \GL_{\lambda_r}(q^\delta).
\end{equation}
Here, $\delta = 2$ if $G = \GU_n(q)$, and $\delta = 1$, otherwise.
If~$G$ is not one of the Clifford groups, an explicit description of the 
embedding of the right hand side of~(\ref{L_J}) as a standard Levi subgroup 
of~$G$ can be found, e.g., in \cite[4.4]{HiKe}.
For the Clifford groups use duality as in \cite[8.5]{luszbuch}. The group
$G_{n-2k}(q)$ is of the same type as~$G$ and corresponds to the component 
of~$J$ of type~$B$,~$C$ or~$D$, respectively; the index~$n-2k$ indicates the 
dimension of the natural (projective) module of $G_{n-2k}(q)$.
It may happen that the component of type~$B$ (respectively~$D$) of~$J$
is empty, in which case $n - 2k \in \{ 0, 1 \}$; in this case the group 
$G_{n-2k}(q)$ is a cyclic torus.

Let~$\psi$ be an irreducible cuspidal character of~$L := L_J$. Then
\begin{equation}
\label{psi}
\psi = \psi_0 \otimes \psi_1 \otimes \psi_2 \otimes \cdots
\otimes \psi_r,
\end{equation}
where $\psi_0$ is an irreducible cuspidal character of $G_{n-2k}(q)$, and 
$\psi_i$ is an irreducible cuspidal character of $\GL_{\lambda_i}(q^\delta)$
for $1 \leq i \leq r$.

Suppose first that $s_1$ and $s_2$ are contained in~$J$ if~$W$ is of 
type~$D_m$. Then there is, for every $1 \leq i \leq r$, an involution
$\sigma_i \in W_G(L)$
which centralizes every direct factor $\GL_{\lambda_j}(q^\delta)$, $1 \leq j 
\leq r$, of the decomposition~(\ref{L_J}) except $\GL_{\lambda_i}(q^\delta)$.
If $G = \GU_n(q)$, the map induced by $\sigma_i$ on $\GL_{\lambda_i}(q^2)$ is
the Frobenius morphism~$F$ of Example~\ref{ExamplesForAsymptotics}(b) 
(restricted from $\GL_{\lambda_i}(\bar{\mathbb{F}}_q)$ to $\GL_{\lambda_i}(q^2)$).
In the other cases~$\sigma_i$ induces (modulo an inner automorphism), the
transpose inverse automorphism on $\GL_{\lambda_i}(q)$. With respect to a
suitable basis of~$V$ (recall that~$V$ is the natural (projective) module 
of~$G$), some preimage of $\sigma_i$ in~$N$ is represented by the matrix
\begin{equation}
\label{si}
\left( \begin{array}{ccccc}
{I}_{a} & 0 & 0 & 0 & 0 \\
0 & 0 & 0 & I_{\lambda_i} & 0 \\
0 & 0 & \tilde{I} & 0 & 0 \\
0 & I_{\lambda_i} & 0 & 0 & 0 \\
0 & 0 & 0 & 0 & {I}_{a}
\end{array}\right),
\end{equation}
with $a = \sum_{j = i+1}^r \lambda_j$. Also, $\tilde{I}^2 = I_b$
with $b = n - 2k + 2\sum_{j=1}^{i-1} \lambda_j$. Note that~$\sigma_i$ need 
not normalize $G_{n-2k}(q)$ (see Example~\ref{ConditionsOnPsi}(a)).

Now suppose that $W$ is of type~$D_m$ and either~$s_1$ or~$s_2$
is not contained in~$J$. Without loss of generality we may assume that
$s_1 \not\in J$. (Indeed, if $s_1 \in J$ but $s_2 \not\in J$, consider
$J' := (J \cup \{ s_2 \}) \setminus \{ s_1 \}$; then~$J'$ is conjugate 
to~$J$ in~$W$, and thus $L_{J}$ and $L_{J'}$ give rise to the same 
Harish-Chandra induced characters of~$G$ (see \cite[(70.11)]{cr1}).)
In this case the Levi subgroup $L_J$ is as in
(\ref{L_J}) with $\delta = 1$, $k = m$ and $G_0(q)$ a cyclic group of
order~$q-1$. Fix an~$i$, $1 \leq i \leq r$. If $\lambda_i$ is even, there
still is an involution $\sigma_i \in W_G(L)$ as in the general case above.
It is represented on~$V$ by the matrix~(\ref{si}) with identity
matrix~$\tilde{I}$. If $\lambda_i$ is odd, then there is, for every
$1 \leq j \neq i \leq r$ with $\lambda_j$ odd, an involution $\sigma_{ij}
\in W_G(L)$ inducing the inverse transpose automorphism (modulo an inner
automorphism) simultaneously on $\GL_{\lambda_i}(q)$ and on 
$\GL_{\lambda_j}(q)$, and centralizing the other direct factors 
$\GL_{\lambda_h}( q )$ of~$L_J$. The matrix of~$\sigma_{ij}$ on~$V$ is the 
product of the matrices for $\sigma_i$ and $\sigma_j$ (which do not lie 
in~$G$, but their product does).

Let us now return to the general case.
For every pair $(i,j)$, $1 \leq i < j \leq r$ such that $\lambda_i = 
\lambda_j$, there is an element $\tau_{ij} \in W_G(L)$ swapping 
$\GL_{\lambda_i}(q^\delta)$ with $\GL_{\lambda_j}(q^\delta)$ and centralizing
the other factors of~$L$. The elements $\sigma_i$ and $\tau_{ij}$ 
(respectively $\sigma_i$, $\sigma_{ij}$, and $\tau_{ij}$) generate
$W_G(L)$ by the results of Howlett \cite{How}.

The $\sigma_i$ and the $\tau_{ij}$ (respectively $\sigma_i$, $\sigma_{ij}$,
and $\tau_{ij}$) arise as the non-trivial elements of subgroups of $W_G(L) = 
(N_{G}(L) \cap N)L/L$ of the form $(N_{L_1}(L) \cap N)L/L$ for Levi 
subgroups $L_1 \geq L$ of~$G$ such that $(N_{L_1}(L) \cap N)L/L$ has order~$2$.
Such Levi subgroups~$L_1$ are groups of $F$-fixed points of admissible
subgroups in the sense of Lusztig \cite[p.~255]{luszbuch}. The $\sigma_i$
(respectively the $\sigma_i$ and $\sigma_{ij}$) are distinguished from the
$\tau_{ij}$ by the fact that the former arise from admissible Levi subgroups
$L_1$ for which the rank of the component of the root system of type~$B$,~$C$ 
or~$D$, respectively is larger than the corresponding rank for~$L$.
Note that the $\sigma_i$ (respectively the $\sigma_i$ and $\sigma_{ij}$)
and the $\tau_{ij}$ are also distinguished by the fact that the $\sigma_i$
(respectively the $\sigma_i$ and $\sigma_{ij}$) normalize the direct factors
$\GL_{\lambda_1}(q^\delta), \ldots , \GL_{\lambda_r}(q^\delta)$ of~$L$. 

Again by \cite[(8.5.13)]{luszbuch}), the group $W_G(L,\psi)$ is generated 
by the elements $\sigma_i$ and $\tau_{ij}$ (respectively $\sigma_i$,
$\sigma_{ij}$, and $\tau_{ij}$) that leave $\psi$ invariant.

\begin{prop}
\label{PrimitivesInClassicalII}
Let~$G$ be one of the classical groups 
from~{\rm \ref{TheClassicalGroups}(b)--(e)},
let~$L_0 = L_J$ be a proper Levi subgroup of~$G$ as in~{\rm (\ref{L_J})},
and let~$\psi$ be a cuspidal irreducible character of~$L_0$ as
in~{\rm (\ref{psi})}.

{\rm (a)} In case~$G$ is an orthogonal group of $+$~type, assume
in addition that $n - 2k \geq 4$ (this is equivalent to the assumption that
$s_1$ and $s_2$ are contained in~$J$).

Then if $\sigma_i \in W_G(L_0,\psi)$ for all $1 \leq i \leq r$, every
irreducible constituent of $R_{L_0}^G(\psi)$ is Harish-Chandra primitive.
Otherwise every such constituent is Harish-Chandra imprimitive.

{\rm (b)} Let~$G$ be an orthogonal group of $+$~type, and 
assume that $n-2k = 0$ (this is equivalent to the assumption that $s_1$ is not
contained in~$J$; recall that we do not need to consider the case $s_1 \in J$,
$s_2 \not\in J$).

Then if $\sigma_i \in W_G(L_0,\psi)$ for all $1 \leq i \leq r$ with
$\lambda_i$ even, and if for all $1 \leq i \leq r$ with $\lambda_i$ odd,
there is some $1\leq j \neq i \leq r$ with $\lambda_j$ odd such that
$\sigma_{ij} \in W_G(L_0,\psi)$, every irreducible constituent of
$R_{L_0}^G(\psi)$ is Harish-Chandra primitive. Otherwise every such
constituent is Harish-Chandra imprimitive.

\end{prop}
\begin{prf}
We only prove~(a); the proof of~(b) is similar.
Suppose first that $\sigma_i \in W_G(L_0,\psi)$ for all $1 \leq i \leq r$.
Let~$L$ be a proper Levi subgroup of~$G$ such that~$L_0$ is a Levi subgroup
of~$L$. Then for some~$i$, no preimage of $\sigma_i$ is contained in $L$.
In particular, $W_G(L_0,\psi) \neq W_{L}(L_0,\psi)$. It follows as in the
proof of Proposition~\ref{PrimitivesInClassicalI} that no irreducible
constituent of $R_{L_0}^G(\psi)$ is Harish-Chandra imprimitive.

Suppose now that $\sigma_i \not \in W_G(L_0,\psi)$ for some $1 \leq i
\leq r$. By conjugating with a suitable element of~$W$, we may assume 
that $\lambda_i = \lambda_{i+1} = \cdots = \lambda_r$, and that $\tau_{ij} 
\in W_G(L_0,\psi)$ for all $j > i$, but $\tau_{i',i} \not\in W_G(L_0,\psi)$ 
for all $i' < i$ with $\lambda_{i'} = \lambda_i$. Then $\sigma_j \not\in 
W_G(L_0,\psi)$ for all $j > i$, since $\sigma_i = \tau_{ij} \sigma_j 
\tau_{ij}$. Put $k' = \sum_{j=i}^r \lambda_j$, and let $L = G_{n-2k'}(q) 
\times \GL_{k'}(q^\delta)$. Then $L \neq G$ and $W_G(L_0,\psi) = 
W_{L}(L_0,\psi)$. The assertion now follows from Theorem~\ref{HCresult}.
\end{prf}

Let us look at an example to see what the conditions in the two propositions
really mean.
\begin{exl}
\label{ConditionsOnPsi}
{\rm
For ${a} \in \GL_n(q)$ write $F'({a}) := J_n {a}^{-T} J_n$. 
(For the notation see Subsection~\ref{NotationFromLinearAlgebra}.)

(a) Let $G = \CSO_{16}^+(q)$ with odd~$q$. Take $J = \{ s_2, \ldots , s_5, 
s_7, s_8\}$. Then $L := L_J = G_0(q) \times \GL_5(q) \times \GL_3(q)$
(so that $r = 2$, $\lambda_1 = 5$ and $\lambda_2 = 3$). We may choose an 
embedding of the direct product above into~$G$ in such a way
that $(z, {a}, {b} )$ with $z \in \mathbb{F}_q^* = G_0(q)$,
${a} \in \GL_5(q)$, ${b} \in \GL_3(q)$ is represented by
the matrix 
$$
\left(
\begin{array}{cccc}
{b} & 0 & 0 & 0 \\
0 & {a} & 0 & 0 \\
0 & 0 & z F'({a}) & 0 \\
0 & 0 & 0 & z F'({b})
\end{array}\right).
$$
Thus $\sigma_{12}$
maps $(z, {a}, {b} )$ to $(z, z F'({a}), 
z F'({b}) )$. Now let $\psi_0$, $\psi_1$, and $\psi_2$
be irreducible cuspidal characters of $G_0(q)$, $\GL_{\lambda_1}(q) 
= \GL_5(q)$, and $\GL_{\lambda_2}(q) = \GL_3(q)$, respectively. Suppose 
that the linear characters of the centers of $\GL_5(q)$ and $\GL_3(q)$
induced by $\psi_1$ and $\psi_2$ are $\xi_1$ and $\xi_2$, respectively.

It is easy to check that $\psi = \psi_0 \otimes \psi_1 \otimes \psi_2$
is invariant under $\sigma_{12}$, if and only if $\xi_1 = \xi_2^{-1}$
and $\psi_i( {c} ) =  \psi_i( {c}^{-T} )$
for all ${c} \in \GL_{\lambda_i}(q)$ and $i = 1, 2$.

(b) Next let $G = \CSO_8^+(q)$ with odd~$q$. We wish to determine all
irreducible Harish-Chandra imprimitive characters of~$G$ corresponding
to the Levi subgroup $L_0 = \CSO_4^+(q) \times \GL_2(q)$ (i.e., $L_0 = L_J$
with $J = \{ s_1, s_2, s_4 \})$.
We may choose an embedding of the direct product into~$G$ in such a way
that $({a}, {b} )$ with ${a} \in \CSO^+_4(q)$, 
${b} \in \GL_2(q)$ is represented by the matrix
$$
\left(\begin{array}{ccc}
{b} & 0 & 0 \\
0 & {a} & 0 \\
0 & 0 & \mu_a F'({b})
\end{array}\right),
$$
where $\mu_a \in \mathbb{F}_q^*$ is the multiplier of~${a}$.

Again, let $\psi_0$ and $\psi_1$ be irreducible cuspidal characters
of $G_4(q) = \CSO_4^+(q)$ and $\GL_{\lambda_1}(q) = \GL_2(q)$, respectively.
We have to determine the exact conditions under which $\psi :=
\psi_0 \otimes \psi_1$ is invariant under $\sigma_1$. Now~$\sigma_1$ maps
$( {a}, {b} )$ to $( {a}, \mu_a F'({b}) )$.
It is easy to check that
$\psi$ is invariant under $\sigma_1$ if and only if the center of
$\GL_2(q)$ is in the kernel of~$\psi_1$ and $\psi_1( {b} ) =
\psi_1( {b}^{-T} )$ for all ${b} \in \GL_2(q)$.
}
\end{exl}

\section{Lusztig series}
Let~${\mathbf G}$ be a connected reductive algebraic classical group 
over $\bar{\mathbb F}_q$ with connected center, and let~$F$ be a 
Frobenius morphism of~${\mathbf G}$, such that $G = {\mathbf G}^F$ 
is one of the groups introduced in~\ref{TheClassicalGroups}(a)--(e). 
As in Section~\ref{LusztigSeries}, we let $(\mathbf{G}^*,F)$ be the 
pair dual to $(\mathbf{G},F)$, and we write $G^* = {\mathbf{G}^*}^F$.
The natural modules for~$\mathbf{G}$ and~$\mathbf{G}^*$ are denoted
by~$\mathbf{V}$ and~$\mathbf{V}^*$, respectively.

Let $s \in G^*$ be semisimple. Since the center of~$\mathbf{G}$ is 
connected, the ordinary irreducible characters of~$G$ lying in 
$\mathcal{E}(G,[s])$ can be labelled as $\chi_{s,\lambda}$,
where~$\lambda$ is a unipotent character of $C_{G^*}(s)$.

\subsection{General linear groups}
Again we begin with the general linear groups $G = \GL_n(q)$ where we
may and will identify~$\mathbf{G}$ with~$\mathbf{G}^*$. An element of
$\GL_n(q)$ is called irreducible if and only if it acts irreducibly on
the natural module of $\GL_n(q)$.

\begin{prop}
\label{PrimitivesInGLn}
Let $G = \GL_n(q)$. Let $s \in G$ be semisimple, and let $\lambda$ be a 
unipotent character of $C_G(s)$.
Then the following assertions are equivalent:

{\rm (a)} The irreducible character $\chi_{s,\lambda}$ 
of~$G$ is Harish-Chandra primitive.

{\rm (b)} The minimal polynomial of~$s$ is irreducible.

{\rm (c)} There are $m, d \in {\mathbb N}$ with $md = n$ such
that $s$ is conjugate in~$G$ to a block diagonal matrix 
with~$m$ blocks, all equal to some irreducible matrix.

{\rm (d)} There are $m, d \in {\mathbb N}$ with $md = n$ such that
$C_G(s) \cong \GL_m(q^d)$.
\end{prop}
\begin{prf}
Let $\mathbf{G} = \GL_n(\bar{\mathbb{F}}_q)$ with standard Frobenius
map~$F$ so that $G = \mathbf{G}^F$. Now $\chi_{s,\lambda} \in
\mathcal{E}(G,[s])$. Thus $\chi_{s,\lambda}$ is Harish-Chandra 
imprimitive if and only if $C_{\mathbf{G}}(s)$ is contained in
a proper split Levi subgroup of $\mathbf{G}$ 
(Theorem~\ref{FirstLusztigSeriesResult}). This is the case if and only 
if $C_{\mathbf{G}}(s)$ stabilizes a non-trivial proper subspace of~$\mathbf{V}$,
i.e., if and only if $C_{\mathbf{G}}(s)$ acts reducibly on~$\mathbf{V}$.

Let $\mu$ denote the minimal polynomial of~$s$, and let~$\mu_1$ be an
irreducible factor of~$\mu$. Then $\mbox{\rm Ker}(\mu_1(s))$ is invariant 
under $C_{\mathbf{G}}(s)$. Hence if $\chi_{s,\lambda}$ is Harish-Chandra
primitive then $\mu = \mu_1$ is irreducible. Thus (b) follows from~(a).
The implications from (b) to (c) and from (c) to (d) are clear (with~$d$ 
the degree of~$\mu$). Also, if $C_G(s) \cong \GL_m(q^d)$ for $m, d \in 
{\mathbb N}$ with $md = n$, then a Coxeter torus of $C_G(s)$ is a
Coxeter torus of~$G$, and so~$C_G(s)$ acts irreducibly on~$V$.
It follows that $C_{\mathbf{G}}(s)$ is not contained in
a proper split Levi subgroup of $\mathbf{G}$ so that $\chi_{s,\lambda}$ 
is Harish-Chandra primitive.
\end{prf}

Two particular cases are worth noting. The case $d = 1$ corresponds to 
the central elements $s \in G$, the Harish-Chandra primitive characters
arising in this way are the unipotent characters of~$G$ (for $s = 1$),
and multiples of unipotent characters with linear characters
(for $s \neq 1$). The case $d = n$ gives rise to the cuspidal irreducible
characters of~$G$.

Note also that Proposition~\ref{PrimitivesInGLn}(d) gives the degrees of
the irreducible Harish-Chandra primitive characters of~$\GL_n(q)$, via
Lusztig's Jordan decomposition of characters.

\subsection{Other classical groups}
We now prove a similar result for the other classical groups.
If~$\mu$ is a monic irreducible polynomial over $\mathbb{F}_{q^2}$, we write
$\mu^*$ for the monic irreducible polynomial over $\mathbb{F}_{q^2}$ whose
roots are the $(-q)$ths powers of the roots of~$\mu$, i.e., if $\mu = 
\prod_{i=1}^d (X - \alpha_i)$ with $\alpha_i \in \bar{\mathbb{F}}_q$, then
$\mu^* = \prod_{i=1}^d  (X - \alpha_i^{-q})$.
Note that if $\mu \in \mathbb{F}_q[X]$, then $\mu^*$ is the monic polynomial
whose roots are the inverses of the roots of~$\mu$. If~$G^*$ is one of
the Clifford groups, i.e.,
a group as in \ref{TheClassicalGroups}~(d) or~(f), we have surjective 
homomorphisms $\mathbf{G}^* \rightarrow \bar{\mathbf{G}}^*$, defined over 
$\mathbb{F}_q$, where $\bar{\mathbf{G}}^* = \SO_{n}(\bar{\mathbb{F}}_q)$.
%in Case~(d), and $\bar{G}^* = \CSO_{n}(\bar{\mathbb{F}}_q)$ in Case~(f).
In these cases we write $\bar{s}$ for the image of $s \in G^*$ in 
$\bar{G}^* := \bar{\mathbf{G}}^{*^F}$. In the other cases we let 
$\bar{\mathbf{G}}^* = \mathbf{G}^*$ and $\bar{s} = s$ for $s \in G^*$. 
Note that the natural module~$\mathbf{V}^*$ of $\mathbf{G}^*$ is also 
the natural module of $\bar{\mathbf{G}}^*$.

\begin{prop}
\label{PrimitivesInOtherClassicals}
Let~$G$ be one of the classical groups 
of~{\rm \ref{TheClassicalGroups}(b)--(e)}.
Let $s \in G^*$ be semisimple and let $\lambda$ be a unipotent character
of $C_{G^*}(s)$. Then the following conditions are equivalent:

{\rm (a)} The irreducible character $\chi_{s,\lambda}$
of $G$ is Harish-Chandra primitive.

{\rm (b)} Every irreducible factor~$\mu$ of the  minimal polynomial 
of~$\bar{s}$ (acting on the natural module~$\mathbf{V}^*$ 
of~$\bar{\mathbf{G}}^*$) satisfies $\mu = \mu^*$.

%{\rm (c)} The element~$s$ does not lie in any proper split
%Levi subgroup of~$G^*$.
%
{\rm (c)} If~$G$ is a unitary group, there are nonzero $m_i, d_i
\in {\mathbb N}$ with $\sum_{i=1}^r m_id_i = n$ such that $C_G(s) \cong
\GU_{m_1}(q^{d_1}) \times \cdots \times \GU_{m_r}(q^{d_r})$.

If~$G$ is one of the other classical groups, there are nonzero $m_i, d_i 
\in {\mathbb N}$ with $\sum_{i=1}^r m_id_i \leq n$ such that $C_{-1} \times
C_1 \times \GU_{m_1}(q^{d_1}) \times \cdots \times \GU_{m_r}(q^{d_r})$ is
isomorphic to a subgroup of $C_{\bar{G}^*}(\bar{s})$ of index at 
most~$2$. Here, for $\zeta = -1, 1$, the group $C_\zeta$ is defined by
$C_\zeta = C_{\bar{\mathbf{G}}^*(\mathbf{V}^*_\zeta)}(\bar{s}_\zeta)^F$,
where $\mathbf{V}^*_\zeta$ denotes the $\zeta$-eigenspace of~$\bar{s}$
on $\mathbf{V}^*$, $\bar{s}_\zeta$ is the map induced by $\bar{s}$
on~$\mathbf{V}^*_\zeta$, and $\bar{\mathbf{G}}^*(\mathbf{V}^*_\zeta)$
is the same type of group as $\bar{\mathbf{G}}^*$, however acting on
$\mathbf{V}^*_\zeta$.
\end{prop}
\begin{prf}
Suppose that~$\mu$ is a monic irreducible factor of the irreducible 
polynomial of~$\bar{s}$ with $\mu \neq \mu^*$. If $\bar{\mathbf{G}}^*
= \mathbf{G}^*$, i.e., if $\mathbf{G}^*$ is not a Clifford group, then 
$\mbox{\rm Ker}_{\mathbf{V}^*}(\mu(\bar{s}))$ is a totally isotropic 
subspace of~$\mathbf{V}^*$. Since
$\mbox{\rm Ker}_{\mathbf{V}^*}(\mu(\bar{s}))$ is invariant under
$C_{\mathbf{G}^*}(s)$, it follows that $C_{\mathbf{G}^*}(s)$ is contained
in a proper split Levi subgroup of~$\mathbf{G}^*$, and thus~(a) does not 
hold.  In case $\mathbf{G}^*$ is a Clifford group, let $\prod_i \bar{s}_i$
denote the primary decomposition of $\bar{s}$, the $\bar{s}_i$ acting on
the subspaces $\mathbf{V}^*_i$ of $\mathbf{V}^*$ (cf.\ \cite[(1.10)]{fosri2}).
Put $\mathbf{C}_i := C_{\SO(\mathbf{V}^*_i)}(\bar{s}_i)$. By what we have
already proved above, there is an~$i$ such that $\mathbf{C}_i$ is contained
in a proper split Levi subgroup of $\SO(\mathbf{V}^*_i)$. Thus the product
$\prod_i \mathbf{C}_i$, embedded in the natural way into $\bar{\mathbf{G}}^*$,
is contained in a proper split Levi subgroup of $\bar{\mathbf{G}}^*$.
By \cite[{\bf (2B)}]{fosri2}, the centralizer $C_{\mathbf{G}^*}(s)$ is
the inverse image of $\prod_i \mathbf{C}_i$ under the surjection from
$\mathbf{G}^*$ to $\bar{\mathbf{G}}^*$. It follows that $C_{\mathbf{G}^*}(s)$
is contained in a proper split Levi subgroup of $\mathbf{G}^*$. 
Hence~(a) implies~(b).

If~(a) does not hold then $C_{\mathbf{G}^*}(s)$ is contained in a proper
split Levi subgroup of~$\mathbf{G}^*$. This implies that $\bar{s}$ is 
contained in a proper split Levi subgroup of~$\bar{\mathbf{G}}^*$, and thus
fixes a non-trivial totally isotropic
subspace of~$\mathbf{V}^*$. The restriction of~$\bar{s}$ to this subspace
has a minimal polynomial all of whose irreducible factors $\mu$ satisfy 
$\mu \neq \mu^*$. Thus~(b) implies~(a).

The equivalence of~(b) and~(c) follows from the description of the
centralizers of semisimple elements in the classical groups as 
given in ~\cite[Proposition (1A)]{fosri1} and \cite[(1.13)]{fosri2}.
%
%Now~(b) is certainly equivalent to the condition that~$\bar{s}$ does not
%lie in any proper split Levi subgroup of~$\bar{\mathbf{G}}^*$. This in
%turn is equivalent to the condition in~(c).
\end{prf}

\section{Examples for the restriction to commutator subgroups}

The results above are valid for finite groups of Lie type arising
from an algebraic group with connected center. In general, the groups
obtained this way are not quasisimple. We therefore have to investigate
the descent from a group as above to its commutator subgroup. We include
an example, a special case of a much more general example shown to us
by C{\'e}dric Bonnaf{\'e}, to give a flavor of the phenomena that can
occur. 

\begin{exl} \label{Bonnafe}
{\rm
Let
$$s = \left(\begin{array}{cc} I_2 & 0 \\ 0 & -I_2 \end{array}\right),$$
where~$I_2$ denotes the identity matrix of dimension~$2$. We view $s$ as 
element of $G = \GL_4(q)$ for some odd~$q$. We then have $C_G( s ) = 
\GL_2(q) \times \GL_2(q)$. Hence $\mathcal{E}(G,[s])$ contains the four 
characters $\chi_{s, 1 \otimes 1}$, $\chi_{s, 1 \otimes q}$, 
$\chi_{s, q \otimes1}$, and $\chi_{s, q \otimes q}$, corresponding to 
the four unipotent characters of $\GL_2(q) \times \GL_2(q)$ (where we have 
denoted the two unipotent characters of $\GL_2(q)$ by their degrees).
Since $C_{\mathbf{G}}( s )$ is a split Levi subgroup of ${\mathbf{G}} = 
\GL_4( \bar{\mathbb{F}}_q )$, the characters 
in $\mathcal{E}(G,[s])$ are Harish-Chandra induced from the characters of 
$\mathcal{E}(C_G( s ), [s])$ by Theorem~\ref{FirstLusztigSeriesResult}.

Let~$z$ be a generator of $\mathbb{F}_q^* = \GL_1(q)$. We identify the
maximally split torus~$L_0$ of~$G$ with the direct product $\GL_1(q) 
\times \cdots \times
\GL_1(q)$, so that $L_0 = \{( z^i, z^j, z^k, z^l ) \mid 1 \leq i, j, k, l 
\leq q - 1\}$. Let $\psi$ be the linear complex character of~$L_0$ defined
by $\psi(z^i, z^j, z^k, z^l) = (-1)^{k+l}$. Then the Harish-Chandra series
of $(L_0,\psi)$ is exactly $\mathcal{E}(G,[s])$. Note that $W_G(L_0,\psi)$
is generated by the two permutation matrices corresponding to the 
transpositions $(1,2)$ and $(3,4)$.

Now let $\tilde{G} = \SL_4(q)$ and $\tilde{L}_0 = \tilde{G} \cap L_0$.
The restriction~$\tilde{\psi}$ of~$\psi$ to~$\tilde{L}_0$ is cuspidal
and irreducible. However, $W_{\tilde{G}}(\tilde{L}_0,\tilde{\psi})$
is a dihedral group of order~$8$, since~$\tilde{\psi}$ is also invariant
under the simultaneous exchange of the first two with the last two
coordinates. Thus $R_{\tilde{L}_0}^{\tilde{G}}(\tilde{\psi})$ has
exactly five irreducible constituents, one occurring with multiplicity~$2$.
These characters also occur as the irreducible constituents of the
restrictions to $\tilde{G}$ of the elements of $\mathcal{E}(G,[s])$,
since the restriction of $R_{L_0}^G( \psi )$ to~$\tilde{G}$ equals
$R_{\tilde{L}_0}^{\tilde{G}}(\tilde{\psi})$ (see, e.g., 
\cite[Proposition 13.22]{DiMi2}).

It follows from Clifford theory that the two characters 
$\chi_{s, 1 \otimes q}$ and $\chi_{s, q \otimes 1}$
of $\mathcal{E}(G,[s])$ restrict to the
same irreducible character of~$\tilde{G}$, where\-as each of the
other two elements of $\mathcal{E}(G,[s])$ splits into two
conjugate irreducible characters of~$\tilde{G}$. 

Since only the $2$-dimensional irreducible character of 
$W_{\tilde{G}}(\tilde{L}_0,\tilde{\psi})$ is induced from a 
proper subgroup, there is at most one irreducible constituent
of $R_{\tilde{L}_0}^{\tilde{G}}(\tilde{\psi})$ which is
Harish-Chandra imprimitive. In fact, the restriction of 
$\chi_{s, 1 \otimes q}$ to~$\tilde{G}$ is imprimitive, being 
Harish-Chandra induced from an irreducible character of degree~$q$ 
of the Levi subgroup $C_{\tilde{G}}(s)$ of~$\tilde{G}$.
Thus the other four irreducible constituents of 
$R_{\tilde{L}_0}^{\tilde{G}}(\tilde{\psi})$ are Harish-Chandra
primitive (hence primitive by Theorem~\ref{ReductionToParabolic}).

Let~$\bar{s}$ denote the image of~$s$ in $\PGL_4(q)$. The irreducible
constituents of the restriction of the elements of $\mathcal{E}(G,[s])$
to~$\tilde{G}$ are contained in $\mathcal{E}(\tilde{G},[\bar{s}])$ (see,
e.g., \cite[Proposition 11.7]{CeBo2}). 
Thus, unlike for groups arising from algebraic groups with connected
center, a rational Lusztig series can contain primitive and imprimitive
characters at the same time.
}
\end{exl}

\chapter{Exceptional groups}
\label{ExceptionalGroups}

In addition to the classical groups introduced in 
Section~\ref{TheClassicalGroups}, we here investigate exceptional groups of 
Lie type. The ground field~$K$, assumed to be algebraically closed as always, 
is now allowed to have any characteristic $\ell \geq 0$, but different from 
the defining characteristic of the group considered. By a character we always
mean a $K$-character.

\section{The exceptional groups of type~$E$ and $F$}\label{TheExceptionalGroups}

Assume that $\mbox{\rm char}(K) = 0$.
Let $\mathbf{G}$ be a simple adjoint algebraic group with a Dynkin diagram
of type~$E$ or~$F_4$.
Moreover, we assume that~$F$ is a Frobenius morphism of $\mathbf{G}$
such that, in the respective cases, $G := \mathbf{G}^F$ is a finite group of 
type $E_i(q)$, $i = 6, 7, 8$, ${^2\!E}_6(q)$ or $F_4(q)$.
In the cases $E_8(q)$ and $F_4(q)$, the group $G$ is simple. In the other 
cases, the commutator subgroup $G'$ of~$G$ is simple.

Since~$\mathbf{G}$ is adjoint, its center is trivial and we may thus apply
Theorem~\ref{FirstLusztigSeriesResultConverse}. By this and 
Theorem~\ref{FirstLusztigSeriesResult}, in order to find the Harish-Chandra
imprimitive irreducible characters of~$G$, we have to determine the semisimple
elements $s \in G^*$ such that $C_{\mathbf{G}^*}( s )$ is not contained in any
proper split $F$-stable Levi subgroup of~$\mathbf{G}^*$. This property only 
depends on the
$G^*$-conjugacy class of $C_{\mathbf{G}^*}( s )$, hence is a property of the
semisimple class type of~$s$. (Two semisimple elements $s_1, s_2 \in G^*$ 
belong to the same {\em class type}, if and only if $C_{\mathbf{G}^*}( s_1 )$
and $C_{\mathbf{G}^*}( s_2 )$ are conjugate in $G^*$; see 
\cite[Subsection~$4.2$]{chevie}.) The semisimple conjugacy classes of the
groups considered here have been determined explicitly by Shoji \cite{shoji}
and Shinoda \cite{shino1} for type $F_4$, and by Fleischmann and Janiszczak 
\cite{FleJaE6uE7,FleJaE8} for type~$E$. 

The semisimple elements in question can be read off from the tables by 
Fleischmann and Janiszczak as follows. These tables give, for each semisimple 
element $s \in G^*$, the order of $Z^\circ( C_{\mathbf{G}^*}( s ) )^F$ as a 
polynomial $f_s(q)$ in~$q$. (As in the proof of Theorem~\ref{Asymptotic}, we
write ${Z^\circ( C_{\mathbf{G}^*}( s ) )}$ for the connected component of the
centre of $C_{\mathbf{G}^*}( s )$.) Then $C_{\mathbf{G}^*}( s )$ is contained 
in a proper split $F$-stable Levi subgroup of $\mathbf{G}^*$, if and only if 
$f_s(q)$ is divisible by $q - 1$ (as a polynomial in~$q$). To see this, first 
note that the split $F$-stable Levi subgroups of $\mathbf{G}^*$ are exactly 
the centralizers of split tori, i.e., $F$-stable tori of~$\mathbf{G}^*$ whose 
order polynomial is of the form $(X - 1)^a$ for some non-negative integer~$a$ 
(see \cite[Section 3.D]{BrouMa}; for the concept of order polynomial see 
\cite[D{\'e}finition~$1.9$]{BrouMa}.) If $Z^\circ( C_{\mathbf{G}^*}( s ) )$ 
contains a split torus $\mathbf{S}^*$, then $C_{\mathbf{G}^*}( s )$ is contained 
in the split $F$-stable Levi subgroup $C_{\mathbf{G}^*}( \mathbf{S}^* )$. On the 
other hand, if $C_{\mathbf{G}^*}( s )$ is contained in the split $F$-stable 
Levi subgroup $\mathbf{L}^*$, then $Z^\circ( \mathbf{L}^* )$ is contained in 
$Z^\circ( C_{\mathbf{G}^*}( s ) )$ and so the latter group contains a split 
torus.

\section{Explicit results on some exceptional groups}
\label{ExplicitExceptional}

Here we give explicit results on those exceptional groups of Lie type, for 
which (partial) generic ordinary character tables are available in the 
literature. These are the groups ${^2\!B}_2(q)$, ${^2\!G}_2(q)$, ${^2\!F}_4(q)$,
$G_2(q)$, and ${^3\!D}_4(q)$. We obtain a complete classification of the
imprimitive irreducible $K$-characters of these groups. The following 
observation essentially reduces to the case of $\mbox{\rm char}(K) = 0$.

\begin{lem}
\label{SmallRank}
Suppose that $\mbox{\rm char}(K) = \ell > 0$, and let~$G$ be a finite group
with a split $BN$-pair of characteristic $p \neq \ell$.
If the irreducible $K$-characters of every proper Levi subgroup~$L$ of~$G$, 
are liftable, then the irreducible Harish-Chandra imprimitive $K$-characters
of~$G$ are liftable to ordinary Harish-Chandra imprimitive characters of~$G$.
\end{lem}
\begin{prf}
The claim is a reformulation of Lemma~\ref{Modp}(2).
\end{prf}

\noindent We note that the assumption of the Lemma above is satisfied
for the groups ${^2\!B}_2(q)$, ${^2\!G}_2(q)$, %${^2\!F}_4(q)$,
$G_2(q)$, and ${^3\!D}_4(q)$. The proper split Levi subgroups of the listed
groups are abelian or of type~$A_1$. By the results of Burkhardt \cite{Bu}
and James \cite{James}, the irreducible characters (in cross-characteristics) 
of groups of the latter type are all liftable. The group ${^2\!F}_4(q)$
has a Levi subgroup ${^2\!B}_2(q) \times (q - 1)$, for which the 
assumption on liftability is not satisfied. The exceptions can be treated 
with Theorem~\ref{HCresult}.

If $G = G_2(q)$ or $G = {^3\!D}_4(q)$, then $G = \mathbf{G}^F$, where the 
algebraic group $\mathbf{G}$ satisfies the hypothesis of Section~\ref{Section73} 
(see \cite[1.19]{cart}). In these cases we can use 
Theorems~\ref{FirstLusztigSeriesResult} 
and~\ref{FirstLusztigSeriesResultConverse}. Occasionally we sketch alternative 
arguments for some of the results. For the groups ${^2\!B}_2(q)$ and 
${^2\!G}_2(q)$ we use Theorem~\ref{HCresult} and its converse, which is clear 
for these groups, since they have only one proper split Levi subgroup.

For the irreducible characters we use the notation of the original
sources for the character tables. We also write $\Phi_m$ for the
value at~$q$ of the $m$th cyclotomic polynomial, $m \in \mathbb{Z}$,
$m > 0$.

We begin with the Suzuki groups.
\begin{prop}
\label{2B2}
Suppose that $\mbox{\rm char}(K) = 0$.
Let $G = {^2\!B}_2( q )$ with $q = 2^{2m+1}$, $m \geq 1$
be a simple Suzuki group. Then exactly the irreducible characters 
of degree $q^2 + 1$ are imprimitive. They are Harish-Chandra induced
from certain linear characters of a Borel subgroup of~$G$.
We have $|\Irr(G)| = q + 3$ and $q/2 - 1$ of these are
imprimitive. 
\end{prop}

We now consider the simple Ree groups of characteristic~$3$.
\begin{prop}
\label{2G2}
Suppose that $\mbox{\rm char}(K) = 0$.
Let $G = {^2\!G}_2( q )$ with $q = 3^{2m+1}$, $m \geq 1$ be a simple Ree 
group (for the group ${^2\!G}_2( 3 )' \cong \SL_2(8)$ see 
Chapter~\ref{ExceptionalMultipliers}). Then exactly the irreducible
characters of degree $q^3 + 1$ are imprimitive. They are Harish-Chandra
induced from certain linear characters of a Borel subgroup of~$G$.
We have $|\Irr(G)| = q + 8$ and $(q - 3)/2$ of these are imprimitive.
\end{prop}

We next consider the simple Ree groups of characteristic~$2$. The unipotent
characters of these groups have been computed by Malle~\cite{malleuni};
the complete character table is available in CHEVIE~\cite{chevie}.
\begin{prop}
\label{2F4}
Let $G = {^2\!F}_4( q )$ with $q = 2^{2m+1}$, $m \geq 1$
be a simple Ree group (of characteristic~$2$).

{\rm (1)} Suppose that $\mbox{\rm char}(K) = 0$.
Using the {CHEVIE} notation for irreducible characters of~$G$,
and denoting the maximally split torus of~$G$ by~$T_1$, 
Table~\ref{ReeImprim} lists the imprimitive
irreducible characters of~$G$. The column headed ``Levi'' indicates
the smallest Levi subgroup of~$G$ of which the irreducible characters
are Harish-Chandra induced.

We have $|\Irr(G)| = q^2 + 4q + 17$, and $9/16\,q^2 + 11/8\,q - 5$
of these are imprimitive.

{\rm (2)} Suppose that $\mbox{\rm char}(K) = \ell > 2$. If $\ell \nmid
q - \sqrt{2q} + 1$, every imprimitive irreducible $K$-character of~$G$
is liftable. Otherwise, there is a series of $(q-2)/2$ non-liftable
imprimitive irreducible $K$-characters of degree $\Phi_1 \Phi_2^3 
\Phi_4 \Phi_6 \Phi_{12}$. These characters are $\ell$-modular constituents 
of the ordinary characters $\chi_{25}$ (notation from Table~\ref{ReeImprim}),
and arise by Harish-Chandra induction from irreducible $K$-characters of
degree $q^2 - 1$ of the Levi subgroup ${^2\!B}_2(q) \times (q - 1)$.
\end{prop}
\begin{prf}
(1) We use the notation of Shinoda~\cite{shino2} and Malle~\cite{mallemax}
for the maximal tori of~$G$. Four of them, $T_1, \ldots , T_4$ lie in
proper Levi subgroups of~$G$. The corresponding irreducible 
Deligne-Lusztig characters are the characters $\chi_{28}$, $\chi_{34}$,
$\chi_{36}$, and $\chi_{35}$, respectively.

Multiplying the trivial character of $\SL_2(q)$ with a non-trivial
character of the torus $q - 1$ we obtain an irreducible character
of the  Levi subgroup $\SL_2(q) \times (q - 1)$ whose Harish-Chandra
induction is a character of type~$\chi_{26}$. Similarly, the Steinberg 
character of $\SL_2(q)$ gives rise to the characters of type~$\chi_{27}$.
In exactly the same manner we obtain the characters $\chi_{22}$,
$\chi_{25}$, $\chi_{23}$ and $\chi_{24}$ from the unipotent characters
of the Suzuki group ${^2\!B}_2(q)$.

The other irreducible characters of~$G$ are primitive, since their
degrees do not match any product of the index of a proper parabolic
subgroup of~$G$ with the degree of an irreducible character of the
corresponding Levi subgroup. This has been checked using {CHEVIE}.

(2) The decomposition numbers of the Suzuki groups ${^2\!B}_2(q)$ in odd
characteristic have been determined by Burkhardt~\cite{BuSz}. 
All irreducible $K$-characters of ${^2\!B}_2(q)$ are liftable unless
$\ell \mid q - \sqrt{2q} + 1$. 

In the latter case, there is a unique non-liftable irreducible $K$-character, 
which is cuspidal and has degree $q^2 - 1$. Multiplying this with a 
non-trivial character of the cyclic factor $q - 1$ of ${^2\!B}_2(q) 
\times ( q - 1)$, we obtain an irreducible, cuspidal $K$-character 
of~$G$ through Harish-Chandra induction (see Theorem~\ref{HCresult}). 
%\marginpar{Why are these not liftable?}
This yields $(q - 2)/2$ distinct imprimitive irreducible $K$-characters
of~$G$ of the asserted degree; these characters are non-liftable, as
$(q + 1)^3$ does not divide the order of~$G$.
\end{prf}

The above result is also true for the smallest of the characteristic~$2$
Ree groups ${^2\!F}_4(2)$, in the sense that none of its irreducible
characters is Harish-Chandra induced. This fact can easily be checked
with the Atlas \cite{ATLAS}. There are of course imprimitive irreducible
characters of ${^2\!F}_4(2)$, whose block stabilizer is the Tits group
${^2\!F}_4(2)'$. For the latter group see Proposition~\ref{TitsGroup}.

\begin{table}[htb]
\caption{\label{ReeImprim}Imprimitive, ordinary irreducible characters
of Ree's groups ${^2\!F}_4(q)$, $q = 2^{2m+1}$, $m \geq 1$}
$$\begin{array}{llll}\hline
\mbox{\rm Char.} & \mbox{\rm Degree}
& \mbox{\rm Number} & \mbox{\rm Levi}
\rule[- 7pt]{0pt}{ 20pt} \\ \hline
\chi_{28} & \Phi_2^2\Phi_4^2\Phi_{6}\Phi_{12} &
\frac{1}{16}(q - 2)(q - 8) & T_1 \rule[ 0pt]{0pt}{ 13pt}\\
\chi_{34} & \Phi_1\Phi_2\Phi_4^2\Phi_{6}\Phi_{12} &
\frac{1}{4}q(q-2) & \SL_2(q) \times (q - 1) \\
\chi_{35} & \Phi_1\Phi_2^2\Phi_6'\Phi_4\Phi_6\Phi_{12} &
\frac{1}{8}(q-2)(q - \sqrt{2q}) & {^2\!B}_2(q) \times (q - 1) \\
\chi_{36} & \Phi_1\Phi_2^2\Phi_6''\Phi_4\Phi_6\Phi_{12} &
\frac{1}{8}(q-2)(q + \sqrt{2q}) & {^2\!B}_2(q) \times (q - 1) \\
\chi_{26} & \Phi_2\Phi_4^2\Phi_6\Phi_{12} &
\frac{1}{2}(q-2) & \SL_2(q) \times (q - 1) \\
\chi_{27} & q\Phi_2\Phi_4^2\Phi_6\Phi_{12} &
\frac{1}{2}(q-2) & \SL_2(q) \times (q - 1) \\
\chi_{22} & \Phi_2^2\Phi_4\Phi_6\Phi_{12} &
\frac{1}{2}(q-2) & {^2\!B}_2(q) \times (q - 1) \\
\chi_{25} & q^2\Phi_2^2\Phi_4\Phi_6\Phi_{12} &
\frac{1}{2}(q-2) & {^2\!B}_2(q) \times (q - 1) \\
\chi_{23} & \frac{1}{2}\sqrt{2q}\Phi_1\Phi_2^2\Phi_4\Phi_6\Phi_{12} &
\frac{1}{2}(q-2) & {^2\!B}_2(q) \times (q - 1) \\
\chi_{24} & \frac{1}{2}\sqrt{2q}\Phi_1\Phi_2^2\Phi_4\Phi_6\Phi_{12} &
\frac{1}{2}(q-2) & {^2\!B}_2(q) \times (q - 1)
\rule[- 7pt]{0pt}{ 5pt} \\ \hline
\end{array}
$$
$\Phi_6' = q + \sqrt{2q} + 1$, $\Phi_6'' = q - \sqrt{2q} + 1$
\end{table}

We now investigate the Chevalley groups $G_2(q)$.
\begin{prop}
\label{G2}
Suppose that $\mbox{\rm char}(K) = 0$.
Let $G = G_2( q )$ be a simple Chevalley of type~$G_2$,
where~$q$ is any prime power larger than~$2$ (for the group
$G_2(2)' \cong \SU_3(3)$ see Chapter~\ref{ExceptionalMultipliers}).

Using the notation for irreducible characters of~$G$ introduced
by Chang and Ree in~\cite{chare}, Tables~\ref{G2odd} and~\ref{G2even}
list the imprimitive
irreducible characters of~$G$. The column headed ``Levi'' indicates
the smallest Levi subgroup of~$G$ of which the irreducible characters
are Harish-Chandra induced.

If $q$ is odd and not divisible by~$3$, we have 
$|\Irr(G)| = q^2 + 2q + 9$. For $q \equiv 1(\mbox{\rm mod\ } 3)$,
there are $7/12\,q^2 + 1/3\,q - 71/12$ imprimitive irreducible characters,
whereas for $q \equiv 2(\mbox{\rm mod\ } 3)$, there are 
$7/12\,q^2 + 1/3\,q - 17/4$
of them.

If~$q$ is a power of~$3$, we have $|\Irr(G)| = q^2 + 2q + 8$,
and there are $7/12\,q^2 + 1/3\,q - 17/4$ imprimitive irreducible
characters.

If~$q$ is even, we have $|\Irr(G)| = q^2 + 2q + 8$. If 
$q \equiv 1(\mbox{\rm mod\ } 3)$, there are $7/12\,q^2 + 1/3\,q - 14/4$ 
imprimitive irreducible characters, wheras if $q \equiv 2(\mbox{\rm mod\ } 3)$,
there are $7/12\,q^2 + 1/3\,q - 3$ of them.
\end{prop}
\begin{prf}
The irreducible characters of~$G$ can be found in~\cite{enya}, if~$q$
is even, in~\cite{enomo}, if~$q$ is a power of~$3$, and in~\cite{chare},
in the general case. Explicit tables are also given in~\cite{HiHa}, as
well as in {CHEVIE}, where a dictionary is given to compare the
various notations in these papers.
\end{prf}

\begin{table}[htb]
\caption{\label{G2odd}Imprimitive, ordinary irreducible characters
of the Chevalley groups $G_2(q)$, $q$ odd}
$$\begin{array}{lllll}\hline
\multicolumn{1}{c}{\mbox{\rm Char.}} &
\multicolumn{1}{c}{\mbox{\rm Degree}} &
\multicolumn{1}{c}{
\begin{array}{c} \mbox{\rm Number} \\ q \equiv 1(\mbox{\rm mod\ } 3) \end{array}} &
\multicolumn{1}{c}{
\begin{array}{c} \mbox{\rm Number} \\ q \not\equiv 1(\mbox{\rm mod\ } 3) \end{array}} &
\multicolumn{1}{c}{\mbox{\rm Levi}}
\rule[- 7pt]{0pt}{ 20pt} \\ \hline
X_1 & \Phi_2^2 \Phi_3 \Phi_6 &
\frac{1}{12}(q^2 - 8q + 19) &
\frac{1}{12}(q^2 - 8q + 15) & (q - 1)^2
\rule[ 0pt]{0pt}{ 13pt}\\
X_a, X_b & \Phi_1\Phi_2\Phi_3\Phi_6 &
\frac{1}{4}(q^2 - 2q + 1) & \frac{1}{4}(q^4 - 2q + 1) & \GL_2(q) \\
X_{1a}' & \Phi_2\Phi_3\Phi_6 &
\frac{1}{2}(q - 5) & \frac{1}{2}(q - 3) & \GL_2(q) \\
X_{1a} & q\Phi_2\Phi_3\Phi_6 &
\frac{1}{2}(q - 5) & \frac{1}{2}(q - 3) & \GL_2(q) \\
X_{1b}' & \Phi_2\Phi_3\Phi_6 &
\frac{1}{2}(q - 3) & \frac{1}{2}(q - 3) & \GL_2(q) \\
X_{1b} & q\Phi_2\Phi_3\Phi_6 & 
\frac{1}{2}(q - 3) & \frac{1}{2}(q - 3) & \GL_2(q) 
\rule[- 7pt]{0pt}{ 5pt} \\ \hline
\end{array}
$$
\end{table}

\begin{table}[htb]
\caption{\label{G2even}Imprimitive, ordinary irreducible characters
of the Chevalley groups of type~$G_2(q)$, $q > 2$ even}
$$\begin{array}{lllll}\hline
\multicolumn{1}{c}{\mbox{\rm Char.}} &
\multicolumn{1}{c}{\mbox{\rm Degree}} &
\multicolumn{1}{c}{
\begin{array}{c} \mbox{\rm Number} \\ q \equiv 1(\mbox{\rm mod\ } 3) \end{array}} &
\multicolumn{1}{c}{
\begin{array}{c} \mbox{\rm Number} \\ q \equiv 2(\mbox{\rm mod\ } 3) \end{array}} &
\multicolumn{1}{c}{\mbox{\rm Levi}}
\rule[- 7pt]{0pt}{ 20pt} \\ \hline
X_1 & \Phi_2^2 \Phi_3 \Phi_6 &
\frac{1}{12}(q^2 - 8q + 16) &
\frac{1}{12}(q^2 - 8q + 12) & (q - 1)^2
\rule[ 0pt]{0pt}{ 13pt}\\
X_a, X_b & \Phi_1\Phi_2\Phi_3\Phi_6 &
\frac{1}{4}(q^2 - 2q) & \frac{1}{4}(q^4 - 2q) & \GL_2(q) \\
X_{1a}' & \Phi_2\Phi_3\Phi_6 &
\frac{1}{2}(q - 4) & \frac{1}{2}(q - 2) & \GL_2(q) \\
X_{1a} & q\Phi_2\Phi_3\Phi_6 &
\frac{1}{2}(q - 4) & \frac{1}{2}(q - 2) & \GL_2(q) \\
X_{1b}' & \Phi_2\Phi_3\Phi_6 &
\frac{1}{2}(q - 2) & \frac{1}{2}(q - 2) & \GL_2(q) \\
X_{1b} & q\Phi_2\Phi_3\Phi_6 &
\frac{1}{2}(q - 2) & \frac{1}{2}(q - 2) & \GL_2(q)
\rule[- 7pt]{0pt}{ 5pt} \\ \hline
\end{array}
$$
\end{table}

Finally we consider Steinberg's triality groups ${^3\!D}_4( q )$.
\begin{prop}
\label{3D4}
Suppose that $\mbox{\rm char}(K) = 0$.
Let $G = {^3\!D}_4( q )$ be a Steinberg triality group,
where~$q$ is any prime power.
Using the notation for irreducible characters of~$G$ introduced
in~\cite{DeMi}, Table~\ref{TrialityImprim} lists the imprimitive
irreducible characters of~$G$. The column headed ``Levi'' indicates
the smallest Levi subgroup of~$G$ of which the irreducible characters
are Harish-Chandra induced.

If $q$ is even, we have $|\Irr(G)| = q^4 + q^3 + q^2 + q + 5$, and
of these,
$$\frac{7}{12} q^4 + \frac{1}{6} q^3 - \frac{5}{6} q^2 - \frac{2}{3} q
- 3$$
are imprimitive.
If $q$ is odd, we have $|\Irr(G)| = q^4 + q^3 + q^2 + q + 6$, and
of these,
$$\frac{7}{12} q^4 + \frac{1}{6} q^3 - \frac{5}{6} q^2 - \frac{2}{3} q
- \frac{17}{4}$$
are imprimitive.
\end{prop}
\begin{prf}
There are three maximal tori lying in proper Levi subgroups, namely
$T_0$, $T_1$ and $T_2$ (in the notation of~\cite{DeMi}).
The corresponding irreducible Deligne-Lusztig characters
are the characters named $\chi_6$, $\chi_8$, and $\chi_{11}$,
respectively.

The characters of type $\chi_{3,1}$ and $\chi_{3,\mathrm{St}}$
are obtained by Harish-Chandra inducing suitable characters of~$T_0$,
the maximally split torus of~$G$. If~$\lambda$ is such a linear
character of~$T_0$, we have $R_{T_0}^G(\lambda) = \chi_1 + \chi_\mathrm{St}$,
with~$\chi_1$ of type $\chi_{3,1}$ and $\chi_\mathrm{St}$ of type
$\chi_{3,\mathrm{St}}$.
If $L$ denotes the Levi subgroup $\SL_2(q^3) \times (q - 1)$, then
$R_{T_0}^L(\lambda) = \psi_1 + \psi_\mathrm{St}$, with $\psi_1$
of degree~$1$ and $\psi_\mathrm{St}$ of degree~$q^3$. It follows
that $\chi_1 = R_L^G(\psi_1)$ and $\chi_\mathrm{St} =
R_L^G(\psi_\mathrm{St})$.
The assertion for the characters of type $\chi_{5,1}$ and
$\chi_{5,\mathrm{St}}$ is proved similarly.

The numbers of the irreducible characters of each type can
be found in~\cite{DeMi}.

The other irreducible characters of~$G$ are primitive, since their
degrees do not match any product of the index of a proper parabolic
subgroup of~$G$ with the degree of an irreducible character of the
corresponding Levi subgroup. This has been checked using {CHEVIE}
\cite{chevie}.
\end{prf}

\begin{table}[htb]
\caption{\label{TrialityImprim}Imprimitive, ordinary irreducible characters
of Steinberg's triality groups ${^3\!D}_4(q)$}
$$\begin{array}{llp{2.70cm}p{2.70cm}l}\hline
\multicolumn{1}{c}{\mbox{\rm Char.}} &
\multicolumn{1}{c}{\mbox{\rm Degree}} &
\multicolumn{1}{c}{
\begin{array}{c} \mbox{\rm Number} \\ q \mbox{\rm\ even} \end{array}} &
\multicolumn{1}{c}{
\begin{array}{c} \mbox{\rm Number} \\ q \mbox{\rm\ odd} \end{array}}  &
\multicolumn{1}{c}{\mbox{\rm Levi}}
\rule[- 7pt]{0pt}{ 20pt} \\ \hline
\chi_6 & \Phi_2^2 \Phi_3 \Phi_6^2 \Phi_{12} &
$\frac{1}{12}(q^4 - 4q^3 + 2q^2 -2q + 12)$ &
$\frac{1}{12}(q^4 - 4q^3 + 2q^2 -2q + 15)$ & (q^3 - 1) \times ( q - 1)
\rule[ 0pt]{0pt}{ 13pt}\\
\chi_8 & \Phi_1\Phi_2\Phi_3\Phi_6^2\Phi_{12} &
$\frac{1}{4}(q^4 - 2q)$ & $\frac{1}{4}(q^4 - 2q + 1)$ & \SL_2(q) \times (q^3 - 1) \\
\chi_{11} & \Phi_1\Phi_2\Phi_3^2\Phi_6\Phi_{12} &
$\frac{1}{4}(q^4 - 2q^3)$ & $\frac{1}{4}(q^4 - 2q^3 + 1)$ & \SL_2(q^3) \times (q - 1)\\
\chi_{3,1} & \Phi_2\Phi_3\Phi_6\Phi_{12} &
$\frac{1}{2}(q - 2)$ & $\frac{1}{2}(q - 3)$ & \SL_2(q^3) \times ( q - 1)\\
\chi_{3,\mathrm{St}} & q^3\Phi_2\Phi_3\Phi_6\Phi_{12} &
$\frac{1}{2}(q - 2)$ & $\frac{1}{2}(q - 3)$ & \SL_2(q^3) \times ( q - 1)\\
\chi_{5,1} & \Phi_2\Phi_3\Phi_6^2\Phi_{12} &
$\frac{1}{2}(q^3-q^2-q-2)$ & $\frac{1}{2}(q^3-q^2-q-3)$ & \SL_2(q) \times (q^3 - 1) \\
\chi_{5,\mathrm{St}} & q\Phi_2\Phi_3\Phi_6^2\Phi_{12} &
$\frac{1}{2}(q^3-q^2-q-2)$ & $\frac{1}{2}(q^3-q^2-q-3)$ & \SL_2(q) \times (q^3 - 1)
\rule[- 7pt]{0pt}{ 5pt} \\ \hline
\end{array}
$$
\end{table}

%\appendix
%    Include appendix "chapters" here.
%\include{}
%
\backmatter
%    Bibliography styles amsplain or harvard are also acceptable.
%\bibliographystyle{amsalpha}
%\bibliography{}

%    See note above about multiple indexes.
%\printindex

\end{document}